\documentclass[12pt]{report}
\usepackage{amssymb}
\usepackage{amsmath,amsfonts,amsthm}
\usepackage{srcltx}
\usepackage{makeidx}
\makeindex
\newcommand{\jump}{\marginpar{\vspace*{-.1in}$$\mid\!\! \mid\!\!\mid$$\vspace*{-.45in}
     $$\blacktriangledown$$}}
\newcommand{\fab}[2]{\langle {#1}_1,{#1}_2,\ldots ,{#1}_{#2} \rangle}
\newcommand{\abs}[1]{\:|#1|}    \newcommand{\Ker}{{\rm Ker\,}}
     \newcommand{\Z}{ Z\!\!\!Z}
\newcommand{\C}{ l\!\!\!C}       
\newcommand{\R}{ I\!\!R}         \newcommand{\N}{ I\!\!N}
\newcommand{\Q}{ l\!\!\!Q}       \newcommand{\notsubset}{ /\!\!\!\!\!\!\subset }
\newcommand{\id}{\mbox{id}}     \newcommand{\notdivide}{\!\not |\,}
\newcommand{\pf}{{\bf Proof: }}    \newcommand{\F}{ I\!\!F}
 
\newcommand{\st}{\; \vline\;}   \newcommand{\norm}[1]{\:\parallel #1 \parallel}
\newcommand{\orb}{\mbox{ \bf Orb}}  \newcommand{\stab}{\mbox{ \bf Stab}}
\newcommand{\rhdn}{\rhd_{_{\!\!\!\!\!\!\neq }}}
\newcommand{\ba}{$$\begin{array}}\newcommand{\ea}{\end{array}$$}
\newcommand{\bea}{\begin{eqnarray*}}\newcommand{\eea}{\end{eqnarray*}}
\newcommand{\be}{\begin{equation}}\newcommand{\ee}{\end{equation}}
\newcommand{\bd}{\begin{definition}\bf} \newcommand{\ed}{\end{definition}}
\newcommand{\brs}{\begin{remarks}\rm}   \newcommand{\ers}{\end{remarks}}
\newcommand{\br}{\begin{remark}\rm}     \newcommand{\er}{\end{remark}}
\newcommand{\bt}{\begin{theorem}}       \newcommand{\et}{\end{theorem}}
\newcommand{\bl}{\begin{lemma}}         \newcommand{\el}{\end{lemma}}
\newcommand{\bco}{\begin{corollary}}    \newcommand{\eco}{\end{corollary}}
\newcommand{\bp}{\begin{proposition}}   \newcommand{\ep}{\end{proposition}}
\newcommand{\bo}{\begin{observation}\rm}\newcommand{\eo}{\end{observation}}
\newcommand{\bex}{\begin{examples}\rm}   \newcommand{\eex}{\end{examples}}
\newcommand{\bos}{\begin{observations}\rm}\newcommand{\eos}{\end{observations}}
\newcommand{\bx}{\begin{example}\rm}   \newcommand{\ex}{\end{example}}
\newcommand{\bexe}{\begin{exercise}\rm}   \newcommand{\eexe}{\end{exercise}}

\newtheorem{definition}{Definition}[chapter]
\newtheorem{theorem}[definition]{Theorem}
\newtheorem{lemma}[definition]{Lemma}
\newtheorem{corollary}[definition]{Corollary}
\newtheorem{proposition}[definition]{Proposition}
\newtheorem{remarks}[definition]{Remarks}
\newtheorem{remark}[definition]{Remark}
\newtheorem{observation}[definition]{Observation}
\newtheorem{examples}[definition]{Examples}
\newtheorem{observations}[definition]{Observations}
\newtheorem{example}[definition]{Example}
\newtheorem{exercise}[definition]{Exercise}
\reversemarginpar 
\begin{document}
\title{ \Huge \bf Excursions in Group Theory}
\author{\large \bf Pramod K Sharma\thanks{e-mail: pksharma1944@yahoo.com}\\ School of Mathematics\\ Devi Ahilya University, Indore(India)}
\date{}
\maketitle
 \tableofcontents
\chapter{ Groups And Subgroups }
\footnote{\it contents group1.tex }
  We encounter integers very early in our mathematical training.
Let $m,n,k$ $ \in$ $ \Z$. Then $$n + m = m + n \in \Z,$$ $$(m + n) + k = m
+ (n + k),$$ and $$0 + n = n = n + 0$$ for all  $n  \in \Z.$ Further, for
every $n \in \Z$ there is $-n \in \Z$ such that $n + (-n) = 0.$ There are
seeds of group theory in these observations. Let us note another such
structure  that  we have encountered in Linear Algebra  i.e., $Gl_n(\R)$,
the set of all n x n- matrices over $\R$ with non-zero determinant. Let
$A,B,C \in Gl_n (\R)$. Then $det(AB) = det A \cdot det B \neq 0$. Hence
$AB \in Gl_n(\R)$.We also have $(AB)C = A(BC)$ and for the identity matrix
$I_n$ of order $n$, $A \cdot I_n = I_n \cdot A = A$.Further, for any $A
\in Gl_n(\R)$, the inverse matrix of $A$ i.e.,$A^{-1}$ satisfies $AA^{-1}
= A^{-1}A = I_n$. One must note that for $m,n \in \Z,  m + n = n + m$,
however, for $A,B \in Gl_n(\R), AB \neq BA$ in general. We, now, define,
\bd\label{d1.1} \bf Let G be a non-empty set provided with a map,
$$\begin{array}{crcl}
  f : & G \times G & \to & G \\
   & (x, y) & \mapsto & x \ast y
\end{array}$$ called multiplication (a
binary operation) satisfying: \\ (i)For all $x,y,z \in G$,
 $(x\ast y)\ast z = x\ast(y\ast z)$ (associativity).\\ (ii)
 There exists an element $e \in G$ such that, $x\ast e = e\ast
x = x$  for all $x \in G$. The element $e$ is called an identity
element of $G$. \\ (iii)
 For each $x \in G$ there is an element $w \in G$(depending on $x$) such that
$x\ast w = w\ast x = e$. The element $w$ is called an inverse of
$x$.\\  Then G is called a \index{group}group with respect to the
binary
operation ($\ast$). \\

If, moreover, $x\ast y = y\ast x$ for all $x, y \in G$, then $G$ is
called an \bf \index{abelian group}abelian group.\ed \rm For a group
$(G,\ast )$, we shall simply write $G$ is a group i.e., the binary
operation will not be mentioned unless  required and we shall write
$xy$ for $x\ast y$ for $x,y \in G$ when $\ast$ is well understood.
Some times we may need to specify the map $f$, in that case we shall
write that $(G, f)$ is a group. \brs\label{r1.2} (i) The identity
element is unique ( If $e_1 ,e_2 \in G$ are identities then $e_1 =
e_1 \ast e_2 = e_2)$. We shall normally denote this by $e$, and
shall write $e_G$ for $e$ if it is necessary to emphasise $G$.\\
(ii) Inverse of an element $x \in G$ is uniquely determined. If
$x_1, x_2$ are inverses of $x \in G$, then $x_1  = x_1\ast e =
x_1\ast (x*x_2) = (x_1\ast x)*x_2 = e\ast x_2 = x_2$,The inverse
of the element $x$ in G is denoted by $x^{-1}$. \ers
\bexe\label{ec1.3} Show that subtraction is not an associative
binary operation on $\Z$. \eexe \bexe\label{ec1.4} Prove that the
map, \bea
  f:\Z\times \Z & \to     & \Z \\
(m,n)           & \mapsto & m\ast n = m + n - mn \eea gives an
associative binary operation, but $(\Z,\ast)$ is not a group.
\eexe  \bexe\label{ec1.5} Let $G$ be a group, and let $b \in G$ be a
fixed element of $G$. Prove that $G$ is group with respect to the
binary operation : \bea
    f_b : G\times G & \to     & G \\
     (x,y)          & \mapsto & xby
\eea Find the identity element of $(G,f_b)$ and also the inverse
of an element  $x \in G$. \eexe  \bexe Let $G$ be a group and $x,
y \in G$. Prove that  \\ (i) $(x^{-1})^{-1} = x$  \\ (ii)
$(x^{-1}yx)^{n} = x^{-1}y^{n}x$ for all $n \geq 0$.  \\ (iii)
$(xy)^{-1} = y^{-1}x^{-1}$.\eexe Group structure is encountered
quite often. We, now, give some examples:
 \bx\label{ex1.6} $(\R,+)$, $(\, \Q,+)$,
$(\, \C,+)$ and $(\R^*,\cdot)$, $(\, \Q^*,\cdot)$, $(\,
\C^*,\cdot)$ are groups, where $\cdot$ denotes the usual product
and $\R^* = \R - \{0\}$ , $\Q^* = \Q - \{0\}$ , $\C^* = \, \C -
\{0\}$.\ex \bx Let $a$ be a non-zero integer. Then  $$\Q_{a} =
\{m/a^{n} \st m,n \in \Z \} $$ is a group with respect to
addition.\ex \bexe Let $a,b$ be two non-zero integers. Then
$\Q_{a}\subset \Q_{b}$ if and only if $a$ divides $b^{r}$ for some
$r\geq 1$. \eexe
\bx\label{ex1.7} {\bf(Additive Group Of Integers Modulo n)} Take an
integer $n
> 0$. For any $x, y $ in the set $$S = \{ 0, 1, 2, \cdots, n-1
\},$$ define
: $$
\begin{array}{ccc}
  x +_{n} y & = & \left\{ \begin{array}{lcl}
    x + y & \mbox{  if  } & 0 \leq x + y < n \\
    x + y - n & \mbox{  if  } & x + y \geq n
  \end{array}\right.
\end{array}$$ Then $(S, +_{n})$ is an abelian group where $0$ is the
identity element and $n - x$ is the inverse of $x$ for any $x \in S$. This
group is denoted by $\Z_{n}$ and the binary operation is called addition
modulo $n$.\ex \bx\label{ex1.8} {\bf(Multiplicative Group Of Integers
Modulo n)} Take an integer $n
> 1$ and put ${\Z^{*}}_{n} = \{1 \leq a < n \mid (a, n) = 1\}$.
For any $a, b \in {\Z^{*}}_{n}$, define : $$a \cdot_{n} b = r$$
where $ab = nq + r$ such that $1 \leq r < n$. Note that we can
always write $ab = nq + r$ where $0 \leq r < n$. However, as $(a,
n) = (b, n) = 1$,  $1 \leq r < n$. Further, it is clear that $(r,
n) = 1$. It is easy to check that $\cdot_{n}$ is associative and
if $a \in {\Z^{*}}_{n}$, then there exist $1 \leq b < n$ such that
$ab + nq = 1$ for some $q \in \Z$. Thus $a \cdot_{n} b = 1$. Hence
it is immediate that $({\Z^{*}}_{n}, \cdot_{n})$ is an abelian
group. The operation $\cdot_{n}$ is called multiplication mod
$n$.\ex \bx\label{ex1.9} $G = \{1, -1, i, -i\}$ is a group with
respect to multiplication. \ex \bx\label{ex1.10} The set \bea G &
= & \left\{
  {\left( \begin{array}{cc}
  1 & 0 \\
  0 & 1
\end{array} \right),}
 { \left( \begin{array}{cc}
  i & 0 \\
  0 & -i
\end{array}\right),}
 {\left( \begin{array}{cc}
  -1 & 0 \\
  0 & -1
\end{array}\right),}
 {\left( \begin{array}{cc}
  -i & 0 \\
  0 & i
\end{array}\right),}  \right.\\
 &     & \quad  \left.
  {\left( \begin{array}{cc}
  0 & -1 \\
  1 & 0
\end{array}\right),}
 {\left( \begin{array}{cc}
  0 & -i \\
  -i & 0
\end{array}\right),}
 { \left( \begin{array}{cc}
   0 & 1 \\
   -1 & 0\
 \end{array} \right),}
  {\left( \begin{array}{cc}
    0 & i \\
    i & 0 \\
  \end{array}\right)}\right\}\\
\eea
 of $2\times 2$-matrices over complex numbers  is  a
  group  with respect to matrix  multiplication. If we write  \\
 $I = \left( \begin{array}{cc}
  1 & 0 \\
  0 & 1
\end{array} \right)$, $a =   \left( \begin{array}{cc}
  i & 0 \\
  0 & -i
\end{array}\right)$, $b = \left( \begin{array}{cc}
  0 & -1 \\
  1 & 0
\end{array}\right)$, and $c = \left( \begin{array}{cc}
  0 & -i \\
  -i & 0
\end{array}\right)$ , then $$G = \{\pm I, \pm a, \pm b, \pm c\}$$ where $a^{2} = b^{2} = c^{2} =
-I$, $ab = c, \, bc = a$ and $ca = b$.  This  is called the group
of \bf Quaternions. \rm We shall denote this group by $\bf H$.\ex
\bx\label{ex1.11} The set $$K = \left\{ I =\left(
\begin{array}{cc}
  1 & 0 \\
  0 & 1
\end{array} \right), a =   \left( \begin{array}{cc}
  1 & 0 \\
  0 & -1
\end{array}\right), b = \left( \begin{array}{cc}
  0 & -1 \\
  1 & 0
\end{array}\right),c = \left( \begin{array}{cc}
  0 & -1 \\
  -1 & 0
\end{array}\right)\right\}$$ of $2 \times 2$ matrices over $\R$ is a
group with respect to multiplication of matrices. Here $ab = c,\, bc
= a,\, ca = b$ and $a^{2} = b^{2} = c^{2} = I$. This is called \bf
Klein's $4$-group\index{Klein's $4$-group} \rm and is denoted by
$V_{4}$.\ex \bx The Example \ref{ex1.11} can also be realized in a
rather interesting way
\\ Let $A$ and $B$ be two light switches for the bulbs $1$ and $2$
such that if the bulb $1$ is off then pressing $A$ once puts it on
and pressing once more puts it off. Similar is the function of the
switch $B$ for the bulb $2$. Then let the initial state be when
both bulbs are off . Let $$\tau = \mbox{ action of pressing } A$$
$$\sigma = \mbox{ action of pressing }B$$ $$e = \mbox{ not
pressing } A \mbox{ and } B$$ Then $$G = \{e, \tau, \sigma,
\tau\sigma\}$$ is a group with respect to the operation of doing
one action after the other on the initial state when both bulbs
are off. Clearly $\tau^{2} = e,\, \sigma^{2} = e,\, \sigma\tau =
\tau\sigma$ and $(\tau\sigma)^{2} = e$.Thus taking $a = \sigma$,
$b = \tau$ and $c = \sigma\tau$ we again get Klein's $4$-group.
\ex \bx\label{ex1.12} Let $n \geq 1$ be fixed. The set of $n^{th}$ roots
of unity in complex numbers, i.e., $$\C_n = \{ z\in \C \, |\, z^n = 1
\},$$ is a group with respect to multiplication. We have $\C_4 =
\{1,-1,i,-i\}$ (Example \ref{ex1.9}).
\ex

 \bx\label{ex1.13} The set ${\bf \C_\infty} = \{\alpha \in \C \, |\, \alpha, \mbox{
a root of unity}\}$ is a group with respect to multiplication. \ex
\bx\label{ex1.13b} If $m\geq 1$ is fixed, then $ \C_{m^{\infty}} =
\{\alpha \in \C \st \alpha^{m^{k}} = 1 \mbox{ for some } k\geq 1\}$
is a group with respect to multiplication.\ex \bx\label{ex1.14} The
set $Gl_n(\R) = \{A: n\times n-\mbox{matrix over} \R \,|\,detA\neq
0\}$  is a group with respect to multiplication, called the {\bf
general linear group}\index{general linear group} of dimension $n$
over $\R$. \ex \br\label{ec1.15} We can replace $\R$ by $\Q$ or $\C$
in the example \ref{ex1.14} but not by $\Z$. Check why? More
generally, $\R$ can be replaced by any field. We define
\[Gl_n(\Z)=\{A:n\times n-\mbox{ matrix over }\Z \st \det A=\pm 1\}.\]
This is a group with respect to multiplication.\er
 \bexe Let $K$ be a field and let $A$ be an $n
\times n$-matrix over $K$. Prove that $A \in Gl_{n}(K)$
  if and only if columns (rows) of $A$ are linearly independent vectors of the
  $K$-vector space $K^{n}$, equivalently, columns (rows) of $A$ form a basis of the
  vector space $K^{n}$.\eexe

\bx\label{es1} The set $$Sl_{n}(\R) =\{A: n\times n-\mbox{matrix
over} \R \,|\,detA=1\}$$  is a group with respect to matrix
multiplication. This is called the \bf special linear
group\index{special linear group} \rm of dimension $n$ over $\R$.\ex
\br It is easy  to check that $\R$ can be replaced by $\Q$ or $\C$
or $\Z$ in the above example. In general, $\R$ can be replaced by
any field. \er \bx\label{ex6} Let $S$ be a non-empty set and $$A(S)
= \{ f: S \to S \; |\; f: \mbox{ one-one, onto map }\}.$$ Then
$A(S)$ is a group with composition of maps as binary operation. This
is called the \bf group of transformations\index{group of
transformations} \rm of $S$, and the elements of $A(S)$ are called
transformations on $S$. \ex \bx\label{ex7} In the example \ref{ex6}
if we take $S = \{1,2,...,n\}$, then $A(S)$ is denoted by $S_n$ and
is also called the \bf symmetric group\index{symmetric group}\rm \,
or \bf permutation group\index{permutation group} \rm of degree $n$.
For any $\sigma \in S_n$ , we represent $\sigma$ as :
$$\sigma = \left(\begin{array}{cccc}
                    1     & 2          & \cdots & n\\
               \sigma (1) & \sigma (2) & \cdots & \sigma(n)
               \end{array}\right)$$
Thus for $n = 2$
$$S_2 = \left\{\left(\begin{array}{cc}
                    1     & 2 \\  1  & 2
              \end{array}\right),
                \left(\begin{array}{cc}
                    1     & 2 \\  2  & 1
              \end{array}\right)\right\}$$
             \ex
             and for $n=3$,
  \bea
S_3  & = & \left\{
                {\left( \begin{array}{ccc}
                    1 & 2 & 3\\  1  & 2 & 3
                     \end{array} \right),}
                     {\left( \begin{array}{ccc}
                    1 & 2 & 3 \\  2  & 3  & 1
                    \end{array} \right),}
                      {\left( \begin{array}{ccc}
                      1  & 2 & 3\\  3  & 2 & 1
                  \end{array} \right) ,} \right.\\
     &   &  \quad \left. {\left( \begin{array}{ccc}
                    1  & 2 & 3\\  2  & 1 & 3
                   \end{array}\right),}
                    {\left(\begin{array}{ccc}
                    1 & 2 & 3\\  1  & 3 & 2
                   \end{array}\right),}
                    {\left(\begin{array}{ccc}
                    1 & 2 & 3\\  3  & 1 & 2
                   \end{array}\right)}
        \right\}
\eea

\bx\label{ex8} Let $(X,d)$ be a metric space  and $$ I(X) = \{T: X
\to X | T : \mbox{ an isometry} \}$$ Then  $I(X)$ is a group with
respect to composition of maps as the binary operation. \ex \bo
Take $X = \R$ with standard metric. Let $ f \in I(\R)$ and $f(0) =
a$.  If for any $b \in \R$, $\tau_b :\R \to \R$ denotes the
isometry $\tau_b (x) = x + b$ then $\tau_{-a} f(0) = 0$. For any
isometry $g \in I(\R)$, if $g(0) = 0$,
 then,
$|g(x)| = |x|$. Thus $g(x) = x$ or $-x$. Clearly, as $g$ is
distance preserving, if $g(x) = -x$ for some $x( \neq 0) \in \R$,
then $g(y) = -y$ for all $y \in \R$. Thus $g = Id$ or $g(x) = - x$
for all $x \in \R$. Hence if $\theta \in I(\R)$ is such that
$\theta(x) = -x$ for  all $x$, we get $\tau_{-a}  f = \theta$ or
Id. Hence if $\tau_{-a}  f = \theta$, then $f = \tau_a \theta$.
However if $\tau_{-a}  f = Id$, then $f = \tau_a.$Therefore,
$$I(\R) = \{\tau_a \theta^t \;|\; t = 0 \mbox{ or $1$ and } a\in
\R\}.$$ \eo \bx\label{ex9} Let for $0\leq \theta \leq 2\pi$, $$
r_{\theta} : \R^2 \to \R^2$$ denotes the rotation of a point by an
angle $\theta$ around the origin. Then $G = \{r_{\theta} \; |\;
0\leq \theta \leq 2\pi \}$ is a group with composition as the binary
operation. By plane geometry, one can see that $$ G =
\left\{\left(\begin{array}{cc} \cos \theta & -\sin \theta\\ \sin
\theta & \cos \theta \end{array} \right) \st 0\leq \theta \leq 2\pi
\right\}$$ where $$r_{\theta} = \left(\begin{array}{cc} \cos \theta
& -\sin \theta\\ \sin \theta & \cos \theta \end{array} \right)$$ is
the $\R$-linear  transformation over $\R^2$ determined by the 2 x
2-matrix. \ex \bx\label{th11} \bf(Dihedral Group\index{Dihedral
Group}) : \rm Consider the $n$-gon $(n \geq 3)$ in $\R^{2}$ with its
vertices the points $P_{k} = (\cos{( 2\pi k/n )}, \; \sin{( 2\pi
k/n)})$ where $0 \leq k \leq n-1$. Let $\sigma$ denotes the counter
clock-wise rotation in the plane by an angle $2\pi /n$ and let
$\tau$ be the reflection about the axis which makes an angle $\pi/n$
with the $x$-axis. Then $\sigma : \R^{2} \rightarrow \R^{2}$ and
$\tau : \R^{2} \rightarrow \R^{2}$ are given by the matrices :
$$\sigma = \left(\begin{array}{cc}
  \cos ( 2\pi/n ) & -\sin ( 2\pi/n ) \\
  \sin ( 2\pi/n ) & \cos ( 2\pi/n )
\end{array}\right)$$ and $$\tau = \left(\begin{array}{cc}
  \cos ( 2\pi/n)  & \sin ( 2\pi/n ) \\
  \sin ( 2\pi/n ) & -\cos ( 2\pi/n )
\end{array}\right)$$  The set $D_{n} = \{Id, \sigma, \sigma^{2},
\cdots,\sigma^{n-1},
 \tau, \sigma\tau, \cdots
, \sigma^{n-1}\tau\}$ forms a group of order $2n$ with respect to
the binary operation, composition of maps. Here we have $\sigma^{n}
= Id = \tau^{2}$ and $\sigma\tau = \tau\sigma^{n-1}$. Check that
$D_{n}$ is the group of all rotational symmetries and reflections of
the regular $n$-gon in the plane.\ex To give further examples, we
define : \bd\label{d2}\bf A collineation\index{collineation} is a
bijection $T : \R^{2} \to \R^{2}$ such that a triple $a,b$ and $c$
of distinct points in $\R^2$ is co-linear if only if $T(a)$, $T(b)$
and $T(c)$ are co-linear. \ed \bd\label{d3}\bf A map $F : \R^{2} \to
\R^{2}$ is called an affine transformation\index{affine
transformation} if there exists an invertible $2 \times 2$ - matrix
$Q$ over $\R$ and a vector $\underline a \in \R^{2}$ such that for
all $\underline x \in \R^2$,
$$ F(\underline x) = Q \underline x + \underline a$$ \ed
\bx\label{ex10} Let $G = \{ T : \R^2 \to \R^2 \;|\; T : \mbox{ a
collineation} \}.$ Then  G is a group with respect to the
composition of maps as binary operation. \ex \bx\label{ex11} The set
$H = \{f:\R^2 \to \R^2 \;| \; f: \mbox{ affine transformation} \}$
is a group with composition of maps as the binary operation. \ex
\bexe\label{ex5e} For any $n \geq 1$, let $\xi = e^{{2 \pi i}/ n}$,
the $n^{th}$ root of unity in $\C$. If $$ a =
\left(\begin{array}{cc}
  0 & 1 \\
  1 & 0
\end{array}\right),
  b = \left(\begin{array}{cc}
  \xi & 0 \\
  0 & \xi^{-1}
\end{array}\right), e = \left(\begin{array}{cc}
  1 & 0 \\
  0 & 1
\end{array}\right),$$
 then show that $$G = \{e, a, b, b^2, \cdots, b^{n-1}, ab, ab^{2}, \cdots, ab^{n-1}\}$$
 is a group with respect to matrix multiplication.  \eexe
\bd\label{d4}\bf For a group $G$, the cardinality of the set $G$
 is called the order of\index{order of group}  $G$ and is denoted by $\circ
(G)$ or $|G|$. In case $G$ is finite, order of $G$ is the number of
elements in G. \ed  \bexe\label{ec25} Find the order of the group in
the example \ref{ex1.12}. \eexe  \bd\label{d5} \bf A non empty
subset $H$ of a group $G$ is called a subgroup\index{subgroup} of
$G$ if $H$ is a group with respect to the binary operation of G. \ed
\bexe Let $K$ be a field. Show that the set $D_{n}(K)$ of all
diagonal matrices in $Gl_{n}(K)$ is a subgroup of $Gl_{n}(K)$. (The
group $D_{n}(K)$ is called the \bf diagonal group \rm of dimension
$n$ over $K$).\eexe
 \bd \bf Let $A$ and $B$ be any two
non-empty subsets of a group $G$ Then $$AB = \{ ab \mid a \in A, b
\in B\} \mbox { and } A^{-1} = \{a^{-1} \mid a \in A\}$$ \ed
\bexe Prove that for any two non-empty subsets $A, B$ of a group
$G$, $(AB)^{-1} = B^{-1}A^{-1}$.\eexe \bexe Prove that if $H$ is a
subgroup of a group $G$, then $H^{-1} = H$  and $HH = H$.\eexe
 \bt\label{t1} Let $G$ be a group. A non empty subset
$H$ of $G$ is a subgroup of $G$ if and only if $a b^{-1} \in H$
for all $a,b \in H$. \et \pf ($\Rightarrow$) Let $a, b \in H$. As
$H$ is a subgroup of $G$, $b^{-1} \in H$. Hence as $a, b^{-1} \in
H$, $a b^{-1} \in H.$
\\($\Leftarrow$) Take $a, b \in H$. By assumption $bb^{-1} = e \in
H$. Further as $e, b \in H$, we conclude $eb^{-1} = b^{-1} \in H$
i.e., $H$ contains identity and has inverse for each of its
elements. Now, as seen above, for $a,b \in H, a,b^{-1} \in H$.
Hence by our assumption $a(b^{-1})^{-1} = ab \in H$. Thus $H$ is a
subgroup of $G$.     $\blacksquare$ \bco\label{co28} Let $H_1,
H_2,...,H_k$ be subgroups of a group $G.$ Then $H_1\cap H_2\cap
\ldots \cap H_k$ is a subgroup of $G.$ \eco \pf Let $x, y \in
H_1\cap H_2\cap \ldots \cap H_k$. Then as each $H_i$ is a subgroup
of $G$, $xy^{-1} \in H_i$ for all $ 1\leq i\leq k.$ Hence $xy^{-1}
\in H_1\cap H_2\cap \ldots \cap H_k$. Hence the result follows by
the theorem. $\blacksquare$ \bco\label{co1} Let $G$ be a group and
$H$ $(\neq \emptyset ) \subset G$ a finite subset. Let for any
$a,b \in H$, $ab \in H$ i.e., $H H \subseteq H$. Then $H$ is a
subgroup of $G$. \eco \pf Take $a,b \in H$. By our assumption on
$H$, $S = \{ b^i \st i \geq 1 \}$ is a subset of $H$. As $H$ is
finite there exist +ve integers $m < n$ such that $b^n = b^m$.
Hence $b^{n-m} = e$. This proves $b^{-1} = b^{n-m-1} \in H$. Now,
as $a,b^{-1} \in H$, $ab^{-1} \in H$.  Therefore by the theorem
 $H$ is a subgroup of $G$. $\blacksquare$
\bco\label{co111} Let $S$ be a non-empty subset of a group $G$.
Then $gp(S) = \{{a_{1}}^{e_{1}}{a_{2}}^{e_{2}}\cdots
{a_{m}}^{e_{m}} | m \geq 1, a_{i} \in S, \mbox{ and } e_{i} = \pm
1 \mbox{ for all } 1 \leq i \leq m\}$ is a subgroup of $G$. \eco
\pf Note that for any $x = {a_{1}}^{e_{1}}{a_{2}}^{e_{2}}\cdots
{a_{t}}^{e_{t}}$ such that $e_{i} = \pm 1$, $a_{i} \in S$ and $t
\geq 1$, we have $x^{-1} = a_{t}^{-e_{t}}a_{t-1}^{-e_{t-1}}\cdots
a_{1}^{-e_{1}} \in gp\{S\}$. Hence the proof is immediate from the
theorem. $\blacksquare$ \bos (i) Let $G$ be a group with identity
$e$, then $G$ and $\{e \}$ are subgroups of G.These are called \bf
trivial subgroups \rm of G.\\ (ii) Any subgroup of a subgroup is a
subgroup. \eos \bt\label{t29} Let $H, K$ be subgroups of a group
$G$.  Then $HK$ is a subgroup of $G$ if and only if $HK = KH$. \et
\pf $(\Rightarrow)$ As $HK$ is a subgroup of $G$, we have \ba{rl}
                  &  (HK)^{-1} = HK\\
\Rightarrow       & K^{-1}H^{-1} = HK\\ \Rightarrow       & KH =
HK \quad (\mbox{ since } K^{-1} = K, H^{-1} = H) \ea
$(\Leftarrow)$ Let $hk, h_1k_1 \in HK$, where $h, h_1\in H$ and
$k, k_1\in K.$ Then, \bea
       (hk)(h_1k_1)^{-1} & = & hk k_1^{-1}h_1^{-1}\\
                 &  = & hk_2h_1^{-1} \mbox { where } k_2 = kk_1^{-1} \in K
\eea Note that $k_2h_1^{-1} \in KH = HK.$  Hence we can write
$k_2h_1^{-1} = h_3k_3$ for some $h_3 \in H$, $k_3\in K$.
Therefore, $$ (hk)(h_1k_1)^{-1} = hh_3k_3\in HK \; (hh_3\in H)$$
Consequently $HK$ is a subgroup of $G.$  $\blacksquare$ \\We give
below an example of a group $G$ and subgroups $H,K$ of $G$ for
which $HK$ is not a subgroup.
 \bx\label{ex30} Let $G = S_3$. Take
$$I = \left(\begin{array}{ccc}
  1 & 2 & 3 \\
  1 & 2 & 3
\end{array}\right),  a  =  \left(\begin{array}{ccc}
                            1 & 2 & 3\\  1  & 3 & 2
                            \end{array}\right)
               b  = \left(\begin{array}{ccc}
                            1 & 2 & 3\\  3  & 2 & 1
                            \end{array}\right) $$
 in $S_3$. We have $a^{2} = b^{2} = c^{2} =I$. Hence $H = \{I,a\}$ and $K = \{I,b\}
 $ are subgroups of $G$ and \bea
        HK & = & \{I,a,b,ab\}\\
\eea Note that $$ab = \left(\begin{array}{ccc}
  1 & 2 & 3 \\
  2 & 3 & 1
\end{array}\right) \mbox{ and } ba = \left(\begin{array}{ccc}
  1 & 2 & 3 \\
  3 & 1 & 2
\end{array}\right) $$ Thus $HK$ is not a subgroup of $G$ since $(ab)^{-1} = ba$ and $ba$
 is not in $HK$. \ex

 \bx\label{ex21} In the example \ref{ex1.6}
of groups $$( \, \C,+) \supset (\, \R,+) \supset (\, \Q,+) \supset (\,
\Z,+)$$ and $$(\, \C^*,\cdot) \supset ( \,\R^*,.) \supset (\, \Q^*,.)$$
are chains of groups, where each member is a subgroup of the preceding
member. Further, $\R^{*}_{+} = \{\alpha \in \R^{*} \st \alpha > 0\}$ and
$\Q^{*}_{+} = \{\beta \in \Q^{*} \st \beta > 0\}$ are subgroups of
$(\R^{*}, \cdot)$ and $(\Q^{*}, \cdot)$ respectively. \ex
\bx The set $S^1 = \{z \in \C \st \abs{z} =1\}$ is a subgroup of $(\C^{*},
\cdot)$.\ex
\bx\label{ex22} The groups in the examples \ref{ex1.12}, \ref{ex1.13} and
\ref{ex1.13b} i.e., $\C_n = \{z\in \C \st z^n = 1 \}$, $\C_{\infty}$ and
$\C_{m^{\infty}}$ are subgroups of $(\, \C^{*}, \cdot).$ Further
$\C_{m^{\infty}}$ is a subgroup of $\C_{\infty}$.
\ex
\bx\label{ex23} In the example \ref{ex1.14}, $Sl_n(\R)$ is a subgroup
of $Gl_n(\R)$ and $Sl_n(\Z)$ is a subgroup of $Sl_n(\R).$ The
following is a diagram of subgroups of $Gl_n(\R).$\bea
\begin{array}{ccccc}
  Gl_{n}(\R) & \supset & Sl_{n}(\R) &  \\
  \cup &  & \cup &  &  \\
 Gl_{n}(\Q)  & \supset & Sl_{n}(\Q) & \supset & Sl_{n}(\Z)
\end{array}
\eea \ex
\bx Let $K$ be a field. The sets $$T_n(K) =\{A = (a_{ij})\;:\;  n\times
n-\mbox{matrix over }K \st a_{ij}=0 \mbox{ for all } i>j, \mbox{ and } det
\, A \neq 0 \},$$ and $$UT_n (K)=\{A \in T_n(A)\st a_{ii} =1 \mbox{ for
all } 1 \leq i\leq n \}$$ are subgroups of $Gl_n (K)$. In fact $UT_n (K)$
is a subgroup of $Sl_n(K)$. The group $T_n(K)$ is called the {\bf group of
upper triangular matrices} over $K$ and $UT_n(K)$ is called the {\bf group
of upper uni-triangular matrices} over $K$. Similarly we can define the
{\bf group of lower triangular matrices} and {\bf lower uni-triangular
matrices}.\ex
\bx\label{ex24} Let $G = I(X)$, the group of isometries of the metric
space $X$ (Example \ref{ex8}). For any subset $F \subset X$, the
set$$I(X,F)=  \{ T \in G \st T(f) = f \mbox{for all} f \in F \} =
  H$$ is a subgroup of G. Not all subgroups of $I(X)$ are
obtained in this way e.g., if $X = \R$ then for any non-empty
subset $F\neq \{0\}$ of $\R$, $I(X,F) = Id.$ \ex \jump
\bx\label{ex25} Let $X = \R^n$ be the n-dimensional real space
with usual  metric i.e., $d(x,y) = \norm{x-y}$ for all $x,y \in
\R^n$ where for any $x = (x_{1},x_{2},\cdots,x_{n}) \in \R^{n}$
and $\|x\|^{2} = \sum_{i=1}^{n}{\abs{x_{i}}}^{2}$. The set of all
$\R$-linear maps from $\R^n$ to $\R^n$ which are isometries of the
metric space $\R^n$ is a subgroup of $I(\R^n).$\\
 Let
 $$T : \R^n \to \R^n$$
be an $\R$-linear map.If $T$ is an isometry, then for any $x \in \R^n$ we
have \ba{rl}
            & \norm{T(x)} = \norm{x}
\ea Now,let $x,y\in \R^n $ and  $\langle x,y \rangle$ denotes the usual
inner product/dot product of $x,y$.Then \ba{rl}
             & \norm{x-y}= \norm{T(x-y)} = \norm{T(x)-T(y)}
\\ \Rightarrow & {\norm{T(x)-T(y)}}^2 = {\norm{x-y}}^2\\
\Rightarrow  & \langle T(x) - T(y),T(x) - T(y)\rangle = \langle x - y,x -
y\rangle \\ \Rightarrow & \langle T(x),T(y)\rangle + \langle
T(y),T(x)\rangle
=
\langle x,y\rangle + \langle y,x\rangle \\ \Rightarrow & 2\langle
T(x),T(y)\rangle = 2\langle x,y\rangle \\ \Rightarrow & \langle
T(x),T(y)\rangle = \langle x,y\rangle \ea Conversely, if for any
$x, y \in \R^n,$ $$ \langle T(x),T(y)\rangle
=
\langle x,y\rangle$$ Then \ba{rl}
            & \langle T(x-y),T(x-y)\rangle = \langle x-y,x-y\rangle\\
\Rightarrow & {\norm{T(x)-T(y)}}^2 = {\norm{x-y}}^2\\ \Rightarrow
& {\norm{T(x)-T(y)}} = {\norm{x-y}} \ea Therefore, $T$ is an
isometry iff $\langle T(x),T(y)\rangle = \langle x,y\rangle$ for
all $x,y \in \R^n$. Now, consider the canonical basis $\{ e_i=
(0,...0, 1^{ith}, 0,...0) \st 1 \leq i \leq n \}$ of $\R^n$. Let
$T(e_i) = \sum_{k = 1}^n a_{ki}e_k$. Then, \ba{rl}
            & \langle T(e_i),T(e_j)\rangle = \langle e_i,e_j\rangle\\
\Rightarrow & \langle \displaystyle\sum_{k=1}^n a_{ki}e_k ,
\sum_{k=1}^n a_{kj}e_k\rangle = \langle e_i,e_j\rangle \\
\Rightarrow & \displaystyle\sum_{k=1}^n a_{ki}a_{kj} = \left\{
{\begin{array}{lr}
                                                    0  \mbox{ if } i \neq j\\
                                                    1 \mbox{ if } i=j
                                                    \end{array}}
                                                    \right.
\ea Hence for the matrix  $A = [a_{ij}]$,  $A A^t = I$. Conversely,
given an $n\times n$- matrix A over reals with $A A^t = I$, consider
the $\R$-linear transformation $T$ of $\R^{n}$ determined  by $A$
i.e., $T(x) = A(x_{1},x_{2},\cdots,x_{n})^{t}$ for any $x =
(x_{1},\cdots,x_{n}) \in \R^{n}$ where $(x_{1},\cdots,x_{n})^{t}$ is
the transpose of row-vector $(x_{1},\cdots,x_{n})$. Then $T$
satisfies $\langle T(e_i) , T(e_j)\rangle = \langle e_i ,
e_j\rangle$ reversing the above argument. Hence by linearity
$$\langle T(x),T(y)\rangle = \langle x,y\rangle$$ for all $x, y \in
\R^{n}.$ Thus the subgroup of $\R$-linear maps from $\R^{n}$ to
$\R^{n}$ which are isometries is $\{ A \in Gl_n(\R) \st A A^t = Id
\}.$ This is denoted by $O(n,\R)$, and is called the \bf orthogonal
group\index{orthogonal group} \rm of degree $n$ over $\R$. The set
$$ SO(n,\R) = \{ A \in O(n,\R) \st \abs{A}= 1 \}$$ is a subgroup of
$O(n,\R)$ and is called the \bf special orthogonal
group\index{special orthogonal group} \rm of degree $n$ over $\R$.
\ex \bexe Let $n\geq 1$, and $K$ be a field.Prove that \\ (i) ({\bf
Orthogonal group over $K$})The set $$ O(n,K)=\{A: n\times
n-\mbox{matrix over } K \st AA^t= Id\}$$ is a group, and $\det A=\pm
1$ for any $A\in O(n,K).$
\\ (ii)({\bf Special orthogonal group over $K$}) The set$$SO(n,K)=
\{ A\in O(n,K) \st \abs{A}= 1\}$$is a subgroup of $O(n,K)$.\eexe
\br\label{ex26} The group of rotations of $\R^2$ (example
\ref{ex9}) i.e., $$ G = \left\{\left(\begin{array}{cc} \cos \theta
& -\sin \theta\\ \sin \theta & \cos \theta \end{array} \right) \st
0\leq \theta \leq 2\pi \right\}$$ is a subgroup of $Gl_2(\R)$.
Note that any element of $G$ is an orthogonal matrix, hence this
is a subgroup of $O(2,\R)$, in fact, of $SO(2,\R).$ \er  \bexe Let
$H$ be the set of all matrices in $O(n.\R)$ with integral entries.
Prove that
\\(i) $H = \{(a_{ij}) = A \st A \mbox{ is } n \times n$
 matrix over $\Z$ with $\abs{a_{ij}} = 1 \mbox{ or }
 0$ for all $1 \leq i,j \leq n$ and each row/column
 of $A$ has exactly one non-zero entry \}.  \\
 (ii) $H$ is a subgroup of $O(n,\R)$ with $\circ(H) =
2^{n}\cdot n!$.  \eexe \bexe Let $A$ be an $n\times n$-matrix over $\R$
and let $A^{(i)}$,$1\leq i \leq n$, denotes the $i^{th}$ column of $A$.
Then $A \in O(n,\R)$ if and only if $\langle A^{(i)},A^{(j)}\rangle =
\delta_{ij}$ for all $1\leq i,j\leq n$ where $\delta_{ij}=0$ if $i \neq j$
and is 1 if $i=j$.\eexe
 {\bf{Note}} The above statement is true for rows of $A$ as well.
  \bexe In the Example \ref{ex25}, take $X = \C^{n}$, the
$n$-dimensional complex space with the metric $d(z, w) = \|z - w \|$ where
for any $z = (z_{1},\cdots , z_{n}) \in \C^{n}$, $\|z\|^{2} =
\displaystyle{\sum_{i=1}^{n} \abs{z_{i}}^{2}}$. Prove that the set of all
$\C$- linear maps from $\C^{n}$ to $\C^{n}$ which are isometries of the
metric space $\C^{n}$ is a subgroup (say) $\bf U\rm( n, \C)= {\bf U_n}$ of
$I( \, \C^{n})$. Further, show that $$\begin{array}{ccl}
  {\bf U_n} & = & \{A \in Gl_{n}(\, \C) | A^{*} = A^{-1} \} \\
   & = & \{A \in Gl_{n}(\, \C) | \|Az - Aw\| = \|z - w\| \mbox{ for
   all }z, w \in \C^{n}\}
\end{array}$$ where $A^{*} = {\overline{A}}^{\;t}$, the conjugate transpose of
$A$.  \\  (Hint: Given $z = ( z_{1},\cdots,z_{n})$ and
$w=(w_{1},\cdots,w_{n})$ in $\C^{n}$, we have $\langle z,w \rangle = \sum
z_{i}\overline{w_{i}}.$ Use that $$2\langle z,w \rangle
=
({\norm{z+w}}^{2}-{\norm{z}}^{2}-{\norm{w}}^{2})+i({\norm{z+iw}}^{2}-{\norm{z}}^{2}-
{\norm{w}}^{2}))$$  \eexe
 \bexe Prove that $${\bf SU_n}=\bf SU\rm(n, \C) = \{ A \in \bf
U\rm(n, \C) | det \, A = 1\}$$ is a subgroup of ${\bf U}(n,\C)={\bf
U_n} $.\eexe \bf Note : \rm The group $\bf U_n$ is called the {\bf
unitary group}\index{unitary group} of dimension $n$ over $\C$, and
$\bf SU_n$ is called the \bf special unitary group\index{special
unitary group}.\rm

Most of the groups encountered above are infinite.  However $\C_n$, the
group of $n^{th}$ roots of unity in $\C$, has degree $n$. Thus for each $n
\in N$, there is a group of order $n$. Further $S_n$, the symmetric group
of order $n$, is also finite and has order $n!$ (use: The number of
permutations on $n$ symbols is $n!$, and that any one-one map from
$\{1,2,\cdots,n\}$ to itself is a permutation of this set, and
conversely). \\ We shall, now, prove a result which enlightens the
geometric content in the definition of $O(n,\R)$. It is proved that any
isometry of $\R^{n}$ which fixes origin is  a linear transformation of
$\R^{n}$ determined by an element of $O(n,\R)$.  \bt\label{t1.54t} Let $f$
be any isometry of the metric space $\R^{n}$ (Example \ref{ex25}) which
fixes the origin. Then there exists an element $A \in O(n,\R)$ such that
for any $x = (x_{1},\cdots,x_{n}) \in \R^{n}$ $$f(x) =
A(x_{1},\cdots,x_{n})^{t}$$ where $(x_{1},\cdots,x_{n})^{t}$ is the
transpose of the $1 \times n$-matrix $(x_{1},\cdots,x_{n})$.\et \pf Since
$f$ is an isometry of $\R^{n}$ which fixes origin, \be\label{a.a}
\norm{f(x)}\; = \; \norm{x} \mbox{ for all }x \in \R^{n}. \ee Further, for
any $x,y \in \R^{n}$, $$\begin{array}{crcl}
   & \norm{f(x)-f(y)} & = & \norm{x-y} \\
  \Rightarrow & {\norm{f(x)-f(y)}}^{2} & = & {\norm{x-y}}^{2} \\
  \Rightarrow & {\norm{f(x)}}^{2}+{\norm{f(y)}}^{2}-2\langle x,y\rangle & = & {\norm{x}}^{2}+
  {\norm{y}}^{2}-2\langle x,y\rangle
\end{array}$$ Hence, using \ref{a.a}, we get \be\label{b.b} \langle
f(x),f(y) \rangle \; = \; \langle x,y \rangle \mbox{ for all } x,y \in
\R^{n}. \ee Let $e_{i} = (0,\cdots,0,1^{i^{th}},0,\cdots,0)$;
$i=1,2,\cdots,n$ be the canonical basis of $\R^{n}$. Then, from the
equation \ref{b.b}, we get \be\label{c.c} \langle f(e_{i}),f(e_{i})
\rangle \; = \; \langle e_{i},e_{i} \rangle = 1 \ee for all $1 \leq i \leq
n$ and if $1 \leq i \neq j \leq n$, then \be\label{d.d} \langle
f(e_{i}),f(e_{j}) \rangle \;= \; \langle e_{i},e_{j} \rangle = 0 \ee Put
$f(e_{i}) = a_{i}$. Then for the $n \times n$-matrix $$A =
[a_{1}^{t},\cdots,a_{n}^{t}]$$ with $i^{th}$ column $a_{i}^{t}$, the
transpose of the vector $a_{i} \in \R^{n}$, $A^{t}A =Id$. Thus $A \in
O(n,\R)$. \\ In the Example \ref{ex25} we have seen that if $B \in
O(n,\R)$ and $B$ denotes the $\R$-linear transformation from $\R^{n}$ to
$\R^{n}$ defined by $B(x) = Bx^{t}$ for any $x \in \R^{n}$, then $B$ is an
isometry. As $B$ fixes the origin, as seen above, $$\langle B(x), B(y)
\rangle \; =\; \langle x,y \rangle \mbox{ for all } x,y \in \R^{n}$$
$$\mbox{ i.e.,   }\langle Bx^{t},By^{t} \rangle \; = \; \langle x,y
\rangle $$ Now, let $A^{t}$ denotes the isometry determined by $A^{t}$,
then for any $e_{i}$, $i = 1,\cdots,n$ we have \be\label{e.e}
(A^{t}f)(e_{i}) = A^{t}a_{i} = e_{i}\ee Further, as $A^{t}f$, the
composite of the isometries $A^{t}$ and $f$ fixes origin, as seen above,
we have for all $1\leq i\leq n,$  $$\begin{array}{crclc}
  & \langle A^{t}f(x),A^{t}f(e_{i}) \rangle & = & \langle x,e_{i} \rangle = x_{i} &  \\
  \Rightarrow & \langle A^{t}f(x),e_{i} \rangle & = & x_{i} &  \mbox{use \ref{e.e}}\\
  \Rightarrow & A^{t}f(x) & = & x &
  \end{array}$$
\hspace{1in} $ \Rightarrow \;\;\; f  =  (A^{t})^{-1}=A,$
 the isometry determined by $A$.  \\ Hence the result.
$\blacksquare$
\brs\label{rea.as} (i) Let us note that not all isometries of $\R^{n}$ fix
origin,e.g., if $a\in \R^n$,$a\neq (0,\cdots,0)$, then the transformation,
$\tau_{a}(x) = x+a$ defines an isometry of $\R^{n}$ which does not fix
origin. \\ (ii) If $f$ is any isometry of the metric space $\R^{n}$ and
$f(0) = a$ where $a \neq (0,\cdots,0)$, then $(\tau_{-a} f)(0) = 0$.
Hence by the theorem $\tau_{-a} f = A$ for some $A \in O(n,\R)$. Thus
$f = \tau_{a} A$ i.e., any isometry of $\R^{n}$ is composite of a
translation and an orthogonal transformation.\ers In view of the Remark
\ref{rea.as}(ii), to understand geometrically $I(\R^{2})$, we need to know
the geometric meaning of the isometries of $\R^{2}$ determined by elements
of $O(2,\R)$ (see Example \ref{ex25}). For this purpose we prove:
\bl\label{la,al} Let $A \in O(2,\R)$. Then the isometry ( Example
\ref{ex25}), $$\begin{array}{rcl}
  A : \R^{2} & \rightarrow & \R^{2} \\
  x = \left(\begin{array}{c}
    x_{1} \\
    x_{2} \
  \end{array}\right) & \mapsto & A\left(\begin{array}{c}
    x_{1} \\
    x_{2} \
  \end{array}\right)
\end{array}$$ (Here an element of $\R^{2}$ is represented by a column
vector) is either a rotation by an angle $\theta$ or is a reflection with
respect to the line passing through the origin making an angle $\theta/2$
with the $+ve$ $X$-axis.\el  \pf Let $$A = \left(\begin{array}{cc}
  a & b \\
  c & d
\end{array}\right)$$ Then $$\begin{array}{crcl}
  & Id = AA^{t} & = & \left(\begin{array}{cc}
    a & b \\
    c & d \
  \end{array}\right)\left(\begin{array}{cc}
    a & c \\
    b & d \
  \end{array}\right) \\
  \Rightarrow & a^{2}+b^{2}=1 & = & c^{2}+d^{2} \\
   & ac+bd & = & 0
\end{array}$$Thus $(a,b)\; , \; (c,d)$ are points on the unit circle in
$\R^{2}$. As any vector on the unit circle in $\R^{2}$ is
$(\cos\theta,\sin\theta)^{t}$ where $\theta$ is the angle of the vector
with the $+ve$ $X$-axis. Therefore  $$A = \left(\begin{array}{cc}
  \cos\theta & \sin\theta \\
  \cos\phi & \sin\phi
\end{array}\right)$$ where $\cos\theta\cos\phi + \sin\theta\sin\phi = 0$.
Hence $$\begin{array}{crlc}
   & \cos(\phi-\theta) & = 0 &  \\
  \Rightarrow & \phi & = \theta+n\pi/2 & n : \mbox{ odd integer }
\end{array}$$ Hence, either  $$A = \left(\begin{array}{cc}
  \cos\theta & \sin\theta \\
  -\sin\theta & \cos\theta
\end{array}\right) \mbox{ or }\left(\begin{array}{cc}
  \cos\theta & \sin\theta \\
  \sin\theta & -\cos\theta
\end{array}\right)$$ Now, as the linear transformation of $\R^{2}$
determined by $$\left(\begin{array}{cc}
  \cos\theta & \sin\theta \\
  -\sin\theta & \cos\theta
\end{array}\right)$$ represents anti-clockwise rotation by an angle
$\theta$ around the origin and the linear transformation determined by
$$\left(\begin{array}{cc}
  \cos\theta & \sin\theta \\
  \sin\theta & -\cos\theta
\end{array}\right)$$ represents reflection with respect to the line
through the origin making an angle $\theta/2$ with the $+ve$ $X$-axis. The
result follows.     $\blacksquare$  \bl\label{laa.al} Let $n>1$ be an
integer. \\ (i) If $A \in O(n,\R)$ then for any eigen value $\lambda$ of
$A$, $\abs{\lambda} = 1$.  \\  (ii) If $n$ is odd and $A \in SO(n,\R)$,
then 1 is an eigen-value of $A$.  \\  (iii) Let $A \in O(n,\R)$, $A
\not\in SO(n,\R)$. Then if $n$ is odd $-1$ is an eigen value of $A$,
moreover, if $n$ is even then both $\pm1$ are eigen values of $A$.\el  \pf
(i) Let $\lambda$ be an eigen value of $A$, then there exists a non-zero
vector $x \in \R^{n}$ such that $Ax = \lambda x$. Therefore
$$\begin{array}{crcl}
   & \langle x,x \rangle & = & \langle Ax,Ax \rangle \\
   &  & = & \langle \lambda x,\lambda x\rangle \\
  \Rightarrow & {\abs{\lambda}}^{2} & = & 1 \mbox{ since } \langle x,x \rangle \neq 0 \\
  \Rightarrow & \abs{\lambda} & = & 1
\end{array}$$ Therefore (i) is proved.  \\  (ii) As $A \in SO(n,\R)$,
$det\;A = 1.$  Hence,$\det A^t=\det A=1$. Therefore, we get
$$\begin{array}{crclc}
   & det((A-Id)A^{t}) & = & det(Id-A^{t}) &  \\
   &  & = & det(Id-A)^{t} &  \\
  \Rightarrow & det(A-Id) & = & det(Id-A) &  \\
   &  & = & (-1)^{n}det(A-Id) &  \\
  \Rightarrow & det(A-Id) & = & 0 & \mbox{ since }n \mbox{ is odd }
\end{array}$$ Therefore 1 is an eigen value of $A$. \\ (iii) As in (ii),
we have $$\begin{array}{crcl}
   & det((A+Id)A^{t}) & = & det(Id+A^{t}) \\
   &   &  =  & det(Id+A)^{t} \\
  \Rightarrow & det(A+Id)det\;A & = & det(Id+A) \\
  \Rightarrow & det(A+Id) & = & 0 \mbox{   since }det\;A \neq 1
\end{array}$$ Hence -1 is an eigen-value of $A$. Now, let $n$ be even.
Then as above  $$\begin{array}{rcl}
  det(A-Id)det\;A & = & det(Id-A) \\
   & = & (-1)^{n}det(A-Id) \\
   & = & det(A-Id) \\
  \Rightarrow \;\;\;\;\; det(A-Id) & = & 0 \mbox{\;\;\;\; since }det\;A \neq 1
\end{array}$$ Therefore 1 is also an eigen-value of $A$ in case $n$ is
even.     $\blacksquare$
 \clearpage
 {\bf EXERCISES}
\begin{enumerate}
  \item
  Prove that if $H$ is a subgroup of a group $G$,
then identity element of $H$ is same as the identity element of
$G$.
 \item
 Let $G$ be a non-empty set together with an associative
binary operation  $$\begin{array}{crcl}
  f : & G \times G & \rightarrow & G \\
   & (x, y) & \mapsto & x \ast y
\end{array}$$ show that if  \\  (i) There exist an element $e \in G$
 such that $e \ast x = x$
for all $x \in G$.  \\  (ii) For all $x \in G$ there exist $w \in
G$ such that $w \ast x = e$ ($w$ depends on $x$ ), \\  then $G$ is
a group.
 \item
 Let $x_{1},x_{2},\cdots,x_{n}\;(n\geq1)$, $x$ and $y$ be elements
 in a group $G$. Prove that  \\
 (i)$x^{-1}(\prod_{i=1}^{n}x_i)x$ $=$
 $\prod_{i = 1}^{n}(x^{-1}x_{i}x)$ \\ (ii)
 $(x_{1}\cdots x_{n})^{-1}\; =
 \;x_{n}^{-1}x_{n-1}^{-1}\cdots x_{1}^{-1}$  \\ (iii)
 $(x^{-1}yx)^{k} \;= \;x^{-1}y^{k}x$ for all $k \in \Z$.
 \item
 Let $\Z = (\Z, +)$ be the additive group of integers. Then
show that  \\  (i) For any $m \in \Z$, the set $\Z m = \{am \mid a
\in \Z\}$ is a subgroup of $\Z$.  \\  (ii) For any two non-zero
integers $m,n$, $\Z m \subset \Z n$ if and only if $n$ divides $m$.
\\  (iii) For any subgroup $H \neq (0)$ of $\Z$ there exists a subgroup $K
\neq (0)$, $K \subsetneqq H$.
\\  (iv) There exists a subgroup $A$ of $\Z$ such that if for any subgroup $B$ of $\Z$,
$A \subsetneqq B \subset \Z$, then $B = \Z$.
 \item
 Prove that in any finite
group $G$, the set $\{x \in G \mid x^{2} \neq e\}$ has even number
of elements. Further, show that if order of $G$ is even, then
there exists an element $x (\neq e )$ in $G$ such that $x^{2} =
e$.
 \item
 Prove that for any group $G$,  \\  (i) $H_{t} = \{
{a_{1}}^{t}{a_{2}}^{t}\cdots {a_{k}}^{t} \mid a_{i} \in G, k \geq
1 \}$ is a subgroup of $G$ for all $t \geq 1$. Further, if $G$ is
abelian then $H_t = \{x^t | x \in G\}$.\\
 (ii) $H_{2} = \{
x_{1}x_{2}\cdots x_{n}.x_{n}x_{n-1}\cdots x_{1} \mid x_{i} \in G,
n \geq 1\}.$
 \item
 Let $H, K, L$ be subgroups of a group $G$ such that
$H \subset K$. Prove that: \\ (i) $HL \cap K = H( L \cap K).$
 \\ (ii) If $HL=KL$ and $H\cap L=K\cap L$, then $H=K$.
 \item
Let $H, K$ be two subgroup of a group $G$. Prove that $H \bigcup
K$ is a subgroup of $G$ if and only if either $H \subset K$ or $K
\subset H$.
\item Let $H,K$ be two subgroups of a group $G$. Prove that for
any two elements $x,y \in G$,either $HxK = HyK$ or $HxK \bigcap
HyK = \emptyset$(null set).
 \item
 Let $G$ be a finite group with $\circ(G) =
n$ and let $S$ be a subset of $G$ with $\abs{S} \geq (n/2) +1$
where $(n/2)$ is the largest integer less than equal to $n/2$.
Show that $G = SS$. Further, show that if $\abs{S} = (n/2)$ then
$G$ need not be equal to $SS$ and that $G$ contains no proper
subgroup of order $> (n/2)$.
\item
Calculate the order of $Gl_{2}(\Z_{3})$ by writing down all its
elements.
 \item
\label{ec26} Let $\F$ be a finite field with $q$ elements. Prove that the
order of the group $Gl_n(\F)$ is $\prod_{k=0}^{n-1}(q^n - q^k)$ and the
order of $Sl_n(\F)$ is $q^{n-1}\prod_{k=0}^{n-2}(q^n - q^k).$\\
 (Hint: Use that an $n
\times n$ matrix over a field $K$ has determinant $\neq 0$ if and
only if all its columns(rows) are linearly independent vectors
over $K$ in the vector space $K^n$.)
 \item
 Let $R$ be a commutative ring with identity. Prove that $\R$ can be
replaced by $R$ in the Example \ref{es1}. Further show that $Gl_n(R)=\{A:
n\times n-\mbox{matrix over} R \,|\,detA \mbox{ is unit.}\}$ is a group.

 \item
The set $$G = \{f : \R \rightarrow \R \mid f(x) =
\frac{ax+b}{cx+d}, a, b, c, d \in \R \mbox{ with } ad - bc = 1
\}$$ forms a group under the binary operation of composition of
maps.
\item
Let $x, y$ be two elements of a group $G$. Prove that if $x^{m} =
y^{m}$, $x^{n} = y^{n}$ for $m, n \in \Z$ with $(m, n) = 1$ i.e.,
g.c.d. of $m$ and $n$ is $1$, then $x = y$.
 \item
\label{.004x} Let $G$ be a group. Prove that if for two relatively
prime integers $m$ and $n$, $x^{m}y^{m} = y^{m}x^{m}$ and
$x^{n}y^{n} = y^{n}x^{n}$,for all $x,y \in G$, then $G$ is
abelian.\\ (Hint
: Show that for $k = (m,n), (x^{m}y^{n})^{k} =
(y^{n}x^{m})^{k}$ using $(xy)^{t} = (x(yx)x^{-1})^{t} =
x(yx)^{t}x^{-1}$ for all integers $t$.)
\item
Prove that two groups of same order need not have same number of
subgroups.
\item
Prove that $S_{3}$ has exactly six subgroups.
\item
Prove that $SO(2,\R)$ is abelian, but $SO(3,\R)$ is not abelian.
\end{enumerate}

 \chapter{ Generators Of A Group And  Cyclic Groups }
 \footnote{\it contents group2.tex }
 Let $G$ be a group
and $S(\neq\emptyset)$,a subset of $G$.Then  $$ gp\{S\} = \{
a_1a_2...a_m \st m\geq 1, a_i \mbox{ or } a_i^{-1} \in S \mbox{
for all } 1 \leq i \leq m\}.$$  is a subgroup of $G$ ( Corollary
\ref{co111}).This is called the subgroup of $G$ generated by $S$,
and  $S$ is called a set of generators  of $gp\{S\}$. If $S =
\emptyset$, the null set, then we define $gp\{S\} = e$, the
identity subgroup of $G$. \bd\label{d123} If for a finite subset
$S$ of a group $G$, $G = gp\{S\}$, then $G$ is called a finitely
generated group. Further, if no proper subset of $S$ generates
$G$, then $S$ is called a minimal set of generators of $G$.\ed \br
Two minimal sets of generators of a group need not have same
number of elements e.g., for the additive group of integers $\Z =
(\Z, +)$, we have $\Z = gp\{1\} = gp\{2, 3\}$ and neither $2$ nor
$3$ generates $\Z$ since for any $m \in \Z$, $gp\{m\} = m\Z$. Thus
$\{1\}$ and $\{2,3\}$ are two minimal sets of generators for
$\Z$.\er
\bexe Let $S$ be a finite subset of a group $G$. Prove that
$gp\{S\}$ is either finite or is countably infinite.\eexe \bexe Show that
$(\R, +)$ and $(\R^{*}, .)$ are not finitely generated groups. \eexe
\bt\label{th2.06th} Prove that every subgroup of a group $G$ is finitely
generated if and only if there is no infinite ascending chain of subgroups
of $G$ i.e., if $A_1\subset A_2\subset \cdots \subset A_n\subset \cdots$
is an infinite ascending chain of subgroups of $G$, then there exists
$m\geq 1$ such that $A_{m+k}=A_m$ for all $k\geq 1.$\et \pf Assume that
every subgroup of $G$ is finitely generated and let $$A_1\subset
A_2\subset \cdots \subset A_n\subset \cdots$$ is an infinite chain of
subgroups of $G$. \\ Put $$A=\bigcup_{n=1}^{\infty}A_n$$ Let $x,y \in A.$
Then since $A_n,\:n\geq 1,$ is ascending chain of subgroups, there exists
$l\geq 1$ such that $x,y\in A_l$. Therefore $xy^{-1}\in A_l\subset A.$
Hence A is a subgroup of $G$ (Theorem \ref{t1}). By assumption $A$ is
finitely generated. Let $A=gp\{x_1,\cdots,x_t\}.$ Then there exists $m\geq
1$ such that $x_i\in A_m$ for all $i=1,2,\cdots,t.$ Hence $A=A_m$, and
consequently $A_m=A_{m+k}$ for all $k\geq 1.$ Conversely, let $G$ has no
infinite ascending chain of subgroups. If the converse is not true, then
there exists a subgroup $H$ of $G$ which is not finitely generated. We
shall, now, construct inductively an infinite ascending chain of subgroups
of $G$. Let $A_1=gp\{x_1\}$ for any $x_1\in H.$ As $H$ is not finitely
generated $H\neq A_1$. Choose $x_2\in H-A_1$, and let $A_2=gp\{x_1,x_2\}.$
As above $H\neq A_2$. Choose $x_3\in H-A_2$ and put
$A_3=gp\{x_1,x_2,x_3\}$. Let we have chosen $x_1,x_2,\cdots x_n$ in $H$
such that for $A_t=gp\{x_1,x_2,\cdots,x_t\},\:1\leq t\leq n,\: x_{t+1}\in
H-A_t.$ Then as $H$ is not finitely generated $H\neq A_n.$ Choose
$x_{n+1}\in H-A_n,$ and put $A_{n+1}=gp\{x_1,x_2,\cdots,x_{n+1}\}.$ Thus,
by induction, we get the infinite ascending chain $$A_1\subsetneqq
A_2\subsetneqq \cdots \subsetneqq A_n\subsetneqq \cdots$$ of subgroups of
$G$. This contradicts our assumption on $G$. Hence every subgroup $G$ is
finitely generated.     $\blacksquare$ \bd\label{d2.06d} An element $x$ in
a group $G$ is called a Non-Generator of $G$ if whenever $G=gp\{x,S\}$ for
a subset $S$ of $G$, $G=gp\{S\}.$\ed \br For any group $G$, $e\in G$ is
clearly a non-generator of $G$.\er \bl\label{l2.07l} Let $G$ be a group.
Then \\ (i)The set $\Phi(G)$ of non-generators of $G$ is a subgroup of a
$G$.  \\  (ii) Let $M$ be a maximal subgroup of $G$ that is $M$ is a
subgroup of $G$, $M\neq G$, contained in no subgroup of $G$ other than $G$
and itself. Then $\Phi(G) \subset M.$\el \pf (i) As $e\in G$ is
non-generator, $\Phi(G)$ is non-empty. Let $x,y \in \Phi(G)$, and let for
a subset of $G$, $$\begin{array}{crlc}
   & G & =gp\{xy^{-1},S\} &  \\
   &  & =gp\{x,y,S\} &  \\
   &  & =gp\{S\} & \mbox{ since $x,y$ are non-generators of } G \\
  \Rightarrow & xy^{-1} & \in \Phi(G) &
\end{array}$$ Therefore $\Phi(G)$ is a subgroup of $G$ (Theorem \ref{t1}).
\\  (ii) Let $x \in \Phi(G)$. If $x\not\in M,$ then by maximality of $M$,
$G=gp\{x,M\}.$ However $G\neq M$. This contradicts the fact that $x\in
\Phi(G)$. Hence $\Phi(G) \subset M$.     $\blacksquare$ \bd\label{d2.08d}
Let $G$ be a group. The subgroup $\Phi(G)$ of set of non-generators of $G$
is called Frattini subgroup of $G$.\ed \bt\label{t2.09t} For a group $G$,
the Frattini subgroup $\Phi(G)$ of $G$ is the intersection of all maximal
subgroups of $G$.\et \pf By lemma \ref{l2.07l} (ii), for any maximal
subgroup $M$ of $G$, $\Phi(G) \subset M.$ Next,let $x$ be any element of
$G$ which lie in all maximal subgroups of $G$. If $x\not\in \Phi(G)$, then
there exists a subset $S$ of $G$ such that $G=gp\{x,S\},$ but $G\neq
gp\{S\}.$ Let $B=gp\{S\}.$ Then $B\neq G$ is a subgroup of $G$ such that
$G=gp\{x,B\}.$ Clearly $x\not\in B$. Put $$\mathcal{A}=\{H\subset G \st H:
\mbox{ subgroup, }H \supset B, x\not\in H\}$$ Clearly $\mathcal{A}\neq
\emptyset$ as $B\in \mathcal{A}.$ For any two $H_1,H_2\in \mathcal{A}$,
define: $$H_1\leq H_2 \Leftrightarrow H_1 \subset H_2.$$ Clearly
$(\mathcal{A},\leq)$ is a partially ordered set. Let $\{H_i\}_{i\in I}$ be
any chain in $\mathcal{A}$. Put $$D=\bigcup_{i\in I}H_i$$ If $x,y\in D$,
then as $\{H_i\}_{i\in I}$ is a chain in $\mathcal{A},$ there exists
$i_0\in I$ such that $x,y\in H_{i_0}.$ Therefore $xy^{-1} \in H_{i_0}
\subset D$. Hence $D$ is a subgroup of $G$, $D\supset B$ and $x\not\in D$.
Thus $D$ is an upper bound of the chain $\{H_{i}\}_{i\in I}$ in
$\mathcal{A}.$ Thus any chain in $\mathcal{A}$ has an upper bound in
$\mathcal{A}.$ Hence by Zorn's lemma $\mathcal{A}$ has a maximal element
(say) $N$. Then $x\not\in N$ and $gp\{x,N\}\supset gp\{x,B\}=G.$ Thus
$gp\{x,N\}=G$. We claim that $N$ is a maximal subgroup of $G$. If not then
there exists a subgroup $N_1$ of $G$, $N_1\supsetneqq N, N_1\neq G.$ By
maximality of $N$ in $\mathcal{A}$, $x\in N_1.$ Hence $N_1=G.$ This is a
contradiction to the fact that $N_1\neq G$. Thus $N$ is a maximal subgroup
of $G$ not containing $x$. Consequently $\Phi(G)$ is the intersection of
all maximal subgroups in $G$.        $\blacksquare$
\bd\label{d23} For any element $x$ of a group $G$, order of $x$ is
defined to be $o(gp\{x\})$ and is denoted by $o(x)$.\ed \bos (i) In a
group $G$, the only element of order 1 is $e$. Further for any $x \in G$,
$gp\{x\} = gp\{x^{-1}\}$. Hence $\circ (x) = \circ (x^{-1})$.  \\ (ii)
Every element $(\neq 0)$ in $(\R,+)$ has order $\infty$.  \\ (iii) The
only non-identity element of finite order in $\R^*$ is -1 and it has order
2. Further, an element in \, $\C^*$ is of finite order if and only if it
is a root of unity. Consider the group (Example \ref{ex1.12}) $$\C_n = \{
z \in \, \C \st z^n = 1 \}$$ of $n^{th}$ roots of unity in $\C$. Every
element of $\C_n$ is a power of $e^{\frac{2\pi i}{n}}.$ Thus $\C_{n}$ has
exactly $n$-elements i.e., $\circ( \C_{n} ) = n$, moreover, $\C^{*}$ has
infinite number of elements of finite order. \eos \bd\label{2.1dd} An
element $x$ of a group $G$ is called a torsion element if
$\circ(x)<\infty,$ otherwise, it is called torsion free element.\ed
\bd\label{2.2dd} let $G$ be a group. Then \\ (i) If every element of $G$
is a torsion element, $G$ is called a torsion / periodic group. \\ (ii) If
every non-identity element of $G$ has infinite order, $G$ is called
torsion free.\ed \br The identity group is torsion as well as torsion
free.\er \bd\label{2.3dd} If for a group $G$, the set $\{\circ(x)\st x\in
G\}$ is bounded, then the least common multiple of the integers $\circ(x),
x\in G$, is called exponent of $G$.\ed We shall, now, find generator sets
for $Gl_n(K)$ and $Sl_n(K)$ ($K$: field). To do this, we need to know:
\bx\label{101ex} Let $K$ be a field and let  $n \geq 2$ be a fixed integer. Let for all
$1\leq i \neq j \leq n$, $E_{ij}(\lambda)$ denotes the $n \times n$-matrix
over $K$ with $(i,j)^{th}$ entry $\lambda$ and all other entries $0$. Put
$e_{ij}(\lambda) = Id + E_{ij}(\lambda)$ where Id denotes the $n \times
n$-identity matrix. Then $e_{ij}(\lambda)$ is the matrix with all entries
on the diagonal equal to $1$ and all off diagonal entries  $0$ except the
$(i,j)^{th}$ entry which is equal to $\lambda$. Clearly, $det
\,e_{ij}(\lambda) = 1$. The subgroup $$E_{n}(K) = gp\{e_{ij}(\lambda) | 1
\leq i \neq j \leq n, \lambda \in K \}$$ of $Sl_{n}(K)$ is called the {\bf
elementary group} of dimension n over $K$. The matrices $e_{ij}(\lambda)$
are called {\bf elementary matrices}.\ex
\bt\label{t109t} For any field $K$, we have  \\  (i) $Gl_{n}(K) = gp\{e_{ij}(\alpha
),d(\alpha ) | \alpha (\neq 0) \in K , 1 \leq i \neq j \leq n\}$
\\ where $d(\alpha )$ denotes the diagonal matrix with $(n,n)^{th}$ entry equal to
 $\alpha $ and all
other diagonal entries equal to 1. \\  (ii)$Sl_{n}(K) = gp\{e_{ij}(\alpha
) | \alpha (\neq 0) \in K , 1 \leq i \neq j \leq n\}$\et  \pf  (i) Let
$M_{(i)}(M^{(i)})$ denotes the $i^{th}$ row (column ) of a matrix $M$. If
$A \in Gl_{n}(K)$ and $ B = e_{ij}(\alpha )A$ then $B_{(t)} = A_{(t)}$ for
all $1 \leq t ( \neq i) \leq n$ and $B_{(i)} = A_{(i)} + \alpha A_{(j)}$.
Further for $C = Ae_{ij}(\alpha ), C^{(t)} = A^{(t)}$ for all $1 \leq t (
\neq j) \leq n$ and $C^{(j)} = A^{(j)} + \alpha A^{(i)}$. Thus
multiplication by $e_{ij}(\alpha )$ on left (right ) to $A$ amounts to
performing an elementary row ( column ) operation on $A$. As $det \, A
\neq 0, A_{(1)} \neq 0$ i.e. , there exists a non zero entry in the first
row of $A$. Hence by performing elementary column operations on $A$ we can
transform $A$ to a matrix $B$ with first row $(1,0,0,\cdots ,0)$. Now , by
elementary row operations we can transform $B$ to a matrix $C$ such that
  $$C = \left(\begin{array}{cc}
    1 & 0\cdots 0 \\
    \vdots & A_{1} \\
    0 &
  \end{array}\right)$$
 Thus there exist $E_{1}, E_{2} \in E_{n}(K )$ such that $E_{1}AE_{2} = C$.
 Note that  $$\begin{array}{rcl}
   det A_{1} & = & det C \\
    & = & det E_{1}AE_{2} \\
    & = & det A \
 \end{array}$$  Hence $A_{1} \in Gl_{n-1}(K)$. As seen above, by performing elementary
 row and column operations, we can bring $A_{1}$ to the form $$\left(\begin{array}{cc}
    1 & 0\cdots 0 \\
    \vdots & A_{2} \\
    0 &
  \end{array}\right)$$  where $A_{2} \in Gl_{n-2}(K)$. Note that we can perform elementary
  row (column) operations on $A_{1}$ by doing the same on $C$ without changing the
  first
  column and row  of $C$. Continuing as above we obtain $E, F \in E_{n}(K )$ such that
  $EAF = d( \alpha )$ for some $\alpha ( \neq 0) \in K$. Hence $A = E^{-1}d( \alpha )F^{-1}$
  , and the result follows. \\ (ii)  Let $A \in Sl_{n}(K)$ in (i).
  Then, as $A = E^{-1}d(\alpha)F^{-1}$, we get
  $$\begin{array}{rcl}
    detA & = & det(E^{-1}d( \alpha )F^{-1}) \\
     & = & det d(\alpha ) \\
   \Rightarrow \hspace{.5in} \alpha & = & 1
  \end{array}$$  Hence $d(\alpha ) = Id$, and the proof follows.
  $\blacksquare$

   \bd\label{d24} A group $G$ is called cyclic if there exists an
element $g \in G$ such that $G = gp\{g\}$ i.e., every element of
$G$ is a power of $g$. \ed \br (i) For any element $x$ in a group
$G$, $x^{m}x^{n}=x^{m+n}=x^{n}x^{m}$ for all $m,n\in\Z$. Hence a
cyclic group is clearly abelian. \\  (ii) The group $(\Z,+)$ is
cyclic and clearly $\Z = gp\{1\} = gp\{-1\}.$ Further, $\Z_{n}$ is
also a cyclic group and $\Z_{n} = gp\{1\}.$
\\ (iii)Let $G$ be a finite group of order $m$, then $G$ is cyclic
if and only if there exists an element of order $m$ in $G$. \er
 {\bf Notation}: We shall denote by $C_n$, a cyclic group of order $n$. \bexe
 Prove that the group $\C_{n}$ is cyclic for any $n\geq 1.$\eexe  \bexe
 Let $p$ be a prime. Prove that for any $x,y \in \C_{p^{\infty}},$ either
 $gp\{x\}\subset gp\{y\}$ or $gp\{y\}\supset gp\{x\}.$ Then deduce that
 any finitely generated subgroup of $\C_{p^{\infty}}$ is $C_{p^n}$ for
 some $n\geq 0.$\eexe
\bexe\label{20}Prove that the groups $(\Q,+)$ and $(\Q^{*},\cdot)$
are not cyclic.\eexe \bt\label{t24a} A subgroup of a cyclic group
is cyclic. Further, a non-identity subgroup of an  infinite cyclic
group is infinite cyclic. \et \pf Let $C = gp\{g\}$ be a cyclic
group, and let $H$ be a subgroup of $C$. If $H = \{e\}$, then it
is clearly cyclic and is generated by $e$. Let $H \neq \{e\}$.
Then there exists a +ve integer $k$ such that $g^k \in H$. Choose
$m$ least positive integer such that $g^m \in H.$ We claim $H =
gp\{g^m\}.$ Let for $t > 0, g^t \in H.$ We can write $$ t = m q +
r$$ where $0 \leq r
< m.$ This gives $$  g^r = g^{t-mq} = g^t \cdot g^{-mq} \in H.$$
By minimality of $m$, we conclude that $r = 0$ i.e., $t = mq.$
Hence $g^t = g^{mq} \in gp\{g^m\}.$ Now, note that $g^k \in H$ if
and only if $g^{-k} = (g^k)^{-1} \in H.$ Hence we get  $H =
gp\{g^m\}.$ This proves $H$ is cyclic.  To prove our next
assertion, let $C$ be an infinite cyclic group. As $H\neq \{e\}$,
there exists $m > 0$, such that $H = gp\{g^m\}.$ If $H$ is finite
then there exist integers $s > t$, such that \ba{rl}
                    & (g^m)^s = (g^m)^t\\
               i.e.,&  g^{ms} = g^{mt}\\
    \Rightarrow     &  g^{m(s-t)} = e
\ea Choose $k$ least positive integer such that $g^k = e.$ Then $C
= \{1,g,\ldots ,g^{k-1} \}$ i.e., $C$ is finite. This, however, is
not true. Hence, $H$ is infinite cyclic. $\blacksquare$ \brs (i)
If a group $G$ has elements of finite order $(\neq 1)$ and also of
infinite order then $G$ is not cyclic( use: Any non- identity
subgroup of an infinite cyclic group is infinite cyclic).\\ (ii)
If every subgroup of a group is cyclic, then obviously as a group
is a subgroup of itself, it is cyclic. However, if we assume that
every proper subgroup of a group is cyclic, then the group need
not be cyclic. To show this we give the following: \ers \bx Let
$$V_{4} = \left\{ I = \left(\begin{array}{cc}
  1 & 0 \\
  0 & 1
\end{array}\right), a = \left(\begin{array}{cc}
  0 & 1 \\
  1 & 0
\end{array}\right), b = \left(\begin{array}{cc}
  0 & -1 \\
  -1 & 0
\end{array}\right), c = \left(\begin{array}{cc}
  -1 & 0 \\
  0 & -1
\end{array}\right) \right\}$$ be the Klein's $4$-group. Then
$H_{1} = \{ I, a \}, H_{2} = \{ I, b \}$, and $H_{3} = \{ I, c \}$
are all the proper subgroups of $G$. Hence every proper subgroup
of $G$ is cyclic, but, $V_{4}$ is not cyclic as it has no element
of order $4$.\ex  \bl\label{888} Let $\Z = ( \Z, + )$ be the
additive group of integers. Then for any finite set $\{ m_{1},
m_{2}, \cdots, m_{k} \}$ of non-zero integers, we have
$$\begin{array}{rcl}
  gp\{ m_{1}, m_{2}, \cdots, m_{k}\} & = & \Z m_{1} + \cdots + \Z m_{k} \\
   & = & \Z d
\end{array}$$  where $d$ is the greatest common divisor (g.c.d.) of $m_{i}'s$, $1 \leq i \leq k$.
\el \pf  As $\Z$ is abelian, clearly $$gp\{ m_{1}, m_{2}, \cdots,  m_{k}\} = \Z
m_{1} + \cdots + \Z m_{k}$$ Further, as $\Z$ is a cyclic group,by
the theorem \ref{t24a}, we get, $$\Z m_{1} + \cdots + \Z m_{k} =
\Z d$$ for some $d \in \Z$. Thus for all $1 \leq i \leq k$, $m_{i}
\in \Z d$ and hence $d \mid m_{i}$. Now, let for an integer $e$,
$e | m_{i}$ for all $1 \leq i \leq k$. Then $m_{i} \in \Z e$ and
hence $$\begin{array}{rcl}
  \Z d & = & \Z m_{1} + \cdots + \Z m_{k} \subset \Z e \\
  \Rightarrow &  & e | d
\end{array}$$
  Hence the result follows. $\blacksquare$
\bco\label{.002co} Let $C = gp\{a\}$ be a cyclic group generated
by $a$. Then, given any finite set $\{ m_{1}, m_{2}, \cdots,
m_{k}\}$  of non-zero integers, $$gp\{ a^{m_{1}}, a^{m_{2}},
\cdots, a^{m_{k}} \} = gp\{a^{d}\}$$ where $d$ is the greatest
common divisor of $m_{i}'s, \, 1 \leq i \leq k$. \eco \pf The
proof follows easily using that $d = a_{1}m_{1} + a_{2}m_{2} +
\cdots + a_{k}m_{k}$ for some $a_{i} \in \Z, \, 1 \leq i \leq k$.
\bt\label{.003t} The special linear group of degree $n$ over $\Z$
is generated by elementary matrices i.e., $$Sl_{n}(\Z) =
gp\{e_{ij}(t) \mid \, 1 \leq i \neq j \leq n, t \in \Z \}$$ \et
\pf The result is trivial for $n = 1$. Hence assume $n\geq 2$. Let
$A_{1} \in Sl_{n}(\Z)$. Since determinant $A_{1}$ is $1$, not all
entries in any given row (column) of $A_{1}$ are zero, moreover,
entries of any row (column) of $A_{1}$ have g.c.d. $1$. Now, note
that we can divide any entry in a given row (column) of $A_{1}$ by
any other entry in the same row (column) by performing elementary
column (row) operations on $A_{1}$. Thus as g.c.d. of the entries
in the first row of $A_{1}$ is $1$, as proved in the theorem
\ref{t109t} (i), there exist an elementary $n \times n$- matrix
$F_{1}$ over $\Z$ such that  $$C = A_{1}F_{1} = \left(
\begin{array}{ccc}
  1 & 0 \cdots & 0 \\
  \ast &  &  \\
  \vdots & A_{2} &  \\
  \ast &  &
\end{array}\right)$$ Now, by subtracting multiples of the first row in $C$ from
other rows we can make all entries (second onwards) in the first
column zero. Thus there exist an elementary $n \times n$-matrix
$E_{1}$ over $\Z$ such that $$B = E_{1}A_{1}F_{1} = \left(
\begin{array}{ccc}
  1 & 0 \cdots & 0 \\
  0 &  &  \\
  \vdots & A_{2} &  \\
  0 &  &
\end{array}\right)$$ Now,
$$\begin{array}{cl}
   & det \, (E_{1}A_{1}F_{1})  =   det \, A_{2}   \\
  \Rightarrow & det\;E_{1}\;\; det \, A_{1}\;\; det\;F_{1}  =  det \, A_{2}  \\
  \Rightarrow &  det \, A_{2} = 1,
\end{array}$$  Since the determinant of any elementary matrix is
$1$. Thus $A_{2} \in E_{n-1}(\Z).$ Now, note that we can perform
elementary row and column operations on $A_{2}$ by doing the same
on $B$ and while doing so keeping first row and column of $B$
unchanged. Hence,continuing as a above, we can find $E, F \in
E_{n}(\Z)$ such that $EA_{1}F = Id$,  the identity matrix of order
$n$. This gives that $A_{1} = E^{-1}F^{-1} \in E_{n}(\Z)$. Hence
the result follows. $\blacksquare$
\bco\label{co2.9} The general linear group $Gl_{n}(\Z)$ of degree
$n$ over $\Z$ is generated by the set $\{e_{ij}(t),\;d(-1) \st
1\leq i\neq j < n,\; t \in \Z\}$ where $d(-1)$ is the $n \times
n$-diagonal matrix over $\Z$ with $(i,i)^{th}$ entry $1$ for all
$1 \leq i \leq n-1$ and $(n,n)^{th}$ entry $-1$.\eco \pf If $A \in
Gl_{n}(\Z),$ then  \\ \hspace*{1.25in}$1 = det\; (AA^{-1}) =
det\;A.det\;A^{-1}$
\\ \hspace*{1in} $\Rightarrow \; \; \; det\; A = \pm 1.$  \\ If $det \; A = 1$,
then $A \in Sl_{n}(\Z).$ Let $det\; A = -1$. Then as in the
theorem there exist $E, F \in E_{n}(\Z)$ such that $$EAF = \left(
\begin{array}{ccccc}
  1 & 0 & \cdots & 0 & 0 \\
  0 & 1 & \cdots & 0 & 0 \\
  \vdots &  &  &  & \vdots \\
  0 & 0 & \cdots & 0 & d
\end{array}\right)$$ $$\Rightarrow \; \; det\; EAF = d$$
$$\Rightarrow\; \; det\; A = d$$ $$i.e. \; \; \; d = -1$$
 Now, the result is immediate from the theorem.   $\blacksquare$
  \bco\label{.11co}  The group $Sl_{2}(\Z)$ is generated by the
matrices   $$A = \left(\begin{array}{cc}
  0 & 1 \\
  -1 & 0
\end{array}\right) \; , \; B = \left(\begin{array}{cc}
  1 & 1 \\
  0 & 1
\end{array}\right) $$\eco \pf First of all, by the theorem, $E_{2}(\Z) = Sl_{2}(\Z)$.
Therefore, it suffices to prove that $E_{2}(\Z) = gp\{A, B\}$. We
have for all $m \geq 1$, \be\label{56} B^{m} =
\left(\begin{array}{cc}
  1 & m \\
  0 & 1
\end{array}\right) \mbox{  and  } B^{-m} = \left(\begin{array}{cc}
  1 & -m \\
  0 & 1
\end{array}\right)\ee Further, $$\begin{array}{ccccl}
  C & = & A^{-1}B^{-1}A & = & \left(\begin{array}{cc}
  0 & -1 \\
  1 & 0
\end{array}\right)\left(\begin{array}{cc}
  1 & -1 \\
  0 & 1
  \end{array}\right)\left(\begin{array}{cc}
    0 & 1 \\
    -1 & 0 \
  \end{array}\right)  \\
   &  &  & = & \left(\begin{array}{cc}
     1 & 0 \\
     1 & 1 \\
   \end{array}\right)
\end{array}$$ Hence for all $m \geq 1$, \be\label{44e} C^{m} = \left(\begin{array}{cc}
  1 & 0 \\
  m & 1
\end{array}\right) \mbox{  and  } C^{-m} = \left(\begin{array}{cc}
  1 & 0 \\
  -m & 1
\end{array}\right)\ee By the equations \ref{56} and \ref{44e}, it is clear that $E_{2}(\Z)
 = gp\{B, C\} \subset gp\{A, B\}$. Finally, as  $$\begin{array}{ccccl}
   A & = & \left(\begin{array}{cc}
   0 & 1 \\
   -1 & 0 \
 \end{array}\right) & = & \left(\begin{array}{cc}
   1 & 1 \\
   0 & 1 \
 \end{array}\right)\left(\begin{array}{cc}
   1 & 0 \\
   -1 & 1 \
 \end{array}\right)\left(\begin{array}{cc}
   1 & 1 \\
   0 & 1 \
 \end{array}\right)   \\
    &  &  & = & BC^{-1}B \
 \end{array}$$ It is clear that $gp\{A, B\} = gp\{B, C\} = Sl_{2}(\Z)$. $\blacksquare$
 \bt\label{t25} Let $x
(\neq e)$ be an element of order $n$ in a group $G$. Then  \\
 (i) The integer $n$ is least +ve such that  $x^n = e.$  \\
(ii) If $x^m = e$ for an integer $m$, then $n$ divides $m$. \et
\pf (i) By definition  $\circ (x) = \circ(gp.\{x\}) = n$. Thus
$e, x, x^2,\ldots ,x^n$ ( $(n+1)$- elements) are not all distinct.
Hence there exist $0 \leq k < l\leq n$ such that \ba{rl}
                   &  x^k = x^l \\
\Rightarrow        &  x^{l-k} = e \\
       i.e.,       &  x^s  = e  \quad \mbox{ for } s = l-k \geq 1
\ea Choose $t$ least +ve such that $x^t = e.$ For any $k \in \Z$,
we can write $k = tq+r$ where $o\leq r\leq t$. Hence  $$\begin{array}{ccccc}
  x^{k} & = & x^{tq}x^{r} &  &  \\
   & = & x^{r} &  & \mbox{since } x^{t} = e
\end{array}$$ Therefore
$$gp\{x\} = \{e, x, x^2,\ldots ,x^{t-1}\}$$ If for $0 \leq k
< l
< t, x^l = x^k,$ then $x^{l-k} = e$, where $1 \leq l-k < t.$ This
contradicts the minimality of $t$. Therefore $e, x, x^2,\ldots
,x^{t-1}$ are all distinct. We, however, have $\circ (gp\{x\}) =
n$. Hence $n = t$ i.e., $n$ is the least +ve integer such that
$x^n = e.$\\ (ii) If $m = 0$, there is nothing to prove. Hence,
assume $m \neq 0$. It is sufficient to assume that $m
> 0$
 since if $x^{m} = e$, then $x^{-m} = e$. We can write
\ba{rl}
                    &  m = nq+r ,\quad q, r \in \Z, 0 \leq r < n \\
\Rightarrow         &  e = x^m = x^{nq} \cdot x^r\\ \Rightarrow &
e = x^r \quad (\mbox { since } x^n = e).\hspace*{.25in} by (i) \ea
This, however, gives by minimality of $n$ that $r = 0$. Thus $m =
nq$ i.e., $n$ divides $m.$ \bt\label{t27} Let $x(\neq e)$ be an
element of order $n$ in a group G.  Then for any integer $k > 0$,
$\circ (x^k) = \frac{n}{(n,k)}$ where $(n,k)$ denotes the g.c.d.
of $n$ and $k$. Thus in particular $\circ(x^{k})$ divides $n$.
\et \pf Since $x^n = e$, we have $${(x^k)}^{\frac{n}{(n,k)}} =
{(x^n)}^{\frac{k}{(n,k)}} = e$$ Further, let \ba{rl}
                     {(x^k)}^t = e\\
\Rightarrow        & n\mid kt \mbox{     Theorem \ref{t25}(ii)}
\\ \Rightarrow & \frac{n}{(n,k)} \mbox{ divides }
\frac{k}{(n,k)}.t \ea As $\frac{n}{(n,k)}$ and $\frac{k}{(n,k)}$
are relatively prime, we get $\frac{n}{(n,k)}$ divides $t$. Hence
$\frac{n}{(n,k)}$ is the least +ve integer such that
${(x^k)}^{\frac{n}{(n,k)}} = e.$ Therefore $\circ (x^{k}) =
\frac{n}{(n,k)}$, and the rest is clear. $\blacksquare$
\bco\label{l1.23} Let $x (\neq e)$ be an element of order $n$ in a
group $G$. Then for any $k \geq 1$, $\circ(x^{k}) = n$ if and only
if $(k, n) = 1$.\eco  \pf The proof is immediate from the theorem.
$\blacksquare$
\bt\label{t26} Let $G$ be a group and let $x,y$ be two elements in
$G$. If $\circ (x) = m,$ $\circ (y) = n, xy = yx,$ and $(m,n) =
1$, then $\circ (xy) = mn.$ \et \pf We have $$ (xy) ^{mn} =
(x^m)^n\cdot (y^n)^m  = e  $$ Further, let for some $ t \geq 1$
\ba{rl}
                (xy)^t = e\\  \Rightarrow   &  x^{t}y^{t} = e
                \mbox{    since }xy = yx  \\
\Rightarrow    & x^t = (y^t)^{-1}\\ \Rightarrow & \circ (x^t) =
\circ (y^t) \ea   $$\begin{array}{crclccc}
  \Rightarrow & \frac{m}{(m, t)} & = & \frac{n}{(n, t)} &  &  & \mbox{(use : Theorem \ref{t27})}  \\
  \Rightarrow & m(n, t) & = & n(m, t) &  &  &  \\
  \Rightarrow & m\mid (m, t), &  & n \mid (n, t) &  &  & \mbox{since } (m, n) = 1 \\
  \Rightarrow & m \mid t, &  & n \mid t &  &  &  \\
  \Rightarrow & mn \mid t &  &  &  &  &
\end{array}$$ Hence $mn$ is the least $+ ve$ integer such that $(xy)^{mn} =
e$.Therefore by the theorem \ref{t25}(i), $\circ(xy) = mn$.
$\blacksquare$
\bd\label{d1.11} For any integer $n \geq 1$, $\bf{\Phi}(n)$ denotes the
number of elements in the set $\{k \in \N \; \mid \; 1 \leq k \leq n, \;
(k, n) = 1\}$. The function ${\bf{\Phi}} : \N \rightarrow \N$ defined here
is called Euler's function.\ed \bl\label{l33b} Let $C = gp\{x\}$ be a
cyclic group of order $n$. Then $C$ has $\bf{\Phi}(n)$ distinct
generators.\el \pf As $C$ has order $n$, $\circ(x) = n$ and $$C = \{ e, x,
\cdots, x^{n - 1} \}$$ An element $x^{k}$, $1 \leq k \leq n$ is a
generator of $C$ if and only if $\circ(x^{k}) = n$. Thus $x^{k}$ is
generator of $C$ if and only if $(k, n) = 1$ (Corollary \ref{t26}). Hence
$C$ has $\bf{\Phi}(n)$ generators. $\blacksquare$
\bt\label{t28} Let $C = gp\{x\}$ be a finite cyclic group of order
$n.$ Then for any $d \geq 1, d\mid n,$ $C$ contains a unique
subgroup of order $d$. \et \pf Let $\frac{n}{d} = n_1$. Then from
the theorem \ref{t27} $\circ (x^{n_1}) = d.$ Thus $gp\{x^{n_1}\}$
has order $d$. Further, let $H$ be a subgroup of $C$ with  $\circ
(H) = d.$ As $H$ is cyclic (Theorem \ref{t24a}) $H = gp\{x^k\}$
for some $0 \leq k < n.$  We have \ba{rl}
                & \circ (gp\{x^k\}) = d\\
\Rightarrow  &  \circ(x^{k}) = d   \\    \Rightarrow & x^{kd} =
e\\ \Rightarrow & n \mid kd\\ \Rightarrow & \frac {n}{d} \mid k\\
\Rightarrow & k = \frac {n}{d}k_1 \quad \mbox { for some }k_{1}
\in \Z.\\ \Rightarrow & x^k \in gp\{x^{n_1}\} \quad \mbox { where
} n_1 = \frac {n}{d}\\ \Rightarrow     & gp\{x^k\}  \subset
gp\{x^{n_1}\} \ea As both groups have order $d$, we conclude  $H =
gp\{x^k\} = gp\{x^{n_1}\}$. Thus a subgroup of order $d$ in the
cyclic group $C$  is unique.  $\blacksquare$
\\  We shall, now prove a number theoretic result using group
theory
:
   \bl\label{98l} Let $n \geq 1$ be any integer. Then  $$n = \sum_{d|n} \bf{\Phi}(d)$$\el
   \pf Consider the cyclic group $\Z_{n}$. We know that $\circ(\Z_{n}) = n$.
   By the theorem \ref{t28}, for each $d|n, \Z_{n}$ contains a
   unique subgroup $G_{d}$ of order $d$, which of course is cyclic (Theorem
   \ref{t24a}). Let $S_{d}$ denotes the set of generators of
   $G_{d}$. By the lemma \ref{l33b}, $G_{d}$, has exactly
   $\bf{\Phi}(d)$ generators. Thus, $\abs{S_{d}} = \bf{\Phi}(d)$. Since every
   element of $\Z_{n}$ generates exactly one of the $G_{d}$, the
   theorem \ref{t27} gives $$\begin{array}{ccccc}
     n & = & \circ(\Z_{n}) & = & \sum_{d|n} \abs{S_{d}} \\
      &  &  & = & \sum_{d|n} \bf{\Phi}(d) \
   \end{array}$$ Hence the result follows. $\blacksquare$
   \bco\label{44co} Let $p$ be a prime and $n \geq 1$. Then
   $$\bf{\Phi}(p^{n}) = (p - 1)p^{n - 1}$$\eco \pf We shall prove
   the result by induction on $n$. Clearly $\bf{\Phi}(p) = p - 1$.
   Hence the result is true for $n = 1$. By the lemma, we have
   $$p^{n} = \sum_{0 \leq k \leq n} \bf{\Phi}(p^{k})$$ Hence, by the
   induction hypothesis $$\begin{array}{ccclcl}
      & p^{n} & = & 1 + (p - 1) + (p - 1)p + \cdots + (p - 1)p^{n - 2} + \bf{\Phi}(p^{n}) &  &   \\
     \Rightarrow & \bf{\Phi}(p^{n}) & = & (p^{n} - 1) - (p - 1)(1 + p + \cdots + p^{n - 2}) &  &  \\
      &  & = & (p^{n} - 1) - (p - 1)\frac{(p^{n - 1} - 1)}{(p -
      1)}\\
      & & = & p^{n - 1}(p - 1)
   \end{array}$$ Thus the result is proved. $\blacksquare$
\bl\label{t211} If every element $(\neq e)$ of a group $G$ has
order 2, then $G$ is abelian.  \, \el \pf By assumption, for any
$x \in G, x^2 = e.$ Hence $x^{-1} = x.$ Now, let $x, y \in G$.
Then \ba{rl}
             &           (xy)^2 = e\\
\Rightarrow  &  xy = (xy)^{-1} =  y^{-1}x^{-1} = y x \ea Hence $G$
is Abelian. \clearpage {\bf EXERCISES}
\begin{enumerate}
  \item
 \label{co211} Show that a group of order $4$ is
abelian.
  \item
 Let $G$ be a finite group and $p$, a prime. Prove that if $ G = \{
x^{p} \; \mid \; x \in G \}$, then $G$ contains no element of
order $p$.
  \item
 Let $G$ be a group of even order. Prove
that $G \neq \{ x^{2} \; \mid \; x \in G \}$.
\item
 Let $G = \{ x \in \Q
\; \mid \; 0 \leq x \leq 1 \}$,and let for any $x, y \in G$,
$$\begin{array}{rcl}
  x + y & = & \left\{\begin{array}{lcl}
    x + y & \mbox{if} & 0 \leq x + y < 1   \\
    x + y - 1 & \mbox{if} & x + y \geq 1 \
  \end{array}\right.
\end{array}$$  Prove that $G$ is an abelian group with
respect to the above binary operation and all element of $G$ have
finite order.
 \item
Prove that an infinite cyclic group $C = gp\{x\}$ has exactly two
generators $x$ and $x^{-1}$.
  \item  Prove that a group $G$ is finite if
and only if $G$ has finite number of subgroups.
  \item
   Prove that the
subgroup of $Gl_{2}(\R)$ generated by the subset $S = \left\{
\left(\begin{array}{cc}
  0 & 1 \\
  1 & 0
\end{array}\right)  \right.,$  \\ $ \left. \left(\begin{array}{cc}
  1 & 0 \\
  0 & -1
\end{array}\right)\right\}$ has order $8$.
  \item
   Let $A
(\neq 0)$ be an $n \times n$-matrix over $\R$ such that $A^{m} =
0$ for some $m \geq 2$. Prove that $I + A \in Gl_{n}(\R)$ ( $I$ :
$n \times n$ identity matrix ). Further, show that if $m = 2$,
then $\circ(I + A) = \infty$.
  \item
   Prove that the elements $ a =
\left(\begin{array}{cc}
  0 & -1 \\
  1 & 0
\end{array}\right)$ and $b = \left(\begin{array}{cc}
  0 & 1 \\
  -1 & -1
\end{array}\right)$ in $Gl_{2}(\R)$ have finite order, but,
$\circ(ab) = \infty$.
  \item
   Show that the groups $(\Q,+)$ and $(\Q^{*},
\cdot)$ are not finitely generated.
\item
Prove that any finitely generated subgroup of $(\Q, +)$ is cyclic.
  \item
   Let $C =
gp\{x\}$ be a finite cyclic group of order $n \geq 3$. Prove that
$C$ has an even numbers of generators.
\item
Find the number of subgroups of a cyclic group of order $30$ i.e.,
of $C_{30}.$
  \item  Show that for any
element $x$ of a group $G$, $\circ(x) = \circ(g^{-1}xg)$ for all
$g \in G$. Then deduce that $\circ(xy) = \circ(yx)$ for all $x, y
\in G$.
  \item
   Let $p$ be a prime and $G = \Z_{p}^{*}$, the multiplicative group of
integers modulo $p$. Prove that $\prod_{x \in G} x = p-1$. Then
deduce that $(p-1)! + 1 \equiv 0 \mod (p)$.
  \item
   Let $C$
be a finite cyclic group and $n| \circ(C)$. Prove that the set $\{
x \in C \; | \; x^{n} = e \}$ has exactly $n$-elements.
  \item
   Let $m
\geq 1$, $n \geq 1$ be two relatively prime integers. Show that if
for an element $x \in G$, $\circ(x) = m$, then $x = y^{n}$ for
some $y \in G$. ( Hint : use $am + bn = 1$ for some $a, b \in
\Z$.)
\item
Show that if a group $G$ is union of its proper subgroups, then
$G$ is not cyclic.
\item Let $G$ be a finitely generated group. Prove that there is no infinite
chain $A_1\subset A_2\subset \cdots \subset A_n\subset \cdots$ of proper
subgroups $A_n,\:n\geq 1,$ of $G$ such that $G=\bigcup_{i=1}^{\infty}A_n.$
\item
Let $C=gp\{a\}$ be an infinite cyclic group. Then \\ (i) Prove that any
maximal subgroup of $C$ is of the type $gp\{a^p\},$ where $p$ is prime.
\\ (ii) Show that $\Phi(C)=id.$
\item
Show that $\Phi(S_3)=id$.
\item
Find the Frattini subgroup of a cyclic group of order 12.
\item
Let $p$ be a prime and $\alpha\geq 1.$ Prove that if $C_{p^{\alpha}}$
cyclic group of order $p^{\alpha}$, then $\Phi(C_{p^{\alpha}})$ is cyclic
group of order $p^{\alpha -1}$.
\end{enumerate}

 \chapter{Cosets Of A Subgroup, Normal Subgroups And Quotient Groups }
\footnote{\it contents group3.tex } In this chapter we start with
the concept of cosets of a subgroup which in some sense, measure
the size of a group with respect to a given subgroup. Then we
define a normal subgroup and associate a group with every normal
subgroup of a group called the factor group or quotient group of
the group with respect to the given normal subgroup. A preliminary
study of these ideas is done here.
 \bd\label{d25} Let $G$ be a
group and $H$, a subgroup of $G$. For any $x \in G$, the set $xH =
\{ xh \st h \in H\}$ is called the left coset of $H$ in $G$ with
respect to $x$, and $Hx = \{ hx \st h \in H \}$ is called the
right coset of $H$ in $G$ with respect to $x$. \ed
\brs\label{brs2.1} (i) For any $h \in H$, $hH = H = Hh$   \\ (ii)
The map $h \mapsto xh$ is one-one from $H$ onto $xH.$ Hence
$\abs{H} = \abs{xH}$ \ers
\bl\label{l22} Let $G$ be a group and $H$, a subgroup of $G$.
Then any two distinct left (right) cosets of $H$ in $G$ are
disjoint. \el \pf Let $xH, yH$ be two distinct left cosets of $H$
in $G$. Let \ba{rl}
                 &  xH \cap yH \neq \phi\\
\Rightarrow      &  xh_1 = yh_2 \;\; \mbox{ for some }h_1, h_2 \in
H \\ \Rightarrow & xh_1 h_2^{-1} = y\\ \Rightarrow      &  xh_1
h_2^{-1}H = yH\\ \Rightarrow      &  xH = yH \mbox{ since }h_1
h_2^{-1} \in H \ea A contradiction to our assumption that the two
cosets are distinct. Hence $xH \cap yH = \phi$. The proof for
right cosets is similar. $\blacksquare$   \bt\label{t23} Let $G$
be a group and $H$, a subgroup of $G$. Then there exists a 1-1
correspondence between the left cosets of $H$ in $G$ and the right
cosets of $H$ in $G$. Thus  the cardinality of the set of all
distinct left cosets of $H$ in $G$ is the same as the cardinality
of the set of all distinct right cosets of $H$ in $G.$ \et \pf
Consider the correspondence $$\theta :xH \to Hx^{-1}$$ If for two
left cosets $xH, yH$  of $H$ in $G$, we have \ba{rl}
                   & Hx^{-1} = Hy^{-1}\\
\Rightarrow        &(Hx^{-1})^{-1} = (Hy^{-1})^{-1}\\ \Rightarrow
& xH^{-1} = yH^{-1}\\
 i.e.,             & xH = yH \quad(\mbox{since }H^{-1} = H)
\ea Hence the correspondence is 1-1. Further, for any right coset
$Hy$ of $H$ in $G$, we have $\theta (y^{-1}H) = Hy.$ Thus the
result holds.  $\blacksquare$  \bd\label{d26} Let $H$ be a
subgroup of a group $G$. Then $[G : H]$ denotes the cardinality of
the set of all distinct left(right) cosets of $H$ in $G$.  This is
called the index of $H$ in $G$. \ed

\br  Let $H$ be a subgroup of
a group $G$, and $ G = \bigcup_{i\in I} (x_iH)$.
 Then for any $x \in H, x = x_ih$ for some
$i \in I$ and $h \in H.$ Therefore $xH = x_ihH = x_iH$. Hence if $x_iH
\cap x_jH = \emptyset$ for all $ i \neq j$, then $[G : H] = \abs I$.  \er
 \bexe Let $n\geq 1$ be any integer. Prove that for any subgroup $n\Z$ of
 the group $\Z = (\Z,+)$, $[\Z:n\Z] = n.$ \eexe
\bt\label{co22}\bf(Lagrange's)  \rm  Let $G$ be a finite group and $H$, a
subgroup of $G$. Then $\circ (H)$ divides $\circ (G),$ and $\circ(G) =
\circ(H)[G:H]$ i.e., $[G:H] = \frac{\circ(G)}{\circ(H)}$. \et \pf Let
$x_1H, x_2H,...,x_nH$ be all distinct left cosets of $H$ in $G.$ Thus, as
$x \in xH$, by the lemma \ref{l22}, $$ G = x_1H \cup x_2H \cup ... \cup
x_nH$$ is a decomposition of G into pair-wise disjoint subsets. Therefore
\ba{rl}
              &  \abs G = n \abs H \mbox{   use: Remark \ref{brs2.1}(ii)}  \\
\Rightarrow   &  \circ(H) \mbox{ divides }\circ(G) \\
      \mbox{and }  & \circ(G) = \circ(H)[G:H] \hspace{1in}\blacksquare  \ea
 \bco\label{bco2.2} Any group of prime order is cyclic.\eco
 \pf Let $G$ be a group of order $p$,a prime. Take any $x(\neq e) \in G$.
 Then, by the theorem, for the subgroup $H = gp\{x\}$, $\circ(H)$
 divides $\circ(G) = p$. As $x \neq e$,$\circ(H) \neq 1$. Hence $\circ(H) = p$. Thus
 $G = gp\{x\}$ i.e., $G$ is cyclic.  $\blacksquare$
 \bco\label{co3.9} Let $G$ be a finite group of order $n$. Then for any
 $x \in G$, order of $x$ divides $n$.\eco \pf By the definition
 \ref{d23}, $\circ(x) = \circ(gp\{x\}).$ Hence the result follows
 by the theorem. \bexe Let $H,K$ be two finite subgroups of a group $G$ such
 that $\circ(H)=m, \circ(K)=n,$ and $(m,n)=1.$ Prove that $H\cap K=\{e\}$ the
 identity subgroup.\eexe \bexe Prove that elements of finite order in $(\C^*,\cdot)$
  is an infinite torsion group.\eexe \bt\label{t3.10} Let $G$ be a finite group such that
 for all $d \geq 1$, $G$ has at most $d$ elements satisfying
 $x^{d} = e$. Then $G$ is cyclic. \et  \pf Let $\circ(G) = n$ and
 $N(d)$ denotes the number of elements of order $d$ in $G$. If $a\in G$
  is an element of order $d$ in $G$, then for any element $y
 \in gp\{a\}$, $\circ(y)$ divides $\circ(a) = d$ (Corollary
 \ref{co3.9}). Thus by the theorem \ref{t25} (i), $y^{d} = e$. By assumption $G$
  has at most $d$ elements satisfying $x^{d} = e$. Hence $gp\{a\}
  = \{ e,a,a^{2},\cdots, a^{d-1}\}$ is the set of all such
  elements. By the corollary \ref{l1.23}, $\bf{\Phi}(d)$ is the
  number of elements of order $d$ in $gp\{a\}$. Hence if $G$
  contains an element of order $d$ in $G$, then $N(d) =
  \bf{\Phi}(d)$, otherwise
  $N(d) = 0 < \bf{\Phi}(d).$ As order of an elements in
  $G$ divides $\circ(G) = n$, we have $$\begin{array}{ccccccc}
    \sum_{d\mid n} N(d) & = & n & = & \sum_{d\mid n}{{\bf \Phi}}(d) &  & (Lemma\ref{98l}) \\
    \Rightarrow &  & N(n) & = & \bf{\Phi}(n) \geq  1 &  &  \\
  \end{array}$$ Hence $G$ is cyclic (Remark \ref{20} (iii)).
  $\blacksquare$      \bco\label{co3.10} Let $G$ be a group of order $n$
  such that for every $d\geq 1$, $d| n$, there is exactly one subgroup
  of order $d$ in $G$. Prove that $G$ is cyclic.\eco  \pf Let for $d \geq
  1$, $d | n$, $N(d)$ denotes the numbers of elements of order $d$ in
  $G$. By our assumption, if for $x,y \in G$, $\circ(x)=\circ(y)=d$, then
  $gp\{x\}=gp\{y\}$. hence by corollary \ref{44co}, $N(d) \leq
  {\bf\Phi}(d)$. Now the proof is immediate from the theorem.
  $\blacksquare$

 \bexe Prove that for the subgroup $\R = (\R,+)$ of $\C =(\C,+)$,
  $[\C:\R] = \infty.$ \eexe
  \bexe Let for a subgroup $H$ of a group $G$, $[G:H] = n < \infty .$
  Prove that for any $x \in G$ $x^{t} \in H$ for some $1 \leq t \leq n.$\eexe
    \bd\label{d20} Let $H$ be a subgroup of a group $G$. A left
transversal of $H$ in $G$ is a subset $T$ of $G$ such that
\\ (i) $e \in T$ \\ (ii) $xH \cap yH = \phi$ for all $x, y \in T,
x \neq y$ \\ (iii) $ \displaystyle{\bigcup_{x \in T} xH} = G $.
\ed \brs (i) We can define a right transversal of $H$ in $G$ by
taking right cosets of $H$ instead of left cosets. \\ (ii) Any
left transversal of $H$ in $G$  is obtained by taking one element
from each left cosets of $H$ in $G$ and taking the identity
element from the coset $H$.
\\(iii) If $\phi : G \rightarrow T$ is the function $\phi(x) =
x_{t} \in T$ such that $xH =  x_{t}H$ , then for any $x,y \in G$
$$\phi(xy) = \phi(x\phi(y))$$ \hspace{1in} and $$\phi(h) =e \mbox{
for all }h \in H.$$
\\ This function $\phi$ is referred as coset representative
function .\ers \bt\label{t24} Let $G$ be a group and let $H,K$ be
two subgroups of $G$. If $H \supset K,$ and $[G : K] < \infty$,
then $[G: H] < \infty, [H : K] < \infty$, moreover,    \\  $[G :
K] = [G:H].[H
: K].$
\et \pf Let $[G : K] = n$, and let $x_1K, x_2K,...,x_nK$ be all
distinct left cosets of $K$ in $G$. Then $$ G = \bigcup_{i=1}^n
(x_iK)$$ \be\label{e21} \Rightarrow G = \bigcup_{i=1}^n (x_iH)
\mbox{ since } K \subset H. \ee Then for any $g \in G$, $g = x_ih$
for some $1 \leq i \leq n$ and $h \in H.$ Hence $gH = x_{i}H$ and
$[G: H] \leq n$. Further, it is clear by definition of cosets that
$[H
: K] < \infty$, since $[G : K] < \infty$. Let $[H
: K] = t$, and let $h_1K, h_2K,...,h_tK$ be all distinct left cosets
of $K$ in $H.$ Then \be\label{e22} H = \bigcup_{i=1}^t (h_iK) \ee
From the equations \ref{e21} and \ref{e22}, we get
\be
G = \bigcup_{i=1}^{n}x_i \left( \bigcup_{j=1}^{t}h_j K\right) =
\bigcup_{\begin{array}{c}
          \scriptstyle 1\leq i\leq n\\
          \scriptstyle1\leq j \leq t
        \end{array}}(x_i h_j K)
\ee Now, if $x_ih_jK = x_{i_1}h_{j_1}K$, then \ba{rl}
                & x_ih_jKH = x_{i_1}h_{j_1}KH\\
\Rightarrow     & x_iH = x_{i_1}H   \mbox{ since } K \subset
H,\mbox{ and } h_j,h_{j_1} \in H\\ \Rightarrow     & i = i_1
\mbox{ by choice of } x_i\mbox{'s.}\\ \Rightarrow     & x_ih_jK =
x_ih_{j_1}K\\ \Rightarrow     &h_jK = h_{j_1}K\\ \Rightarrow     &
j = j_1 \mbox{ by choice of } h_j\mbox{'s.} \ea Hence $x_ih_jK =
x_{i_1}h_{j_1}K$  implies  $i = i_1 \mbox{ and } j = j_1.$
 Therefore $[G : K] = nt = [G : H][H : K]$ This completes the
proof. $\blacksquare$    \bt\label{t210} Let $H, K$ be any  two
 finite subgroups of a group $G.$ Let $\circ (H) = m$, $\circ (K) = n$,
and $\circ (H \cap K) = t.$ Then $$ | HK | = \frac {mn}{t}.$$\et
\pf We know that $H\cap K$ is a subgroup of $H$ as well as of $K.$
Let
\be\label{e23}(H\cap K)y_1\cup (H\cap K)y_2\cup \ldots \cup (H\cap
K)y_p = K \ee be a right coset decomposition of $K$  with respect
to $H\cap K.$ Then,
\begin{eqnarray}
  HK & = & H(H\cap K)y_1\cup H(H\cap K)y_2\cup \ldots \cup H(H\cap K)y_p\\
\label{eq3.1}     & = & Hy_1\cup Hy_2\cup \ldots \cup Hy_p \quad
(\mbox{use: } H(H\cap K) = H)
\end{eqnarray}
We claim that $Hy_1, Hy_2, \ldots ,Hy_p$ are all distinct cosets.
Let for $i\neq j,$ we have \ba{rl}
               &         Hy_i  =  Hy_j\\
\Rightarrow    &     H  = Hy_jy_i^{-1}\\ \Rightarrow    &
y_jy_i^{-1}  \in  H\cap K\\ \Rightarrow   & (H\cap K)y_i = (H\cap
K)y_j \ea This contradicts that \ref{e23} is coset decomposition
of $K$ with respect to $H\cap K.$  Hence $Hy_1, Hy_2, \ldots Hy_p$
are all distinct. We have $\abs{Hy} = \abs H.$ Therefore from
(\ref{eq3.1}), we get \bea
   \abs{HK}= p\abs{H} & = & \abs{H}\cdot [K : H\cap K]\\
                    & = & \frac{\circ (H) \cdot \circ (K)}{\circ (H\cap K)}
                    \mbox{   (use: Theorem \ref{co22}) }\\
                    & = & \frac {mn}{t}\hspace{1in} \blacksquare
\eea
 \bt\label{t212} Let $H, K$ be two
subgroups of a group $G$ such that $[G:H] = m,$ and $[G:K] = n.$
Then $[G:H\cap K] \leq mn.$ \et \pf Let $a_1H, a_2H, \ldots ,a_mH$
be all distinct left cosets of $H$ in $G$ and let $b_1k,b_2K,
\ldots ,b_nK$ be all distinct left cosets of $K$ in G. Take $x, y
\in a_iH \cap b_jK$. Then we can write $$x = a_ih = b_jk  \qquad
(h \in H, k \in K)$$ $$ y = a_ih^{\prime} = b_jk^{\prime}
\qquad(h^{\prime} \in H, k^{\prime} \in
                                                                        K)$$
Hence, \ba{rl}
        &  y^{-1}x = {h^{\prime}}^{-1}h = {k^{\prime}}^{-1}k \in H\cap K\\
\Rightarrow & x(H\cap K) = y(H\cap K) \ea Further for any $g\in
G$, $g\in (a_{i_1}H\cap b_{j_1}K)$ for some $1\leq i_1\leq m$,
$1\leq j_1\leq n.$  Hence the number of distinct left cosets of
$H\cap K$ in $G$ is less than or equal to the number of members in
the set $ \{(a_iH\cap b_jK \st 1\leq i\leq m, 1\leq j\leq n\},
i.e.,mn.$ Thus $$[G:H\cap K] \leq [G:H]\cdot [G:K] =
mn.\hspace{.5in} \blacksquare$$
\bco\label{co212} If we assume that $(m,n) = 1$,in the theorem,
then $$ [G:H\cap K] = [G:H] [G:K] = mn$$ \eco \pf We have
$$[G:H\cap K] = [G:H] [K:H\cap K] = [G:K] [H:H\cap K]$$ \ba{rl}
\Rightarrow & m\mid [G:H\cap K], n\mid [G:H\cap K]\\ \Rightarrow &
mn \mid [G:H\cap K]\mbox{ since } (m,n) = 1. \ea Hence, by the
theorem $$ [G: H\cap K] = [G:H] [G:K] = mn \hspace{.5in}
\blacksquare$$
\bco\label{co212b} Let $H_1, H_2, \ldots H_k$ be subgroups of a
group $G$ with $[G:H_i] = t_i ; i= 1,2,..,k.$ Then
$$[G:\bigcap_{i=1}^n H_i] \leq t_1t_2\cdots t_k.$$
\eco \pf The proof follows by induction on {k}. $\blacksquare$ \bexe\label{b}
Show that in the theorem \ref{t212}, $[G : H \cap K] \neq mn$ in
general. \eexe
\bt\label{t215a} Let $G$ be a finitely generated group and $H$, a
subgroup of finite index in $G$.
 Then $H$ is finitely generated.\et
 \pf  Let $M$ be a finite subset of $G$ such that $G = gp\{M\}$,
 and let $T$ be a left transversal of $H$ in $G$. As $[G : H] <
 \infty$, $T$ is finite. To prove the result we shall show that
 $$H = gp\{\phi(x^{e}t)^{-1}x^{e}t |\; t\in T, x\in M  \mbox{ and } e = \pm1 \}$$
 where $\phi : G \rightarrow T$ is the coset representative
 function.  \\   For any $h \in H$, we can write \be\label{e25} h =
 x_{1}^{e_{1}}x_{2}^{e_{2}}\cdots x_{k}^{e_{k}} \ee
 where $x_{i} \in M$ and $e_{i}=\pm 1$ for all $i = 1,\cdots
 ,k.$ Now, put $u_{0}=e$, $u_{1}=x_{k}^{e_{k}}$ and $u_{i+1}=x_{k-i}^{e_{k-i}}\cdots
 x_{k}^{e_{k}}$ for all $i = 1,2,\cdots,k-1$. Clearly, $\phi(u_{k}) = \phi(h) = e = \phi(u_{0}).$
 Hence, from the equation
 \ref{e25}, we get $$h =
 \phi(u_{k})^{-1}x_{1}^{e_{1}}\phi(u_{k - 1})\phi(u_{k - 1})^{-1}x_{2}^{e_{2}}\phi(u_{k - 2})
 \phi(u_{k - 2})^{-1}\cdots$$
 $$ \cdots \phi(u_{1})\phi(u_{1})^{-1}x_{k}^{e_{k}}\phi(u_{0})$$
 Now, note that $$\begin{array}{rcl}
   \phi(u_{i+1})  & = & \phi(x_{k -i}^{e_{k - i}}u_{i}) \\
    & = & \phi(x_{k - i}^{e_{k - i}}\phi(u_{i})) \
 \end{array}$$ for all $i = 0, 1, \cdots, k - 1.$ Therefore
 $\phi(u_{i+1})^{-1}x_{k - i}^{e_{k - i}}\phi(u_{i})=\phi(x_{k-i}^{e_{k-i}}\phi(u_i))^{-1}
 x_{k-i}^{e_{k-i}}\phi(u_i)$ is an element of the form
 $\phi(x^{e}t)^{-1}x^{e}t , (t \in T, x \in M, e = \pm 1)$, for all $i = 0, 1, \cdots, k - 1$.
 Hence the result follows. $\blacksquare$  \brs\label{d23.02f} (i) Let $M$ be a
  generator set of a group $G$ and let
    $H$ be a subgroup of $G$. If $T$ is a left transversal of $H$
    in $G$, then by the proof of the theorem it follows that the set $S =
    \{\phi(x^{e}t)^{-1}x^{e}t | t \in T, x \in M, e = \pm 1\}$ is a
    set of generators of $H$. Thus finiteness of $M$ and $[G :
    H] < \infty$ is used to ascertain that $S$ is finite.   \\
     (ii)  We have $$\begin{array}{rrcl}
       & \phi(x^{e}t)H & = & x^{e}tH \\
      \Rightarrow & x^{-e}\phi(x^{e}t)H & = & tH \\
      \Rightarrow & \phi(x^{-e}\phi(x^{e}t)) & = & t \\
    \end{array}$$  Hence   $$\begin{array}{rcl}
      (\phi(x^{e}t)^{-1}x^{e}t)^{-1} & = & t^{-1}x^{-e}\phi(x^{e}t) \\
       & = & \phi(x^{-e}\phi(x^{e}t))^{-1}x^{-e}\phi(x^{e}t) \\
       & = & \phi(x^{-e}\widehat{t})^{-1}x^{-e}\widehat{t} \\
    \end{array}$$  $$\Rightarrow\;\;\; \phi(x^et)^{-1}x^et=(\phi(x^{-e}\widehat{t}\;)
    x^{-e}\widehat{t}\;)^{-1}$$  where $\widehat{t} = \phi(x^{e}t)$.  \\
    Therefore, we have
    $$H = gp\{ \phi(xt)^{-1}xt | t \in T, x \in M \}$$
    (iii) From the proof of the theorem and (ii), we get $$
    H = gp\{ T^{-1}MT \cap H \} \mbox{ since } \phi(xt)^{-1}xt \in H.$$\ers

   \bd\label{d27} A subgroup $H$ of a group $G$ is
called a normal subgroup if $xH = Hx$ for all $x \in G.$ We denote
this fact by the notation $H \lhd G$ or $G \rhd H$.
\ed\brs\label{a} (i)  Any subgroup of an abelian group is normal.
\\ (ii) In a group, the identity subgroup \{e\} and $G$ are always
normal subgroups.
\\ (iii) Any normal subgroup of a group $G$ other then $G$ and
identity subgroup is called a \bf proper normal subgroup \rm of
$G$.
\\ (iv) Let $H$ be a normal subgroup of $G$,$H \neq G$. Then $H$ is
called a {\bf maximal normal subgroup} of $G$ if there exists no normal
subgroup $K$ of $G$ such that $H \varsubsetneqq K \varsubsetneqq G$. \\
 (v) A normal subgroup $H(\neq e)$ of $G$ is called a {\bf minimal normal
 subgroup}
  of $G$ if for any normal subgroup $K(\neq e)$ of $G$, $K\subset H$ implies $K=H.$\ers
\bl\label{w}Let $H,K$ be two subgroups of a group $G$. Then \\ (a)
If $H \lhd G$,$HK$ is a subgroup of $G$. Further, if $K$ is also
normal in $G$, then  $HK \lhd G$.
\\ (b) If $H \lhd G$, $K \lhd G$ and $H \cap K$ = \{e\},
then $hk = kh$ for all $h \in H$, $k \in K$.\el \pf (a) As $H \lhd
G$, $HK = KH$. Hence by the theorem \ref{t29}, $HK$ is a subgroup
of $G$. Further, let $K \lhd G$. Then for any $x \in G$,
$$\begin{array}{cccccc}
  xHK & = & HxK &  &  & (H \lhd G) \\
   & = & HKx &  &  & (K \lhd G)
\end{array}$$
   Hence $HK \lhd G$.    \\  (b) Let $h \in H$, $k \in K$. Then as
   $H \lhd G$, $K \lhd G$,  $$hkh^{-1}k^{-1} \in H \cap K =
   \{e\}$$
$$\begin{array}{crcl}
  \Rightarrow & hkh^{-1}k^{-1} & = & e \\
  \Rightarrow & hk & = & kh
\end{array}$$
Hence the result follows.   $\blacksquare$\bt\label{t213} Let H be
a normal subgroup of a group G.Then the set $$ G/ H = \{xH \st x
\in G\} $$ is a group with respect to the product of subsets in G.
\et \pf Let $xH, yH \in G/ H.$ Then \bea
                (xH)(yH) & = & x(Hy)H\\
                         & = & x(yH)H  \qquad \mbox{ since } H \lhd G \\
                         & = & xy(HH) = xyH.
\eea Thus $G/ H$ is closed with respect to set product.\\ For $eH
= H \in G/ H$ , we have $$(xH)(H) = xH = H(xH)$$ for all $xH \in
G/ H$. Hence $H$ is the identity element in $G/ H.$ As product in
$G$ is associative, the set product is associative in $G$.
Finally, for any  $xH \in G/ H$, $x^{-1}H \in G/ H$ and
$$(xH)(x^{-1}H) = H = (x^{-1}H)(xH),$$ i.e., every element of
$G/H$ admits an inverse with respect to $H.$ Thus $G/ H$ is a
group.  $\blacksquare$  \bco\label{co213a} Let $G$ be a finite
group and let $H$ be a normal subgroup of $G.$ Then $G/ H$ is a
group of order $\frac{\circ (G)}{\circ (H)}$. \eco  \pf : By the
theorem \ref{co22} $[G:H] = \frac{\circ (G)}{\circ (H)}$.
Hence$G/H$ is a group of order $\frac{\circ (G)}{\circ (H)}$.
$\blacksquare$\bd\label{d28} Let $H$ be a subgroup of a group $G$
such that $H\lhd G.$ Then the group $G/ H$ defined above is called
the factor group or quotient group of $G$ by $H.$ \ed
\bd\label{27a} Let $G$ be a group with $\circ(G) > 1$. The group
$G$ is called a simple group if there exists no proper normal
subgroup of $G$.\ed
\bx\label{e213} Let $G$ be a group of prime order. As order of a
subgroup divides the order of the group, for any subgroup $H$ of
$G$, either $H = \{e\}$ or $H = G$. Hence a group of prime order
is simple.\ex  We shall give more examples of simple groups in due
course and shall show that any abelian simple group is of prime
order. \bt\label{t214} Let $H$ be a subgroup of a group $G$. Then
$H$ is a normal subgroup of $G$ if and only if $x^{-1}Hx \subset
H$  for all $x \in G.$ \et \pf If $H\lhd G$, then for any $x \in
G$, \ba{rl}
                &  Hx = xH\\
\Rightarrow     &  x^{-1}Hx = H \ea Conversely, if for all $x \in
G$, $$x^{-1}Hx \subset H$$ Then $$Hx \subset xH$$ Taking $x^{-1}$
for $x$, we get $$Hx^{-1}\subset x^{-1}H $$ $$ \Rightarrow \qquad
xH \subset Hx$$ Consequently, $Hx = xH$  for all $x \in G.$ Thus
$H \lhd G.$  $\blacksquare$ \bco\label{co214} For a subgroup $H$
of a group $G$, $H\lhd G$ if and only if $xH = yH$ implies $Hx =
Hy$ for all $x, y \in G.$ \eco \pf If $H\lhd G$, then $xH = yH$
implies $Hx = Hy$ by the definition of the normal subgroup.
Conversely, let the condition holds. Take $x\in G$, $h\in H.$
Then \ba{rl}
                       & xH = xhH\\
\Rightarrow            & Hx = Hxh\\ \Rightarrow            & H =
Hxhx^{-1}\\ \Rightarrow            & xhx^{-1}\in H \mbox{ for all
}h\in H\\ \Rightarrow            & xHx^{-1}\subset H \mbox { for
all }x\in G \ea Hence $H\lhd G.$   $\blacksquare$
 \bexe  Let $K$ be a field. Prove that $Sl_{n}(K)$ is a normal subgroup of $Gl_{n}(K)$.
 (Hint: Apply the theorem \ref{t214} ) \eexe
 \bexe Let $\{H_i\}_{i\in I}$ be a family of normal subgroups of a group
 $G$. Prove that $\bigcap_{i\in I}H_i$ is a normal subgroup of $G$.\eexe
 \bd For any subset $S$ of a group $G$, the intersection of all normal
 subgroups of $G$ containing $S$ is called the normal closure of $S$ in
 $G$, and is denoted by $\overline{N}_G(S)$.\ed \brs (i) $\overline{N}_G(S)$
 is the smallest normal subgroup of $G$ containing $S$. \\ (ii) A subgroup
 $H$ of $G$ is normal if and only if $\overline{N}_G(H)=H.$\ers
 \bt\label{t215} Let $H$ be a subgroup of a group $G.$ For any $x,y \in G$, define a
relation $R$ as: $$x R y  \Leftrightarrow  xy^{-1}\in H $$ Then
$R$  is an equivalence relation over $G$ and $Hx$ is the
equivalence class of $x$ with respect to $R.$ \et \pf For any
$x\in G, xx^{-1} = e \in H.$  Thus $x R x$ holds. Further, for any
$x, y \in G$, $$xy^{-1}\in H \Leftrightarrow yx^{-1} =
(xy^{-1})^{-1} \in H$$ $$ \mbox{Hence, } x R y \Leftrightarrow y R
x$$ Next,let  $x, y, z \in G$  be such that $x R y, y R z$ hold.
Then $$xy^{-1}\in H, yz^{-1}\in H$$ $$\Rightarrow \quad
(xy^{-1})(yz^{-1}) = xz^{-1}\in H$$ Therefore $x R z$   holds.
Thus $R$ is reflexive, symmetric and transitive i.e., $R$ is an
equivalence relation over $G$. Finally, for $x, y\in G$ \bea
                   x R y  & \Leftrightarrow & xy^{-1}\in H\\
                          & \Leftrightarrow & yx^{-1}\in H\\
                          & \Leftrightarrow & y\in Hx
\eea Therefore, $Hx$ is the equivalence class of $x$ in $G$ with
respect to the relation $R.$  \\ \clearpage
 {\bf EXERCISES}
\begin{enumerate}
\item
Let $H$ be a normal subgroup of a group $G$ with $[G:H] =n <
\infty.$ Show that $x^{n} \in H$ for all $x \in G.$
\item
 Show that $[\C^{*}:\R^{*}] = \infty.$
 \item
 Prove that if $H$ is a subgroup of $\C^{*}$ with $[\C^{*}:H] <
 \infty$, then $H = \C^{*}.$
 \item
 Prove that the statement of the exercise $3$ is not true if we
 replace $\C$ by $\Q$ or $\R$.
 \item
 Let $H$ be a subgroup of $\Q = (\Q,+)$ such that $[\Q:H]<\infty.
 $ Prove that $H = \Q$.
\item
Show that in the exercise 5, $\Q$ can be replaced by $\R =(\R,+)$
or by $\C = (\C,+).$
\item
Let for an abelian group $G$, $G = \{x^{n} \st x \in G\}$ for all
$n \geq 1.$ Prove that for any subgroup $H$ of finite index in
$G$, \;$G = H$.
\item
For the subgroup $\Z = (\Z,+)$ of the group $\Q = (\Q,+)$, prove
that  \\ (a) Every element of $\Q/\Z$ has finite order.  \\  (b)
For any element $\overline{q} = q+\Z$ in $\Q/\Z$ and $n(\neq 0)
\in \Z$, the equation $nX = \overline{q}$ has a solution in $\Q /
\Z$.  \\  (c) Any finitely generated subgroup of $\Q /\Z$ is
cyclic.  \\  (d) For any $m / n (\neq 0) \in \Q$, order of the
element $m /n +\Z$ in $\Q / \Z$ is $n$ if and only if $(m,n) = 1.$
\item
Prove that a normal subgroup of a normal subgroup need not be
normal.
\item
Let for a subgroup $H$ of a group $G$, $x^{2} \in H$ for all $x
\in G$. Prove that $H\lhd G$.
\item
Prove that any subgroup of the group of Quaternions is a normal
subgroup.
\item
Show that if a normal subgroup $H$ of a group $G$ is cyclic, then
every subgroup of $H$ is normal in $G$.
\item
Let $G$ be an abelian group and let $S = \{ x \in G \st x^{2} =
e\}$. Prove that  $$\prod_{g \in G}g = \prod_{x \in S}x =
\left\{\begin{array}{ccl}
  e & if & \abs{S} \neq 2 \\
  x & if & S =\{e,x\},x\neq e
\end{array}\right.$$  (Hint : Note that $S$ is a subgroup of $G$.
If $\abs{S} > 2$, then for any $y(\neq e) \in S$, consider the
subgroup $A =\{e,y\}$ of $S$ and show that if $[S:A] = m$, then
$\prod_{x \in G}x = y^{m}$.)
\item
Let $H$ be a subgroup of a group $G$ with finite index. Prove that
for any $x \in G$, there exist $a_{1},a_{2},\cdots,a_{n}$ in $G$
such that $$HxH = \bigcup_{i = 1}^{n}a_{i}H = \bigcup_{i =
1}^{n}Ha_{i}.$$
\item
Let $H$ be a subgroup of finite index in $G$. Show that there
exist $a_{1},a_{2},\cdots,a_{n}$ in $G$ such that $$G = \bigcup_{i
= 1}^{n}a_{i}H = \bigcup_{i = 1}^{n}Ha_{i}.$$
 \item
 Let $H_i; i=1,2,\cdots,n$ be finite normal subgroups of a group $G$.
 Prove that $H=H_1 H_2\cdots H_n$ is a finite normal subgroup of $G$ and
 $\circ(H)$ divides $\circ(H_1)\cdot \circ(H_2)\cdots \circ(H_n).$
 \item
 Let $H,K$ be two subgroups of a group $G$. Prove that if $Hx=Ky$ for some
 $x,y\in G$ then $H=K.$
 \item
 Prove that for any subgroup $H$ of the cyclic group $G=C_n,$ there exists
 a subgroup $K$ of $G$ such that $K$ is isomorphic to $G/H.$
 \item
 Let $H$ be a normal subgroup of a group $G$. Prove that if $H$ and $G/H$
 are finitely generated, then $G$ is finitely generated.
 \item
 Let $C_n=gp\{a\}$ be a cyclic group of order $n$. Prove that
 $\Phi(C_n)=gp\{a^m\}$ where $m$ is the product of all distinct prime
 divisors of $n$.
 \item
 Prove that for any group $G$, and $n\geq 0,$ the subgroup
 $P_n(G)=gp\{x^n\st x\in G\}$ is a normal subgroup of $G$.
 \item
 Let $S$ and $T$ be two subsets of a group $G$. Prove that \\ (i)
 $\overline{N}_G(\overline{N}_G(S))=\overline{N}_G(S).$ \\ (ii) $\overline{N}_G(S\bigcup
 T)=\overline{N}_G(S)\overline{N}_G(T).$
 \item
 Prove that for any subset $S$ of a group $G$, $\overline{N}_G(S)=gp\{g^{-1}xg\st
 g\in G, x\in S\}$.
 \item
 Let $A$ be an abelian normal subgroup of a group $G$. Prove that if for a
 subgroup $H$ of $G$, $AH=G$ then $H\cap A$ is normal in $G$.
\end{enumerate}

 \chapter {The Commutator Subgroup And The Center Of A Group }
 \footnote{\it contents group4.tex }
One of the most important properties by which general groups differ from
$(\Z, +), (\R, +)$ etc. is commutativity. We shall define concepts of
commutator subgroup and center of a group which give measures of non -
commutativity of a group. To make matters precise, let us define:
\bd\label{d101z} For any two elements $a, b$ of a group $G$, the element
$$[a,b] = aba^{-1}b^{-1}$$ is called the commutator of the ordered pair
$\{a,b\}$.\ed \brs (i) We have $$[a,b]^{-1} = bab^{-1}a^{-1} = [b,a]$$
(ii) Note that $$\begin{array}{crcl}
   & [a,b] & = & e, \mbox{    the identity element}\\
  \Leftrightarrow & aba^{-1}b^{-1}& = & e \\
  \Leftrightarrow & ab & = & ba
\end{array}$$
Thus the commutator of the pair $\{a,b\}$ is in some sense the
measure of non-commutativity of $a$ and $b$.\ers \bd\label{d102z}
let $G$ be a group. Then $$\begin{array}{rclc}
  G^{'} & = & [G, G] &  \\
   & = & gp\{[a,b] | a, b \in G\} &
\end{array}$$
is called the (first) derived group or commutator subgroup of $G$.
\ed \brs (i) For any two subgroups $A,B$ of a group $G$,we write
$$[A,B] = gp\{[a, b] | a \in A, b \in B \}$$ (ii) The group $G$ is
abelian if and only if $G' = \{e\}$.\ers \bexe Let $G$ be a group.
Prove that  \\ (i) If $x,a,b \in G$, then  $[x^{-1}ax, x^{-1}bx] =
x^{-1}[a,b]x$. Thus there exist $c,d \in G$ such that $[a,b]x =
x[c,d].$  \\  (ii) If $a_{i},b_{i},x \in G$, $1\leq i\leq n$, then
$$x^{-1}(\prod_{i = 1}^{n}[a_{i},b_{i}])x = \prod_{i =
1}^{n}[x^{-1}a_{i}x, x^{-1}b_{i}x]$$ \\ (iii) The derived group
$G'$ of $G$ is normal in $G$.
\eexe
\bexe Let $G$ be a group. Prove that for any $x,y,z \in G$, \\ (i)
$[xy,z]=x[y,z]x^{-1}[x,z]$ \\ (ii) $[x,yz]=[x,y] y[x,z]y^{-1}$
\\ (iii) $y^{-1}[[y,x^{-1}],z^{-1}]y z^{-1}[[z,y^{-1}]x^{-1}]z
x^{-1}[[x,z^{-1}],y^{-1}]x=Id$  \\ hold. \eexe
 A more direct measure of the commutativity of a group is given by the
following
:
\bd\label{d103z} For any group $G$, $$ Z(G) = \{x \in G | xg = gx
\mbox{ for all } g \in G\}$$ is called the center of $G$ and is
denoted by $Z(G)$.\ed    \br  Clearly $Z(G) = G$ if and only if
$G$ is abelian. \er  \bl\label{l101l} Let $H$ be a subgroup of a
group $G$ such that $G'\subset H$. Then $H$ is normal in $G$. In
particular $G'\lhd G.$\el \pf Let $g \in G$ and $h \in H$. Then
$$\begin{array}{crclcc}
  & ghg^{-1}h^{-1} & \in & H &  &  \\
  \Rightarrow & ghg^{-1} & \in & Hh = H &  &  \\
  \Rightarrow & gHg^{-1} & \subset & H &  & \mbox{for all } g \in
  G
\end{array}$$ Hence by the theorem \ref{t214}, $H$ is a normal in
$G$.      $\blacksquare$
\bco\label{b101z} The factor group $G / G'$ is abelian. \eco \pf
For any $x, y \in G$ , $$\begin{array}{rcl}
  (xG')(yG') & = & xyG' \\
   & = & xy[y^{-1},x^{-1}]G' \\
   & = & xyy^{-1}x^{-1}yxG' \\
   & = & yxG' \\
   & = & (yG')(xG')
\end{array}$$
Hence $G/ G'$ is abelian.       $\blacksquare$ \bl\label{l102l}
For any normal subgroup $H$ of a group $G$, $G/ H$ is abelian if
and only if $G' \subset H$.\el    \pf  Let $G/ H$ be abelian. Then
for any $x, y \in G$ , $$\begin{array}{crclcc}
   & (xH)(yH) & = & (yH)(xH) &  &  \\
  \Rightarrow & xyH & = & yxH &  &  \\
  \Rightarrow & Hxy & = & Hyx &  &  \\
  \Rightarrow & Hxyx^{-1}y^{-1} & = & H &  &  \\
  \Rightarrow & [x, y] & = & xyx^{-1}y^{-1} & \in & H \\
  \Rightarrow & G' \subset H &  &  &  &
\end{array}$$  Conversely if $G' \subset H$, then for any $x, y
\in G$, $$\begin{array}{crclcc}
  & [x, y] & = & xyx^{-1}y^{-1} \in  H  & \\
 \Rightarrow & Hxyx^{-1}y^{-1} & = & H &  &  \\
  \Rightarrow & Hxy & = & Hyx &  &  \\
  \Rightarrow & (Hx)(Hy) & = & (Hy)(Hx) &  &
\end{array}$$
Hence $G/ H$ is abelian.   $\blacksquare$ \bt\label{4.44t} For any
group $G$, $$G' = \{(x_{1}\cdots x_{n})(x_{n}\cdots x_{1})^{-1}
\st x_{i} \in G, n\geq 1\}$$ \et  \pf Put $$K = \{(x_{1}\cdots
 x_{n})(x_{n}\cdots x_{1})^{-1} \st x_{i} \in G, n\geq 1\}$$ We
shall first show that $K$ is a subgroup of $G$. Let $a =(x_{1}\cdots
x_{n})(x_{n}\cdots x_{1})^{-1}$ and $b = (y_{1}\cdots y_{k})(y_{k}\cdots
y_{1})^{-1}$ be two elements of $K$. Then taking, $c = (y_{k}y_{k-1}\cdots
y_{1})$, we have $$\begin{array}{ccl}
  ab & = & (x_{1}\cdots x_{n}\cdot x_{1}^{-1}x_{2}^{-1}\cdots x_{n}^{-1})(y_{1}y_{2}\cdots y_{k}y_{1}^{-1}\cdots y_{k}^{-1}) \\
   & = & x_{1}x_{2}\cdots x_{n}y_{1}^{-1}y_{2}^{-1}\cdots y_{k}^{-1}cx_{1}^{-1}x_{2}^{-1}\cdots x_{n}^{-1}y_{1}y_{2}\cdots y_{k}c^{-1}\\
   & = & (x_{1}x_{2}\cdots x_{n}y_{1}^{-1}y_{2}^{-1}\cdots y_{k}^{-1}c)(cy^{-1}_{k}\cdots y^{-1}_{2}y_{1}^{-1}x_{n}\cdots x_{2}x_{1})^{-1}\\
   \end{array}$$ and $$a^{-1} = (x_{n}\cdots x_{1})(x_{1}\cdots
x_{n})^{-1}$$ Hence $ab,  a^{-1}\in K$ for all $a, b \in K$. Thus $K$ is a
subgroup of $G$. Now, note that for any $x,y \in G$, $$\begin{array}{ccl}
  (xy)(yx)^{-1} & = & xyx^{-1}y^{-1} \\
   & = & [x,y]
\end{array}$$ Thus $[x,y] \in K$ for all $x,y \in G.$
Consequently $G'\subset K.$ To complete the proof, we have to show
that $K \subset G'$. To do this we shall show by induction on $n$,
that for any $x_{i} \in G$, $1\leq i \leq n$, $(x_{1}\cdots
x_{n})(x_{n}\cdots x_{1})^{-1} \in G'.$ The result is clear for $n
= 1,2$. Now, let $n\geq 3$. Then $$\begin{array}{ccl}
  (x_{1}\cdots x_{n})(x_{n}\cdots x_{1})^{-1} & = & (x_{1}\cdots x_{n})(x_{1}^{-1}\cdots
x_{n}^{-1}) \\
   & = & ([x_{1},x_{2}]x_{2}x_{1}x_{3}\cdots x_{n})(x_{1}^{-1}\cdots x_{n}^{-1}) \\
   & = & ([x_{1},x_{2}]x_{2}[x_{1},x_{3}]x_{3}x_{1}x_{4}\cdots x_{n})(x_{1}^{-1}\cdots x_{n}^{-1})
\end{array}$$ Continuing as above, as $xG' = G'x$ for all $x\in G$, we
get $$\begin{array}{ccl}
  (x_{1}\cdots x_{n})(x_{n}\cdots x_{1})^{-1} & = & ([x_{1},x_{2}]x_{2}[x_{1},x_{3}]x_{3}\cdots [x_{1},x_{n}]x_{n}x_{1})(x_{1}^{-1}\cdots x_{n}^{-1}) \\
   & = & z(x_{2}\cdots x_{n})(x_{2}^{-1}\cdots x_{n}^{-1}) \;\;\;\;\;\; (z\in G')\\
   & = & (x_{2}\cdots x_{n})(x_{n}\cdots x_{2})^{-1}
\end{array}$$ As $(x_{2}\cdots x_{n})(x_{n}\cdots x_{2})^{-1} \in
G'$ by induction, we get that $(x_{1}\cdots x_{n})(x_{n}\cdots
x_{1})^{-1} \in G'.$ Consequently $$G' = \{(x_{1}\cdots
x_{n})(x_{n}\cdots x_{1})^{-1} \st x_{i} \in G, n\geq 1\}.
\hspace{.25in} \blacksquare$$
\bl\label{l103l} Let $G$ be a group. Then $Z(G) \lhd G$ i.e. ,
$Z(G)$ is a normal subgroup of $G$. \el \pf First of all $Z(G)
\neq \emptyset $, since $e \in Z(G)$. Now, let $x, y \in Z(G)$ and
$g \in G$. Then $$\begin{array}{crcl}
  & x^{-1}gx & = & g \\
  \Rightarrow & x^{-1}g & = & gx^{-1} \\
  \Rightarrow & x^{-1}gy & = & gx^{-1}y \\
  \Rightarrow & x^{-1}yg & = & gx^{-1}y
\end{array}$$
Hence $x^{-1}y \in Z(G)$. Therefore $Z(G)$ is a subgroup of $G$ (Theorem
\ref{t1}). The fact that $Z(G)$ is normal is clear by definition. Hence
the result follows.  $\blacksquare$ \bt\label{t4.14tt} Let $G$ be a group
with $[G:Z(G)]=n$. Then $\circ(G')\leq (n^2-3n+3)^{n(n^2-3n+3)}\leq
n^{2n^3}.$ \et  \pf We shall prove that the result in steps. \\ {\bf Step
I: } $G$ has atmost $n^2-3n+3$ commutators. \\ Let $Z=Z(G)$ and let
$a_1Z=Z,a_2Z,\cdots,a_nZ$ be all distinct cosets of $Z$ in $G$. If $x,y\in
G,$ then $x=a_ic,y=a_jd$ for some $1\leq i,j\leq n$, and $c,d\in Z.$
Therefore $$[x,y]=[a_i,a_j]$$ Note that $a_1\in Z$, and $[y,b]=e$ for all
$b\in Z$ and $y\in G.$ Therefore the number of commutators in $G$ is less
than or equal to $$2((n-2)+(n-3)+\cdots +1)+1=n^2-3n+3.$$ {\bf Step II:}
For any $x,y\in G, [x,y]^{n+1}=[x,y^2][yxy^{-1},y]^{n-1}.$ \\ By assumption
$\circ(G/Z)=n.$ Hence for any $x,y\in G, [x,y]^n\in Z.$ Therefore
$$\begin{array}{cl}
  [x,y]^{n+1} & =xyx^{-1}[x,y]^ny^{-1} \\
   & =xy^2x^{-1}y^{-2}y[x,y]^{n-1}y^{-1} \\
   & =[x,y^2][yxy^{-1},y]^{n-1}
\end{array}$$ {\bf Step III:} Each element of $G'$ is a product of at
most
$n(n^2-3n+3)$ commutators. \\ Let us observe for $a,b,c,d$ in $G$,
$$\begin{array}{cl}
  [c,d][a,b] & =[a,b][a,b]^{-1}[c,d][a,b] \\
   & =[a,b][e,f]
\end{array}$$ where $e=x^{-1}cx, f=x^{-1}dx$ for $x=[a,b].$ Hence, by step
I, any element in $G'$ is a product of powers of distinct commutators in
$G$ which are atmost $n^2-3n+3$ in number. Using step II, we can assume
that the power of each commutator is $\leq n.$ This proves the assertion . \\
 Now, the theorem follows easily from steps I to III.        $\blacksquare$
\bexe\label{ex101} For the group $G = S_{3}$ (Example \ref{ex7}), prove that
$Z(S_{3}) = \{e\}$, and $$G' = \left\{\left(\begin{array}{ccc}
    1 & 2 & 3 \\
    1 & 2 & 3 \
  \end{array}\right),
\left(\begin{array}{ccc}
  1 & 2 & 3 \\
  3 & 1 & 2
\end{array}\right),
\left(\begin{array}{ccc}
  1 & 2 & 3 \\
  2 & 3 & 1
\end{array}\right)\right\}$$
    \eexe
\bl\label{l104l} Let $K$ be a field . Then the center of $Gl_{n}(K)$ is the set of all
diagonal matrices in $Gl_{n}(K)$ with constant entry i.e. , $$Z(Gl_{n}(
K)) = \{ \alpha \, Id | \alpha ( \neq 0) \in K \}.$$ where $Id$ is the
identity matrix of order $n$.\el \pf Clearly, for any $A = (a_{ij}) \in
Gl_{n}(K)$, $$\begin{array}{rcl}
 ( \alpha Id )A & = & \alpha (a_{ij}) \\
   & = & A( \alpha Id )
\end{array}$$
Further, if $A \in Z(Gl_{n}(K))$, then for any elementary matrix
$e_{ij}(r)$, $r( \neq 0) \in  K$, $$\begin{array}{crclc}
  & Ae_{ij}(r) & = & e_{ij}(r)A &  \\
  \Rightarrow & AE_{ij}(r) & = & E_{ij}(r)A &  \\
  \Rightarrow & ra_{pi} & = & 0  & \mbox{ for all } p \neq i, \\
  \mbox{ and } & ra_{ii} & = & ra_{jj} &  \\
  \Rightarrow & a_{ii} & = & a_{jj} & \mbox{ for all } i \neq j, \\
  \mbox{ and } & a_{pi} & = & 0 & \mbox{ for all } p \neq i
\end{array}$$
Hence $A = \alpha Id$ where $\alpha = a_{ii}$ for all $1\leq i\leq
n$. Further, clearly $\alpha \neq 0$, since $detA \neq
 0$. Thus the result follows.       $\blacksquare$
\bco\label{c101z}  $Z(Sl_{n}(K)) = \{ \xi Id | \xi :n^{th} \mbox{
root of 1 in }K \}$ \eco \bl\label{l105l} Let $K$ be a field and $n \geq
3$. Then the derived group of $E_{n}(K)$ is equal to $E_{n}(K)$ i.e.,
$(E_{n}(K))' = E_{n}(K)$.\el \pf We know that  $$E_{n}(K) = gp\{ e_{ij}(
\lambda ) = Id + E_{ij}( \lambda ) | \lambda ( \neq 0) \in K , 1 \leq i
\neq j \leq n\}$$ in $Gl_{n}(K)$,where $E_{ij}( \lambda )$ denote the $n
\times n $ matrix over $K$ with $(i,j)^{th}$ entry $\lambda $ and all
other entries $0$. Let us observe that if $j \neq k$, then $E_{ij}(\lambda
)E_{kl}(\mu ) = 0$ and $E_{ij}(\lambda )E_{jk}(\mu ) = E_{ik}( \lambda \mu
)$ for all $k \neq i$. Therefore, for $i \neq j$, $j \neq k$, $k \neq i$,
$$\begin{array}{rcl}
  & &[e_{ij}(\lambda ), e_{jk}( \mu )]\\
   & = & (Id + E_{ij}(\lambda ))(Id + E_{jk}(\mu ))(Id + E_{ij}
   (- \lambda ))(Id + E_{jk}(-\mu)) \\
   & = & e_{ik}( \lambda \mu )
\end{array}$$ Thus for any $e_{ij}( \lambda )$, taking $k \neq
i,j$ $$[e_{ik}(\lambda ),e_{kj}(1)] = e_{ij}(\lambda ) \in (E_{n}(K))'$$
Hence  $(E_{n}(K))' = E_{n}(K)$. $\blacksquare$
\bl\label{l106l} For any field $K$, the subgroup $Sl_{n}(K)$ of $Gl_{n}(K)$ is its
derived group for all $n\geq 3.$\el \pf  For any $A, B \in Gl_{n}(K)$, we
have $$\begin{array}{rcl}
  det[A, B] & = & det(ABA^{-1}B^{-1}) \\
   & = & 1
\end{array}$$   Hence $[Gl_{n}(K), Gl_{n}(K)] \subset Sl_{n}(K
)$. Further, by theorem \ref{t109t}(ii), $$Sl_{n}(K) = gp\{e_{ij}( \alpha )
| \alpha ( \neq 0) \in K \}.$$ Now, as in the lemma \ref{l105l}, for any
$i \neq j, k \neq j,$ and $ k \neq i,$
 $$\begin{array}{rl}
   &[e_{ik}( \alpha ), e_{kj}( \beta )]\\

   = & e_{ij}(\alpha \beta )
\end{array}$$
Hence $Sl_{n}(K) \subset [Gl_{n}(K), Gl_{n}(K)]$ and the result follows.
$\blacksquare$ \br All the results from the lemma \ref{l104l} to the lemma
\ref{l106l} are true for $\Z$ as well. To see the validity of the lemma
\ref{l106l} in the case of $\Z$ use the theorem \ref{.003t} instead of
theorem\ref{t109t}.\er \clearpage
 {\bf EXERCISES}
\begin{enumerate}
  \item
  Let $H$ be a normal subgroup of a group $G$. Prove that  \\
(i) $\displaystyle{\left(\frac{G}{H}\right)' = \frac{G'H}{H}}$  \\ (ii)
$\displaystyle{Z\left(\frac{G}{H}\right) \neq \frac{Z(G)H}{H}}.$
\item
Let $H,K,L$ be normal subgroup of a group. Prove that \\ $[HK,L] =
[H,L][K,L].$ \\ and (ii) $[[H,K],L]\subset [[K,L],H][[L,H],K].$
\item
Let for a group $G$, $G'\subset Z(G).$ Prove that for any $x,y \in G$,
$$(xy)^{n} = x^{n}y^{n}[y,x]^{n(n-1)/2}.$$
\item
Show that for the group of Quaternions $\bf H$, $Z(\bf H) = {\bf
H}'.$
\item
Prove that for the symmetric group $S_{3}$, $$S_{3}' = \left\{ Id,
\left(\begin{array}{ccc}
  1 & 2 & 3 \\
  1 & 3 & 2
\end{array}\right),\left(\begin{array}{ccc}
  1 & 2 & 3 \\
  2 & 3 & 1
\end{array}\right)\right\}.$$
\item
Let $p,q$ be two distinct primes and let $G$ be a group such that for any
$x(\neq e) \in G$, $\circ(x) = p$ or $q$. If for any $x,y \in G$ with
$\circ(x) = \circ(y) = p$, $\circ(xy) = q$ or $1$, and elements of order
$q$ in $G$ together with identity form a cyclic group, then $G'$ is cyclic
of order $q^{m}(m \geq 0).$
\item
Prove that the derived group of the dihedral group $$D_{n} = \{
Id, \tau, \sigma, \sigma^{2}, \cdots, \sigma^{n-1}, \tau\sigma,
\tau\sigma^{2}, \cdots,\tau\sigma^{n-1} \;\st \; \tau^{2} =
\sigma^{n} = Id, \sigma\tau = \tau\sigma^{n-1}\} $$ is the cyclic
group $gp\{\sigma^{2}\}.$
\item
Prove that \\ (i) $Z(O(n.\R))=\{\pm Id\}$ \\ (ii)
$Z(SO(n,\R))=\left\{\begin{array}{cl}
  Id & \mbox{ if $n$ is odd }\\
  \pm Id & \mbox{ if $n$ is even }
\end{array}\right.$  \\ (iii) $Z(U(n,\C))=\{diag(z,z,\cdots,z)\st
\abs{z}=1\}$  \\ (iv) $Z(SU(n,\C))=\{diag(w,w,\cdots,w)\st w^n =1\}$
\item
Let $\F$ be a finite field with $q$-elements. Prove that $Z(Sl_n(\F))$ has
$d = g.c.d.(n.q-1)$ elements.
\item
Let $D_n$ be the Dihedral group. Prove that  \\ (i) If $n$ is odd,
$\circ(Z(D_n))=1$  \\ (ii) If $n$ is even, $\circ(Z(D_n))=2$.
\item
Let $G$ be a group. Prove that every commutator is a product of squares in
$G$.
\item
Let $S$ be a set of generators of a group $G$. Prove that $G'$ is the
smallest normal subgroup of $G$ containing $\{[a,b]\st a,b\in S\}.$
\item
Let $H,K$ be two subgroups of a group $G$. Prove that $[H,K]$ is a normal
subgroup of $gp\{H\cup K\}.$
\item
Let $H$ be a normal subgroup of a group $G$. Prove that if $[H,G']=id,$
then $[H',G]=id.$
\item
Prove that the exercise 14 is true even when $H$ is any subgroup of $G$.
\item
Let $G$ be a group and $x,y\in G$ be such that $[x,y]$ commutes with $x$
and $y$. Prove that if $\circ(x)=m$, then $\circ([x,y])$ divides $m.$
\item
Let for a group $G$, $G=gp\{x,y\}.$ Prove that if $[x,y]\in Z(G),$ then
$G'\subset Z(G).$
\item
Prove that for a subgroup $H$ of a group $G$, $[G,H]H$ is the smallest
normal subgroup of $G$ containing $H.$
\item
Let $N$ be a normal subgroup of a group $G$ with $N\cap G'=\{e\}$. Prove
that $N\subset Z(G).$
\item
Let for a group $G$,$x^2y^2=y^2x^2$ and $x^3y^3=y^3x^3$ for all $x,y\in
G.$ Prove that $[x,y]\in Z(G)$ for all $x,y\in G.$
\item
Let $A,B$ be abelian subgroups of a group $G$ such that $G=AB.$ Prove that
$Z(G)=(A\cap Z(G))(B\cap Z(G)).$
\item
Prove that if $A,B$ are abelian subgroups of a group $G$ and $G=AB$, then
$[A,B]$ is an abelian group. Further, show that $G'=[A,B].$
\end{enumerate}

\chapter{Homomorphisms And Some Applications}
 \footnote{\it contents group5.tex }
 Here, we define maps between groups called homomorphisms and some variants.
 These maps help in connecting different group structures. Image of a homomorphism
 is in a sense an algebraic deformation of the group into another group. We shall study
 these ideas below.
 \bd\label{d31}
Let $G, H$ be two groups.  A mapping $\varphi: G \to H$ is called
a homomorphism if it satisfies : $$\varphi (xy) = \varphi (x)
\varphi (y) \mbox { for all } x, y \in G.$$ \ed
\bd\label{d32} Let
$G,H$ be two groups and $\varphi: G \to H$,  a homomorphism.  We
shall say that $\varphi$ is a monomorphism if $\varphi$ is
one-one. Further $\varphi$ is called an epimorphism if $\varphi
(G) = H$ i.e., $\varphi$ is onto.
\ed
\bd\label{d33} A
homomorphism $\varphi$ from a group $G$ to a group H  is called an
isomorphism if $\varphi$ is one-one and onto. \ed \br Some authors
define a monomorphism as (into) isomorphism. \er \bd\label{d34} A
homomorphism from a group $G$ to itself is called an endomorphism,
and an isomorphism from $G$ to $G$ is called an automorphism. \ed
\bos Let $\varphi: G \to H$  be a homomorphism of groups.  Then\\
(i)  For the identity  $e_G$ of $G$,  $\varphi (e_G) = e_H,$ the
identity of $H.$  \\ We have \ba{rl}
            & \varphi (e_G) = \varphi (e_G\cdot e_G) = \varphi (e_G)\cdot \varphi (e_G)\\
\Rightarrow & \varphi (e_G) = e_H \ea (ii) For any $x\in G,
\varphi (x^{-1}) = (\varphi (x))^{-1}$\\ We have \ba{rl}
            & \varphi (x\cdot x^{-1}) = \varphi (e_G) = e_H\\
\Rightarrow &              \varphi (x)\cdot \varphi (x^{-1}) =
e_H\\ \Rightarrow &               \varphi (x^{-1}) = (\varphi
(x))^{-1} \ea (iii) Identity mapping from a group to itself is an
automorphism.\\ (iv)  Any two isomorphic groups have the same
cardinality. \eos \bexe\label{ec35} (i) Is a monomorphism of a
group into itself an automorphism? \\ (ii) Is a homomorphism of a
group onto itself an automorphism ? \eexe \bl\label{l35} If
$\varphi: G \to H$ and $\psi: H \to K$ are two group
homomorphisms, then $\psi\varphi: G \to K$ is a group
homomorphism. Further, if $\varphi$ and $\psi$ are monomorphism (epimorphism or isomorphisms),
then so is $\varphi \psi $. \el \pf For $x,y \in G$, \bea
\psi\varphi (xy) & = & \psi (\varphi (xy))\\
                  & = &  \psi (\varphi (x)\varphi (y))\\
                  & = &  \psi\varphi (x)\psi\varphi (y)
\eea Hence $\psi\varphi$ is a group homomorphism from $G$ to $K.$
Further, if $\psi,\varphi$ are one-one (onto), then clearly  $\psi\varphi$
 is one-one (onto). Hence the rest of the statement is immediate.  $\blacksquare$
 \bx\label{ex31} Let $G = (\Z,+)$
and $H = (\C^* \cdot)$.  Then for the map \bea
 \varphi : G & \to & H\\
    n & \mapsto & e^{\pi in}
\eea \bea
           \varphi (m+n) & = & e^{\pi i(m+n)}\\
                      & = & e^{\pi im}\cdot e^{\pi in}\\
                      & = & \varphi (m)\cdot \varphi (n)
\eea Hence $\varphi$ is a homomorphism. Further $\abs{\varphi (m)}
= 1$ for all $m\in \Z$, hence $\varphi$ is not onto. The image of
$\varphi$ is the subgroup $\{1, -1\}$ of $\C.$ \ex \bx\label{ex32}
Let $G = (\Z,+)$, and $H = (k\Z,+)$, where $k$ is a fixed
integer. Then \bea
 \varphi : G & \to & H\\
    n & \mapsto & kn
\eea is an isomorphism. \ex \bx\label{ex33} Let $G$  be the
multiplicative group $\{1,-1\}$ and $H = \Z_{2}$ $=
\{0,1\}$ be  the additive group of integers $mod\; 2$.
Then \bea
 \varphi : H & \to & G\\
      0 & \mapsto & 1 \\
     1 & \mapsto & -1
\eea is an isomorphism. \ex  \bx Let $G =\Z_{n}$, the additive group of integers
modulo $n$ (Example \ref{ex1.7}) and $H =\Z/n\Z$ the factor group of $\Z = (\Z,+)$ by
 the subgroup $n\Z$. Then the map $$\begin{array}{ccccl}
   \varphi & : & G & \rightarrow & H \\
    &  & k & \mapsto & k+n\Z \
 \end{array}$$ is an isomorphism.  \\  If $k,l \in \Z_{n}$ and $k+l < n$, then
   $$\begin{array}{ccl}
     \varphi(k+_{n}l) & = & k+l+n\Z \\
      & = & (k+n\Z)+(l+n\Z) \\
      & = & \varphi(k)+\varphi(l). \\
   \end{array}$$ Moreover, if $k+l \geq n$, then  $$\begin{array}{ccl}
     \varphi(k+_{n}l) & = & \varphi(k+l-n) \\
      & = & k+l-n+n\Z \\
      & = & (k+n\Z)+(l+n\Z) \\
      & = & \varphi(k)+\varphi(l) \\
   \end{array}$$  Thus $\varphi$ is a homomorphism. Further, for any $m \in \Z$, by
    division algorithm, we can write  $$\begin{array}{rclcl}
      m & = & nq+r, &  & 0\leq r\leq n \\
      \Rightarrow \; \; \; \;\; m+n\Z & = & r+n\Z &  & \mbox{since }nq \in n\Z. \\
           &  = &  \varphi(r) &  &  \\
\end{array}$$ Hence $\varphi$ is onto. Now, as $\circ(G) = n = \circ(H).$ The map
$\varphi$ is clearly one-one and hence is an isomorphism. \ex
 \bx\label{ex34} We have a general form
of the example \ref{ex33}. Let $G = \C_{n}$ be the multiplicative
group of $n^{th}$ roots of unity in $\C$, and $H = \Z_n$, then
\bea
 \varphi : \Z_n & \to & G\\
     m &  \mapsto & e^{\frac{2\pi mi}{n}} \quad (0\leq m < n)
\eea is an isomorphism. \ex  \pf If for $k,l \in\Z_{n}$, $k+l<n$, then
$$\begin{array}{ccl}
  \varphi(k+_{n}l) & = & e^{2\pi i(k+l)/n} \\
   & = & e^{2\pi ik/n}\cdot e^{2\pi il/n} \\
   & = & \varphi(k)\varphi(l)
\end{array}$$ Further, if $k+l\geq n$,  $$ \begin{array}{ccl}
  \varphi(k+_{n}l) & = & \varphi(k+l-n) \\
   & = & e^{2\pi ik/n}\cdot e^{2\pi il/n} \\
   & = & \varphi(k)\varphi(l)
\end{array}$$   Hence $\varphi$ is a homomorphism. It is easy to see that $\varphi$
 is onto and hence, since $\circ(G) = \circ(H) = n$, $\varphi$ is one-one.
  Thus $\varphi$ is an isomorphism.  \bx\label{ex35} Let $G = \C^*$ be the
multiplicative group of non-zero complex numbers and $H = S^1 =
\{z\in \C \st \abs{z} = 1\}$, the subgroup of $G$ of complex
numbers with modules 1.  Consider \bea
 \alpha : \C^*  & \to & S^1\\
        z & \mapsto & \frac{z}{\bar{z}}
\eea \bea
 \beta : \C^* & \to & S^1\\
       z & \mapsto & \frac{z}{\abs {z}}
\eea Clearly $\alpha$ and $\beta$ are homomorphisms. Moreover ,for
any $z = \cos\theta + i \sin\theta$ in $S^1$, we have $$ \beta(z)
= z,\mbox{ and } \alpha(\cos\frac{\theta}{2} + i
sin\frac{\theta}{2}) = z$$ Thus $\alpha$ and $\beta$ are
surjective. \ex \bx\label{ex36} Let $G$ be an Abelian group. Then
\bea
 \theta : G & \to & G\\
      x & \mapsto & x^{-1}
\eea is an automorphism  of G. \ex \bx\label{ex37} Let $H$ be a
normal subgroup of a group $G$. Then for the map \bea
 \eta : G & \to & G/H\\
      g & \mapsto & gH
\eea $\eta(g_1g_2) = g_1g_2H = (g_1H)(g_2H) = \eta(g_1)\eta(g_2)$
i.e., $\eta$ is a homomorphism from  $G$ to $G/H$. This
homomorphism is normally referred to as {\bf{natural}} or {\bf {canonical
homomorphism}}. \ex \bx\label{ex38} Given any two groups $G$ and $H$
, the mapping \bea
 \tau : G & \to & H\\
      g & \mapsto & e_H
\eea which maps all elements of $G$ to the identity element of H
is a homomorphism. This is called the {\bf {trivial homomorphism}}. \ex
\bx\label{ex39} Let $G = (\C^*,\cdot)$ and $H = Gl_2(\R)$. Then
$\theta: G \to H$ defined by $$ a+ib \mapsto
\left[\begin{array}{cc}
                          a & b\\
                          -b & a
                       \end{array}\right]$$
is a monomorphism.
 \ex
\bexe\label{ec35a} Let $G =S_{n}$ and $H = Gl_{n}(\R)$. Then the
map $$\begin{array}{ccccc}
  \theta & : & S_{n} & \rightarrow & Gl_{n}(\R) \\
   &  & \sigma & \mapsto & [e_{\sigma (1)},e_{\sigma (2)},\cdots,e_{\sigma (n)}]
\end{array}$$
where $e_{i}$ is the column vector $(0,0,\cdots,1^{i
th},0,\cdots,0)^{t}$ is a monomorphism. \eexe    \bx\label{ex310}
Let $G = (\R,+)$ and $H = (\C^{*},\cdot)$. Then $t \mapsto e^{it}$
is a homomorphism from $G$ into $H.$ \ex \bx\label{ex311} Let $G =
(\R^*,\cdot)$ and $H = \{\pm 1\}$ with multiplication. Then $$ r
\mapsto \left\{
\begin{array}{cc}
                      +1 & \mbox { if } r > 0\\
                      -1 & \mbox { if } r < 0
                      \end{array}\right. $$
is a homomorphism from $G$ onto $H.$ \ex \bexe\label{ec32} Show
that there exists no non-trivial homomorphism  from $(\R,+)$ to
$(\Z,+)$. (Hint : If there exists a non-trivial homomorphism then
there exists an epimorphism.) \eexe \bexe\label{ec33} Show that
the groups $(\R,+)$ and the multiplicative group $(\R^+,.)$ are isomorphic where $\R^+$ is
the set of positive reals \eexe \bexe\label{ec34} Show that
$(\R,+)$ is not isomorphic to $(\R^*,\cdot)$ \eexe \bx Let $K$ be a field, and
$G=Gl_n(K)$. Then $$\begin{array}{rl}
  \phi \; :\;G\;\rightarrow & G \\
  P\;\mapsto & (P^t)^{-1}
\end{array}$$ is an automorphism, where $P^t$ denotes the transpose of $P$.\ex
\bexe Let $K$ be a field, and $\alpha :Gl_n(K)\rightarrow K^{*}$ (the multiplicative group
of non-zero elements of $K$), be any homomorphism. Prove that  $$\begin{array}{rcl}
  \psi\;:\;Gl_n(K) & \rightarrow & Gl_n(K) \\
  P & \mapsto & \alpha(P)P
\end{array}$$ is an automorphism. \eexe
 \bd\label{d35}Let $\varphi : G \to H$  be a homomorphism of groups.Then the set
$\{x\in G\st \varphi (x) = e_H \}$ is called the kernel of
$\varphi$ and is denoted by $\ker\varphi$. \ed \bt\label{t31} For
a group homomorphism  $\varphi : G \to H$, the $\ker\varphi$ is a
normal subgroup of $G$. Further $\ker\varphi = \{e\}$ if and only
if $\varphi$ is one-one. \et \pf For the identity element $e_G\in
G$, $\varphi (e_G) = e_H$, the identity of $H$. Hence $\ker\varphi
\neq \emptyset$. Let $x, y \in \ker\varphi$. Then, \bea
           \varphi (xy^{-1}) & = & \varphi (x) \varphi (y^{-1})\\
                         & = & e_H \quad (\mbox{ since }\varphi (x) = e_H = \varphi (y)) )\,
\eea $\Rightarrow  \quad xy^{-1}\in \ker\varphi$ for all $x, y \in
\ker\varphi$.   \\  Therefore, $\ker\varphi$ is a subgroup of $G$.
Further for any $g\in \ker\varphi$, and $x\in G$, \ba{rl}
     & {\begin{array}{rcl}
          \varphi (x^{-1}gx) & = & \varphi (x^{-1})\varphi (g)\varphi (x)\\
                    & = & \varphi (x^{-1})\varphi (x) \quad \mbox{ since } \varphi(g) = e_H\\
                    & = & {\varphi(x)}^{-1}\varphi(x) \\
                    & = & e_{H}

       \end{array}}\\
\Rightarrow  & x^{-1}gx \in \ker\varphi\\ \Rightarrow  &
x^{-1}(\ker\varphi) x \subset \ker\varphi \mbox { for all }x\in G.
\ea Hence $\ker\varphi \lhd G$. Next, for $x, y \in G$, we have
\bea
          \varphi (x) = \varphi (y) & \Leftrightarrow  & \varphi (x)\varphi (y)^{-1} = e_H\\
                              & \Leftrightarrow  & \varphi (x)\varphi (y^{-1}) = e_H\\
                              & \Leftrightarrow  & \varphi (xy^{-1}) = e_H\\
                              & \Leftrightarrow  & xy^{-1}\in \ker\varphi
\eea Hence if $\ker\varphi= e_G,$ then, \ba{rl}
                   &  \varphi (x) = \varphi (y)\\
\Rightarrow        &  xy^{-1} = e_G\\ \Rightarrow        &   x = y
\ea Therefore $\varphi$ is one-one.\\ Conversely, let $\varphi$ be
one-one. If $x, y \in \ker\varphi$, then $xy^{-1}\in \ker\varphi$ as
$\ker\varphi$ is a subgroup. Hence \ba{rl}
              &  \varphi (xy^{-1}) = e_H\\
\Rightarrow   & \varphi (x)\varphi (y)^{-1} = e_H\\ \Rightarrow
&  \varphi (x) = \varphi (y)\\ \Rightarrow   &     x = y  \quad
\mbox { as }\varphi \mbox { is one-one.} \ea Therefore
$\ker\varphi = e_G$ (use: $e_G \in \ker\varphi$) Hence the result
follows.     $\blacksquare$ \bt\label{t32} Let $G,H$ be two groups and let  $\varphi
: G \to H$ be a homomorphism of groups. Then $\varphi$ is an
isomorphism if and only if there exists a homomorphism $ \psi : H
\to G$ such that  $\psi\varphi = Id_G$ and $\varphi\psi = Id_H$.
\et \pf Let $\varphi$ be an isomorphism. Then it is
one-one and onto. Hence  $\varphi$ admits inverse map and $$
\varphi^{-1}(y) = x \mbox{ if } \varphi (x) = y$$ Let $\varphi
(x_1) = y_1$,  $\varphi (x_2) = y_2.$  Then \ba{rl}
               & \varphi (x_1x_2) = \varphi (x_1)\varphi (x_2) = y_1y_2\\
\Rightarrow    & \varphi^{-1}(y_1y_2) = x_1x_2 = \varphi^{-1}(y_1)
\varphi^{-1}(y_2) \ea Therefore $\varphi^{-1}$ is a homomorphism
of groups and clearly $$\varphi \varphi^{-1} = Id_H \mbox { and }
\varphi^{-1} \varphi = Id_G.$$ Conversely, let there
exists a homomorphism  $\psi : H \to G$
 such that $\varphi \psi = Id_H \mbox { and } \psi \varphi = Id_G$. Then
\ba{rl}
               &    \varphi (x_1) = \varphi (x_2) \quad (x_1, x_2\in G)\\
\Rightarrow    & \psi\varphi (x_1) = \psi\varphi (x_2)\\
\Rightarrow    & x_1 = x_2 \ea Hence $\varphi$ is one-one. Next,
if $y \in H$, then \ba{rl}
                &   y = \varphi\psi (y)  \quad(\varphi \psi= Id_H)\\
\Rightarrow     & \varphi \mbox { is onto.} \ea Thus, as $\varphi$
is one-one, onto, homomorphism, it is an isomorphism.     $\blacksquare$
\bt\label{t301a} If $\varphi : G \rightarrow H$ is a group
 homomorphism, then  \\  (i) For any subgroup $A$ of $G$,
 $\varphi(A)$ is a subgroup of $H$. Further if  $\varphi$ is onto and $A
 \lhd G$, then $\varphi(A)$ is a normal subgroup of $H$.  \\  (ii)
 Inverse image of a subgroup of $H$ is a subgroup of $G$. Moreover, inverse image of a normal subgroup
 of $H$ is normal in $G$.\et   \pf (i) Let $x,y \in A$, then
 $xy^{-1} \in A$, and
 $$\begin{array}{ccccc}
   \varphi(x)(\varphi(y))^{-1} & = & \varphi(x)\varphi(y^{-1}) &  &  \\
    & = & \varphi(xy^{-1}) & \in & \varphi(A) \
 \end{array}$$
 Hence $\varphi(A)$ is a subgroup of $H$.  \\  Next, let $\varphi$ is onto and $A \lhd G$.
 As $\varphi$ is onto, any element of $H$ is of the from $\varphi(g)$
 for some $g \in G$. Hence for any $a \in A$,
 $$\begin{array}{ccc}
  \varphi(g)\varphi(a)(\varphi(g))^{-1} & = & \varphi(g)\varphi(a)\varphi(g^{-1}) \\
    & = & \varphi(gag^{-1}) \in\varphi(A) \
 \end{array}$$   Therefore $\varphi(A) \lhd H$.  \\
 (ii) Let $B$ be a subgroup of $H$ and let $A = \varphi^{-1}(B)$. Clearly $A \neq \emptyset$ since $e \in A$ If
 $x,y \in A$, then $\varphi(x),\varphi(y) \in B$. As $B$ is a subgroup
 of $H$,we have
 $$\begin{array}{ccc}
   \varphi(xy^{-1}) & = & \varphi(x)(\varphi(y))^{-1} \in B \\
   \Rightarrow  & xy^{-1} \in A &  \\
   i.e. ,   &xy^{-1} \in A & \forall  x, y \in A \
 \end{array}$$

Hence $A$ is a subgroup of $G$.Further, let $B \lhd H$. If $x \in
A$ and $g \in G$, then $$\varphi(gxg^{-1}) =
\varphi(g)\varphi(x)(\varphi(g))^{-1} \in B$$   $$\Rightarrow
  gxg^{-1} \in A \mbox{ for all } g \in G, x \in A$$   Hence $A \lhd H$.     $\blacksquare$
\bco\label{c301a} Let $K$ be a normal subgroup of a group $G$.
Then any subgroup of $G /K$ is of the form $A /K$ where $A$ is a
subgroup of $G$ containing $K$. Further $A / K  \lhd  G / K$ if
and only if $A \lhd G$. \eco     \pf Consider the canonical
homomorphism $$\begin{array}{ccccc}
  \eta & : & G & \rightarrow & G / K\\
   &  & x & \mapsto & xK
\end{array}$$
Clearly $\eta$ is an onto homomorphism. Further, for any $x \in K$,
$$\begin{array}{rrcc}
   & \eta(x) = K &  & (\mbox{identity element of }G/K) \\
  \Leftrightarrow & xK = K &  &  \\
  \Leftrightarrow & x \in K &  &  \\
  \mbox{Thus } & K = \ker \eta . &  &
\end{array}$$ Hence by the
theorem, for any subgroup $B$ of $G /K$, $A =
\eta^{-1}(B)$ is a subgroup of $G$ containing $K.$ As $\eta$ is
onto, $$\begin{array}{ccrcc}
 \eta(A) & = & \eta(\eta^{-1}(B)) & = & B\\
 \Rightarrow &  &\{ xK | x \in A \} & = & B \\
  i.e.,  & & A / K & = & B
\end{array}$$
Further, if $A$ is a subgroup of $G$ containing $K$ then by the
theorem $$\eta(A) = \{ xK | x \in A \} = A / K$$ is a
subgroup of $G / K$. Finally by the theorem it is immediate
 that $A / K \lhd G / K$ if and only if $A \lhd G$.     $\blacksquare$
\bl\label{l301a} Let $H$ be a normal subgroup of a group $G$. Then
$G / H$ is simple if and only if $H$ is a maximal normal subgroup
of $G$.\el \pf  Let $G / H$ be a simple group. Then $\circ(G / H)
> 1$ and $G / H$ has no proper normal subgroup. Hence by the
theorem \ref{t301a},there dose not exist any $A \lhd G$ such that
$H \subset A$, $A \neq H$. Therefore $H$ is a maximal normal
subgroup of $G$. Converse is also immediate from the theorem
\ref{t301a} as any normal subgroup of $G / H$ is of the form $A /
H$ where $A \lhd G$, $H \subset A$.     $\blacksquare$  \bl\label{l302a} Any
homomorphic image of a simple group $G$ is either identity group or is
isomorphic to $G$.\el \pf  Let
  $$\begin{array}{ccccc}
    \varphi & : & G & \rightarrow & H \
  \end{array}$$ be an onto homomorphism and let $K = ker\varphi$.
  then as $K \lhd G$, $K = G$ or $K = \{e\}$. Hence either
  $\varphi(G) = \{e\}$ or $\varphi$ is an isomorphism onto $H$.     $\blacksquare$

\bt\label{t33} {\bf(FIRST ISOMORPHISM THEOREM)} Let $G$ and $H$ be two
groups and let $\varphi: G \to H$ be a homomorphism. If
$\ker\varphi = K$, then $K$ is a normal subgroup of $G$ and there
exists a unique monomorphism $\bar{\varphi} : G/K \to H$ such that
$\bar{\varphi}\eta = \varphi$ where $\eta : G \to G/K$ is
the  canonical homomorphism $\eta(x) = xK$ for all $x\in G.$ \et
\pf By the theorem \ref{t31}, $K$ is a normal subgroup of $G$.
Consider the diagram \ba{cll}
                      G & \stackrel{\varphi}{\to} & H\\
                      \eta \downarrow &  &\\
                     G/K   &  &
\ea of groups and homomorphisms where $\eta$ is the canonical homomorphism.
As $\eta$ is onto there is
precisely one  map (say) $\bar{\varphi}$ from $G/K$ to $H$ i.e.,
$\bar{\varphi}(xK) = \varphi (x)$ for which $\bar{\varphi}\eta =
\varphi.$ Thus to complete the proof it is sufficient to show that
$\bar{\varphi}$ is a monomorphism. Let $xK, yK \in G/K$. Then \bea
         \bar{\varphi}(xK.yK) & = & \bar{\varphi}(xyK)\\
                           & = & \varphi (xy)\\
                           & = & \varphi (x)\varphi (y)\\
                           & = & \bar{\varphi}(xK)\bar{\varphi}(yK)
\eea Therefore $\overline{\varphi}$ is a homomorphism. Next, let \ba{rl}
                   & \bar{\varphi}(xK) = \bar{\varphi}(yK)\\
\Rightarrow        & \varphi (x) = \varphi (y)\\ \Rightarrow &
\varphi (xy^{-1}) = e\\ \Rightarrow        & xy^{-1}\in
\ker\varphi = K\\ \Rightarrow        &  xK = yK. \ea Thus
$\overline{\varphi}$ is one-one and hence a monomorphism.     $\blacksquare$
 \bco\label{c43a}
Let $$\begin{array}{ccccc}
  \varphi & : & G & \rightarrow & H
\end{array}$$
be an epimorphism and let $K = \ker \varphi$. then $G / K$ is
isomorphic to $H$.\eco  \pf  It is clear that the map
$\overline{\varphi}$ in the theorem \ref{t33} is the isomorphism.     $\blacksquare$
\bexe Let $$\phi:G \rightarrow H$$ be a group epimorphism and let $B$ be a normal subgroup
 of $H$. Prove that for $\phi^{-1}(B)=A$, $G/A$ is isomorphic to $H/B$.\eexe
\bt\label{t34.01} Let $H,K,L$ be subgroup of a group $G$ such that $K\lhd H$ and $L\lhd G$,
 then $HL/KL$ is isomorphic to $H/K(L\cap H)$. Thus in particular if $K = \{e\}$,
 $HL/L$ is isomorphic to $H/L\cap H$. \et \pf Since $L\lhd G$, $HL = LH$ and $KL = LK.$
 Thus by the theorem \ref{t29}, $HL$ and $KL$ are subgroup of $G$. Further, as $K\subset H$
 $KL \subset HL$. Since $L\lhd G$ and $K\lhd H$, for any $x \in H$, $y \in L$, we
have $$\begin{array}{cclc}
  (xy)KL & = & x(yL)K &  \\
   & = & xLK & \mbox{ since }yL = L \\
   & = & xKLy &  \\
   & = & KL(xy) &
\end{array}$$ Hence $KL\lhd HL$. Now, for $H\subset HL$, consider the map:
  $$\begin{array}{rcl}
    \varphi : H & \rightarrow & HL/KL \\
    x & \mapsto & xKL \
  \end{array}$$  For any $x_{1},x_{2} \in H$, $$\begin{array}{ccl}
    \varphi(x_{1}x_{2}) & = & x_{1}x_{2}KL \\
     & = & (x_{1}KL)(x_{2}KL) \\
     & = & \varphi(x_{1})\varphi(x_{2}) \
  \end{array}$$ Hence $\varphi$ is a homomorphism. Further, as above, for any
   $x \in H$,$y \in L$, $$\begin{array}{ccl}
     xyKL & = & xyLK \\
      & = & xLK \\
      & = & xK \\
      & = & \varphi(x) \\
   \end{array}$$  Thus $\varphi$ is an epimorphism. Now, for any $x \in H$,
     $$\begin{array}{crclc}
        & \varphi(x) & = & KL &  \\
       \Rightarrow & xKL & = & KL &  \\
       \Rightarrow & x & = & ab & \;\;\;\;\;\; (a \in K, b \in l) \\
       \Rightarrow & a^{-1}x & = & b \in H\cap l &  \\
       \Rightarrow & x & \in & K(H\cap L) &  \\
\end{array}$$ Conversely if $a \in K$ and $c \in H\cap L$, then $ac \in HK = H$,
 and $$\begin{array}{ccc}
   \varphi(ac) & = & acKL \\
    & = & (aK)(cL) \\
    & = & KL \
 \end{array}$$ Consequently $\ker \varphi = K(H\cap L)$ and by the theorem \ref{t33}
 $H/K(H\cap L)$ is isomorphic to $HL/KL$. The rest of the statement is now immediate.     $\blacksquare$
 \bl\label{l34.01l} Let $C = gp\{a\}$ be a multiplicative cyclic group. Then $C$ is isomorphic
 to the factor group $\Z/n\Z$ of the group $\Z = (\Z,+)$ where $n = \circ(C)$ if $\circ(C) = n< \infty$
  and 0 otherwise. In particular any two cyclic group of same order are isomorphic.\el
  \pf Consider the map:  $$\begin{array}{rcl}
    \varphi : \Z & \rightarrow & C \\
    m & \mapsto & a^{m} \\
  \end{array}$$  Then  $$\begin{array}{ccl}
    \varphi(m+n) & = & a^{m+n} \\
     & = & \varphi(m)\varphi(n) \\
  \end{array}$$ Hence $\varphi$ is a homomorphism. Clearly $\varphi$ is onto. Now, let
  $\circ(C) = n < \infty$. Then $\circ(a) = n.$  \\  For any $m \in \Z$,  $$\begin{array}{crcc}
     & \varphi(m) = e &  &  \\
    \Leftrightarrow & a^{m} = e &  &  \\
    \Leftrightarrow & n|m &  & \mbox{ Theorem \ref{t25} }(ii) \\
  \end{array}$$  Therefore $\ker \varphi = n\Z$ and by the theorem \ref{t33}, $\Z/n\Z$ is isomorphic
   to $C$. Next, let $\circ(C) = \infty$. Then $\circ(a) = \infty$. Hence  $a^{m} = \varphi(m) = e$
    if and only if $m = 0$. Thus $\ker \varphi = (0)$ and as above $\Z$ is isomorphic to $C$. Then
     rest of the statement follows by the lemma \ref{l35}.     $\blacksquare$
     \bt\label{t34} Let $\varphi: G \to H$ be a group homomorphism.
      If $x\in G$ is an element of
order $n$, then order of  $\varphi (x)$ divides order of $x.$
Further if $\varphi$ is an isomorphism then $\circ (x) = \circ
(\varphi (x)).$ \et \pf As $\circ (x) = n$, we have \ba{rl}
               &   \varphi (x^n) = \varphi (e_G) = e_H\\
\Rightarrow    &  (\varphi (x))^n = e_H\\ \Rightarrow    & \circ
(\varphi (x)) \mid n = \circ (x) \;\;\;\;\;\;\;\; \mbox{ Theorem \ref{t25} }(ii). \ea
If $\varphi$ is an
isomorphism then $\varphi^{-1}$ exists and is a homomorphism from
$H$ to $G.$ Thus, as above  \bea
            \circ (\varphi^{-1}(\varphi (x)) & \mid & \circ (\varphi (x))\\
    i.e.,\quad \circ (x) & \mid & \circ (\varphi (x)).
\eea Hence, if $\varphi$ is an isomorphism $$ \circ (\varphi (x)) = \circ
(x).     \blacksquare$$ \bco\label{l000l} Let $\varphi : G \rightarrow H$ be a
homomorphism of groups. If $x \in G$ is an element of finite order
then $o(\varphi (x)) < \infty$. \eco  \pf It is clear from the proof
of the theorem.     $\blacksquare$ \bl\label{l31} Let $G$ be a finite group and
$\varphi
: G \to H$, an epimorphism. Then $\circ (H) \mid
\circ (G).$ \el \pf Let $K$ be the kernel of $\varphi$. Since
$\varphi$ is onto, by the corollary\ref{c43a}, $G/K$ is isomorphic
to $H$. Hence $\circ (G/K) = \circ (H).$ Further, by the corollary
\ref{co213a}, we have \ba{rl}
              &   \circ (H) = \circ (G/K) = \frac{\circ (G)}{\circ (K)}\\
\Rightarrow   &   \circ (H).\circ (K) = \circ (G). \ea Hence the
result follows.     $\blacksquare$ \bt\label{t35} The set $AutG$ of all automorphisms
of a group $G$ is a group  with respect to composition operation.
\et \pf Clearly $I$, the identity map over $G$, is an automorphism
of $G$ and $\varphi  I = I\varphi =\varphi$ for all
$\varphi \in AutG.$ Further, for any $\varphi \in AutG$,
$\varphi^{-1}$ is an automorphism (Theorem \ref{t32}) and $\varphi
\varphi^{-1} = \varphi^{-1} \varphi = I.$ Now, let $\psi, \varphi
\in AutG$. Then as $\varphi$ and $\psi$ are one-one, onto, so is
$\varphi  \psi$.  Moreover, by the lemma \ref{l35}, $\varphi
 \psi$ is a homomorphism. Hence composite of two automorphisms is
an automorphism. As composition of maps is associative, we
conclude $AutG$ is a group.     $\blacksquare$ \br The group $Aut(G)$ in the theorem
is called the automorphism group of $G.$ \er \bt\label{t000t} Let
$G$ be an additive group and $\varphi, \psi$ be two endomorphisms
of $G$. Consider the maps : $$f :  G  \rightarrow  G \mbox{ and }
         h :  G  \rightarrow  G$$
       such that $f(x) = x - \varphi \psi (x)$, $h(x) = x - \psi
       \varphi (x)$. Then $f$ is onto (one-one) if and only if $h$
       is onto (one-one).\et
\pf First of all, note that for any $x,y \in G$,$h(x) - h(y) = h(x
- y)$. Hence $h$ is an endomorphism. Similarly $f$ is an
endomorphism.Now, let $f$ be onto $G$. Then for any $x \in G$,
there exist $y \in G$  such that $$\begin{array}{cccl}
  & x & = & y  -  \varphi \psi (y) \\
  \Rightarrow & \psi(x) & = & \psi(y)  -  \psi \varphi \psi(y) \\
   &  & = & h(\psi(y))  \in  Im(h) \\
  \Rightarrow & \psi(G) & \subset & Im(h)
\end{array}$$ Further, for any $x \in G$,$$\begin{array}{ccl}
  h(x) & = & x - \psi \varphi(x)  \in  Im(h) \\
  \Rightarrow & x & \in  Im(h)
\end{array}$$ \hspace{2in} since $\psi(\varphi(x)) \in Im(h)$ and $Im(h)$ is a subgroup
(Theorem \ref{t32} (i))
 $$\begin{array}{cccc}
   \Rightarrow & G & = & Im(h) \
 \end{array}$$ Thus $h$ is onto. Similarly we can prove the converse. \\
 For the other part, let $h$ be one-one, and let for some $x \in G$,
 $$\begin{array}{crclcc}
   & h(x) & = & 0 &  &  \\
   \Rightarrow & x & = & \varphi \psi(x) &  &  \\
   \Rightarrow & \psi(x) & = & \psi \varphi \psi(x) &  &  \\
   \Rightarrow & h(\psi(x)) & = & 0 &  &  \\
   \Rightarrow & \psi(x) & = & 0 &  & \mbox{since h is one-one}    \\
   \Rightarrow & x & = & \varphi \psi(x) = 0 &  &  \mbox{ homomorphism}  \
 \end{array}$$ Hence $f$ is one-one.The rest of the statement is clear.
      \bt\label{t36} Let $G$ be a group . For any $x\in G$ , the map
\bea
               \tau_x : G & \to &G\\
                        g & \mapsto &xgx^{-1}
\eea is an automorphism. The set $I_G = \{ \tau_x \st x\in G \}$
is a normal subgroup of $Aut (G).$ \et \pf For $x,g_1,g_2\in G$ ,
we have \bea
       \tau_x (g_1g_2) & = & x(g_1g_2)x^{-1}\\
                       & = & x(g_1)x^{-1} x(g_2)x^{-1}\\
                       & = & \tau_x(g_1) \tau_x(g_2)
\eea Hence $\tau_{x}$ is an endomorphism of $G.$ Further for any
$g\in G$, \bea
         (\tau_x \ \tau_{x^{-1}})(g) & = & \tau_{x} (x^{-1}gx)\\
                                          & = & x(x^{-1}gx)x^{-1} = g
                                          \eea
Similarly, we can prove that
   \bea
       (\tau_{x^{-1}}  \tau_x)(g) & = & g \eea
       Hence \bea
        (\tau_x \ \tau_{x^{-1}}) & = & Id =
(\tau_{x^{-1}} \tau_x) \eea Therefore, by the theorem \ref{t32} ,
$\tau_x$ is an automorphism of $G$. We also note that
$(\tau_x)^{-1} = \tau_{x^{-1}}$. Now for $\tau_x, \tau_y\in I_G$,
we have for any $g\in G$, \bea
 (\tau_x \tau_{y^{-1}})(g) & = & \tau_x (y^{-1}gy)\\
                                 & = & x(y^{-1}gy)x^{-1}\\
                                 & = & {(xy^{-1})}^{-1}g(xy^{-1})^{-1}\\
                                 & = & \tau_{xy^{-1}} (g)
\eea Hence $$\tau_{x}  \tau_{y^{-1}} = \tau_{xy^{-1}} \in I_{G}$$  Thus by
the theorem \ref{t1}, $I_G$ is a subgroup of $Aut (G).$ Finally,
let $\psi \in Aut (G)$, and $\tau_x \in I_G.$ Then for any $g\in
G,$ \bea
   (\psi  \tau_x \psi^{-1})(g) & = & \psi \tau_x(\psi^{-1} (g))\\
                                          & = & \psi (x \psi^{-1} (g)x^{-1})\\
                                          & = & (\psi(x))g\psi(x)^{-1}\\
                                          & = & \tau_{\psi(x)}(g)
\eea Thus, \ba{rl}
            & \psi \tau_x \psi^{-1} = \tau_{\psi (x)} \in I_G\\
\Rightarrow & \psi I_G \psi^{-1} \subset I_G  \mbox{    for all }\psi \in Aut G
\ea
Hence, by the theorem \ref{t214},  $ I_G \lhd Aut(G).$     $\blacksquare$
\br The subgroup $I_G$ is called the group of
{\bf inner automorphism} \rm of $G.$ \er
\bl\label{l37.05l} For any group $G$, the map  $$\begin{array}{rcl}
  \tau : G & \rightarrow & I_{G} \\
  x & \mapsto & \tau_{x}
\end{array}$$ is homomorphism with $\ker\tau = Z(G)$, the center of $G$. \el
 \pf For any $x,y,g \in G$,  $$\begin{array}{crcl}
    & \tau_{xy}(g) & = & (xy)g(xy)^{-1} \\
    &  & = & x(ygy^{-1})x^{-1} \\
    &  & = & \tau_{x}\tau_{y}(g) \\
   \Rightarrow & \tau_{(xy)} & = & \tau_{x}\tau_{y} \\
   \Rightarrow & \tau(xy) & = & \tau(x)\tau(y) \
 \end{array}$$  Hence $\tau$ is a homomorphism. Now, let for $x \in G$, $\tau_{x}$
  be identity automorphism. Then for every $g \in G$,$\tau_{x}(g) = g$. Note that
    $$\begin{array}{crcl}
       & \tau_{x}(g) & = & g \\
      \Leftrightarrow & xgx^{-1} & = & g \\
      \Leftrightarrow & xg & = & gx \\
\end{array}$$  hence, $x \in \ker\tau$ if and only if $x \in Z(G)$. Thus the result
 is proved.     $\blacksquare$  \bt\label{t37} Let $C_n = gp\{x\}$
be a finite cyclic group of order of $n.$ Then the group
$Aut(C_n)$ has order ${\bf{\Phi}} (n)$, where ${\bf{\Phi}}$ is the
 Euler's function.
\et \pf For any integer $k\geq 1$, define \bea
               \alpha_k: C_n & \to & C_n\\
                       x^m & \mapsto & x^{mk}
\eea For $x^a,x^b\in C_n$, we have  \ba{cl}
        \alpha_k(x^a.x^b)  & = \alpha_k(x^{a+b})    \\
                           & = x^{(a+b)k}  \\
                           & = x^{ak}x^{bk}  \\
                           & =\alpha_k(x^a)\alpha_k(x^b) \
               \ea
 Hence $\alpha_k$ is an endomorphism and
$$\alpha_n (x^a) = x^{na} = e$$ since $\circ (x) = n.$ Thus $\alpha_n$
is the trivial endomorphism. Note that for $k, l\in \Z$, $k\geq
1$, $l\geq 1$, \ba {rl}        &  \alpha_k = \alpha_l\\
\Leftrightarrow &  \alpha_k (x) = \alpha_l (x)\mbox{  for all }x \in C_{n}\\ \Leftrightarrow & x^k
= x^l\\ \Leftrightarrow & x^{k-l} = e\\ \Leftrightarrow & n \mid
(k-l) \ea Therefore $\alpha_1, \alpha_2,...,\alpha_n$ are all distinct
endomorphisms of $G$. Clearly image of $\alpha_k =
gp\{x^k\}.$ Hence $\alpha_k$ is onto (equivalently, one-one) if and only if $C_n =
gp\{x^k\}$  i.e., $\circ (x^k) = n.$ Now ,by the corollary \ref{l1.23}
 $\circ (x^k) = n$ if and only if $(n,k)
= 1.$ Therefore $\alpha_k$, $1\leq k \leq n$, $(k,n) = 1$, are
automorphisms of $C_n$. Finally, let $\beta$ be any automorphism of
$C_n$. Then \ba{rl}
                  & \beta(x) = x^k \mbox { for some} 1\leq k\leq n\\
\Rightarrow        &  \beta (x^a) = x^{ak} \mbox { for all }a\in
\Z\\ \Rightarrow       &  \beta = \alpha_k \mbox { for some } 1\leq
k \leq n. \ea Hence $Aut(C_n) = \{ \alpha_k \st 1\leq k\leq n, (n,k)
= 1\}.$ This proves order of $Aut(C_n)$ is ${\bf{\Phi}} (n).$     $\blacksquare$
\bl\label{l36}
Let $C = gp\{x\}$ be an infinite cyclic group. Then Aut(C) has
 order $2.$
\el
\pf Let $k\in \Z$. Then the map \bea
                 \alpha_k : C & \to & C\\
                          x^m & \mapsto & x^{mk}
\eea is an endomorphism of $C$ as in the theorem\ref{t37}.
Further, if $\beta$ is an endomorphism of $C$, then if \ba{rl}
                  & \beta (x) = x^k \;\;\;\;\;\;\;\;(k \in \Z)   \\
\Rightarrow       & \beta (x^a) = x^{ak} \mbox{ for all }a\in
\Z\\ \Rightarrow       & \beta = \alpha_k \ea Now, if $\alpha_m$ is
an automorphism, there exists $n\in \Z$ such that \ba{rl}
                    & \alpha_m \alpha_n = \alpha_n \alpha_m = Id.\\
\Rightarrow         & \alpha_n \alpha_m (x) = x\\ \Rightarrow
& {(x^m)}^{n} = x\\ \Rightarrow         & x^{mn-1} = e\\
\Rightarrow         & mn-1= 0\\ \Rightarrow         & mn = 1\\
\Rightarrow         & m = +1 \mbox{ or } -1. \ea Hence $\alpha_1 =
Id$ and $\alpha_{-1}$ are the only two automorphisms of $C.$ As $x\neq x^{-1}, \alpha_{1}$ and
$\alpha_{-1}$ are distinct. Hence the result is proved.     $\blacksquare$

\bt\label{t38} Let $G = S_3$ be the symmetric group of degree 3.
Then all the automorphisms of $S_3$ are inner and $Aut S_3$ is
isomorphic to $S_3$ \et \pf In the symmetric group (Example \ref{ex7}),

  \bea
S_3  & = & \left\{
                Id={\left( \begin{array}{ccc}
                    1 & 2 & 3\\  1  & 2 & 3
                     \end{array} \right),}
                    a_{1} = {\left( \begin{array}{ccc}
                    1 & 2 & 3 \\  2  & 3  & 1
                    \end{array} \right),}
                     a_{2} = {\left( \begin{array}{ccc}
                      1  & 2 & 3\\  3  & 2 & 1
                  \end{array} \right) ,} \right.\\
     &   &  \quad \left.a_{3} = {\left( \begin{array}{ccc}
                    1  & 2 & 3\\  2  & 1 & 3
                   \end{array}\right),}
                   a_{4} = {\left(\begin{array}{ccc}
                    1 & 2 & 3\\  1  & 3 & 2
                   \end{array}\right),}
                    a_{5} = {\left(\begin{array}{ccc}
                    1 & 2 & 3\\  3  & 1 & 2
                   \end{array}\right)}
        \right\}
\eea
 on the three symbols \{1,2,3\}, the elements $a_{1},a_{2},a_{3}$ are precisely the elements of
 order 2. We have $a_{1}a_{2} = a_{4}$ and $a_{1}a_{3} = a_{5}$. Hence $S_{3} = gp\{a_{1},a_{2},a_{3}\}$.
By the theorem \ref{t34} any $\sigma \in Aut\; S_{3}$ maps the set $a_{1},a_{2},a_{3}$
onto itself. Further, as $S_{3} = gp\{a_{1},a_{2},a_{3}\}$, any
 $\sigma \in Aut S_3$
has unique representation  $$\theta(\sigma) = \left(\begin{array}{ccc}
  a_{1} & a_{2} & a_{3} \\
  \sigma(a_{1}) & \sigma(a_{2}) & \sigma(a_{3})
\end{array}\right) $$ as a permutation on the three symbols $\{a_{1},a_{2},a_{3}\}$.
 Moreover, if $\sigma,\rho \in Aut\; S_{3}$, then  $$\begin{array}{ccl}
   \theta(\sigma\rho) & = & \left(\begin{array}{ccl}
     a_{1} & a_{2} & a_{3} \\
     \sigma\rho(a_{1}) & \sigma\rho(qa_{2}) & \sigma\rho(a_{3}) \\
   \end{array}\right) \\
    & = & \left(\begin{array}{ccc}
      \rho(a_{1}) & \rho(a_{2}) & \rho(a_{3}) \\
      \sigma\rho(a_{1}) & \sigma\rho(a_{2}) & \sigma\rho(a_{3}) \\
\end{array}\right)\left( \begin{array}{ccc}
  a_{1} & a_{2} & a_{3} \\
  \rho(a_{1}) & \rho(a_{2}) & \rho(a_{3})
\end{array}\right)\\
    & = & \theta(\sigma)\theta(\rho) \
 \end{array}$$
 Hence, if we denote the group of permutations on the three symbols $\{a_{1},a_{2},a_{3}\}$
  also by the $S_{3}$, then $$\begin{array}{rcl}
    \theta : Aut\; S_{3} & \rightarrow & S_{3} \\
    \sigma & \mapsto & \theta(\sigma) \\
  \end{array}$$
gives a monomorphism of groups. Now let for any $x\in S_3$,
$\tau_x$ denotes the inner automorphism determined by $x$. Then
for $a_{1} \in S_3$, $$ \tau_{a_{1}} (a_{2}) = a_{1}a_{2}a_{1}^{-1} = a_{1}a_{2}a_{1} = a_{3}$$
 and by the lemma \ref{l37.05l}, $(\tau_{a_{1}})^{2} = \tau_{a_{1}^{2}}$ is the identity
 map since $a_{1}^{2} = Id$. Hence $\tau_{a_{1}} \neq Identity$, and $\circ
(\tau_{a_{1}}) = 2$. Further, for $a_{4} \in S_{3}$,   \ba{rcl}
           \tau_{a_{4}}(a_{1}) & = & a_{2}\\
      {(\tau_{a_{4}})}^2 (a_{1}) & = & a_{3}\\
 \mbox{and    }\;\;\; {(\tau_{a_{4}})}^3 (x) & = & \tau_{{a_{4}}^{3}}
\ea is identity automorphism since $a_{4}^{3} = Id$.
 Hence $Aut S_3$ contains an element of order 2  i.e., $\tau_{a_{1}}$
and an element of order 3 i.e., $\tau_{a_{4}}$. Therefore, by the corollary \ref{co3.9},
 $\circ (Aut S_3) \geq 6.$ Further, as $\theta$ is a monomorphism and
$\circ (S_3) = 6$, it follows $\theta$ is an isomorphism.
Moreover, as $\tau_{a_{1}}, \tau_{a_{4}} \in I_{(S_3)}$, order of $I_{(S_3)}$
is $\geq 6.$ Thus as $I_{(S_3)}$ is a subgroup of $Aut(S_3)$, we
conclude again by the corollary \ref{co3.9} $I_{(S_3)} = Aut(S_3)$, i.e., all the automorphisms of
$S_3$ are inner. Hence the result follows.     $\blacksquare$ \bt\label{t39} Let
$\varphi$ be an automorphism of $(\R,+)$, the additive group of
real numbers. If $\varphi (xy) = \varphi (x)\varphi (y)$ for all
$x,y\in \R$ then $\varphi = Id$, the identity map. \et \pf As
$\varphi$ is an automorphism, $\varphi(0) = 0$, and
$\varphi(-\alpha) = - \varphi (\alpha)$ for all $\alpha \in \R.$
Further, we have \ba{rl}
             &  \varphi (1) = \varphi (1.1) = {\varphi (1)}^2\\
\Rightarrow  & \varphi(1) = 1  \mbox{ since }\varphi (t)
                            \neq 0 \mbox { for any } t(\neq 0)\in \R.
\ea Now, as $\varphi (-1) = -\varphi (1) = -1$, by the additivity
of $\varphi$, we get $\varphi (m) = m$  for $m\in \Z.$ Further, if
$n(>0)\in \Z$, then \bea
     n \varphi (\frac{1}{n}) & = & \varphi (n\cdot \frac{1}{n}) = \varphi (1)\\
     \Rightarrow \quad \varphi(\frac{1}{n}) = \frac{1}{n}
\eea Hence, for any rational number $\frac{m}{n}$, we have $$
\varphi(\frac{m}{n}) = \frac{m}{n}$$ Now, let $\alpha > 0$ be any
real number. Then there is $\beta (\neq 0) \in \R$ such that
$\beta^2 = \alpha$. Hence $$ \varphi (\alpha) = \varphi (\beta^2)
= {\varphi (\beta)}^2 > 0$$ Thus if $\alpha > \beta$ $(\alpha,
\beta \in \R)$, then \ba{rl}
              &  \alpha - \beta > 0\\
\Rightarrow   & \varphi (\alpha - \beta) = \varphi (\alpha) -
\varphi (\beta) > 0\\ \ea Therefore $$\varphi (\alpha)  >  \varphi
(\beta)$$ Now, let for $\alpha \in \R$, $$ \varphi (\alpha) >
\alpha$$ Choose a rational number $\frac{m}{n}$ such that
\ba{rl}
               &  \varphi (\alpha) > \frac{m}{n}  > \alpha\\
\Rightarrow    & {\frac{m}{n}} = {\varphi (\frac{m}{n})} >
\varphi(\alpha) \ea This contradicts the fact that $\varphi(\alpha) > \frac{m}{n}.$
Hence $\varphi (\alpha) \leq\ \alpha.$ As above, we can show that
$\varphi (\alpha) < \alpha$ is not possible. Hence $\varphi
(\alpha) = \alpha$ for all $\alpha \in \R,$ i.e., $\varphi = Id.$

 We, now, prove a result which has some geometric interest.     $\blacksquare$
\bt\label{t310} The group of collineations and the group of affine
transformations on $\R^{2}$ are same. \et \pf We shall prove the
result in steps. \\ {\bf{ Step}} I. Every affine transformation is a
collineation.  \\   Let $Q = \left[\begin{array}{cc}
           a & b\\ c & d
           \end{array} \right]$
be an invertible  $2\times 2$ -matrix over reals. Consider the
linear map $$ Q : \R^2  \to  \R^2 $$ determined by $Q$.  Let $A =
(x_1,y_1)$, $B = (x_2,y_2)$ and $C = (x_3,y_3)$ be three collinear
points in $\R^2$.  Then the area of the triangle formed by these
points is $0$ i.e., \ba{rl}
            &   {det \left[\begin{array}{ccc}
                  x_1 & x_2 & x_3\\
                  y_1 & y_2 & y_3\\
                     1  & 1 & 1
                    \end{array} \right]} = 0 \\
\mbox{Therefore   } &          \\
            &   {det \left[\begin{array}{ccc}
                  a & b & 0\\
                  c & d & 0\\
                  0 & 0 & 1
                    \end{array} \right]} \cdot
                {det \left[\begin{array}{ccc}
                  x_1 & x_2 & x_3\\
                  y_1 & y_2 & y_3\\
                     1  & 1 & 1
                    \end{array} \right]} = 0\\
\Rightarrow
             &  {det \left[\begin{array}{ccc}
                 ax_1 + by_1 & ax_2+by_2 & ax_3+by_3\\
                 cx_1+dy_1 & cx_2+dy_2 & cx_3+dy_3\\
                     1  & 1 & 1
                    \end{array} \right]} = 0
\ea Hence the points  $Q(A), Q(B)$ and $Q(C)$ are collinear.
Further for any point $x = (a,b) \in \R$,  we have $$  {det
\left[\begin{array}{ccc}
         x_1+a & x_2+a & x_3+a\\
         y_1+b & y_2+b & y_3+b\\
          1  & 1 & 1
       \end{array} \right]} =
       {det \left[\begin{array}{ccc}
          x_1 & x_2 & x_3\\
          y_1 & y_2 & y_3\\
             1  & 1 & 1
         \end{array} \right]} = 0$$
Therefore, it follows that any affine transformation is a
collineation.   \\  {\bf{ Step}} II.  Let $f : \R^2 \to \R^2$ be a
collineation such that $f(0)=0.$ Then for any two non-proportional
points $\nu, w \in \R^2$, $f(\nu + w) = f(\nu) + f(w).$   \\
 Let $l_1 = \nu + [w], l_2 = w + [\nu]$ be two lines in $\R^{2}.$ As $\nu$
and $w$ are non-proportional $l_1\cap l_2 = \{\nu + w \}.$ Thus
$\{ f(\nu + w) \} = f(l_1)\cap f(l_2).$  Further, note that as $f$
is a collineation such that $f(0) = 0$ , $f([w]) = [f(w)]$
(use:$f([w])$, and $[f(w)]$ are straight lines passing through $0$
and $f(w)$ ). We have \ba{rl}
               &  (\nu + [w])\cap [w] = \emptyset\\
\Rightarrow    & f(l_1)\cap [f(w)] = \emptyset \ea Hence $f(l_1)$
is parallel to $[f(w)]$ and passes through $f(\nu).$ Similarly
$f(l_2)$ is parallel to $[f(\nu)]$ and passes through $f(w).$ Thus
$f(\nu) + f(w)$ is a point on both $f(l_1)$ and $f(l_2).$ We,
however, have $$ f(l_1)\cap f(l_2) = \{f(\nu +w) \}$$ Hence
$$f(\nu)+f(w) = f(\nu + w).$$

{\bf{Step}} III. Let $f : \R^2 \to \R^2$ be a collineation such that
$f(0) = 0$, $f(e_1) = e_1$ and $f(e_2) = e_2.$  Then $f = Id$, the
identity map.   \\  We are given that $f(0) = 0$ and $f(e_1) = e_1$.
Hence $f$ being collineation this maps the $x$-axis onto itself.
For similar reason, $f$ maps the $y$-axis onto itself. Let for
$t\in \R$, we have $$ f(te_1) = \varphi (t)e_1 \mbox { and }
f(te_2) = \psi (t)e_2$$ Then clearly  $\varphi (0) = 0 = \psi(0)$
and $\varphi (1) = \psi (1) = 1.$ Hence using the Step II for any
$t, s\in \R$, we have \ba{rl}
   f(te_1+se_2) & =  f(te_1)+f(se_2)\\
                & =  \varphi (t)e_1+\psi (s)e_2\\
\Rightarrow    & f(e_1+e_2) = e_1+e_2 \ea Now, as $f(0)= 0$ and
$f(e_1+e_2) = e_1+e_2$, we conclude that $f$ maps the line $x = y$
onto itself. Thus, as \ba {rl}
   f(te_1+te_2) & =  \varphi (t)e_1+\psi (t)e_2\\
\Rightarrow     & \varphi (t) = \psi (t) \mbox { for all } t\in
\R. \ea As $\varphi (0) = 0$, it is clear that if $t, s\in \R$ and
$t$ or $s$ is $0$, then $$ \varphi (t)\varphi (s) = \varphi (ts)$$
We shall prove that this is true for all $t, s\in \R.$ Assume
$t\neq 0, s\neq 0.$  Then by the step II, \bea
   f(e_1+se_2) & = & f(e_1)+f(se_2)\\
                & = & e_1+\varphi (s)e_2
\eea and \bea
   f(te_1+tse_2) & =  f(te_1)+f(tse_2)\\
                & =  \varphi (t)e_1+\varphi (ts)e_2
\eea As $(0,0), (1,s)$ and $(t,ts)$ are colinear points so are
their  images under $f.$ Hence $ \varphi (t)\varphi (s) = \varphi
(ts)$ for all $t,s \in \R$.   Taking $t = s = - 1$, we get \ba{rl}
              &  \varphi (1) = \varphi (-1)^{2}\\
\Rightarrow   &  \varphi (-1) = -1  (use :\varphi (1) = 1,
                           \mbox{ and }\varphi \mbox { is one-one })
\ea Further, for any $t, s\in \R$, \bea
       f(te_1+se_2) & = & \varphi (t)e_1+\varphi (s)e_2\\
                    & = & (\varphi (t)+\varphi (s))e_1+\varphi (s)(e_2-e_1)
\eea Also \bea
     f(te_1+se_2) & = & f((t+s)e_1+s(e_2-e_1))\\
                 & = & \varphi (t+s)e_1+f(se_2-se_1)\\
                 & = & \varphi (t + s)e_1+\varphi (s)(e_2-e_1)
\eea Consequently \be\label{e32} \varphi (t+s) = \varphi (t) +
\varphi (s)  \mbox { for all }t, s\in \R \ee Now by the theorem \ref{t39},
 we conclude that $\varphi = Id.$ Hence \bea
          f(te_1+se_2) & = & \varphi (t)e_1+\varphi (s)e_2\\
                       & = & te_1+se_2
\eea Thus, $f = Id.$ \\  {\bf{Step}} IV: A collineation is an affine
transformation.  \\  Let $T : \R^2 \to \R^2$ be a collineation, and let
$T(0) = \underline{\lambda }$ , $T(e_1) = \underline{a } , T(e_2)
= \underline{b}.$ Then define a linear map $F : \R^2 \to \R^2$,
determined by $$ F(e_1) = e_1+\underline{a}- \underline{\lambda},
 F(e_2) = e_2+\underline{b}- \underline{\lambda}$$
Now, for the affine transformation\, \bea
        A : \R^2 & \to & \R\\
    \underline{x} & \mapsto & F \underline{x}+\underline{\lambda}
\eea we have \bea
                (A-T)(0) & = & 0\\
              (A-T)(e_1) & = & Ae_1-Te_1 = e_1\\
              (A-T)(e_2) & = & e_2
\eea As an affine transformation is a collineation $A-T$ is a
collineation. Hence by the Step III, we get \ba{rl}
                &     A-T = I\\
\Rightarrow     & T\underline{x} =
(F-I)\underline{x}+\underline{\lambda} \ea Therefore, $T$ is an affine
transformation, and the proof is complete.     $\blacksquare$
\clearpage
 {\bf EXERCISES}
\begin{enumerate}
  \item
Let $G$ be a group.Prove that   \\ (i) Any homomorphism $\varphi$ from $(\Z,+)$, the additive
group of integers, to $G$ is of the from $\varphi(n) = x^{n}$ where $x \in G$ is a fixed
element.  \\ (ii) For any $x\in G$, the map $$\begin{array}{ccc}
  \tau_{x} : (\Z,+) & \rightarrow & G \\
  m & \mapsto & x^{m}
\end{array}$$ is a homomorphism.
 \item
Show that there exists no non-trivial homomorphism from $(\Q,+)$ to $(\Z,+)$.
\item
Let $G,H$ be the two groups. Prove that $f : G \rightarrow H$ is a homomorphism if and only
if $\{(x,f(x) \st x \in G\}$ is a group with respect to co-ordinate wise product i.e.,
with respect to the binary operation $$(x,f(x))(y,f(y)) = (xy,f(x)f(y)).$$
\item
Let $H$ be a subgroup of a group $G$, and $[G:H] < \infty$. Then let $K =\bigcap_{g \in G}\tau_{g}(H)$,
 where $\tau_{g} : G\rightarrow G$ is the inner automorphism determined by $g$, is a normal subgroup of $G$
 and $[G : K] < \infty$. Thus any subgroup of finite index in $G$ contains a normal subgroup
  of finite index.
  \item
  Let $G,W$ be two groups and let $\varphi : G\rightarrow W$ be an epimorphism. Prove that
   for any subgroup $H$ of $G$ with $\ker f \subset H$,$[G:H] =$ \\ $[W:\varphi(H)]$.
   \item
   Let $\varphi:G\rightarrow H$ be a group isomorphism. Prove that $\varphi(Z(G)) = Z(H)$
   and $\varphi(G') =H'.$
   \item
   Find all automorphisms of $(\Q,+)$, the additive group of rationals.
   \item
   Prove that any group with more then two elements has a non-trivial automorphism.
   \item
   Find the number of automorphisms of $\Z_{176}$.
   \item
   let $f_{i} : G \rightarrow H$,$i = 1,2$ be two group homomorphisms. Show that if $gp\{x\}
   = G$ then $f_{1} = f_{2}$ if and only if $f_{1}(x) = f_{2}(x)$ for all $x \in G$.
   \item
   Let $f : G \rightarrow H$ be a group homomorphism where $\circ(H) = n < \infty$. Prove that if for an element
    $x \in G$,$(\circ(x),n) =1$, then $x \in \ker f$.
    \item
   Let $G = gp\{a\}$ and $H$ be two groups. Prove that for any $b \in H$, there exists a
   homomorphism $\varphi : G \rightarrow H$ such that $\varphi (a) = b$ if and only if $\circ(b)$
    divides $\circ(a)$.
    \item
    Let $G$ be a finite abelian group with $\circ(G) = n$. Then for any $n \Z$,
    $\varphi : G \rightarrow G$ such that $\varphi(x) = x^{n}$ for all $x \in G$, is an automorphism
     if and only if $(m.n) = 1$.
     \item
     Prove that if a group $G$ has order $n$, then $\circ(Aut\; G)$ divides $(n-1)!$.
     \item
     Let $\sigma$ be an endomorphism of a group $G$ such that identity is the only element of
      $G$ fixed by $\sigma$. Prove that the map $$\begin{array}{rcl}
        \theta : G & \rightarrow & G \\
        x & \mapsto & x^{-1}\sigma(x) \\
        \end{array}$$ is one-one. Further, show that if $\theta$ is onto and $\sigma^{2} = Id$,
         then $G$ is abelian.
     \item
     Let $G$ be an infinite group. Prove that if $G$ is isomorphic to every non-identity
      subgroup of itself. Then $G$ is infinite cyclic.
     \item
     Find all proper subgroup of the quotient group $\Z/154\Z$.
     \item
     Let $S =\{s_{1},s_{2},\cdots,s_{n}\}$ and $U =\{u_{1},\cdots,u_{n}\}$ be two sets each with $n$ elements.
     Prove that the group $A(S)$ (Example \ref{ex6}) is isomorphic to $A(U)$.
  \item
  For the subgroup $(\Z,+)=\Z$ of $(\Q,+)$, prove that $\Q/\Z$ is isomorphic to $\C_{\infty}.$
  \item
  Let $\psi$ be an automorphism of a group $G$. Prove that if $x$ is a non-generator of $G$,
   then so is $\psi(x)$. Then deduce that  for the Frattini subgroup $\Phi(G)$ of $G$,
   $\psi(\Phi(G))=\Phi(G).$
   \item
   Let $n\geq 3$ and let
   $D_n=\{e,\sigma,\sigma^2,\cdots,\sigma^{n-1},\tau,\sigma\tau,\cdots,\sigma^{n-1}\tau\st
   \sigma^n=e=\tau^2, \sigma\tau=\tau\sigma^{-1}\}$ (Example \ref{ex9}) be
   the Dihedral group. Prove that \\ (i) For any fixed $0\leq i\leq n-1,$
   the map \[\alpha : D_n\rightarrow D_n\] satisfying
   $\alpha(\sigma^l\tau)=\sigma^{l+i}\tau$ for all $0\leq l\leq n-1$ and
   $\alpha(\sigma^k)=\sigma^k$ for all $k\geq 0,$ is an automorphism of
   $D_n$. \\ (ii) For any fixed $1\leq t\leq n$, $(t,n)=1,$ the map
   \[\beta: D_n\Rightarrow D_n\] satisfying $\beta(\sigma^k)=\sigma^{kt}$
   and $\beta(\sigma^k\tau)=\sigma^{kt}\tau$ for all $0\leq k\leq n-1$ is
   an automorphism of $D_n.$
   \item
   Let $H$ be a subgroup of group $G$ and $\alpha\in Aut\,G.$ Prove that
   \[\alpha(\overline{N}_G(H))=\overline{N}_G(\alpha(H)).\]
  \end{enumerate}

\chapter{Direct Product, Semi Direct Product And Wreath Product}
 \footnote{\it contents group6.tex }
In this chapter we shall describe some methods of constructing new groups
from a given collection of groups and also analyze when a given group is
constructed by a collection of groups (up to isomorphism).  \\   Let
$\{G_i \st i\in I\}$ be a family of groups. Then there is a natural way of
constructing a new group with the help of these groups $G_i$'s. Let us
denote by $G = \prod_{i\in I}G_i$, the set of functions : $$ \left\{ f:I
\to \coprod G_i \mid f(i)\in G_i\right\} $$ where $\coprod G_i$ denotes
the disjoint union of the $G_i$'s. For $f,g$ in $G$, define $$ (f\cdot
g)(i) = f(i) g(i)$$ Then it is easy to check that $G$ is a group with
respect to the above binary operation . The function $\Theta \in G$
defined by $\Theta (i) = e_i,$ the identity element of the group $G_i$, is
the identity element of $G$. For any $f\in G, \bar{f} \in G$ defined by
$\bar{f}(i) = {(f(i))}^{-1}$ is the inverse of $f$ in $G.$ The group
$\prod_{i\in I}G_i$, defined above, is called the cartesian product of the
groups $G_i's$. Some times we write, $$\prod_{i\in I}G_i =
\{\underline{x}=(x_i) \mid x_i \in G_i\}$$ where $\underline{x} : I
\rightarrow \coprod G_i$ is the function $\underline{x}(i)=x_i$. Clearly
$(x_i)(y_i)=(x_i y_i), (x_i)^{-1}=(x_{i}^{-1})$ and $\underline{e}=(e_i)$
is the identity element in $G$.
 The map, $$\begin{array}{cccc}
  p_{i} : & G & \rightarrow & G_{i} \\
   & f & \mapsto & f(i)
\end{array}$$ is a homomorphism and is called the $i^{th}$ projection of
$G$. For any $f \in G=\prod G_i$, the set $$supp(f) = \{ i \in I | f(i)
\neq e_{i}\}$$ is called the support of $f$. Clearly the  set of functions
in $G$ with finite support is a subgroup of $G$. We shall denote this by
${\Large X}_{_{_{_{\!\!\!\!\!\!\!\! i\in I}}}}G_i$. Any group isomorphic to
the subgroup ${\Large X}_{_{_{_{\!\!\!\!\!\!\!\! i\in I}}}}G_i$ is called
{\bf external direct product} of $G_{i'^{s}}$ and the groups $G_{i'^{s}}$ are
called its {\bf direct factors}. Thus direct product of $G_{i'^{s}}$ is
unique upto isomorphism.
\\ If $\abs{I}
< \infty$, then clearly direct product and cartesian product
coincide. In case $G_{i'^{s}}$ are additive groups we call direct sum instead
of direct product and $G_{i'^{s}}$ are called its direct summands instead of
factors. We shall write $$G = {{\bigoplus_{i \in I} G_{i}}}$$for the
direct sum of the groups $G_{i'^s}$
\\  We shall first develop the general properties of direct product for
finite index set $I$. If $I=\{1,2,\cdots,n\}$, then for the cartesian
product $G_1\times G_2\times \cdots \times G_n =
\{\underline{x}=(x_1,x_2,\cdots,x_n) \st x_i \in G_i \mbox{ for
}i=1,2,\cdots,n\}$ the map $$\begin{array}{rcl}
  \tau : G & \rightarrow & G_1\times G_2\times \cdots \times G_n =
\{\underline{x}=(x_1,x_2,\cdots,x_n)\} \\
  f & \mapsto & (f(1),f(2),\cdots,f(n))
\end{array}$$ is a bijection. Identifying $G$ with $G_1\times G_2\times \cdots \times
G_n$ under $\tau$ transfers in a natural way the group structure
on $G$ to a group structure on $G_1\times G_2\times \cdots \times
G_n$, where for $\underline{x}=(x_1,x_2,\cdots,x_n)$ and
$\underline{w}=(w_1,w_2,\cdots,w_n)$ in $G_1\times G_2\times
\cdots \times G_n$, the product $\underline{x}\cdot\underline{w}$
of $\underline{x}$ and $\underline{w}$ is
$$\underline{x}\,\underline{w}=(x_1w_1,\cdots,x_n w_n).$$ The two
groups are of course isomorphic. We shall normally understand the
group $G_1\times G_2\times \cdots \times G_n$ defined above as the
external direct product of the groups $G_1,G_2\cdots,G_n.$ For any
$1\leq i\leq n$, put $$\widehat{G_{i}} =
\{(e_1,\cdots,e_{i-1},a,e_{i+1},\cdots,e_n) \st a \in G_i\}$$ It
is easy to check that $\widehat{G_i}$ is a normal subgroup of
$G=G_1\times \cdots \times G_n,$ moreover, we have
\\ (i) $\underline{a}\,\underline{b} = \underline{b}\,\underline{a}$ for
all $\underline{a}\in \widehat{G_i}, \underline{b}\in \widehat{G_j}, i\neq
j$.
\\  (ii)  $\widehat{G_i}\cap (\widehat{G_1}\cdots
\widehat{G_{i-1}}\widehat{G_{i+1}}\cdots \widehat{G_n}) = \{e\}$ for all
$i=1,2,\cdots,n$. \\ and \\ (iii) $G = \widehat{G_1}\cdots \widehat{G_n}.$
\\
 If for a group $G$, there exist subgroups $H_1,H_2,\cdots ,H_n$
such that the map \bea
   \mu : H = H_1\times H_2\times \cdots \times H_n & \to & G\\
   (x_1,x_2,\cdots ,x_n) & \mapsto & x_1 x_2 \cdots  x_n
\eea where $x_1 x_2 \cdots  x_n$ is the product of $x_i$, $i=1,2,\cdots,n,$ in
$G$, is an isomorphism from the external direct product of $H_1,H_2,\cdots
,H_n$ onto $G$, then $G$ is called the {\bf internal direct product} of
$H_1,H_2,\cdots ,H_n$ As noted above, we have
\\ (i) $\underline{a}\,\underline{b} = \underline{b}\,\underline{a}$ for
all $\underline{a}\in \widehat{H_i}, \underline{b}\in \widehat{H_j}, i\neq
j$.
\\  (ii) $\widehat{H_i}\cap (\widehat{H_1}\cdots
\widehat{H_{i-1}}\widehat{H_{i+1}}\cdots \widehat{H_n}) = \{e\}$ for all
$i=1,2,\cdots,n$. \\ and \\ (iii) $H = \widehat{H_1}\cdots \widehat{H_n}.$
\\ As $\widehat{H_i}\lhd H$ for all $i=1,2,\cdots,n,$ clearly $\mu(\widehat{H_i})
 = H_i$, is normal in $G$ and
\\ (i) $ab = ba$ for all $a \in H_i, b\in H_j, i\neq j.$ \\ (ii) $H_i \cap
(H_1\cdots H_{i-1}H_{i+1}\cdots H_n) = id.$  \\  (iii) $G = H_1 \cdots
H_n.$ \\ We, now, prove a result which shows that these condition on
$H_i's$ are infact sufficient for $G$ to be internal direct product of
$H_i's.$
\bt\label{t41} Let $H_1,H_2,\cdots ,H_n$ be subgroups of a group $G$ such that
  \\ (i)$ab = ba $ for all $a\in H_i ,b\in H_j \mbox{ and }i\neq j.$\\
(ii)$ H_i\cap (H_1 H_2 \cdots H_{i-1} H_{i+1} H_{i+2}\cdots H_n) = id$ for
$i = 1,2,\cdots ,n.$\\ and \\ (iii) $G = H_1 H_2 \cdots  H_n.$\\  Then $G$
is the internal direct product of the subgroups
 $H_1,H_2,\cdots ,H_n.$
\et \pf Consider the map \bea
    \Phi : H_1\times H_2\times \cdots \times H_n & \to G\\
\underline{x} = (x_1,x_2,\cdots ,x_n) & \mapsto & x_1x_2 \cdots x_n. \eea
For $\underline{x} = (x_1,x_2,\cdots ,x_n)$ , $\underline{y} =
(y_1,y_2,\cdots ,y_n)$ in $\prod H_i = H_1\times H_2\times \cdots \times
H_n$, we have \bea
     \Phi(\underline{x} \underline {y}) & = &
             \Phi (x_1y_1,x_2y_2,\cdots ,x_ny_n)\\
             & = & x_1y_1x_2y_2 \cdots x_ny_n\\
             & = & x_1x_2\cdots x_n y_1y_2,\cdots y_n \mbox { (use (i)) }\\
             & = & \Phi(\underline{x}) \Phi(\underline{y})
\eea Hence $\Phi$ is a homomorphism from $\prod H_i$ to $G$.
Further, let $
               \Phi(\underline{x}) = e,$ the identity of $G$. Then \ba{rl}
\Rightarrow    & x_1x_2\cdots x_n = e\\ \Rightarrow    & x_1x_2\cdots
x_{i-1}x_{i+1}x_{i+2}\cdots x_nx_i = e
                                   \;\;\;\;\mbox{(use(i)) }\\
\Rightarrow    & x_1\cdots x_{i-1}x_{i+1}\cdots x_n = {x_i}^{-1} \in
H_i\cap (H_1 H_2\cdots H_{i-1} H_{i+1} H_{i+2}\cdots H_n)\\ \Rightarrow &
{x_i}^{-1} = e\;\;\;\; \mbox{ (use (ii)) }\\ \Rightarrow & x_i = e. \ea As
$i$ is arbitrary, we get $\underline{x} = \underline{e}$, the identity
 element of $\prod H_i$. Hence $\Phi$ is one-one. Finally $\Phi$ is onto
 follows from the assumption (iii). Hence $\Phi$ is an isomorphism.
\bos (1) It is clear from the isomorphism described above that
every element of $G$ can be expressed uniquely as product of elements of
$H_i$'s.\\ (2) If we assume in the theorem that (i) holds and every
element of $G$ can be expressed uniquely as product of elements of
$H_i$'s, then (ii) and (iii) also hold.  \\ (3) For the group $(\Z, +)$
for any two non trivial subgroups $m\Z, n\Z$ of $\Z$, $mn (\neq 0) \in m\Z
\cap n\Z$. Hence $\Z$ can not be decomposed as direct sum of its
subgroups.  \eos
\bco\label{co6.3co} Let $H_1,H_2,\cdots H_n$ be normal subgroups of a
group $G$ such that \\ (i) $H_i\cap (H_1 H_2\cdots H_{i-1} H_{i+1}
H_{i+2}\cdots H_n) = \{e\}$ for all $i=1,2,\cdots,n.$ \\ and  \\  (ii) $G
= H_1 H_2 \cdots H_n.$ \\ Then $G$ is internal direct product of $H_i's.$
\eco \pf By (i), it is clear that $H_i \cap H_j = \{e\}$ for all $i\neq
j.$ If $x \in H_i, y \in H_j$ for $i \neq j,$ then as $H_i's$ are normal,
$xyx^{-1}y^{-1} \in H_i \cap H_j.$ Thus $$\begin{array}{crl}
   & xyx^{-1}y^{-1} & = e \\
  \Rightarrow & xy & = yx
\end{array}$$ Hence the result follows from the theorem.
$\blacksquare$  \bco\label{co6.4.1} Let a group $G$ be (internal) direct
product of its two subgroups $A$ and $B$. Then for any subgroup $H$ of
$G$, $H \supset A, H$ is internal direct product of $A$ and $H\cap B$.
\eco \pf As $G$ is internal direct product of $A$ and $B$, $A \lhd G,
B\lhd G, A\cap B = id$, and $G = AB$. Therefore $A\lhd H, B_1 = B\cap H
\lhd H$ and $A\cap B_1 = id.$ Thus to complete the proof it is sufficient
to prove that $H = AB_1$. Since $G= AB$, for any $h \in H,$
$$\begin{array}{cll}
   & h=ab & \mbox{ where }a\in A, b\in B \\
  \Rightarrow & a^{-1}h = b \in B\cap H = B_1 & \mbox{ since }A \subset H  \\
  \Rightarrow & h = a(a^{-1}h) \in AB_1 &  \\
  \Rightarrow & H = AB_1 &
\end{array}$$ Hence the result follows.       $\blacksquare$
\bexe\label{e401e} Show that $(\Q, +)$ can not be expressed as
direct sum of two proper subgroups. \eexe \bexe\label{e401.1e} Let a group
$G$ be internal direct product of its subgroups $A$ and $B$. Prove that
$G/A$ is isomorphic to $B$. \eexe
\bt\label{t42} Let a group $G$ be isomorphic to the external direct
product $G_1\times G_2\times \cdots \times G_n$ of the groups $G_i
; i = 1,2\cdots ,n.$ Then there exist subgroups $H_i ; i =
1,2\cdots ,n.$ of $G$ such that $G_i$ is isomorphic to $H_i$ and\\ (i)$ab
= ba $ for all $a\in H_i ,b\in H_j, i\neq j.$\\ (ii)$ H_i\cap (H_1
H_2,\cdots H_{i-1} H_{i+1} H_{i+2}\cdots H_n) = id$\;\; for all $i =
1,2,\cdots ,n.$\\ (iii) $G = H_1 H_2 \cdots  H_n.$ \et \pf Let $$
\hat{G_i} = \{(e_1, \cdots ,e_{i-1},x,e_{i+1} \cdots ,e_n)\st x\in G_i ,
                   e_j = \mbox { identity of } G_j \mbox{, for all } j \neq i\}$$
It is easy to see that $ \hat{G_i}$ is a subgroup of $\prod G_i$ and we
have   \\  (a) $\underline{a}\, \underline{b} = \underline{b}\,
\underline{a}$ for all
    $\underline{a} \in \hat{G_i} ,\underline{b} \in \hat{G_j}$ and $i\neq
    j.$  \\
(b)$ \hat{G_i}\cap (\hat{G_1}\hat{G_2}\cdots \hat{G_{i-1}} \hat{G_{i+1}}
\hat{G_{i+2}}\cdots \hat{G_n}) = id$ for $i = 1,2,\cdots ,n.$   \\  (c)
$\prod G_i = \hat{G_1} \hat{G_2} \cdots  \hat{G_n}$. \\ If $\alpha
: \prod G_i \to G$ is an isomorphism. Put $\alpha \widehat{(G_i)} = H_i$
for $i = 1,2,\cdots ,n.$ Then from (a),(b),(c) above we get the required
conditions for $H_i$'s.\bexe\label{ex100} Given a finite set $\{
G_{1},\cdots, G_{n}\}$ of groups and $\sigma \in S_{n}$, the map$$G_{1}
\times \cdots \times G_{n} \rightarrow G_{\sigma(1)}\times \cdots \times
G_{\sigma(n)}$$ $$(x_{1}, \cdots,x_{n}) \mapsto (x_{\sigma(1)}, \cdots,
x_{\sigma(n)})$$ is an isomorphism of groups.
\eexe  \bd A group $G$ with $\circ(G) > 1$ is called indecomposable if $G$ is not
a direct product of its proper subgroups.\ed \bx The group $(\Q,+)$ and
$(\Z,+)$ are indecomposable.\ex
\bd A subgroup $H$ of a group $G$ is called a direct factor of $G$ if for
a subgroup $K$ of $G$, $G$ is (internal) direct product of $H$ and $K$.\ed
  \brs (i)
An indecomposable group has no direct factor.  \\  (ii) Any simple group
is indecomposable. \\ (iii) For a group $G$, the group $G$ and its
identity subgroup are always direct factors of $G$.
\ers
\bt\label{t43} Let $C_n$ denotes the cyclic group of order $n$. Then  \\
 $C_n \simeq C_m \times C_k $ if and only if $n = mk$ and $(m,k) = 1.$
\et \pf $(\Rightarrow)$  We have $\abs{C_m \times C_k} = mk.$
Hence $n = mk.$ Further, let $ (m,k) = d.$ Then, for any element
$(x,y)\in C_m \times C_k $, we have \ba{rl}
             & {(x,y)}^{\frac{mk}{d}} =
                            ({x}^{\frac{mk}{d}},{y}{\frac{mk}{d}}) = id\\
\Rightarrow  &  \circ (x,y) \leq {\frac{mk}{d}}. \ea Now, note
that the image of the generator of $C_n$ is an element of order $
n = mk.$  Hence $(m,k) = d = 1.$  \\  $(\Leftarrow)$  Let $C_m =
gp\{a\},  C_k = gp\{b\}.$ Then, $ {(a,b)}^{mk} = (a^{mk},b^{mk}) =
(e,e),$ the identity of $C_m \times C_k$. Hence $\circ ((a,b))
\leq mk = n.$ Further, if \ba{rl}
                & {(a,b)}^l = (e,e)\\
\Rightarrow     & (a^l,b^l) = (e,e)\\ \Rightarrow     & a^l= e,
b^l = e\\ \Rightarrow     & m \mid l \mbox{ and }k \mid l\\
\Rightarrow     & n = mk \mbox { divides } l \mbox { since }(m,k)
= 1 \\ \Rightarrow   & l\geq n.\ea Hence $\circ (a,b) = n.$ Now,
as $\abs{C_m \times C_k} = mk = n$, it follows that $C_m \times
C_k $ is cyclic of order $n$. Hence $C_m\times C_k$ is isomorphic
to $C_n.$ \br Let $C_n$ denotes a cyclic group of order $n > 1.$
Then $C_n$ is indecomposable if and only if $n = p^k$ where $p$ is
a prime and $k \geq 1.$ \er \bt\label{t4.50t} Let $f$ be an
endomorphism of a group $G$ such that $f^2=f$. If $H$, the image
of $f$, is a normal subgroup of $G$ and $K$ is the Kernel of $f$,
 then $G$ is the internal direct product of $K$ and $H.$ \et
\pf Let $x\in K\cap H.$ As H is the image of $f$ , $x = f(a)$ for some
$a\in G.$ Further,
 as $x\in K = \ker f$
\bea
            e & = & f(x)\\
              & = & f(f(a))\\
              & = & f(a) \mbox { (use $f^2 = f$) }\\
              & = & x. \\
\Rightarrow\;\;\; K\cap H  & = & {e}.
\eea
 Now, if $x\in H, y\in K,$ then, since $H$ and $K$ are normal subgroups of $G$,
$$xyx^{-1}y^{-1} \in H\cap K = {e}$$ $$\Rightarrow  xy = yx.$$ Further,
for any $g\in G$, $g = f(g)f(g^{-1})g$, and $$ \begin{array}{cl}
  &  f(f(g^{-1})g)
= f^2 (g^{-1})f(g) = f(g)f(g^{-1}) = e \\ \Rightarrow  & f(g^{-1})g \in
\ker f = K. \end{array}$$ Hence, $G=HK$. Consequently, by theorem
\ref{t41}, $G$ is internal direct product of the subgroups $K$ and $H.$
\bt\label{t45} Let $G$ be a finite abelian group in which every
non-identity element has order $2$. Then $G$ is direct product of cyclic
groups of order $2$. \et \pf Let $\{a_1,a_2, \cdots ,a_n \}$ be a minimal
set of generators for $G$. Clearly $a_i\neq e$ for any $i$. Put $C_i =
gp\{a_i \}$. Then each $C_i$ is a cyclic group of order $2$. As
$\{a_1,a_2, \cdots ,a_n \}$ is a set of generators for $G$, every element
$g$ in $G$ can be expressed as $$ g = {(a_1)}^{\alpha_1}
{(a_2)}^{\alpha_2} \cdots {(a_n)}^{\alpha_n}$$ where $\alpha_i = 0$ or
$1$. As $G$ is abelian, it is easy to check that the map
\bea
        \Phi : C_1\times C_2\times \cdots \times C_n & \to & G\\
  \underline{x} = (x_1,x_2, \cdots ,x_n) & \mapsto & x_1x_2 \cdots x_n
\eea is an onto homomorphism . To complete the proof we have to show that $\Phi$
 is one-one. Let for $x_i =
{(a_i)}^{\alpha_i}$, we have $ {(a_1)}^{\alpha_1} {(a_2)}^{\alpha_2}
\cdots {(a_n)}^{\alpha_n} = e.$
             If $\underline{x} \neq id,$ then $\alpha_i =1$ for some $i$,
             and hence
              then
$ a_i =  {(a_1)}^{\alpha_1}\cdots {(a_{i-1})}^{\alpha_{i-1}}
              {(a_{i+1})}^{\alpha_{i+1}} \cdots {(a_n)}^{\alpha_n},$
since ${a_i}^{-1} = a_i$ for all $i$. This contradicts that $\{a_1,a_2,
\cdots ,a_n \}$ is a minimal set of generators for $G$. Hence $\Phi$ is
one-one and is an isomorphism.      \br
 We have $\abs{G}= \abs{ C_1\times C_2\times \cdots \times C_n} = 2^n.$
\er
\bl\label{l45}
Let $H$ denotes the multiplicative group of real numbers\\ $\{+1,-1 \}$,
and $G = H \times H$ be the direct product of $H$ with itself. Then the
group $Aut G$ is isomorphic to $S_3$.
\el
\pf We have $\circ (G) = 4$, and $a = (1,-1) , b = (-1,1)$, $c = (-1,-1)$
are precisely the elements of order $2$ in $G$. Note that $ab = c$, $bc =
a$, $ca = b.$  Hence for any permutation $\tau$ of
 the subset $\{a,b,c\}$, if $\hat{\tau} (a) = \tau (a)$,
$\hat{\tau} (b) = \tau (b)$ , $\hat{\tau} (c) = \tau (c)$ and $\hat{\tau}
(1,1) =  (1,1)$, then $\hat\tau \in Aut G$. Further, for any $\sigma \in
Aut G$, $\sigma$ restricted to ${a,b,c}$ gives a permutation of
$\{a,b,c\}$. Therefore $\sigma = \hat{\tau}$ for some permutation $\tau$
of the three symbols ${a,b,c}$. Clearly, the map $$\begin{array}{ccc}
  S_3 & \rightarrow & Aut\,G \\
  \tau & \mapsto & \widehat{\tau}
\end{array}$$ from the group of permutations on $\{a,b,c\}$ to $Aut\, G$
is one-one, and for any permutations $\tau_1,\tau_2$ of $\{a,b,c\}$
$\widehat{\tau_1 \tau_2} = \widehat{\tau_1} \widehat{\tau_2}$. Hence $S_3$
is isomorphic to $Aut\, G$.
\br We have proved that $Aut S_3$ is isomorphic to $S_3$ (Theorem
\ref{t38} ). Hence, as $S_3$ is not an abelian group, non-isomorphic
groups may have isomorphic automorphism groups.
\er  \bt\label{t47} Let $H$ be a normal subgroup
of a group $G$ such that $Z(H) ={e}$ and $Aut H = I-Aut H$. Then $H$ is a
direct factor of $G$ i.e., there exists a normal subgroup $A$ of $G$ such
that $G \cong A \times H$.\et \pf Consider the homomorphism :
$$\begin{array}{cccc}
  \tau : & G & \rightarrow & Aut H \\
   & g & \mapsto & \tau_{g}
\end{array}$$ where $\tau_{g}(h) = ghg^{-1}$ for all $h \in H$. By
assumption $Aut H = I-Aut H$, hence $\tau$ is onto. Let $A = \ker \tau$,
then $A \lhd G$. Note that$$\begin{array}{ccccccc}
  & \tau_{g} & = & id &  &  &  \\
  \Leftrightarrow & \tau_{g}(h) & = & h &  &  & \forall h \in H \\
  \Leftrightarrow & gh & = & hg &  &  & \forall h \in H
\end{array}$$  Hence, as $Z(H)=e$, $A \bigcap H = id$. Now, as $\tau$ is onto, for any
$x \in G$, $\tau_{x} = \tau_{h}$ for some $h \in H$. Therefore
$$\begin{array}{cccccc}
  & id & = & \tau_{x}\tau_{h}^{-1} & = & \tau_{xh^{-1}} \\
  \Rightarrow & xh^{-1} & \in & A &  &  \\
  \Rightarrow & x & \in & AH &  &  \\
  \Rightarrow & G & = & AH &  &  \\
   &  &  &  &  &
\end{array}$$ Hence by corollary \ref{co6.3co} $G \cong A \times H$.
$\blacksquare$ \bl\label{l6.36l} Let $N_i ; i=1,2,\cdots,k,\,
k\geq 2$, be distinct minimal normal subgroups of a group $G$ such
that $G=N_1 N_2\cdots N_k$, the product of the subgroups $N_i's$.
Then $G$ is the direct product of a subset of
$\{N_1,N_2,\cdots,N_k\}$.\el \pf First of all, since $N_i\lhd G$
for all $i=1,2,\cdots,k;$ $N_1 N_2\cdots N_k = N_{\sigma(1)}
N_{\sigma(2)}\cdots N_{\sigma(k)}$ for any $\sigma \in S_k$, the
symmetric group on $\{1,2,\cdots,k\}$. Let $t$ be the smallest
integer such that $$G=N_{i_{1}}N_{i_{2}}\cdots N_{i_{t}}$$ where
$1\leq i_l\leq k$ for all $l=1,2,\cdots,t$. Without any loss of
generality, by change of notation, we can assume that $G=N_1
N_2\cdots N_t$. Clearly $t\geq 2$. As $N_i's$ are distinct minimal
normal subgroups of $G$, $N_i\cap N_j=e$ for all $i\neq j$.
Therefore $xy=yx$ for all $x \in N_i,$ and $y\in N_j, 1\leq i\neq
j \leq k$ (Lemma \ref{w} (b)). Further, by induction , using lemma
\ref{w} (a), any finite product of normal subgroups of $G$ is
normal. Hence for $1\leq i\leq t$, if $a \in N_i\cap N_1 N_2\cdots
N_{i-1} N_{i+1}\cdots N_t$, then as $N_i$ is minimal normal
subgroup of $G$, either $a=e$ or $N_i\subset N_1 N_2\cdots
N_{i-1}N_{i+1}\cdots N_t.$ However, if $N_i \subset N_i\cap N_1
N_2\cdots N_{i-1}N_{i+1}\cdots N_t$, then $$G = N_1N_2\cdots
N_{i-1}N_{i+1}\cdots N_t,$$ a contradiction to the choice of $t$.
Consequently, $$ N_i\cap N_1N_2\cdots N_{i-1}N_{i+1}\cdots N_t =
\{e\}$$ for all $1 \leq i \leq t$. Therefore $G$ is isomorphic to
the direct product $N_1\times N_2 \times \cdots \times N_t.$
(Theorem \ref{t41})     $\blacksquare$ \bco\label{co6.37} Let for
a group $G$, $G=N_1 N_2 \cdots N_k, k\geq 2$, where $N_i's$ are
distinct minimal normal subgroups of $G$ for all $i=1,2,\cdots,k$.
Then any proper normal subgroup $N$ of $G$ is a direct factor of
$G$.\eco  \pf Choose a subset $J=\{j_1,j_2,\cdots,j_t\}$ of
$\{1,2,\cdots,k\}$ with $\abs{J}=t$ minimal such that $G=NP$ for
$P=N_{j_1}N_{j_2}\cdots N_{j_t},$ the product of the subgroups
$N_{j_k}, 1\leq k\leq t,$ in $G.$ As seen in the theorem $P\lhd
G$. By change of notation, if necessary, we can assume
$J=\{1,2,\cdots,t\}$. As $N\lhd G$ and $N_i$ a is minimal normal
subgroup $G$, $N\cap N_i=N_i$ or $\{e\}$. Let for some $1\leq
i\leq t$, $N_i\cap N=N_i$. Then $N_i\subset N$, and for $Q=N_1
N_2\cdots N_{i-1}N_{i+1}\cdots N_t$, $G=NQ$. This contradicts the
minimality of $\abs{J}$. Consequently $N\cap N_i =
\{e\}\,\;\forall 1\leq i \leq t$. Now, let $N\cap P \neq id$. Then
there exists $x(\neq e) \in N$ and $a_i \in N_i, \, i=1,\cdots t$,
such that $$x=a_1 a_2 \cdots a_t.$$ As $x \neq e$, one of the $a_i
\neq e$. As $a_ia_j=a_ja_i$ for all $1\leq i\neq j\leq t,$ if
necessary, by change by notation we can assume that $a_1 \neq e$.
Hence $$\begin{array}{rrl}
   & xa_{t}^{-1}\cdots a_{2}^{-1} & = a_1 \neq e  \\
  \Rightarrow & N_1 \cap NQ_1 & \neq e
\end{array}$$ where $Q_1 = N_2 \cdots N_t$. Note that $NQ_1$ is a normal
subgroup of $G$. Hence, as above, $N_1\cap NQ_1 \neq e$ implies $N_1
\subset NQ_1$ Then as $NQ=G$, $NQ_1 = G$. A contradiction to minimality of
$\abs{J}$. Consequently $N \cap P = \{e\}.$ Now, by theorem \ref{t41}, $G$
is internal direct product of $N$ and $P.$ Hence the result follows.
$\blacksquare$
 \\  We shall, now, define a class of groups such
that any member of the class if contained in an abelian group is
its direct factor.  \bd\label{d26.06d} A group $G$ is called
divisible if for every $n\geq 1,$ and $g \in G$, there exists
at-least one solution of the equation $x^n = g (nx=g$ in additive
notation) in $G$.\ed \br A group $G$ is divisible if and only if
$G=\{x^p\mid x\in G\}$ for all primes $p$.\er \bx The group
$(\Q,+), (\R,+)$ and $(\C,+)$ are divisible.\ex \bx The group
$(\C^*,\cdot)$ is divisible.\ex \bx For any prime $p$,
$\C_{p^\infty}$ is divisible.\ex \bexe Show that the
multiplicative group $(\R_{+}^{*}, \cdot)$ of non-zero positive
real numbers, is divisible.\eexe \bexe\label{p12} Show that any
homomorphic image of a divisible group is divisible. Further show
that direct product of divisible groups is divisible. \eexe
\bt\label{t6.40} Let $G$ be an abelian group, and $H$, a proper
subgroup of $G$. If $H$ is divisible, $H$ is a direct factor of
$G$ i.e., $G$ is internal direct product of $H$ and a subgroup $K$
of $G$.\et \pf Let $$\mathcal{S} = \{A \subset G \st A: \mbox{
subgroup, }A \cap H = \{e\}\}$$ Define an ordering $\leq$ in
$\mathcal{S}$ as follows :
\\ For $A_1,A_2 \in \mathcal{S},$  $A_1 \leq A_2$ if and
only if $A_1 \subset A_2.$ \\ Then $\leq$ is a partial ordering on
$\mathcal{S}.$ If $\{A_i\}_{i\in I}$ is a chain in $(\mathcal{S},
\leq)$, then put, $$ B = \bigcup_{i\in I}A_i$$  Let $x,y \in B$.
Then, as $\{A_i\}_{i \in I}$ is a chain, there exists $i_0 \in I$
such that $x,y \in A_{i_{0}}.$ hence $xy^{-1} \in A_{i_{0}}
\subset B.$ Thus $B$ is a subgroup of $G$. Further $$B\cap H =
\bigcup_{i\in I}(A_i \cap H) = e$$ Hence $B \in \mathcal{S}.$
Clearly, $B$ is an upper bound of the chain $\{A_i\}_{i\in I}$ in
$\mathcal{S}$. Therefore, by Zorn's lemma $(\mathcal{S}, \leq)$
has a maximal element (say) $K$. We shall prove that $G$ is
internal direct product of $K$ and $H$. To prove this, it suffices
to show that $HK=G$ (Theorem \ref{t41}). If $HK \neq G,$ then
choose $g \in G$, $g \not\in HK.$ By choice of $K$,  $$
\begin{array}{clc}
   & gp\{g,K\}\cap H \neq e &  \\
  \Rightarrow & g^m k = h (\neq e) & \mbox{ for some }h\in H, k\in K, m>1 \\
  \Rightarrow & g^m = hk^{-1} \in HK &
\end{array}$$ Let $m$ be the smallest $+ve$ integer such that $g^m \in HK.$ We
can assume that $m$ is prime. If not, then for any prime $p$ dividing $m$,
put $n=\frac{m}{p}$ and $g_1=g^n$. Clearly $g_{1}^{p}=g^m \in HK$, and
$g_1 \not\in HK$ by minimality of $m$. Now, as $H$ is divisible, we can
write $$\begin{array}{clc}
   & g^m=h_{1}^{m}k^{-1} & \mbox{ for some }h_{1}\in H \\
  \Rightarrow & (h_{1}^{-1}g)^m=k^{-1} &
\end{array}$$ Note that $h_{1}^{-1}g\not\in HK$, since $g\not\in HK.$ Thus
 for $g_2=h_{1}^{-1}g,\: g_{2}^{m}\in K$ where $m$ is prime, $g_2\not\in HK.$
 By choice of $K$, as above $$\begin{array}{cc}
    & gp\{g_2,K\}\bigcap H\neq e \\
   \Rightarrow & g_{2}^{t}k_2=h_2(\neq e) \
 \end{array}$$ for some $k_2\in K$ and $h_2\in H$, where $t>1.$ As
 $g_{2}^{m}\in K,$ we can assume $t<m.$ since $H\bigcap K=e.$
   As $g_{2}^{t},g_{2}^{m} \in HK$ for some $1<t<m,$ and $m$ is prime, it follows that
$g_2\in HK.$ (use
: $am+bt = 1$ for some $a,b \in \Z$). This contradicts the fact that $g_2
\not\in HK$. Hence $HK = G$, and the result follows.     $\blacksquare$
  \bd A group $G$ is called semi-simple (or
completely reducible) if every normal subgroup of $G$ is a direct
factor of $G$.\ed \bx (i) Any simple group is semi-simple.  \\
(ii) The Klein's 4-group $V_4$ (Example \ref{ex1.11}) is
semi-simple, but is not simple.\ex \bexe Prove that in a
semi-simple group $G$, normal subgroup of  normal subgroup is
normal.\eexe  \bexe Give an example of a non-abelian semi-simple
group which is not simple.\eexe \bl\label{l6.41l} Every normal
subgroup of a semi-simple group is semi-simple.\el \pf Let $N$ be
a normal subgroup of a semi-simple group. Then there exists a
normal subgroup $K$ of $G$ such that $G$ is internal direct
product of $N$ and $K$. Hence if $A\lhd N,$ then $A\lhd G.$ Thus
as $G$ is semi-simple $$G=AB \mbox{ (direct product) }$$ for a
normal subgroup $B$ of $G$. Therefore by corollary \ref{co6.4.1}
$N=A(N\cap B)$ (direct product). Hence $N$ is semi-simple.
$\blacksquare$ \bl\label{l6.41.1l} Homomorphic image of a
semi-simple group is semi-simple.\el \pf Let for a semi-simple
group $G$,  $$\varphi : G \rightarrow H$$ be an epimorphism of
groups. Let $N=\ker\varphi.$ Then $N\lhd G$ (Theorem \ref{t31}).
Since $G$ is semi-simple $G=NA$ (direct product) for a normal
subgroup $A$ of $G$. The factor group $G/N$ is isomorphic to $H$
(Corollary \ref{c43a}). Further, $G/N$ is also isomorphic to $A$,
which is semi-simple (Lemma \ref{l6.41l}). Hence $H$ is isomorphic
to $A$, hence is semi-simple.     $\blacksquare$ \bl\label{l6.42l}
Let $G$ be a semi-simple group. Then any normal subgroup $N(\neq
id)$ of $G$ contains a minimal normal subgroup of $G$.\el \pf If
$G$ is simple there is nothing to prove. Hence, let $G$ be
non-simple. Let $N(\neq id)$ be a normal subgroup of $G$. Choose
$x(\neq e) \in N,$ and let $H=gp\{g^{-1}xg \st g \in G\}.$ Then
$H\lhd G$ and $H \subset N.$ Put $$\mathcal{S}=\{K\varsubsetneqq H
\st K\lhd G\}$$ Define an ordering $\leq$ over $\mathcal{S}$ as
follows : \\ For $K_1,K_2 \in \mathcal{S}$,
\\ \hspace{.2in} $K_1\leq K_2 \Leftrightarrow K_1 \subset K_2$ \\ Then
$(\mathcal{S},\leq )$ is a partially ordered set. Take a chain
$\{K_{\alpha}\}_{\alpha \in A}$ in $\mathcal{S}$. Then $L=\cup_{\alpha \in
A}K_{\alpha}$ is normal subgroup of $H$. Note that $L\neq H,$ since
otherwise $x\in L=\cup K_{\alpha}$ which implies $x \in K_{\alpha}$ for
some $\alpha \in A$. As $K_{\alpha} \lhd G$, clearly $x \in K_{\alpha}$
implies $H=K_{\alpha}.$ This contradicts our assumption that $K_{\alpha}
\varsubsetneqq H.$ Hence $L\neq H$, and $L \lhd G.$ Thus $L$ is an upper
bound of the chain $\{K_{\alpha}\}_{\alpha \in A}$ in $\mathcal{S}$. Now,
by Zorn's lemma $\mathcal{S}$ has a maximal element (say) $S$. By lemma
\ref{l6.41l}, $H$ is semi-simple. Hence $H=ST$ (direct product) for a
normal subgroup $T$ of $H$. As any normal subgroup of $H$ is normal in
$G$, $S$ is a maximal normal subgroup of $H$, and hence $H/S \cong T$ is
simple. Thus $T$ is a minimal normal subgroup of $H$ and hence of $G$.
$\blacksquare$   \\ To Prove the next result we need to define the product
of an arbitrary family $\{A_{\alpha}\}_{\alpha \in \Lambda}$ of subsets of
a group $G$. This is given by $$\Pi^s \{a_1\cdots a_n\mid a_i \in
A_{\alpha_i}; i=1,2,\cdots,n\}$$, the set of all finite products of
elements in $\displaystyle{\bigcup_{\alpha\in\Lambda}A_{\alpha}}.$ if
$A_{\alpha}\lhd G$ for each $\alpha \in \Lambda,$ then $\Pi^s$ is a normal
subgroup of $G$. Since any finite product of normal subgroups in $G$ is a
normal subgroup of $G$.
\bl\label{l6.43l} Every semi-simple group is a product of its minimal
normal subgroups.\el \pf Let $G(\neq e)$ be a semi-simple group. By lemma
\ref{l6.42l}, $G$ contains a minimal normal subgroup. Let
$(N_{\alpha})_{\alpha \in A}$ be the set of all minimal normal subgroups
of $G$. Put $N=\prod_{\alpha \in A}N_{\alpha}$. If $N=G$, there is nothing
to prove. If not, then as $G$ is semi-simple and $N\lhd G,$ there exists
$K\lhd G,$ $K\neq id$, such that $G=NK$ (direct product). By lemma
\ref{l6.42l}, $K$ contains a minimal normal subgroup $L$ of $G$. Thus
$L=N_{\alpha}$ for some $\alpha \in A$. This is not possible as $L \cap N
= id$. Hence $K = \{e\}$ and $N=G$.     $\blacksquare$
\section{ Direct product (General Case ) }
For an arbitrary family of groups $\{G_{i}\}_{i \in I}$, we have defined
$G = \displaystyle{{{\Large X}_{_{_{_{\!\!\!\!\!\!\!\! i\in I}}}}
G_{i}}}$, the external direct product of the groups $G_{i}'s$, to be the
set of all functions $$f : I \rightarrow {\bigsqcup G_{i}}_{i \in I}$$
such that $\abs{supp(f)} < \infty$. Let $$\widehat{G_{i}} = \{f \in G \mid
f(j) = e_j, \mbox{ where $e_j$ is the identity of }G_j, \mbox{ for all } j
\neq i\}$$ It is easy to check that $\widehat{G_{i}}$ is a normal subgroup
of $G$ and is isomorphic to the group $G_{i}$ under the projection map
$p_{i}$. We identify $G_{i}$ to $\widehat{G_{i}}$ via $p_{i}$, then it is
clear that  \\ (i) Any element of $G$ is finite product of elements of
$G_{i}'s$.  \\  (ii) For all $x \in G_{i}, y \in G_{j}, i \neq j, xy =
yx$.  \\ (iii) For any finite subset $J = \{j_{1}, j_{2}, \cdots, j_{k}\}$
of $I$ and $i \in I - J$,$$G_{i} \bigcap G_{j_{1}}G_{j_{2}} \cdots
G_{j_{k}} = id.$$  As defined in the case $\abs{I} < \infty$, we again
define :  \bd\label{dooo3} Let $\{A_{i}\}_{i \in I}$ be a family of
subgroups of a group $A$. We shall say that $A$ is an internal direct
product of the groups $A_{i}'s$ if the map : $$\begin{array}{rcl}
  {\displaystyle{{\Large X}_{_{_{_{\!\!\!\!\!\!\!\! i\in I}}}} A_{i}}} & \rightarrow & A \\
  f & \mapsto & {\displaystyle{\prod_{i\in supp(f)} f(i)=\prod_{i\in I}f(i)}}
\end{array}$$ is defined and is an isomorphism of groups.\ed  {\bf Note} :  In
the above definition it assumed that the product is independent of the
order of $f(i)'s$.
\\ As seen above if $A$ is an internal direct product of $A_{i}'s$ then as
$\widehat{A_{i}} = \{f \in \displaystyle{{\Large
X}_{_{_{_{\!\!\!\!\!\!\!\! i\in I}}}} A_{i}} \mid f(j) = id \mbox{ for all
} j \neq i\}$ is a normal subgroup of $\displaystyle{{\Large
X}_{_{_{_{\!\!\!\!\!\!\!\! i\in I}}}}A_{i}}$, each $A_{i}$ is normal
subgroup of $A$, and we have  \\ (i) Any element of $A$ is a finite
product of elements of $A_{i}'s$.  \\ (ii) For all $x \in A_{i}, y \in
A_{j}, i \neq j, xy = yx$.  \\  (iii) For any finite subset $J = \{j_{1},
j_{2}, \cdots, j_{k}\}$ in $I$ and $i \in I - J$,$$A_{i} \bigcap
A_{j_{1}}A_{j_{2}} \cdots A_{j_{k}} = id.$$ Conversely, if we have a
family $\{A_{i}\}_{i \in I}$ of subgroups of a group $A$, satisfying the
properties $(i), (ii)$ and $(iii)$ above then one can easily check that
the map : $$\begin{array}{cccc}
  \theta : &{\displaystyle {\prod_{i\in I} A_{i}}}  & \rightarrow & A \\
   & f & \mapsto & {\displaystyle{\prod_{i\in I} f_{i}}}
\end{array}$$ is an isomorphism of groups i.e., $A$ is internal
direct product of $A_{i}'s$ in this case.  \bx  For the group $\Q^{\star}$
of non-zero rationals with multiplication, $\Q^{\star}$ is internal direct
product of the subgroups $\{ \pm 1\}$ and $gp\{p\}$, where $p$ is a prime
i.e., $$\Q^{\star} = \{\pm 1\}\,\, \mbox{\LARGE
x}\!\!\!\!\!\!\!\!\!\!_{_{_{p
: \mbox{\scriptsize{ prime}}}}}\!\!\normalsize{gp \{p\}}$$
\ex
\section{Extensions Of A Group And Semi-Direct Products :}
\bd\label{d0000} A group $G$ is called an extension of a group $A$
by a group $B$ if there exists an epimorphism $\alpha : G \rightarrow B$
such that $\ker \alpha$ is isomorphic to $A$. Further, if there exists a
homomorphism $\beta : B \rightarrow G$ such that $\alpha \beta = id$, then
$G$ is called a split extension of $A$ by $B$.\ed \brs\label{r0000} (i) If
$H$ is a normal subgroup of a group $G$, then the natural map :
$$\begin{array}{cccc}
  \eta : & G & \rightarrow & G/H \\
   & x & \mapsto & xH
\end{array}$$ is an epimorphism and $\ker \eta = H$. Hence $G$ is
an extension of $H$ by $G/H$.  \\  (ii) An extension of a group
$A$ by a group $B$ is not unique (up to
isomorphism)e.g.,$$\begin{array}{rcl}
    \Z_{2} \bigoplus {\Z}_{2} & \stackrel{p_{1}}{\rightarrow } & {\Z}_{2} \\
    (x, y) & \mapsto & x
   \end{array}$$
   and $$\begin{array}{rcl}
    \Z_{4} & \stackrel{\theta }{\rightarrow } & {\Z}_{2} \\
   x & \mapsto & 2x
\end{array}$$ are group epimorphisms and $\ker \theta \cong \ker
p_{1} \cong {\Z}_{2}$,but ${\Z}_{4}$ is not isomorphic to $\Z_{2}
\bigoplus {\Z}_{2}$. Thus, both ${\Z}_{4}$ and $\Z_{2} \bigoplus {\Z}_{2}$
are non-isomorphic extensions of ${\Z}_{2}$ by ${\Z}_{2}$. In fact,
$\Z_{2} \bigoplus {\Z}_{2}$ is a split extension.  \\  (iii) If $G$ is a
split extension of $A$ by $B$, then $\beta$ is a monomorphism , since if
$\beta(b) = id$, then $b = \alpha \beta (b) = id$.\ers \bexe Let $H$ be a
normal subgroup of a group $G$ such that $G/H$ is a cyclic group of order
$m$. Prove that if $G$ contains an element $x$ of order $m$ such that $o(xH)=m$,
then $G$ is a
split extension of $H$ by cyclic group of order $m$.\eexe \bexe Prove that
$Gl_{n}(\R)$ is a split extension of $Sl_{n}(\R)$ by ${\R}^{*}$.\eexe
\bexe Show that $S_{3}$ is a split extension of $C_{3}$ by $C_{2}$, and
use this to show that two split extensions of a group $A$ by a group $B$
need not be isomorphic.\eexe \bexe Let $m > 1$, $n > 1$ be two co-prime
integers. Prove that any extension of $C_{m}$ by $C_{n}$ is a split
extension.\eexe
\bl\label{l0000} Let $G$ be a group. Then $G$ is a split extension
of a group $A$ by a group $B$ if and only if there exists a normal
subgroup $H$ of $G$ and a subgroup $K$ such that $H \cong A$, $K \cong B$,
$G = HK$, and $H \bigcap K = \{e\}$. \el \pf By definition of split
extension there exists an epimorphism $$\alpha
: G \rightarrow B$$ such that $A$ is isomorphic to $\ker\alpha$ and there is
 a homomorphism $$\beta : B \rightarrow G$$
such that $\alpha \beta  = id$. As noted above $\beta$ is a monomorphism.
If $\beta(B) = K$. Then $B$ is isomorphic to $K$. Let $H = \ker \alpha$.
Then $H \lhd G$. For any $x \in G$, $$\begin{array}{crclcc}
   & \alpha(x\beta \alpha(x^{-1})) & = & e &  &\mbox{since }\alpha\beta = id  \\
  \Rightarrow & x\beta \alpha(x^{-1}) & \in & \ker \alpha  = H &  &  \\
\end{array}$$ If $$\begin{array}{crl}
   & x\beta\alpha(x^{-1}) & =h\in H \\
  \Rightarrow & x & =h\beta\alpha(x) \in HK \\
  &  \mbox{since  }  K & =im\beta
\end{array}$$ Now, to complete the proof of direct part, we have to show that $H
\bigcap K = \{e\}$. Let $$\begin{array}{crclclccl}
   & y & \in & H \bigcap K &  &  &  &  &  \\
  \Rightarrow & y & = & \beta(b) & = & h &  &  & \mbox{where }h \in H, b \in B \\
  \Rightarrow & \alpha(y) & = & \alpha\beta(b) & = & \alpha(h) & = & e &  \\
  \Rightarrow & b & = & e &  &  &  &  & \mbox{since }\alpha\beta = id \\
  \Rightarrow & y & = & \beta(b) & = & e &  &  &
\end{array}$$ Hence direct part is proved. Conversely, let for a normal
subgroup $H$ of $G$ and a subgroup $K$, $G = HK,$ and $H \bigcap K =
\{e\}$. Then every element $x$ of $G$ can be uniquely expressed as $x =
hk$, where $h \in H, k \in K$. Now, define the map : $$\begin{array}{crcl}
   & \alpha\; : G & \rightarrow & K \\
   & x = hk & \mapsto & k
\end{array}$$ If $x =hk, y = h_{1}k_{1} (h, h_{1} \in H, k,k_{1}
\in K)$ are any two elements of $G$, then$$\begin{array}{cclc}
  \alpha(xy) & = & \alpha(hkh_{1}k_{1}) &      \\
   & = & \alpha(hkh_{1}k^{-1}kk_{1}) &    \\
   & = & kk_1 &    \mbox{since }hkh_{1}k^{-1} \in H = \ker\alpha, \mbox{ as } H\lhd G \\
   & = & \alpha(x)\alpha(y)   &
\end{array}$$ Therefore $\alpha$ is an
epimorphism and $\ker\alpha = H$. As for the inclusion map $\beta : K
\hookrightarrow G$, $\alpha\beta = id$,  the result follows. \br In the
converse of above lemma,  if $K$ is also normal in $G$, then $G \cong H
\times K$. Further, if $G$ is direct product of its subgroups $H$ and $K$
then $G$ is a split extension of $H$ by $K$. \er
\bd\label{d0001} Let $H$ and $K$ be two subgroups of a group $G$.
Then $G$ is called semi-direct product of $H$ and $K$ if $H \lhd G, G =
HK$ and $H \bigcap K = \{e\}$.\ed  {\bf Notation } \rm : If $G$ is
semi-direct product of its subgroups $H$ and $K$, then we write $G = H
\propto K$.
\bexe Show that the symmetric group $S_{3}$ is semi-direct product of the
subgroups $$H = \left\{ id, \left(\begin{array}{ccc}
  1 & 2 & 3 \\
  3 & 1 & 2
\end{array}\right), \left(\begin{array}{ccc}
  1 & 2 & 3 \\
  2 & 3 & 1
\end{array}\right)\right\}$$ and$$K = \left\{ id, \left(\begin{array}{ccc}
  1 & 2 & 3 \\
  2 & 1 & 3
\end{array}\right)\right\}$$ of $S_{3}$.\eexe
We shall, now, consider a special extension of a group $G$. For the group
of automorphisms of $G$ i.e., $ Aut G,$ consider the cartesian product $$G
\times Aut G = \{(x, \alpha) \mid x \in G, \alpha \in Aut G \},$$  and
define a binary operation on this set as below : $$(x, \alpha )(y,
\beta) = (x\alpha(y), \alpha\beta)$$ We have, $$(x, \alpha)(e, id) = (x,
\alpha) = (e, id)(x, \alpha)$$,$$(x, \alpha)(\alpha^{-1}(x^{-1}),
\alpha^{-1}) = (e, id) = (\alpha^{-1}(x^{-1}),\alpha^{-1})(x,\alpha)$$
and$$\begin{array}{rcl}
  ((x, \alpha)(y, \beta))(z, \gamma) & = & (x\alpha(y)\alpha\beta(z), (\alpha\beta)\gamma) \\
   & = & (x, \alpha)((y, \beta)(z, \gamma))
\end{array}$$ Hence $G \times Aut G$ is a group with respect to
the binary operation. This is called the {\bf Holomorph} of $G$, and is
denoted by$Hol G$. It is easy to check that the maps:
$$\begin{array}{crcl}
  \psi : & Aut G & \rightarrow & Hol G \\
   & \alpha & \mapsto & (e, \alpha)
\end{array}$$ and  $$\begin{array}{crcl}
  \theta : &  G & \rightarrow & Hol G \\
   & x & \mapsto & (x, id)
\end{array}$$  are group monomorphism (i.e.,embeddings) and for
any $\beta \in Aut G,$  $$\begin{array}{ccl}
  (e, \beta)(x, id)(e, \beta)^{-1} & = & (\beta(x), \beta)(e, \beta^{-1}) \\
   & = & (\beta(x), id)
\end{array}$$ Thus if we identify $G$ and $Aut G$ by their images
under $\theta$ and $\psi$ to subgroups of $Hol G$, then any
automorphism of $G$ is obtained as a restriction of an inner
automorphisms of $Hol G$. Further, as $$\begin{array}{ccl}
  (y, \beta)(x, id)(y, \beta)^{-1} & = & (y \beta(x), \beta)(\beta^{-1}(y^{-1}), \beta^{-1}) \\
   & = & (y \beta(x)y^{-1}, id)
\end{array}$$  and $$(y, \beta) = (y, id)(e, \beta)$$ We get
that$G \lhd Hol G$, $Hol G = G.Aut G$.  \\  As it is clear that $G
\bigcap Aut G = id$, we conclude  $$Hol G = G \propto Aut G$$ Let
us note that for any subgroup $A$ of $Aut G$, $GA$ is a subgroup
of $Hol G$ and is, of course, a split extension of $G$ by $A$. \bx
Let $G = (\Z, +)$, then $Aut \Z = \{\pm id\}$. Consider the map
: $$\begin{array}{crcl}
  \psi : & Hol \Z & \rightarrow & Gl_{2}(\Z) \\
   & (m, \pm id) & \mapsto & \left(\begin{array}{cc}
     1 & 0 \\
     m & \pm1 \\
   \end{array}\right)
\end{array}$$ We have $$\begin{array}{cclcr}
  \psi((m, id)(n, id)) & = & \psi(m+n, id)  &  &  \\
   & = & \left(\begin{array}{cc}
     1 & 0 \\
     m+n & 1 \\
   \end{array}\right) & = & \left(\begin{array}{cc}
     1 & 0 \\
     m & 1 \\
   \end{array}\right)\left(\begin{array}{cc}
     1 & 0 \\
     n & 1 \\
   \end{array}\right) \\
   & = & \psi(m, id)\psi(n, id) &  &
\end{array}$$     $$\begin{array}{cclcc}
  \psi((m, id)(n, -id)) & = & \psi(m+n, -id) &  &  \\
   & = & \left(\begin{array}{cc}
     1 & 0 \\
     m+n & -1 \\
   \end{array}\right) & = & \left(\begin{array}{cc}
     1 & 0 \\
     m & 1 \\
   \end{array}\right)\left(\begin{array}{cc}
     1 & 0 \\
     n & -1 \\
   \end{array}\right) \\
   & = & \psi(m,id)\psi(n, -id) &  &
\end{array}$$ and  $$\begin{array}{cclcc}
  \psi((m, -id)(n, -id)) & = & \psi(m-n, id) &  &  \\
   & = & \left(\begin{array}{cc}
     1 & 0 \\
     m-n & 1
   \end{array}\right) \\
   & = & \left(\begin{array}{cc}
     1 & 0 \\
     m & -1
   \end{array}\right)\left(\begin{array}{cc}
     1 & 0 \\
     n & -1
   \end{array}\right) \\
   & = & \psi(m,-id)\psi(n, -id) &  &
\end{array}$$  Hence $Hol \Z$ is isomorphic to the subgroup
$$H = \left\{\left(\begin{array}{cc}
  1 & 0 \\
  m & e
\end{array}\right) \mid e = \pm 1, m \in \Z \right\}$$ of $Gl_{2}(\Z)$.
\ex The following example gives a general formulation of the
 above example.  \bx\label{x00} Let $A$ be an additive abelian group. Then put
 $$K = \left\{\left(\begin{array}{cc}
  id & 0 \\
  x & \alpha
\end{array}\right) | x \in A, \alpha \in Aut A, id :
\mbox{ identity automorphism of } A \right\}$$ Define a binary
operation over $K$ as below :  $$\left(\begin{array}{cc}
  id & 0 \\
  x & \alpha
\end{array}\right)\ast\left(\begin{array}{cc}
  id & 0 \\
  y & \beta
\end{array}\right) = \left(\begin{array}{cc}
  id & 0 \\
  x+\alpha(y) & \alpha\beta
\end{array}\right)$$ Then $(K, \ast)$ is a group , and the map : $$\begin{array}{crcl}
  \psi : & Hol A & \rightarrow & K \\
  & (x,\alpha) & \mapsto & \left(\begin{array}{cc}
    id & 0 \\
    x & \alpha \
  \end{array}\right)
\end{array}$$ is an isomorphism.\ex  \bx If $\Q = (\Q, +)$, the
additive abelian group of rational numbers, then $Aut \Q$
isomorphic to $\Q^{\ast}$, the multiplicative group of non-zero
rationals.Hence by example \ref{x00} $$Hol{ \Q} \cong
\left\{\left(\begin{array}{cc}
  1 & 0 \\
  s & r
\end{array}\right) | s \in \Q, r \in \Q^{\ast}\right\},$$  a
subgroup of $Gl_{2}( \Q)$.\ex
\section{ Wreath Product }
Let $G$, $H$ be groups. Put $$G^{[H]} = \{f : H \rightarrow G\}$$ the set
of all set maps from $H$ to $G$, and $$G^{(H)} = \{f : H \rightarrow G
\mid \abs{supp(f) }
< \infty \}$$ Clearly, $G^{[H]}$ is the cartesian product of copies of
the group $G$ indexed by $H$ and $G^{(H)}$ is the direct product of copies
of $G$ indexed by $H$. For each $f \in G^{[H]}$ and $h \in H$, define :
$$f^{h}(x) = f(xh) \mbox{ for all } x \in H$$ Then for the function
$$\begin{array}{crcl}
  \widehat{h} : & G^{[H]} & \rightarrow & G^{[H]} \\
   & f & \mapsto & f^{h}
\end{array}$$ for any $f,g \in G^{[H]}$, and $x \in H$, $$\begin{array}{ccl}
  (fg)^{h}(x) & = & (fg)(xh) \\
   & = & f(xh)g(xh) \\
   & = & f^{h}(x)g^{h}(x) \\
   & = & (f^{h}g^{h})(x)
\end{array}$$ $$\begin{array}{cccc}
  \Rightarrow & (fg)^{h} & = & f^{h}g^{h}
\end{array}$$ Hence $\widehat{h}$ is a homomorphism.  \\  Clearly
$\widehat{h}\widehat{h^{-1}} = id = \widehat{h^{-1}}\widehat{h}$.
Therefore $\widehat{h}$ is an automorphism of $G^{[H]}$, moreover,
$\widehat{h}(G^{(H)}) \subset G^{(H)}$. Thus restriction of $\widehat{h}$
to $G^{(H)}$ gives as automorphism of $G^{(H)}$. Now, consider the map :
$$\begin{array}{crcl}
  \tau : & H & \rightarrow & Aut G^{[H]} \\
   & h & \mapsto & \widehat{h}
\end{array}$$ For any $x \in H$, $h,h_{1} \in H$ and $f \in Aut
G^{[H]}$,  $$\begin{array}{ccl}
  \widehat{hh_{1}}(f)(x) & = & f^{hh_{1}}(x) \\
   & = & f(xhh_{1}) \\
   & = & (f^{h_{1}})^{h}(x) \\
   & = & \widehat{h}(\widehat{h_{1}}(f))(x)
\end{array}$$ $$\begin{array}{rrcl}
  \Rightarrow & \widehat{hh_{1}} & = & \widehat{h}\widehat{h_{1}} \\
  i.e., & \tau(hh_{1}) & = & \tau(h)\tau(h_{1})
\end{array}$$ thus $\tau$ is a homomorphism from $H$ to $Aut
G^{[H]}$. It is easy to check that if $\widehat{h} = id$, then $h = e$.
Therefore $\tau$ is a monomorphism (i.e., embedding ). Identify $H$ under
$\tau$ to the subgroup $\tau(H)$ of $Aut G^{[H]}$. Then the subgroup
$G^{[H]}.H$ of $Hol G^{[H]}$ is called the unrestricted wreath product of
the group $G$ by the group $H$. It is denoted by $G Wr H$. In fact, this
is the set product $G^{[H]} \times H$ with the binary operation  $$(f,
h)(g,h_1) = (fg^{h}, hh_{1})$$ The group $G^{[H]}$ is called base group of
$G Wr H$ and $H$ is called its top group. We have noted that for any $h
\in H$, the restriction of $\widehat{h}$ to $G^{(H)}$ is an automorphism
of $G^{(H)}$. It is easy to see that the set product $G^{(H)} \times H$ is
a subgroup of $G Wr H$. We call this subgroup, the restricted wreath
product of $G$ by $H$ and is denoted by $Gwr H.$ Here again $G^{(H)}$ is
called base group of $G wr H$ and $H$ is  called its top group. By the
definition of $G Wr H$ and $G wr H$, it is obvious that $G Wr H = G wr H$
if and only if either $H$ is finite or $G = \{e\}$.  \\  The diagonal
subgroup of $G^{[H]}$ i.e., $$Diag(G^{[H]}) = \{f \in G^{[H]} | f(x) = a
(\mbox{ fixed }) \mbox{ for all } x \in H\}$$ is also called the diagonal
subgroup of $G Wr H$. \bexe Show that if $G$ and $H$ are torsion groups
i.e., all elements in $G$ and $H$ have finite order, then $G wr H$ is a
torsion group.\eexe  \bl\label{l0001} Let $G$ and $H$ be two groups such
that $\circ(G) > 1$ and $\circ(H) = \infty$. Then $Z(G wr H) = id$.\el \pf
We know that $G wr H$ is the cartesian product $G^{(H)} \times H$ where
for any $(f,h),(g, h_1)$ in $G^{(H)} \times H$, the binary operation is
defined as : $$(f, h)(g, h_1) = (fg^{h}, hh_{1})$$ Hence if $(f, h) \in
Z(G wr H)$, then for any $g \in G^{(H)}$ and $h_{1} \in H$, $$(fg^{h},
hh_{1}) = (gf^{h_{1}}, h_{1}h)$$ $\Rightarrow$ $h \in Z(H)$ and $fg^{h} =
gf^{h_{1}}$ for all $h_{1} \in H$, and $g \in G^{(H)}$.  \\  If $g = id$,
then $g^{h} = id$, and we get $f = f^{h_{1}}$ for all $h_{1} \in H$.
Therefore for any $x \in H$, $$\begin{array}{crclcl}
   & f(x) & = & f(xh_{1}) &  & \forall h_{1} \in H \\
  \Rightarrow & f & \in & Diag G^{(H)} &  &  \\
  \Rightarrow & f & = & id &  & \mbox{since } \circ(H) = \infty
\end{array}$$ Now, as $f = id$, $(id, h) \in Z(G wr H)$ implies
$g^{h} = g$ for all $g \in G^{(H)}$. This can hold only if $h = e$. Hence
$Z(G wr H) = id$. \bco\label{c00} In the lemma if we take $H$ to be any
group i.e., order of $H$ is not necessarily $\infty$, then $$Z(G Wr H) =
Diag (Z(A)^{[H]})$$ \eco   \pf From the proof of the lemma it is clear
that if $(f, h) \in Z(G Wr H)$ then $f \in Diag G^{[H]}$, $h \in Z(H)$ and
$fg^{h} = gf^{h_{1}}$ for all $h_1 \in H$ and $g \in G^{[H]}.$ Let $f(x)=a
\in G$ for all $x \in H$. Then for all $x, h_{1} \in H$ and $g \in
G^{[H]}$, $$\begin{array}{crclcccc}
   & f(x)g^{h}(x) & = & g(x)f^{h_{1}}(x) &  &  &  &  \\
  \Rightarrow & ag(hx) & = & g(x)a &  &  &  &  \\
  \Rightarrow & a & \in & Z(G) & \mbox{and} & g(hx) & = & g(x) \\
  \Rightarrow & a & \in & Z(G) & \mbox{and} & h = e &  &
  \mbox{since } \circ(G) > 1.
\end{array}$$ Hence the result is proved. \\ \clearpage
 {\bf EXERCISES}
\begin{enumerate}
  \item
  Prove that $\C^*$ is isomorphic to $S^1 \times \R^{*}_{+}$.
  \item
  Show that $S_n (n\geq 2)$ is indecomposable group i.e., is not a direct
  product of its proper subgroups.
  \item
  Prove that the group $\C_{\infty}$ is a direct product of its
  proper subgroups.
  \item
  Let $G$ be a torsion free abelian divisible group and $x (\neq e) \in G$.
   Prove that the map :  $$\begin{array}{rcl}
     \varphi : (\Q,+) & \rightarrow & G \\
     \frac{m}{n} & \mapsto & y \\
   \end{array}$$ where $y^n = x^m$, is a monomorphism. Then deduce that any
   uncountable torsion free abelian divisible group is
   decomposable.
   \item
   Show that $(\R^{*}_{+},\cdot)$ is decomposable.
   \item
   Let $G$ be a divisible abelian group and $x(\neq e) \in G$ has finite order.
   Prove that $G$ contains a subgroup isomorphic to $\C_{p^{\infty}}$ for
    some prime $p$.
   \item
   Let $G$ be a abelian group such that $\circ(x)<\infty$ for all $x \in G$
    (such a group is called {\bf periodic}). For any prime $p$;  \\
    put $$G_{p^{\infty}} = \{x \in G \st x^{p^{m}}=e \mbox{ for some }
    m\geq 0\}$$ Show that $G_{p^{\infty}}$ is a subgroup of $G$ and
    $G=\prod_{p:\mbox{ prime }}G_{p^{\infty}}.$ Thus, in particular ,
    $\C_{\infty} = \prod_{p:\mbox{ prime }}\C_{p^{\infty}}.$
    \item
    Let $G$ be a finite group, and let for subgroups $H_1,H_2,\cdots,H_n$
     of $G$, $G$ is internal direct product of $H_i's$. Show that if
     $(\circ(H_i),\circ(H_j)) = 1$ for all $1\leq i\neq j \leq n,$ then
     $Aut\,G \cong Aut\,H_1 \times \cdots \times Aut\, H_n$ i.e., $Aut\, G$
     is external direct product of $Aut\, H_i$.
     \item
     Prove that if $n\geq 3$ is odd, then $O_n(\R)$ is isomorphic to
     $SO_n(\R) \times \Z_2.$
     \item Show that $(\R, +)$ can be expressed as direct sum of two proper
    subgroups.
    \item Prove that $Sl_{n}(\C)$ is a direct factor of
$Gl_{n}(\C)$ i.e., there exists a subgroup $H$ of $Gl_{n}(\C)$ such that
$Sl_{n}(\C) \times H$ is isomorphic to $Gl_{n}(\C)$.
  \item Let $G_{i} , i
= 1, 2, \cdots, n$ be groups and $H_{i}$ be subgroups of $G_{i}$ for all
$i = 1, 2, \cdots, n$. Prove that $H_{1} \times H_{2} \times \cdots \times
H_{n}$ is a normal subgroup of $G_{1} \times \cdots \times G_{n}$ if and
only if $H_{i} \triangleleft G_{i}$.
    \item Let $G_{i} , i = 1, 2, \cdots,
n$ be groups and let $H_{i} \triangleleft G_{i}$ for all $i = 1, 2,
\cdots, n$. Prove that \\ (a) $G_{1}/ H_{1} \times \cdots G_{n} / H_{n}$
is isomorphic to $G / H$ where $G = G_{1} \times \cdots \times G_{n}$ and
$H = H_{1} \times \cdots \times H_{n}$.\\ (b) The direct product $G_{1}
\times \cdots \times G_{n}$ is finite if and only if each $G_{i}$ is
finite. \\ (c) The group $G_{1} \times \cdots \times G_{n}$ is abelian if
and only if each $G_{i}$ is abelian.
  \item Let a group $G$ be internal
direct product of it's subgroups $G_{i}, i = 1, 2, \cdots, n$. If for any
$x \in G, x = x_{1}x_{2} \cdots x_{n} (x_{i} \in G_i)$, then show that  \\
(i) The map $$\begin{array}{cccc}
  p_{i} : & G & \rightarrow & G_{i} \\
   & x & \mapsto & x_{i}
\end{array}$$ is an epimorphism.   \\   (ii) $p_{i}p_{j}(x) = e$
for all $i \neq j$.  \\   (iii) $p_{1}(x)p_{2}(x) \cdots p_{n}(x) = x$ for
all $x \in G$.  \\  (iv) $\bigcap_{i = 1}^{n} \ker{p_{i}} = e$.
  \item Let for a group $G$, we have endomorphisms $$p_{i}
: G \rightarrow G_i\;\; ; i = 1, 2, \cdots, n$$ such that for all $x \in
G$, we have (i) $p_{i}p_{j}(x) = e$ whenever $i \neq j$ (ii) $p_{i}^{2}(x)
= p_{i}(x)$ (iii)  $\prod p_{i}(x) = p_{1}(x)p_{2}(x) \cdots p_{n}(x) = x$
and (iv) $\bigcap_{i = 1}^{n} \ker{p_{i}} = e$.
\\ Prove that if $H_{i} = p_{i}(G)$, then $G$ is internal direct
product of $H_{i}'^{s}$.
\item
 Prove that the group $\C_{p^{\infty}}$ is an indecomposable group.
 \item
 Prove that a finite abelian group is semi-simple if and only if every
 element of the group is of square free order.
 \item
 Show that a cyclic group is semi-simple if and only if it is finite of
 square free order.
 \item
 Let $G$ be an additive abelian group and $x$, an element of finite order
 in $G$. Prove that if $G=gp\{x\}\bigoplus H$ for a subgroup $H$ of $G$,
 then for any $y\in H, \circ(y)$ dividing $\circ(x),
 G=gp\{x+y\}\bigoplus H.$
 \item Let $G$ be a finite, additive abelian group with $\circ(G)=p^n$ for
 a prime $p$ and $n\geq 1.$ If for an element $x\in G$,
 $G=gp\{x\}\bigoplus H$ for a subgroup $H$ of $G$, then for any $y \in H$
 with $\circ(y)\leq \circ(x),$ $G=gp\{x+y\}\bigoplus H.$
 \item
 Let $n\geq 3.$ Prove that for the Dihedral group $D_n$ (Example
 \ref{ex9}), the group $Aut \,D_n$ is semi-direct product of its
 subgroups,\[H=\{\alpha\in Aut\,D_n\st \alpha(\sigma)=\sigma\}\] and
 \[K=\{\alpha\in Aut\,D_n\st \alpha(\tau)=\tau\}.\] Further show that
 $\circ(H)=n$ and $\circ(K)=\Phi(n).$
 \item
 Let $H,K$ be two subgroups and let \[\begin{array}{rcl}
   f:H & \rightarrow & Aut\,K \\
   y & \mapsto & f_y \
 \end{array}\] be a homomorphism. Put \[G=H\times K=\{(x,u)\st x\in H,u\in
 K\}\] If for $(x,u),(y,v)$ in $G$, we define :
 \[(x,u)*(y,v)=(xy,(f_y(u))v),\] then show that $G$ is a group with
 respect to binary operation $*$, and $G$ is semi-direct of its subgroups
 $A=\{(x,e)\st x\in H\}$ and $B=\{(e,y)\st y\in K\}.$ where $A$ is
 isomorphic to $H$ and $B$ is isomorphic to $K$.
 \item
 Let $m>1, N>1$ be two integrals such that $m|n$. Let $C_m=gp\{x\}$, and
 $C_n=gp\{y\}$ be two cyclic groups of orders $m$ and $n$ respectively.
 Prove that for the direct product $G=C_m\times C_n,$ \[\begin{array}{rcl}
   \theta : G & \rightarrow & G \\
   (x^s,y^t) & \mapsto & (x^s,y^{-t}) \
 \end{array}\] and \[\begin{array}{rcl}
   \alpha : G & \rightarrow & G \\
   (x^s,y^t) & \mapsto & (x^{s+t},y^t) \
 \end{array}\] are automorphisms of $G$. Further, if $n\neq 2,$ show that
 $\theta\alpha\neq\alpha\theta$.
 \item
 Let a group $G$ be direct product of its subgroups $H_i ; i=1,\cdots,n$. Prove that
  $\Phi(G)\subseteq\Phi(H_1)\times\cdots\times\Phi(H_n).$
  \item
  Let $G$ be a finitely generated group and let $G$ be direct product of
  its subgroups $H_i ; i=1,\cdots,n$. Prove that $\Phi(G)$ is direct
  product of $\Phi(H_i) ; i=1,\cdots,n$.
\end{enumerate}

\chapter{Action Of A Group}
 \footnote{\it contents group7.tex }
 In many practical situations we encounter groups and their effects
 on certain objects. We shall describe, below, this phenomenon in mathematical terms.
 \bd\label{d51} A group $G$ is said to
act on a non-empty set $X$ or is called a transformation group of $X$ if
there exists a homomorphism from $G$ into $A(X)$, the group of  all
transformations of $X.$ \ed {\bf Note:}   Given an action of a group $G$
on a set $X$, we shall denote the image of any $g \in G$ in $A(X)$ by $g$
itself, unless a specific notation is called for.
\brs (1)For a transformation group $G$ of a set $X$, we have
$$(i) (g_1g_2)(x) = g_1(g_2(x)) \mbox{ for all }g_1,g_2 \in G \mbox{ and }
x\in X$$ $$(ii) e(x) = x \mbox { where }e \mbox{ is the identity of } G
\mbox { and } x\in X.$$ (2)If $G$ is mapped to the Identity subgroup of
$A(X)$,then $G$ is said to have the trivial action on $X$ and in this case
$g(x) = x$ for all $x\in X.$   \\  (3) Clearly a group $G$ can act on $X$
in several ways depending on the homomorphism from $G$ to $A(X)$.  \\
 (4) If a group $G$ is acting on a set $X$, then for any subgroup $H$ of $G$, $H$ acts on $X$
  where action of any elements of $H$ is its action as an element of $G$.  \ers
\bx\label{ex51} Let $G$ be a group and $X = G$. Let for any $g\in
G$ \bea
             \tau_g : G & \to & G\\
                   x & \mapsto & gx
\eea Then \bea
             \tau : G & \to & A(G)\\
                   g & \mapsto & \tau_g
\eea defines a group action of $G$ over $G$. \ex \bx\label{ex52}
In the example \ref{ex51}, if we take $$\tau_g(x) = gxg^{-1}$$ then we get
another action of $G$ over $G$. \ex \bx\label{ex53} Let $G$ be a group and
$H$, a subgroup of $G$. Let $X = \{g_iH\}_{i\in I}$ be the set of all
distinct left cosets of $H$ in $G$. If for any $g\in G$, \bea
             \tau_g : X & \to & X\\
                   g_iH & \mapsto & gg_iH
\eea Then \bea
             \tau : G & \to & A(X)\\
                   g & \mapsto & \tau_g
\eea defines an action of $G$ over $X$. \ex \bx\label{ex54} Let
$G = Gl_n(\R)$ and $X = \R^n$, the vector space of dimension $n$ over
$\R$. Let for any $A\in Gl_n(\R)$, \bea
             \tau_A : X & \to & X \\
 \underline{\lambda} = (\lambda_1, \lambda_2, \ldots ,\lambda_n)^{t} & \mapsto &
 A\underline{\lambda}
\eea Then \bea
             \tau : Gl_n(\R) & \to & A(X)\\
                   A & \mapsto & \tau_A
\eea defines a group action of $Gl_n(\R)$ over $X$. \ex
\bx\label{ex55} Let $G = S_n$, be the symmetric group on $n$-symbols
$\{1,2,\cdots,n\}$ and $X = \R^n$. For any $\sigma \in S_n$, let \bea
             \tau_{\sigma} : X & \to & X\\
\underline{\lambda} = (\lambda_1, \lambda_2, \ldots
,\lambda_n)^{t} & \mapsto & (\lambda_{\sigma(1)},
\lambda_{\sigma(2)}, \ldots ,\lambda_{\sigma(n)})^{t} \eea Then
\bea
             \tau : S_n & \to & A(X)\\
               \sigma & \mapsto & \tau_{\sigma}
\eea defines a group action of $S_n$ over $X$.
 \bx\label{ex7.8x} Let $G$ and $H$ be two groups, and let $\varphi : G \rightarrow H$
  be a homomorphism. If $H$ acts on a set $X$, then for any $g \in G, x \in X$,
  define :  $$\tau_{g}(x) = \varphi(g)(x) $$ the image of $x$ under $\varphi(g).$
  Clearly $\tau_g \in A(X)$ and $$\begin{array}{rcl}
    G & \rightarrow & A(X) \\
    g & \mapsto & \tau_g \\
  \end{array} $$ defines a group action of $G$ on $X$.\ex  \bx\label{ex7.9e} Suppose
   a group $G$ acts on a set $X$. Then for any set $Y$, $G$ acts on $\mathcal{F} = \{f :
   X \rightarrow Y \st f : \mbox{ set function }\}$ as follows : \\
   Let $g \in G,$ and $f \in \mathcal{F}.$ Then for any $x \in X,$ define
   : $$ (\tau_g(f))(x) = f(g^{-1}(x)).$$ \ex
   \bx\label{ex7.10} The orthogonal group $O(n,\R)$ acts in a natural way
   on the unit circle $S^{n-1}$ in $\R^n$. For any $A \in O(n,\R)$, and $x
   = (x_1,x_2,\cdots,x_n) \in S^{n-1}$, we define $\tau_A(x) =
   A(x_1,x_2,\cdots,x_n)^t.$  \\
    (Use : orthogonal transformations preserve length.) \ex
    \bt\label{t52} Let
$G$ be a transformation group of a set $X$. Let for $x_1,x_2\in X,$ $$
x_1R x_2 \Leftrightarrow g(x_1) = x_2 \mbox { for some }g\in G.$$ Then $R$
is an equivalence relation on $X.$ \et \pf (Reflexivity) For any $x\in X,
e(x) = x.$ Hence $xRx.$  \\   (Symmetry) Let for $x_1,x_2 \in X, x_1Rx_2$
holds. Then there exists $g\in G$ such that \ba{rl}
                 & g(x_1) = x_2\\
\Rightarrow      & g^{-1}(g(x_1)) = g^{-1}(x_2)\\ \Rightarrow & e(x_1) =
g^{-1}(x_2)\\ \Rightarrow      & x_1 = g^{-1}(x_2) \ea Therefore $x_2Rx_1$
holds. \\ (Transitivity) Let for $x_1,x_2,x_3 \in X$, $x_1Rx_2$ and
$x_2Rx_3$ hold. Then there exist $g,h\in G$ such that \ba{rl}
                &  g(x_1) = x_2 \mbox { and }h(x_2) = x_3\\
\Rightarrow     & (hg)(x_1) =h(g(x_1)) = x_3\\ \Rightarrow     & x_1Rx_3.
\ea Hence $R$ is an equivalence relation on $X.$      $\blacksquare$
\bt\label{t53} {\bf (Caley's Theorem) } Any group $G$ is isomorphic to a subgroup of the
transformation group of $G$. Thus, in particular, if $\circ(G) = n$, then
$G$ is isomorphic to a subgroup of $S_{n}$.
\et \pf For any $g\in G$, define $\tau_g
: G \to G$ by $\tau_g (x) = gx.$ Then $\tau_g$ is a one-one, onto mapping
i.e., $\tau_g \in A(G).$ Consider the mapping $$ g \mapsto \tau_g$$ from
$G$ to $A(G)$. For $g,h,x \in G$, we have \ba{rl}
               &  \tau_{gh}(x) = (gh)(x)\\
\Rightarrow    &  \tau_{gh}(x) = (g(h(x)) = \tau_g (\tau_h (x)).
\ea Therefore $\tau_{gh} =  \tau_g \circ \tau_h$ for all $g,h \in
G.$ Thus  $g \mapsto \tau_g$ is a homomorphism from $G$ to $A(G).$
Let $\tau_g$ be identity for some $g\in G.$ Then, for any $x\in
G,$ \ba{rl}
                &  \tau_g(x) = x\\
\Rightarrow     & gx = x\\ \Rightarrow     & g = e (\mbox{ taking} x = e).
\ea Hence $g \mapsto \tau_g$ is an isomorphism of $G$ onto its image in
$A(G)$. The final part of the statement is immediate because if $\circ(G)
= n$, then $A(G)$ is isomorphic to $S_{n}$.      $\blacksquare$
\br The isomorphism of $G$ into the transformation group provided in the
above theorem is quite inefficient in the sense that if we have $\circ (G)
= n$, then $\circ (A(G)) = n!.$ Thus $G$ is mapped into a large
transformation group. The situation can be improved sometimes. \er The
following result shall be useful in this direction. \bt\label{t54} Let $G$
be a group and $H$ be a subgroup of $G.$ Let $X = \{g_iH \st i\in I\}$ be
the set of all distinct left cosets of $H$ in $G.$ Then there exists a
homomorphism $\theta$ from $G$ into $A(X)$ with $$ \ker\theta = \cap
\{x^{-1} Hx \st x\in G\}$$ \et \pf Let $g\in G.$ Then for each $g_iH\in X,
g(g_iH) = gg_iH\in X.$ Further, for $g_iH, g_jH \in X$, \ba{rl}
                  &  g(g_iH) = g(g_jH)\\
\Rightarrow       &  g_iH = g_jH \ea We also have $g(g^{-1}g_iH) =
g_iH.$ Hence, the map \bea
                  \theta_g : X & \to & X\\
            g_iH & \mapsto & gg_iH
\eea is in $A(X).$ One can easily see, as in the above theorem, that
the map $\theta (g) = \theta_g$ is a homomorphism from $G$ to $A(X).$
Further, $g\in \ker\theta$ if and only if $\theta_g = Id.$ Thus $g\in
\ker\theta$ if and only if for any $g_iH\in X$, $\theta_g(g_iH) = g_iH$.
Note that \ba{rl}
                &  \theta_g(g_iH) = g_iH\\
\Leftrightarrow & gg_iH = g_iH\\ \Leftrightarrow & g_i^{-1}gg_i \in H\\
\Leftrightarrow & g\in g_iHg_i^{-1} \ea Now, if $x\in G,$ then $x\in g_iH$
for some $i$. Thus $x = g_ih$ for some $h\in H.$ We then have $$ xHx^{-1}
= g_ihHh^{-1}g^{-1} = g_iHg_i^{-1}.$$ Consequently $$ g\in \ker\theta
\Leftrightarrow g\in \bigcap_{x\in G} x^{-1}Hx.$$ Thus  the proof is
complete. $\blacksquare$ \\ We,now,give some applications of this result.
\bt\label{t55p} Let $G$ be a finitely generated group and $n\geq 1$, an
integer. Then the set of all subgroups of $G$ with index $n$ is finite.\et
\pf Let $H$ be a subgroup of index $n$ in $G$ and let $X = \{g_{1}H = H,
g_{2}H, \cdots, g_{n}H\}$ be all distinct left cosets of $H$ in $G$. Then,
by the theorem \ref{t54}, for any $x \in G$, the map
$$\begin{array}{ccrcl}
  \theta_x & : & X & \rightarrow & X \\
   &  & g_{i}H & \mapsto & xg_{i}H
\end{array}$$ is a permutation on $X$ and $$\begin{array}{ccccr}
  \theta & : & G & \rightarrow & A(X) \\
   &  & x & \mapsto & \theta_x
\end{array}$$ is a homomorphism. Writing $i$ for $g_{i}H$ for all $1\leq i\leq n$,
$\theta$ gives the homomorphism $$\begin{array}{ccccc}
  \theta_H & : & G & \rightarrow & S_n
\end{array}$$ where $\theta_{H}(x)(i) = j$ if $xg_{i}H = g_{j}H$.
Now, if $H \neq K$ are two subgroups of $G$ with index $n$, then $H
\not\subset K$, and $K\not\subset H$. Choose $h \in H$, $h \not\in K$.
Then $\theta_{K}(h)(1) \neq 1$ and $\theta_{H}(h)(1) = 1$. Thus
$\theta_{H} \neq \theta_{K}$. As $G$ is finitely generated any
homomorphism from $G$ to $S_{n}$ is completely determined by its action
over the generators of $G$. Thus the number of homomorphism from $G$ to
$S_{n}$ is finite. Hence the result is immediate.  $\blacksquare$
\bco\label{l7.1l} Let $G$ be a finitely generated group and $n\geq 1$, an
integer. Let $$\mathcal{A}_{n} = \{H \subset G : \mbox{ subgroup }|\;
[G:H] = n\}$$ Then for any epimorphism (onto, homomorphism) $f$ from $G$
to $G$ the map
: $$\begin{array}{ccrcl}
  \widehat{f} & : & \mathcal{A}_{n} & \rightarrow & \mathcal{A}_{n} \\
   &  & H & \mapsto & f^{-1}(H)
\end{array}$$ is a bijection (one-one, onto map).\eco \pf If $H \in \mathcal{A}_{n}$, then
 $$G = \bigcup^{n}_{i=1} x_{i}H$$ for some $x_{i} \in G$, $1\leq i\leq n$, where $x_{i}H\bigcap
 x_{j}H = \emptyset$ for all $i\neq j$. Hence $$G = f^{-1}(G) = \bigcup^{n}_{i=1}y_{i}f^{-1}(H)$$
  where $f(y_{i}) = x_{i}$, is a disjoint union. Therefore $[G:f^{-1}(H)] = n$. Thus $\widehat{f}$
  maps $\mathcal{A}$ to $\mathcal{A}$. Further, let for $H,K \in \mathcal{A}_{n}$, $f^{-1}(H) = f^{-1}(K)$. Then as $f$
   is an epimorphism $H = f(f^{-1}(H)) = f(f^{-1}(K)) = K$. Thus $\widehat{f}$ is one-one. Now,
   as $\mathcal{A}_{n}$ is finite by the theorem, $\widehat{f}$ is a bijection.  $\blacksquare$
   \bd\label{d7.1d} A group $G$ is called residually finite if the intersection
    of all subgroups of $G$ of finite index in $G$ is identity. \ed \bx Any finite group is residually
finite. \ex \bx Any infinite cyclic group is residually finite.
\\ Let $C = gp\{a\}$ be an infinite cyclic group. Then any subgroup of $C$ is of the
 from $gp\{a^{n}\}$, $n \geq 1$. Clearly $\bigcap_{n\geq 1} gp\{a^{n}\} = Id$.
 Hence $C$ is residually finite. \ex
  \bt\label{t55} Let $G$ be a simple group of order 60. Then $G$
contains no subgroup of order 30 or 20 or 15. \et \pf We shall show that
there is no subgroup of order 15. Proof for the other cases is similar.
Let $H$ be a subgroup of order 15. Then $[G:H] = 4.$ Hence, by the theorem
\ref{t54} we have a homomorphism $$ \theta
: G \to S_4$$ Let $K$ be the kernel of $\theta$. Then $K$ is a normal subgroup of $G$
 contained in $H$. As $G$
is simple, $K = G$ or $K = id.$ As $K \subseteq H.$, $K = id.$ Thus
$\theta$ is an isomorphism of $G$ into $S_4$. As $\circ (S_4) = 24$, in
view of theorem \ref{co22}, $S_4$ can not contain a subgroup of order 15.
 Thus the result follows by contradiction.  $\blacksquare$
\bd\label{d56} Let $G$ be a group of transformations on a set $X.$
Then for any $x\in X,$ the orbit of $x$ under $G$-action is $$ \mbox {\bf
Orb}(x) = \{g(x) \st g\in G \}$$ \ed \brs (1) If $R$ is the equivalence
relation defined in theorem \ref{t52}, then $\mbox {\bf Orb}(x) = \{y\in X
\st xRy \},$ the set of all elements of $X$ which are related to $x$ under
the action of $G.$\\ (2) As $R$ is an equivalence relation on $X$, any two
orbits are either identical or disjoint.
 \ers \bd\label{d57} Let $G$ be a transformation group of a set
$X$. A subset $Y$ of $X$ is called invariant under $G$ if $g(Y) \subseteq
Y$ for all $g\in G.$ \ed \bd\label{d58} Let $G$ be a transformation group
of a set $X$. For any $x\in X,$ the set $\mbox{ \bf Stab}(x) = \{g\in G
\st gx = x \}$ is called the Stabilizer of $x$ in $G.$ \ed

\bl\label{l510} Let $G$ be a transformation group of a set $X$.
For any $x\in X,$ $\stab(x)$ is a subgroup of $G$ and $\abs{\orb(x)} =
[G:\stab(x)].$ Thus $\abs{\orb(x)}$ divides $\circ(G).$ \el \pf Clearly
$e\in \stab(x)$. Further, if $g_1,g_2\in \stab(x),$ then $$ (g_1g_2)(x) =
g_1(g_2(x)) = g_1(x) = x$$ and $$ g_1^{-1}(x) =  g_{1}^{-1}(g_1(x)) = x.$$
Hence $\stab(x)$ is a subgroup of $G$. To prove $\abs{\orb(x)} =
[G:\stab(x)],$ let $g_1,g_2\in G$ . Then for any $x\in X,$ \ba{rl}
            & g_1(x) = g_2(x)\\
\Leftrightarrow & g_2^{-1}(g_1(x)) = x\\ \Leftrightarrow & g_2^{-1}g_1 \in
\stab(x)\\ \Leftrightarrow & g_1\stab(x) = g_2\stab(x). \ea Hence
$\abs{\orb(x)}= [G:\stab(x)].$ Finally, as $\abs{G} = [G
: \stab(x)]\abs{\stab(x)},$ the result follows.     $\blacksquare$
\bl\label{l59} Let $G$ be a group acting on a set $X$. Let for $x,y \in
X$, $y = g_1(x)$, for some $g_{1} \in G$. Then $\stab(y) =
g_1\stab(x)g_1^{-1}$ and hence $\abs{\stab(x)} = \abs{\stab(y)}$ \el \pf
As $g_1(x) = y$, for any $g \in \stab(x)$, \ba{rl}
               & (g_1gg_1^{-1})(y) = (g_1gg_1^{-1})(g_1(x)) = (g_1(x)) = y\\
\Rightarrow    & g_1gg_1^{-1} \in \stab(y) \ea Hence
\begin{equation}\label{e51}
 g_1\stab(x)g_1^{-1} \subseteq \stab(y)
\end{equation}
Now, as $g_1(x)=y, x=g_{1}^{-1}(y).$ Therefore as above \ba{rl}
              & g_1^{-1}\stab(y)g_1 \subseteq \stab(x) \ea
\be\label{e52}
 \Rightarrow \;\;\;  \stab(y) \subseteq g_1\stab(x)g_1^{-1}
\ee Now, from the equations \ref{e51} and \ref{e52}, we get $\stab(y) =
g_1\stab(x)g_1^{-1}$. Thus \bea
               \stab(x) & \to & \stab(y)\\
               g & \mapsto & g_1gg_1^{-1}
\eea is one-one and onto, and hence $\abs{\stab(x)} =
\abs{\stab(y)}$.  $\blacksquare$ \br If $Y$ is an invariant subset of $X$,
then for any $y\in Y$, $ \orb(y) \subseteq Y$ as $g(y) \in Y$ for all
$g\in G.$ Hence $Y$ is a union of orbits. \er We shall now give some
interesting group actions on $\R^2$. \bx\label{ex56} For each $b\in \R^*$,
define $ \phi_b : \R^2 \to \R^2 \mbox { by }\phi_b((x,y)) = (b^2x,by).$ It
is easy to see that $\phi_b$ gives a transformation of $\R^2$. Further,
$$(\phi_c \phi_b)(x,y) = \phi_c(b^2x,by) = (c^2b^2x,cby) =
\phi_{cb}(x,y).$$ Clearly $\phi_1 = id.$ Thus we get an action of $\R^*$
on $\R^2$ where each $b\in \R^*$ acts via $\phi_b.$ In this case, for any,
$a\in \R^*$, $\orb(a,2a) = \{(b^2a,2ab) \st b\in \R^* \}$ is the parabola
$Y^2 = 4aX$ except the point $(0,0)$.
\ex
\bx\label{ex57} For each  $\Theta \in \R$, define \bea
       \delta_{\Theta} : \R^2 & \to & \R^2\\
       (x,y) & \mapsto & ( x\cos\Theta+y\sin\Theta,y\cos\Theta-x\sin\Theta).
\eea For any $\Theta_1, \Theta_2\in \R$, we have \bea & &
(\delta_{\Theta_1} \delta_{\Theta_2})(x,y)\\ & =
&\delta_{\Theta_1}(x\cos\Theta_2+y\sin\Theta_2 ,
y\cos\Theta_2-x\sin\Theta_2)\\ & = &
(((x\cos\Theta_2+y\sin\Theta_2)\cos\Theta_1+(y\cos\Theta_2-x\sin\Theta_2)\sin\Theta_1),\\
&  &
((y\cos\Theta_2-x\sin\Theta_2)\cos\Theta_1-(x\cos\Theta_2+y\sin\Theta_2)\sin\Theta_1))\\
& = &
(x\cos(\Theta_1+\Theta_2)+y\sin(\Theta_1+\Theta_2),y\cos(\Theta_1+\Theta_2)-x\sin(\Theta_1+\Theta_2))\\
& = & \delta_{\Theta_1+\Theta_2}(x,y). \eea Further, $\delta_0 =
id.$ Thus $\delta : (\R,+) \to A(\R^2)$,
 where $\delta(\Theta) = \delta_{\Theta}$, is a homomorphism i.e., $(\R,+)$ is
 a group of transformations of $\R^2$, where $\Theta \in \R$ acts via
 $\delta_{\Theta}.$  For any $r > 0,$ we have
$$\orb((0,r)) = \{(r\sin\Theta,r\cos\Theta) \st 0 \leq \Theta \leq
2\pi \},$$
 the circle of radius $r$ in $\R^2$ with centre as origin.
\ex \bd\label{d510} Let $G$ be a group acting over a set $X(\neq
\emptyset)$. The action of $G$ is called transitive if for any
$x,y \in X$ there exists $g\in G$ such that $gx = y$. \ed
\bl\label{l511} Let $G$ be a finite group acting over a finite set
$X$ transitively. Then\\ (i) Order of $G$ is divisible by
$\abs{X}$.\\ (ii) For any normal subgroup $N$ of $G$, orbits of
$X$ under $N$-action have same cardinality. Hence the length of an
$N$-orbit divides $\abs{X}$. \el \pf (i) Let $x\in X$ be any
element. Then \ba{rl}
        & \abs{X} = \abs{Gx} = [G:\stab(x)] = \frac{\abs{G}}{\abs{\stab(x)}}\\
\Rightarrow &  \abs{X} \abs{\stab(x)} = \abs{G} \ea Hence (i) is
proved.\\ (ii) Let $x,y\in X$. Then there exists $g\in G$ such
that $y = gx$. Therefore $$ \abs{Ny} = \abs{Ngx} = \abs{(Ng)x} =
\abs{(gN)x} = \abs{g(Nx)} = \abs{Nx}$$ Hence all $N$-orbits of $X$
have the same cardinality (say)$l$. Since $X$ is
 a disjoint union of $N$-orbits of $X$, it is clear that $l$ divides
 $\abs{X}$.      $\blacksquare$
 \bl\label{l512}
 {\bf(Burnside)} Let $G$ be a finite group acting over a finite set $X$. Let
 $B_1, B_2, \ldots , B_k$ be all the distinct orbits of $X$ under $G$-action.
 If $\chi(g)$ denotes the cardinality of the set $\{x \in X \st gx = x \}$, then
 $$ k = \frac{1}{\abs{G}} \sum_{g\in G} \chi(g)$$
 \el
 \pf To prove the lemma, we shall calculate the number of elements of the set
 $$ M = \{(g,x) \st g\in G, x\in X, gx = x \}$$
 in two different ways. Let for any $g\in G$, $x\in X$,
 $$ A_g = \{(g,x) \st x\in X, gx = x \}$$
 $$ B_x = \{(g,x) \st g\in G, gx = x \}$$
 Then
\begin{equation}\label{e53}
             M = \bigcup_{g\in G} A_g = \bigcup_{x\in X} B_x
\end{equation}
are two disjoint decompositions of $M$ into it's subsets. Clearly $\chi(g)
= \abs{A_g}$ and $\abs{B_x} = \abs{\stab(x)}$. Hence from the equation
\ref{e53}, we get
\begin{eqnarray}\label{e54}
       \abs{M} & = & \sum_{g\in G} \abs{A_g} = \sum_{x\in X}
                                                        \abs{B_{x}}\nonumber\\
               & \Rightarrow & \sum_{g\in G} \chi(g)\nonumber = \sum_{x\in X} \abs{stab(x)}\nonumber\\
& = & \sum_{i=1}^k \sum_{x_i\in B_i} \stab(x_i)
\end{eqnarray}
By the lemma \ref{l59}, $\abs{\stab(x)} = \abs{\stab(y)}$ for any $x,y\in
B_i$. Hence, the equation \ref{e54} gives that $$ \sum_{g\in G} \chi(g) =
\sum_{i=1}^k \abs{B_i} \abs{\stab(x_i)},
                                                              (x_i\in B_i)$$
Now, by the lemma \ref{l510} $$ \abs{B_i} =
\frac{\abs{G}}{\abs{\stab(x_i)}}.$$ Hence \bea \sum_{g\in G} \chi(g) & = &
\sum_{i=1}^k \frac{\abs{G}}{\abs{\stab(x_i)}} = k\abs{G} \eea This implies
\bea
     k & = & \frac{1}{\abs{G}} \sum_{g\in G} \chi(g)
\eea Thus the result is proved.  $\blacksquare$\\ We shall, now, give an
application of Burnside's lemma. \bx\label{ex58} Find the number of all
distinct factorizations of 216 into three natural numbers when the order
of factors is ignored. \ex \bf Solution : \rm We have $$ 216 = 2^3 3^3$$
Hence any factor of $216$ is of the form $2^{\alpha} 3^{\beta}$. Let $$
216 = (2^{\alpha_1} 3^{\beta_1}) (2^{\alpha_2} 3^{\beta_2})(2^{\alpha_3}
3^{\beta_3})$$ where $ \alpha_i \geq 0, \beta_i \geq 0$. Then $$
\alpha_1+\alpha_2+\alpha_3 = 3 = \beta_1+\beta_2+\beta_3$$ Clearly, each
factorization of $216$ into three natural numbers correspond to a tuple
$(\alpha_1,\alpha_2,\alpha_3,\beta_1,\beta_2,\beta_3)$ of non-negative
integers such that $\alpha_1+\alpha_2+\alpha_3 = 3 =
\beta_1+\beta_2+\beta_3$, and conversely. Put $$X =
\{(\alpha_1,\alpha_2,\alpha_3,\beta_1,\beta_2,\beta_3)\st \sum_{i=1}^3
\alpha_i = 3 = \sum_{i=1}^3 \beta_i, \alpha_i(\geq 0), \beta_i (\geq 0)
\in \Z \}$$ Define for each $\sigma \in S_3$, $$ \sigma
(\alpha_1,\alpha_2,\alpha_3,\beta_1,\beta_2,\beta_3)
   = (\alpha_{\sigma(1)},\alpha_{\sigma(2)},\alpha_{\sigma(3)},
          \beta_{\sigma(1)},\beta_{\sigma(2)},\beta_{\sigma(3)})$$
This defines an action of $S_3$ over $X$. Further, any two tuples
in $X$ which lie in the same orbit under the action of $S_3$
define same factorization ( upto order of factors ) of $216$, and
conversely. Hence to calculate the number of distinct
factorizations (ignoring the order of
 factors) of $216$ into three natural numbers, we have to calculate the
 number of distinct orbits (say) $k$ of $X$ under the $S_3$-action defined.
 By the Burnside's lemma,
\bea
     k & = & \frac{1}{\abs{S_3}} \sum_{\sigma \in S_3} \chi(\sigma)
\eea We know that the order of $S_3$ is $6$. Hence to obtain $k$
we have to calculate $\chi(\sigma)$ for all $\sigma \in S_3$.
First of all, note that $\abs{X}$ is $100$. For the identity
element $e\in S_3$, $\chi(e) = 100$. For $ \sigma_1 =
\left(\begin{array}{ccc} 1 & 2 & 3\\ 2 & 1 & 3
                 \end{array} \right)$ in $S_3$,
\ba{rl}
     & \sigma_1 (\alpha_1,\alpha_2,\alpha_3,\beta_1,\beta_2,\beta_3)
           = (\alpha_1,\alpha_2,\alpha_3,\beta_1,\beta_2,\beta_3)\\
\Leftrightarrow & \alpha_1 = \alpha_2,  \beta_1 = \beta_2 \ea One
can easily check that such possible tuples are $$ (0,0,3,0,0,3),
(0,0,3,1,1,1), (1,1,1,0,0,3), (1,1,1,1,1,1).$$ Thus
$\chi(\sigma_1) = 4$. Similarly we can check that for \bea
\sigma_2 & = & \left(\begin{array}{ccc} 1 & 2 & 3\\ 3 & 2 & 1
                 \end{array} \right)\\
\sigma_3 & = & \left(\begin{array}{ccc} 1 & 2 & 3\\ 1 & 3 & 2
                 \end{array} \right) \\
\sigma_4 & = & \left(\begin{array}{ccc} 1 & 2 & 3\\ 3 & 1 & 2
                 \end{array} \right)\\
\sigma_5 & = & \left(\begin{array}{ccc} 1 & 2 & 3\\ 2 & 3 & 1
                 \end{array} \right)
\eea We have $$ \chi(\sigma_2) = \chi(\sigma_3) = 4,
\chi(\sigma_4) = \chi(\sigma_5) = 1$$ Hence \bea
               k & = & \frac{1}{6} (100 + 3\times 4 + 2\times 1)\\
                 & = & \frac{114}{6} = 19
\eea \ex \bt\label{t56} Let $G$ be a finite group and $H$, a
subgroup of $G$, with $[G:H] = n > 1$. If $n!$ is not divisible by
$\circ(G)$ then $H$ contains a proper normal subgroup of $G$. \et \pf Let
$X$ be the set of all distinct left cosets of $H$ in $G$. Then $\abs{X} =
n$. For any $x\in G$ and $gH\in X$, define $$ \hat{x}(gH) = xgH$$  Then
\bea
          \theta : G & \to & A(X)\\
                  x & \mapsto & \hat{x}
\eea $\theta$ is a group homomorphism. Thus we get an action of $G$ over $X$.
 As $\abs{X} = n,$ the
order of the group $A(X)$ is $n!$. By our assumption $\circ(G)$ does not
divide $\abs{A(X)} = n!$. Hence $\theta$ is not a monomorphism. Let $K$ be
the Kernel of $\theta$. As $\theta$ is not a monomorphism $ K \neq id$.
Therefore $K$ is a nontrivial normal subgroup of $G$. As $K = \ker
\theta$, we get \ba{rl}
                  &  Kx = x \mbox { for all } x\in X\\
\Rightarrow       & KH = H\\ \Rightarrow       & K \subset H \ea Hence the
result follows.     $\blacksquare$
\bco\label{co7.32co} Let $G$ be a finite group and $H$, a subgroup of $G$.
If $[G:H]=p,$ the smallest prime dividing order of $G$, then $H\lhd
G.$\eco \pf If $\circ(G)=p $, then $H=\{e\}$, and the result is clear.
Further, if $[G:H]=2$, then $H$ is normal. Hence we can assume that
$\circ(G)>p >2.$  \\ Now, as $p$ is the smallest prime dividing
$\circ(G)$, $\circ(G)$ does not divide $p!$. Thus by the proof of the
theorem $H$ contains a proper normal subgroup $K$ of $G$ such that
$\circ(G/K) = [G:K]$ divides $p!$. We have $$\begin{array}{ccl}
   & [G:K] & =[G:H][H:K] \\
   &  & =p[H:K] \\
  \Rightarrow &  & [H:K] \mbox{ divides }p!
\end{array}$$ \;\;\;\;\; $\Rightarrow\;\; [H:K]$ has a prime factor less
than $p$ and divides $\circ(G)$. This contradicts our assumption on $p$,
unless $[H:K]=1$. Thus $[H:K]=1$ i.e., $H=K$. Hence $H\lhd G$.
$\blacksquare$
   \\ \\  \\
      {\bf EXERCISES}
\begin{enumerate}
  \item
  Let a group $G$ acts transitively on a set $X$. Prove that the cosets of
  the subgroup $G_{x} = {\bf{Stab}}(x)$ for some $x \in X$, are exactly
  the sets $S(y) = \{g \in G \st g(x) = y\}$ where $y \in X$.
  \item
  Let $H$ be a finite subgroup of a group $G$. Prove that for any $(h,h')
  \in H \times H$ and $x \in G$, $$(h,h')(x) = hxh'^{-1}$$ defines an
  action of $H \times H$ on $G$. Further, show that $H \lhd G$ if and only
  if every orbit in $G$ with respect to this action of $H \times H$ has
  exactly $\abs{H}$ elements.
  \item
  Let a group $G$ acts on a set S. Prove that for an element $x \in S$,
  ${\bf{Stab}}(x) = {\bf{Stab}} (gx)$ for all elements $g \in G$ if and only
  if ${\bf{Stab}}(x)$ is a normal subgroup of $G$.
  \item
  Let $G$ be a group of order $p^n(n\geq 1)$, where $p$ is a prime. Prove
  that if $G$ acts on a finite set $X$ and $m$ is the number of the fixed
  points in $X$, then $\abs{x}\equiv m(mod\, p)$.
  \item
  Let $G$ be a group with $\circ(G)=55$. Prove that if $G$ acts on a set
  $X$ with $\abs{X}=19$, then there are at-least 3 fixed points.
  \item
  Prove that under the action of the group $Gl_n(\R)$ on $S=\R^n$ (Example
  \ref{ex54}), there are exactly two orbits.
  \item
  Consider the group $G=S_7$, the symmetric group on 7-symbols. Find the
  stabliser of the element $$\sigma = \left(\begin{array}{ccccccc}
    1 & 2 & 3 & 4 & 5 & 6 & 7 \\
    3 & 4 & 2 & 1 & 5 & 6 & 7 \
  \end{array}\right)$$ in $G$ under the action of $G$ over $G$ as in the
  example \ref{ex52}.
  \item
  Prove that in the example \ref{ex7.10}, orbit of any vector $x \in
  S^{n-1}$ is whole of $S^{n-1}.$
  \item
  Consider the subgroup $G=\left\{\left(\begin{array}{cc}
    s & 0 \\
    t & t \
  \end{array}\right) \st s,t \in \R s>0, t>0 \right\}$ of $Gl_2(\R)$. The
  action of $Gl_n(\R)$ on $\R^2$ (Example \ref{ex54} ) gives an action of
  $G$ on $\R^2$. Prove that $G$ acts on the set of all straight lines in
  $\R^2$. Show that $X$-axis and $Y$-axis are the only fixed points under
  the action of $G$. Further, prove that there are twelve orbits of this
  action.
  \item
  Prove that a group $G$ is residually finite if and only if $G$ is
  isomorphic to a subgroup of a direct product of finite groups.
  \end{enumerate}

\chapter{Conjugacy Classes}
 \footnote{\it contents group8.tex }
 We shall consider here a special action of a group over itself. This
 gives a very useful concept in group theory. \\
 For an element $g$ of a group
$G,$ define \bea
             \tau_g : G & \to & G\\
                 x & \mapsto & gxg^{-1}
\eea We can  easily check that $\tau_g$ is a transformation of
$G$, i.e.,it is a one-one map from $G$ onto $G.$ Further, for
$x,g,h \in G,$ we have \bea
                     \tau_{gh}(x) & = & ghx(gh)^{-1}\\
                                 & = & g(hxh^{-1})g^{-1}\\
                                 & = & \tau_g \tau_h(x).
\eea Hence $\tau_{gh} = \tau_g \tau_h$. Thus $G$ can be considered
a transformation group of $G$ where action of $g$ over $G$ is given by
$\tau_g$. As seen in the Chapter 7, this action defines an equivalence
relation over $G.$ This relation has special importance in group. We
define
:
\bd\label{d61} Let $G$ be a group and $x,y \in G.$ We shall say that $y$
is a conjugate of $x$ in $G$ if there exists an element $g\in G$ such that
$gxg^{-1} = y.$ \ed
\brs (i) As seen in Chapter 7, conjugacy is an equivalence
relation
 over $G$ ( Theorem \ref {t52}).  \\
(ii) For every $g \in G$, $\tau_{g}$ is an automorphism of $G$
[Theorem \ref{t36}].   \\  (iii) The map $ g \mapsto \tau_g$ is a
homomorphism from $G$ to $AutG$ since $\tau_{gh} = \tau_g \tau_h$
for any $g,h \in G.$
\ers \bd\label{d62} Let $G$ be a group and $x\in G.$ Then the set
$$ [x] = \{\tau_g(x) = gxg^{-1} \st g\in G\}$$ is called the conjugacy
class of $x$ in $G.$ \ed \brs (i) As noted above, conjugacy is an
equivalence relation over $G.$  Hence any two conjugacy classes are either
identical or disjoint i.e.,for $x,y\in G,$ either $[x] = [y]$ or $[x] \cap
[y] = \emptyset.$  \\  (ii) For the identity element $e$ of $G, [e] =
\{e\}.$
\ers
\bd\label{d63}
 Let $G$ be a group. Then the set
$$ Z(G) = \{x\in G \st [x] = \{x\}\}$$ is called the {\bf centre} of $G.$
\ed {\bf Note}:(i) In case it is not necessary to emphasize the group $G,$
we write $Z$ for $Z(G)$, the centre of $G.$\\ (ii) Note that \ba{rl}
              &  [x] = \{x\}\\
\Leftrightarrow & g^{-1}xg = x \mbox { for all } g\in G\\
\Leftrightarrow & xg = gx \mbox { for all } g\in G \ea Hence $$
Z(G) = \{x\in G \st xg = gx \mbox { for all } g\in G \}.$$
\bd\label{d64} For an element $x$ in a group $G,$ $$ C_G(x) =
\{g\in G \st gx =xg \mbox { i.e., } \tau_g(x) = x \}$$ is called the
centralizer of $x$ in $G.$ \ed \bt\label{t65} Let $G$ be a group. Then\\
(i) $Z(G)$ is a normal subgroup of $G.$\\ (ii) For any $x\in G, C_G(x)$ is
a subgroup of $G.$ \et \pf (i) Clearly $[e] = \{e\}$, hence $e\in Z(G).$
Take $x,y\in Z(G).$ Then for
 any $g\in G,$
\bea
     gxy^{-1}g^{-1} & = & gx g^{-1}g y^{-1}g^{-1}\\
                    & = & x(gyg^{-1})^{-1} \mbox { (since } [x] = \{x\})\\
                    & = & xy^{-1}      \mbox { (since } [y] = \{y\})
\eea Hence $[xy^{-1}] = \{xy^{-1}\}.$ Thus for $x,y\in Z(G),
xy^{-1}\in Z(G).$ Consequently $Z(G)$ is a subgroup of $G.$ Further for
any $x\in Z(G), g\in G,$ we have, \ba{rl}
              & gxg^{-1} = x \in Z(G)\\
\Rightarrow   & gZ(G)g^{-1}  \subseteq Z(G)\\ \Rightarrow   & Z(G)
\mbox{ is a normal subgroup of }G. \ea (ii) The proof is similar
to (i).  $\blacksquare$ \br For a group G, $$ Z(G) = \bigcap_{x\in
G} C_G(x).$$
\er
\bd\label{d66}
 Let $S$ be a subset of a group $G.$ The set
$$C_G(S) = \{g\in G \st gx = xg
                      \mbox{ i.e., }gxg^{-1} = x \mbox{ for all }x\in S \}$$
is called the centralizer of $S$ in $G.$ \ed \bd\label{d67} For a subset
$S$ of a group $G$, the set $$  N_G(S) = \{ g\in G \st gS = Sg \}$$ is
called the normalizer of $S$ in $G.$ \ed
\brs (i) If $S = \{x\}$ is a singleton set, then $$     C_G(S) =
N_G(S).$$ (ii) By definition of $C_G(S),$ we have $$ C_G(S) =
\bigcap_{x\in S} C_G(x).$$ Thus by the theorem \ref{t65} (ii), it follows
that $C_G(S)$ being an intersection of subgroups is a subgroup of $G.$\\
(iii) One can easily check that $N_G(S)$ is a subgroup of $G.$ \\  (iv) If
$H$ is a subgroup of $G$, then $H \subset N_G(H)$. (Use : $xH=H=Hx$ for
all $x \in H.$)\ers
\bexe Give an example of a proper subgroup $H$ of a group $G$ such that $H \subset
N_G(H)$.\eexe
\bd\label{d68}
 Let $A,B$ be two subsets of a group $G$. We shall say that $B$ is a conjugate of
 $A$ if there exists an element $g\in G$, such that $g^{-1}Ag = B$ i.e.,
 $\tau_g(A) = B.$
\ed \brs (i) If $A,B$ are conjugate subsets of $G$, then as
$\tau_g$ is a transformation of $G$ (in fact an automorphism of
$G$), $$\abs{A} = \abs{B}$$ (ii) Let $H$ be a subgroup of $G$.
Then any conjugate of $H$ is $g^{-1}Hg = \tau_g(H)$ for some $g\in
G.$ Thus being the image of H under $\tau_g$, it is a subgroup.
(Use : Image of a subgroup under a homomorphism
 is a subgroup.)
\ers \bt\label{t69}
 Let $G$ be a group and $x\in G.$ Then the cardinality of the conjugate class
of $x$ i.e.,$\abs{[x]}$ is equal to $[G:C_G(x)].$ \et \pf Let $g,
h\in G.$ Then \ba{rl}
                 & gxg^{-1} = hxh^{-1}\\
\Leftrightarrow  & h^{-1}gx = xh^{-1}g\\ \Leftrightarrow  &
h^{-1}g \in C_G(x)\\ \Leftrightarrow  & gC_G(x) = hC_G(x) \ea
Hence, $\abs{[x]} = [G:C_G(x)].$  $\blacksquare$  \\ More
generally, we have
\bt\label{t610} For a subset $S$ of a group $G$, the cardinality
of the set of conjugates of $S$ in $G$ is equal to $[G:C_G(S)].$
\et \pf The proof is similar to theorem \ref{t69}.  $\blacksquare$
\bco\label{c610} For any subgroup $H$ of a group $G$, the number
of  subgroups conjugate
 to $H$ is $[G:N_G(H)].$
\eco \pf (Use : Conjugate of a subgroup is a subgroup and the
theorem \ref{t610}.)  $\blacksquare$
\bl\label{l611}
Let $H$ be a subgroup of a group $G$. Then for any $x\in G$ $$
N_G(x^{-1}Hx) = x^{-1}N_G(H)x.$$
\el
\pf For any $g\in N_G(H)$, we have \bea
     x^{-1}gx \cdot x^{-1}Hx & = & x^{-1}gHx\\
                             & = & x^{-1}Hgx\\
                             & = & x^{-1}Hx \cdot x^{-1}gx
\eea Therefore \be\label{e61}
 x^{-1}N_G(H)x \subseteq N_G(x^{-1}Hx)
\ee Taking $x^{-1}Hx$ for $H$, we get from \ref{e61} $$
xN_G(x^{-1}Hx)x^{-1}  \subseteq  N_G(x(x^{-1}Hx)x^{-1}) =
N_G(H)$$ Therefore \be\label{e62}
  N_G(x^{-1}Hx) \subseteq  x^{-1}N_G(H)x
\ee Now, from \ref{e61} and \ref{e62}, we conclude $$
x^{-1}N_G(H)x = N_G(x^{-1}Hx) \hspace{1.5in} \blacksquare$$ \br The result
is true for any non-empty subset $A$ of $G.$ \er \bl\label{l612} Let $H$
be a subgroup of a group $G$, and $x\in G.$ Then $$[G:H] = [G:x^{-1}Hx].$$
\el \pf Let $\{x_iH\}_{i\in I}$ be the set of all distinct left coset
 representatives of $H$ in $G$. Then
 $$ G = \cup_{i\in I} x_iH$$
is a disjoint union. For any $x\in G$, \ba{rl}
          &  G = x^{-1}Gx  = \cup_{i\in I}x^{-1}x_iHx\\
                      & = \cup_{i\in I}x^{-1}x_ix(x^{-1}Hx)\\
\Rightarrow   & [G:H] \geq [G:x^{-1}Hx] \mbox{ for any }x\in G. \ea Hence,
\bea
           [G:H] & \geq & [G:x^{-1}Hx]\\
                 & \geq &[G:x(x^{-1}Hx)x^{-1}]\\
                 & = & [G:H].
\eea Consequently, $[G:H] = [G:x^{-1}Hx].$          $\blacksquare$ \bt\label{t613} Let $H$ be a
subgroup of finite index in a group $G.$ Then there exists a normal
subgroup $K$ of $G$, $K \subseteq H$, such that $[G:K] < \infty.$ \et \pf
By the definition of $N_G(H), H \subseteq N_G(H).$ Hence as $[G:H] <
\infty$, we get $[G:N_G(H)] < \infty.$ Thus, by corollary \ref{c610}, the
subgroup $H$ has finite number of conjugates in $G.$ Let $ x_1^{-1}Hx_1$,
$x_2^{-1}Hx_2$, $\ldots x_t^{-1}Hx_t$ be all distinct conjugates of $H$ in
$G$. Then $$ \bigcap_{i=1}^t (x_i^{-1}Hx_i) = \bigcap_{x\in G} (x^{-1}Hx)
= K \mbox{ say }$$ Hence, for any $y\in G$, \ba{rl}
               & y^{-1}Ky = \bigcap_{x\in G}y^{-1}(x^{-1}Hx)y  \\
\Rightarrow    & y^{-1}Ky = \bigcap_{x \in G}(xy)^{-1}H(xy) =
\bigcap_{x\in G} (x^{-1}Hx) = K\\ \Rightarrow    & K \mbox{ is a normal
subgroup of }G. \ea By the lemma \ref{l612}, we get $$ [G:H] =
[G:x_i^{-1}Hx_i] \mbox { for all } i = 1,2,\ldots ,t.$$ Hence, by
corollary \ref{co212},
\bea
          [G:K] & = & [G:\cap_{i=1}^t (x_i^{-1}Hx_i)]\\
                & \leq & \prod_{i=1}^t [G:(x_i^{-1}Hx_i)] = [G:H]^t < \infty.
\eea Thus the proof is complete.    $\blacksquare$
  \bt\label{t8.22tt} Let $H$ be a normal subgroup of a group $G$ such that
  $G/H$ is infinite cyclic. If $N$ is a torsion free, normal subgroup of
  $G$ contained in $H$ with $[H:N]<\infty,$ then there exists a torsion
  free, normal subgroup $M$ of $G$ with $[G:M]<\infty.$\et \pf Let
  $G/H=gp\{aH\},$ and $A=gp\{a\}.$ Then $AH=G.$ Since $gp\{aH\}$ is
  infinite cyclic, $A\cap H=id.$ As $N\lhd G,$ $B=AN$ is a subgroup of
  $G$. Let for $g=a^m$ and $x\in N,\circ(a^mx)=t<\infty.$ Then \[\begin{array}{crc}
     & (a^mxH)^t & =H \\
    \Rightarrow & (a^mH)^t & =H \\
    \Rightarrow & (aH)^{mt} & =H \\
    \Rightarrow & m & =0 \
  \end{array}\] Therefore $\circ(x)=t<\infty.$ Now, as $N$ is torsion free
, $x=e.$ Hence if $\circ(a^mx)<\infty, a^mx=e.$ Consequently $B$ is
torsion free. If \[H=\bigcup_{i=1}^{m}h_iN\] is a coset decomposition of
$H$ with respect to $N$, then \[G=AH=HA=\bigcup_{i=1}^{m}h_iN\] Further,
\[\begin{array}{crlc}
   & h_iAN & =h_jAN & \\
  \Leftrightarrow & h_{j}^{-1}h_i & \in AN & \\
  \Leftrightarrow & h_{j}^{-1}h_i & \in N & \mbox{ since }H\cap A=e, \mbox{ and }N\subset H \\
  \Leftrightarrow & h_iN & =h_jN &
\end{array}\] Hence $[G:B]=[H:N]<\infty.$ Now, the result follows from
theorem \ref{t613}.          $\blacksquare$ \br From the proof of the
theorem it is clear that if $G/H$ is cyclic and $N$ is not necessary
torsion then thee exists $M\lhd G,$ and $[G:M]<\infty.$\er
   \br\label{r8.008r} For any group $G$, if
  $$\mathcal{A} = \{H \subset G :\mbox{ subgroup }\st [G:H] < \infty\}$$
  and  $$\mathcal{N} = \{N\lhd G \st [G:N] < \infty \},$$
  then $$\bigcap_{H \in \mathcal{A}}H = \bigcap_{N \in \mathcal{N}}N $$
  Hence $G$ is residually finite if and only if $$\bigcap_{N \in \mathcal{N}}N = \{e\}$$ \er
  \bexe Show that a group $G$ is residually finite if and only if
  $G$ is isomorphic to a subgroup of a direct product of finite
  groups.\eexe
   \bd Let $G$ be a finite group and
let $[x_1],[x_2], \ldots ,[x_k]$ be all the distinct conjugate
classes of $G$ such that $\abs{[x_i]} > 1$ for $i = 1,2,\ldots
,k$. Then
\begin{eqnarray}
   \abs{G} & = & \abs{Z(G)} + \sum_{i=1}^k \abs{[x_i]}\nonumber\\
            & = & \abs{Z(G)} + \sum_{i=1}^k [G:N_G(x_i)]\label{e64}\
\end{eqnarray}
The equation \ref{e64} is called the class equation of the
group $G$. \ed  We shall, now, give some applications of the class
equation.
\bt\label{t614} Let $p$ be a prime and $G$ a finite group with $p^n (n>0)$
elements. Then $\abs{Z(G)} = p^{\alpha}, \alpha > 0.$ \et \pf By theorem
\ref{t69}, for any $x\in G, \abs{[x]} = [G:C_G(x)].$ Hence, if $\abs{[x]}
> 1$, then $C_G(x) \subsetneqq G$. By theorem \ref{co22} for any finite group,
order of a subgroup divides the order of the group. Hence, if
$\abs{[x]}>1$, then $\abs{C_G(x)} = p^m, 1\leq m < n$. This implies
$\abs{(x_i)} = [G:C_G(x_i)] = p^{n-m_i}$ is a multiple of $p$ for all
$i\geq 1.$ Now, consider the class equation for $G$ \ba{rl}
      &  \abs{G} = \abs{Z(G)} + \sum_{i=1}^k \abs{[x_i]}\\
 \ea where $[x_1],[x_2],\cdots,[x_k]$ are all the distinct conjugate
 classes of $G$ with $\abs{[x_i]}>1.$ As noted above we shall have
 $\abs{C_G(x_i)}=p^{m_{i}}, 1\leq m_i <n,$ and hence $\abs{[x_i]} =
 [G:C_G(x_i)]=p^{n-m_i}$ for all $i=1,2,\cdots,k.$ Therefore the class
 equation gives $$\begin{array}{clc}
    & p^n = \abs{Z(G)}+\sum_{i=1}^{k}p^{n-m_i} & \mbox{where }n-m_i >0 \mbox{ for all }
    i=1,2,\cdots,k  \\
   \Rightarrow & p\mid \abs{Z(G)} &  \\
   \Rightarrow & \abs{Z(G)} > 1 &  \
 \end{array}$$ As $Z(G)$ is a subgroup of $G$, by theorem \ref{co22},
 $\circ(Z(G))$ divides $\circ(G)=p^n.$ Hence $\circ(Z(G))= p^{\alpha}, \alpha
 >0.$      $\blacksquare$
\bco\label{c615} Let $p$ be a prime. Then any group of order $p^2$
is abelian. \eco \pf Let $G$ be a group of order $p^2$, by the theorem,
$\abs{Z(G)} = p$ or $p^2$. If $\abs{Z(G)} = p^2$, then $G = Z(G)$. Hence
$G$ is abelian. Next, let $\abs{Z(G)} = p$. As $Z(G)$ is normal in $G,
G/Z(G)$ is a group of order $p$. A group of prime order is cyclic. Hence
$G/Z(G) = gp\{aZ(G)\}$ for some $a\in G$. Thus every element in $G$ is of
the form $a^ih (0\leq i < p), h\in Z(G)$. Let $x = a^ih_1, y = a^jh_2$ be
any two elements of $G$, where $h_1,h_2 \in Z(G).$ Then
\bea
                xy & = &  (a^ih_1)(a^jh_2)\\
                   & = & a^{i+j}h_1h_2\\
                   & = & a^{i+j}h_2h_1\\
                   & = & a^jh_2a^ih_1\\
                   & = & yx
\eea Thus $G$ is abelian.  $\blacksquare$
\bco\label{co8.24co} Let $G$ be a finite group of order $p^n$ where $p$ is
a prime. If $H$ is a proper subgroup of $G$, then $N_G(H) \supsetneqq H.$
Thus if $G$ contains a proper subgroup, then $G$ has a proper normal
subgroup.\eco \pf We shall prove the result by the induction on $n$.
Clearly $n>1$, since a group of prime order has no proper subgroup. If
$n=2$, then by corollary \ref{c615}, $G$ is abelian. Hence $N_G(H)=G,$ and
the result is clear. Therefore assume $\circ(G)=p^n,n\geq 3.$ By the
theorem $\circ(Z(G)) = p^{\alpha}, \alpha>0.$ If $Z=Z(G) \not\subset H,$
then choose any $z\in Z,z\not\in H.$ As $zH=Hz, z \in N_G(H)$ and hence
$H\subsetneqq N_G(H).$ Now, assume $Z \subset H.$ Then as $Z\lhd G,
\overline{H}=H/Z$ is a subgroup of $\overline{G}=G/Z,$ where
$\circ(\overline{G})=p^{n-\alpha} < p^n.$ Hence by induction
$N_{\overline{G}}(\overline{H})\supsetneqq \overline{H}.$ Consider the
natural epimorphism  $$\begin{array}{rcl}
  \eta : G & \rightarrow & G/Z \\
  x & \mapsto & xZ
\end{array}$$  As $Z\subset H, \eta^{-1}(\overline{H}) = H.$ Since
$N_{\overline{G}}(\overline{H})\supsetneqq \overline{H},$ there exist $z
\in G$ such that $\overline{z}=\eta(z) \in N_{\overline{G}}(\overline{H}),
\overline{z} \not\in \overline{H}.$ Then $z \not\in
\eta^{-1}(\overline{H})=H,$ and $$\begin{array}{crl}
   & \eta(z^{-1}Hz) & =\overline{z^{-1}}\;\overline{H}\;\overline{z} \\
   &  & \subset \overline{H} \\
  \Rightarrow & z^{-1}Hz \subset \eta^{-1}(\overline{H}) & = H
\end{array}$$ Thus $z \in N_G(H), z \not\in H.$ Hence $N_G(H)\supsetneqq
H.$ The final part of the statement follows by observing that as $G$ is
finite, the chain $H \subsetneqq N_G(N_G(H)) \subsetneqq \cdots$ has to
terminate at $G$.     $\blacksquare$
\bexe Let $\phi : A \rightarrow B$ be an onto, homomorphism of groups. Prove that for any
subgroup $K$ of $B$, $\phi^{-1}(N_B(K))=N_A(\phi^{-1}(K)).$ \eexe
\bl\label{l616}
If for a group $G$, $G/Z(G)$ is a cyclic group, then $G$ is
abelian.
\el
\pf Let $G/Z(G) = gp\{aZ(G)\}$ for some $a\in G$. Then $$ G = \cup_{n\in
\Z} a^nZ(G)$$ Let $g,h\in G$. Then there exist $m,n\in \Z$ and $x,y\in
Z(G)$ such that $g = a^nx$ and $h = a^my$. Now, as in the proof of
corollary \ref{c615}, $gh = hg$. Hence $G$ is abelian. $\blacksquare$
\bco\label{c617} Let $G$ be a non-abelian group of order $p^3$ ($p$ a
prime). Then $\abs{Z(G)} = p.$ \eco \pf By theorem \ref{t614},
$\circ(Z(G)) = p^{\alpha}, \alpha > 0.$ Hence $\circ(Z(G)) = p$ or $p^2$
or $p^3$. We, however, have that $G$ is non-abelian. Hence $\circ(Z(G))
\neq p^3$. Further, if $\circ(Z(G)) = p^2,$ then $G/Z(G)$ is cyclic of
order $p$. Then by lemma \ref{l616}, $G$ will be abelian. This contradicts
our assumption. Consequently $\circ(Z(G)) = p.$ $\blacksquare$
\bt\label{CAUCHY} (CAUCHY) Let $p$ be a prime dividing the order of a
finite group $G$. Then $G$ has an element of order $p.$ \et \pf
Let$\circ(G) = n.$ We shall prove the result by induction on $n$. As $ p
\mid n, n \geq p.$ If $n = p$, then $G$ is cyclic of order $p$ and hence
every element $(\neq e)$ of $G$ has order $p$. Now, assume $n > p$. We
shall consider two cases:\\ Case I: $G$ is abelian.\\ Assume that for any
proper subgroup $H$ of $G$, $p$ does not divide the order of $H$, since
otherwise by induction assumption $H$ contains an element of order $p$,
which gives required element in $G$. As $p \mid \circ(G)$ and $\circ(G)
\neq p$, $n$ is a composite number. Hence $G$ has a proper subgroup $H$. (
Use: Any group of composite order has a proper subgroup.) Now, note that
$\circ(H)\cdot \circ(G/H) = \circ(G) $, and $p$ does not divide $\circ(H)$
by our assumption. Hence $p$ divides $\circ(G/H).$ Hence by induction
assumption $G/H$ has an element (say) $aH$ of order $p$, i.e., $aH \neq
H$, and $(aH)^p = H$. Thus $a^p\in H$. Let $a^p = b \in H$, and $\circ(b)
= m.$ Then $ {(a^m)}^p = a^{pm} = e.$ If we show that $a^m \neq e$, then
$a^m$ shall be an element of order $p$ in $G$. Let us assume the contrary,
i.e., $a^m = e$. Then $(aH)^m = H$. Hence $\circ(aH) = p$ divides $m$. As
$\circ(b) = m$, $m \mid \circ(H).$ This implies $p\mid \circ(H)$, which is
not true. Hence $a^m \neq e$, and the proof follows in this case.\\ Case
II: $G$ is any group, i.e., not necessarily abelian.\\ As in case I, we
can assume that $p$ does not divide $\circ(H)$ for any proper subgroup $H$
of $G.$ Now, consider the class equation $$ \abs{G} = \abs{Z(G)} +
\sum_{i=1}^m [G:N_G(a_i)]$$ where $[a_1],[a_2], \ldots ,[a_m]$ denote all
distinct conjugate classes of $G$ with $\abs{[a_i]} > 1$. As $\abs{[a_i]}
> 1$, $N_G(a_i)$ is a proper subgroup of $G$. Hence $p$ does not divide
$\circ(N_G(a_i))$.  Thus $p$ divides $[G:N_G(a_i)]$ for all $ i = 1,2,
\ldots , m.$ Now, as $p\mid \circ(G)$, we conclude from the class equation
that $p\mid \circ(Z(G))$. As $Z(G)$ is an abelian subgroup of $G$, the
result follows from the case I.  $\blacksquare$  \\    \\   {\bf SECOND
PROOF OF CAUCHY'S THEOREM.}\\ We shall infact prove the following more
general statement. \\ If a prime divides the order of a finite group $G$,
then $G$ has $kp(k\geq 1)$ solutions to the equation $ x^p = e$.\\ \pf Let
$\circ(G) = n$. Put $$X = \left\{ \underline{a} = (a_1,a_2, \ldots
,a_p)\st
      a_1.a_2. \ldots .a_p = e, a_i\in G \mbox { for }1\leq i \leq p \right\}$$
The set $X$ has $n^{p-1}$ elements. (Use : If first $p-1$ elements are
arbitrarily chosen, then the $p^{th}$-element is the inverse of the
product $a_1.a_2. \ldots .a_{p-1}$.) For any $\tau \in S_p$ and
$\underline{a} = (a_1,a_2, \ldots ,a_p) \in X$, put $$\tau(\underline{a})
= (a_{\tau(1)},a_{\tau(2)}, \ldots ,a_{\tau(p)})$$ This defines an action
of $S_p$ over $X$. Let  \\ $ \sigma = \left(\begin{array}{cccccc}
                       1 & 2 & 3 & \ldots & p-1 & p\\
                       2 & 3 & 4 & \ldots & p   & 1
                 \end{array} \right) \in S_p$, and $A = gp\{\sigma\}$.
  \\  Restricting the action of $S_p$ to $A$ we get an action of $A$
over $X$. We
 shall, now, concentrate on the $A$-orbits of $X$. Let us observe that
 $\circ(\sigma) = p$. Hence $\circ(A) = p$. Clearly any subgroup of $A$ has
 order $p$ or $1$. (Use: Order of a subgroup divides the order of the group.)
 Hence, by lemma \ref {l510}, any $\underline{a} \in X$,
 $\abs{\orb(a)} = p$ or $1$. Note that for the identity element $e$
 of $G$, the element $\underline{e} = (e,e \ldots ,e) \in X$ has only one
 element in it's $A$-orbit. Let $r$ be the number of $A$-orbits of $X$ with
 one element and $l$ be the number of $A$-orbits with $p$ elements. Then
 \ba{rl}
                 & \abs{X} = r.1 + p.l\\
\Rightarrow      & n^{p-1} = r + pl\\ \Rightarrow      & p\mid r \mbox{
since } p\mid n\\ \Rightarrow      & r = kp\;\;\;\;\;\mbox{ for some
}k\geq 1 . \ea Clearly $r \geq 2$ If $\underline{a} = (a_1,a_2, \ldots
,a_p) \neq (e,e, \ldots ,e)$ is an element of $X$ whose $A$-orbit has only
one element, then \ba{rl}
              & \sigma (\underline{a}) = \underline{a}\\
\Rightarrow   & (a_2,a_3, \ldots ,a_p,a_1) = (a_1,a_2, \ldots ,a_p)\\
\Rightarrow   &  a_1 = a_2 = \ldots = a_p = y \mbox { say }\\ \ea i.e.,
$y$ is an element of order $p$ in $G$. As $r=kp\geq 2,$ the equation $x^p
= e$ has $kp$ solutions in $G$. Each of the non-identity solutions has
ofcourse order $p$.  \bco\label{co8.27co} Let $G$ be a finite group and
$p$, a prime. Then $\circ(G) = p^{\alpha}$ if and only if every element of
$G$ has order power of $p$. \eco  \pf It is clear from the theorem.  \\
We, now, define
:
\bd\label{d8.30d} Let $p$ be a prime. A group $G$ is called a $p$-group if
every $x\in G$ has order $p^{\alpha}, \alpha$ depending on $x$.\ed \br A
finite group $G$ is a $p$-group if and only if $\circ(G)=p^n.$\er \bexe
Give an example of an infinite $p$-group. Further show that any two finite
$p$-groups need be isomorphic.\eexe   \clearpage
{\bf EXERCISES}
\begin{enumerate}
  \item
  Write all the conjugates of the element $A=\left(\begin{array}{cc}
    2 & 0 \\
    1 & 2 \
  \end{array}\right)$ in $Gl_2(\Z_3)$ and determine its centralizer.
  \item
  Write all conjugacy classes in $S_4$ and verify the class equation.
  \item
  Write the class equation for the following groups :  \\ (i)
  $Sl_2(\Z_5)$\; (ii)The group of quaternions , and (iii) $D_5$.
  \item
  Prove that if a group $G$ has a conjugacy class with exactly two
  elements, then $G$ is not simple.
  \item
  Let $H$ be a proper subgroup of a finite group $G$. Then $\cup_{x\in
  G}x^{-1}Hx \neq G.$
  \item
  Let $p$ be a prime. prove that a group of order $p^2$ is either a cyclic
  group or is a direct product of two cyclic groups of order $p$.
  \item
  Prove that a subgroup $H$ of a group $G$ is normal if and if $H$ is a
  union of conjugacy classes of $G$.
  \item
  Let $m\geq 1$ and $n=2m.$ Consider the Dihedral group (Example
  \ref{th11}) : $$D_{n} = \{Id, \sigma, \sigma^{2},
\cdots,\sigma^{n-1},
 \tau, \sigma\tau, \cdots
, \sigma^{n-1}\tau,\,\sigma^{n} = Id = \tau^{2},\;\sigma\tau =
\tau\sigma^{n-1}\}.$$  Prove that $\circ(C_G(\sigma^i))=n$ for all $1\leq
i\leq m-1,$ $C_G(\sigma^m)=D_n$ and $\circ(C_G(\tau))=4$. Further show
that :  \\ Conjugacy class of $\sigma^i=\{\sigma^i , \sigma^{-i}\}$ for
all $1\leq i\leq m-1$ \\  and  \\ Conjugacy class of
$\tau=\{\tau,\tau\sigma^2,\tau\sigma^4,\cdots,\tau\sigma^{2(m-1)}\}.$
\item
In the above exercise if $n=2m+1,$ then show that :  \\  Conjugacy class
of $\sigma^i=\{\sigma^i,\sigma^{-i}\}$ for all $1\leq i\leq m$ and  \\
Conjugacy class of $\tau=\{\tau\sigma^i \st 0\leq i\leq 2m\}.$
\item
Let $G$ be a finite group with 3 conjugacy classes. Prove that
$\abs{Z(G)}=1$ or 3.
\item
Let $p$ be a prime. Prove that $G=\Q_p/\Z$ is a $p$-group, and $p^nG=G$
for all $n\geq 0.$
\item
Let $x$ be a non-generator of a group $G$. Prove that for any $g\in G,
gxg^{-1}$ is also a non-generator of $G$. Then deduce that the Frattin
subgroup of $G$ is a normal subgroup.
\item
Let for a group $G$, $\circ(G')=m.$ Prove that any element of $G$ has
atmost $m$ conjugates.
\item
Let $G=gp\{x_1,x_2,\cdots,x_k\}$ and let $\circ(G')=m.$ Prove that
$[G:Z(G)]\leq m^k.$
\item
Let $\phi:G\rightarrow H$ be an epimorphism of groups. Prove that if $A,B$
are two conjugate subgroups of $H$, then $\phi^{-1}(A)$ and $\phi^{-1}(B)$
are conjugate subgroups of $G$.
\item
Let $A,B$ be two proper subgroups of a group $G$ such that $G=AB.$ Prove
that for any $x,y\in G$, $G=(x^{-1}Ax)(y^{-1}By).$
\item
Let $A$ be a proper subgroup of a group $G$. Prove that if $G=AB$ for a
subgroup $B$ of $G$, then $B$ is not a conjugate of $A$.
\item
Let $H$ be a normal subgroup of a group $G$. Prove that $C_G(H)\lhd G$,
and deduce that for $C=C_G(\Phi(G))$, $\Phi(G)C/C\subset \Phi(G/C).$
  \end{enumerate}

\chapter{Symmetric Group}
 \footnote{\it contents group9.tex }
 Let $T = \{1,2,....,n\}$. Then
$A(T) = \{ f:T \to T \st f \mbox{ is one-one }\}$ is a group with respect
to the composition of maps as binary operation (Example \ref{ex7}). We
denote this group by $S_n$.This is called the Permutation group/ Symmetric
group on $n$ symbols. An element of $S_n$ is called a permutation. We have
proved that any finite group of order $n$ is isomorphic to a subgroup of
$S_n$ (Cayley's theorem). Hence it is of particular interest to have a
closer look at $S_n$. We represent an element $\sigma \in S$ by $$
\left(\begin{array}{ccccc}
                 1 & 2 & 3 & \cdots  & n\\
                 i_1 & i_2 & i_3 & \cdots & i_n
         \end{array}\right)$$
where $\sigma (k) = i_k$ for $k = 1,2,\ldots ,n.$ Let us elaborate
the binary
 operation by an example.
\bx\label{ex71}
 Let $\sigma_1,\sigma_2 \in S_4$, where
$$ \sigma_1 = \left(\begin{array}{cccc}
                 1 & 2 & 3 & 4\\
                 3 & 1 & 4 & 2
               \end{array}\right) ,
\sigma_2 = \left(\begin{array}{cccc}
                 1 & 2 & 3 & 4\\
                 2 & 1 & 4 & 3
         \end{array}\right)$$
Then $(\sigma_1 \sigma_2)(1) = \sigma_1 (\sigma_2(1)) = 1,(\sigma_1
\sigma_2)(2) = 3, (\sigma_1 \sigma_2)(3) = 2$, and $(\sigma_1 \sigma_2)(4)
= 4.$ Hence $$ \sigma_1 \sigma_2 = \left(\begin{array}{cccc}
                          1 & 2 & 3 & 4\\
                          1 & 3 & 2 & 4
               \end{array}\right)$$
\ex We now define : \bd\label{d72} An element $\sigma \in S_n$,
is called a cycle of length $k$ if there exist $k$-elements $i_1,i_2,\dots
,i_k$ in $T=\{1,2,\cdots,n\},\, 1\leq k \leq n $, such that $$ \sigma(i_1)
= i_2, \sigma(i_2) = i_3,\ldots ,\sigma(i_{k-1}) = i_k,
                                   \sigma(i_k) = i_1.$$
and $\sigma(j) = j$ for all $j \in T \setminus \{i_1,i_2,\ldots ,i_k \}.$
In this case $\sigma$ is represented by $(i_1,i_2,\ldots ,i_k).$  \ed
\brs (i) Clearly $(i_1,i_2,\cdots,i_k) = (i_2,i_3,\cdots,i_k,i_1)$ etc.
\\  (ii) Any cycle of length 1 is the identity permutation. Thus for any
$1\leq i\leq n$, $(i)$ is identity permutation.  \\  (iii) If $\sigma$ is
a cycle of length $k>1$, then whenever $a \neq \sigma(a)=b,\; \sigma(b)
\neq b.$\ers {\bf Notation : We shall denote the length of a cycle
$\sigma$ in $S_n$ by $l(\sigma)$.}
\bd\label{d73} Two cycles $\sigma = (a_1,a_2,\ldots ,a_k), \tau =
(b_1,b_2,\ldots ,b_l)$ in $S_n$ are called disjoint if
 $\{a_1,a_2,\ldots ,a_k\}\cap \{b_1,b_2,\ldots ,b_l \} = \emptyset$
\ed \br If $\sigma, \tau \in S_n$ are two disjoint cycles, then
$\sigma \tau =  \tau \sigma$ since for any $1\leq j\leq n$, $\tau(j) \neq
j$ implies $\sigma(j) = j$ and $\sigma(j) \neq j$ implies $\tau(j) = j.$
\er  \bexe\label{ex73ex} Let $\tau,\sigma \in S_n, \; n\geq 3$. If for two disjoint subsets
$A,B$ of $T =\{1,2,\cdots,n\}$, $\tau(i) \neq i$ if and only if $i \in A$
and $\sigma(j) \neq j$ if and only if $j \in B$, then $\sigma\tau =
\tau\sigma$. Thus if $\tau$ and $\sigma$ move no common element then
$\tau\sigma =\sigma\tau.$\eexe
\bd\label{d74} A cycle of length $2$ in $S_n$ is called a {\bf{Transposition}}.
\ed \bo Any cycle $(i_1,i_2,\ldots ,i_k)$ of length $k \geq 2$ is a
product of transpositions. In fact $$(i_1,i_2,\ldots ,i_k) =
(i_1,i_k)(i_1,i_{k-1})\ldots (i_1,i_2).$$ If $n \geq 2$, then $id =
(1,2)(1,2)$ i.e., a cycle of length one is also a product of
transpositions.\\ A decomposition of an element $\sigma \in S_n$ into
transpositions is not unique. For example if $1\neq a,b$, then $(a,b) =
(1,a)(1,b)(1,a)$ as well. \eo
\bexe\label{ec71} Prove that the number of cycles of length $r$ in $S_n$ is $(r-1)!
\cdot {}^nC_r = \frac{1}{r}\cdot {}^np_r.$
\eexe
\bl\label{l76}
Any conjugate of a cycle of length $k$ in $S_n$ is a cycle of length $k$.
Further, any two cycles of length $k$ are conjugate to each other.
\el
\pf Let $\sigma, \tau \in S_n$, and let $$ \tau =
\left(\begin{array}{ccccc}
                 1 & 2 & 3 & \cdots  & n\\
                 i_1 & i_2 & i_3 & \cdots & i_n
         \end{array}\right)$$
Then $(\sigma \tau \sigma^{-1})((\sigma(k)) = \sigma(i_k).$ Hence $$
\sigma\tau\sigma^{-1} = \left(\begin{array}{ccccc}
             \sigma(1) & \sigma(2) & \sigma(3) & \cdots  & \sigma(n)\\
             \sigma(i_1) & \sigma(i_2) & \sigma(i_3) & \cdots & \sigma(i_n)
         \end{array}\right)$$
i.e., the conjugate $\sigma \tau \sigma^{-1}$ of $\tau$ by $\sigma$ is
obtained by applying $\sigma$ on the symbols representing $\tau.$ Now, if
$\tau = (i_1,i_2,\ldots ,i_k)$ is a cycle of length $k$, then clearly
$(\sigma \tau \sigma^{-1}) =  (\sigma(i_1) \sigma(i_2),\ldots
,\sigma(i_k))$ i.e., $(\sigma \tau \sigma^{-1})$ is a cycle of length
$k.$\\ Next, let $\sigma = (i_1,i_2,\ldots ,i_k), \tau = (j_1,j_2,\ldots
,j_k)$ be two cycles of length $k$ in $S_n$. Then for any $\theta \in S_n$
such that $\theta(i_p) = j_p$ for $p = 1,2,\ldots ,k,$ we have $\theta
\sigma \theta^{-1} = (\theta(i_1), \theta(i_2),\ldots ,\theta(i_k))
                           = (j_1,j_2,\ldots ,j_k).$
Hence the result follows. \bt\label{t77} Any $\sigma (\neq Id)\in S_n, n
\geq 2,$ is a product of disjoint cycles
 $(\neq id)$ uniquely upto order of cycles.
\et \pf Let us define a relation $\sim$ on $T = \{1,2,\ldots ,n\}$
as below : \\  For $i,j \in T$, $i \sim j$ iff $\sigma^k(i) = j$ for some
$k\in \Z.$  \\  One can easily check that $\sim$ is an equivalence
relation on $T.$ Let for $k\in T$, $[ k ]$ denote the equivalence class of
$k$ with respect to $\sim$ and let $[ k_1 ], [ k_2 ],\ldots ,[ k_m ]$ be
all the distinct equivalence classes of $T$ with more than one element.
Note that for any $k\in T,$ $[ k ] = \{ k, \sigma (k),\ldots
,\sigma^{t-1}(k) \}$ where $t\geq 1$ is the smallest the integer such that
$\sigma^t(k) = k.$ Moreover $\abs{[k]} = t$. Now, let  $$  [ k_1 ]  = \{
k_1, \sigma (k_1),\ldots ,{\sigma}^{t_1 -1}(k_1) \}$$ $$  [ k_2 ] =  \{
k_2, \sigma(k_2),\ldots ,{\sigma}^{t_2 -1}(k_2) \}$$ $$ \cdots \qquad
\cdots$$ $$ [ k_m ]  =  \{ k_m, \sigma(k_m),\ldots ,{\sigma}^{t_m -1}(k_m)
\}$$ where $\abs{[k_{i}]} = t_i$ for all $1\leq i\leq m$.  Then, consider
the cycles
\bea
        \tau_1 & = & (k_1, \sigma(k_1),\ldots ,{\sigma}^{t_1 -1}(k_1)\\
        \tau_2 & = & (k_2, \sigma(k_2),\ldots ,{\sigma}^{t_2 -1}(k_2)\\
               & \vdots &\\
        \tau_m & = & (k_m, \sigma(k_m),\ldots ,{\sigma}^{t_m -1}(k_m) )
\eea in $S_n$ determined by equivalence classes of $T$ with
respect to $\sim$ having more than one element. As any two equivalence
classes are disjoint, the cycles $\tau_1,\tau_2,\ldots ,\tau_m$ are
pairwise disjoint and hence commute with each other. Further, note that if
$k \in [k_i]$ for some $1 \leq i\leq m$, then $\tau_i(k) = \sigma(k)$. We
claim $$ \sigma = \tau_1\tau_2 \ldots \tau_m.$$ If $k \in T$, and $k
\not\in [k_i]$ for any $i$, then $\sigma(k) = k = \tau_p(k)$ for all $p =
1,2,\ldots ,m.$ Hence $$\sigma(k) = (\tau_1\tau_2 \ldots \tau_m)(k) = k$$
Further, if $k\in [k_j]$ for some $1\leq j\leq m$, then $\tau_p(k) = k$
for all $1\leq p\neq j\leq m.$ Therefore $$(\tau_1\tau_2 \ldots \tau_m)(k)
= \tau_j(k) = \sigma(k).$$ Consequently, $$\sigma = \tau_1\tau_2 \ldots
\tau_m$$ Thus any $\sigma (\neq id) \in S_n$ is a product of disjoint
cycles
 $(\neq id).$
We shall, now, prove the uniqueness. Let $$ \sigma = \omega_1 \omega_2
\ldots \omega_s$$ where $\omega_i$ for $i = 1,2,\ldots ,s$ are disjoint
cycles $(\neq id).$ For any $k\in T$, $\sigma(k) = k$ if and only if
$k\not\in \omega_i$ for any $i$. Further, if $k\in \omega_j$ for some
$1\leq j \leq s$, then clearly $\sigma(k)=w_j(k)$ and  $\omega_j = [k].$
i.e., the set of elements in $w_j$ is precisely the set of elements in
$[k]$. Hence $\omega_1, \omega_2, \ldots ,\omega_s$ are formed as above
precisely from all distinct equivalence classes of T with respect to
$\sim$ having more than one element. Thus the uniqueness follows.
$\blacksquare$
\br
 As any cycle is a product of transpositions, any permutation is a
product of transpositions. Further, for any transposition $(a b)$, where
$a\neq 1,b \neq 1$, we have $$(a b) = (1 a)(1 b)(1 a).$$ Hence $S_n$ is
generated by $(1 2),(1 3),.....(1 n).$ \er \bl Let $n\geq 3$. Then center
of $S_{n}$ is trivial i.e., $Z(S_n) = Id$.\el \pf Let $f (\neq Id) \in
S_n$. Then there exists an $1\leq i \leq n$ such that $f(i)=j, j\neq i$.
Choose any $g \in S_n$ such that $g(i) =i$, and $g(j)=k$, $k \neq i,j$
(such a $g$ exists as $n\geq 3$). Then $fg(i)=j$ and $gf(i) = k$. Thus $fg
\neq gf$. Hence $Z(S_n) = Id$.     $\blacksquare$  \bexe Prove that for
any $n\geq 2$, $S_n$ is generated by $(12)$ and $(12\cdots n)$. Further,
show that this is a minimal set of generators.\eexe  Let for $\sigma \in
S_n$,
\be\label{e71}
  \sigma = \tau_1. \tau_2 \ldots .\tau_m
\ee where $\tau_i$ for $i = 1,2,\ldots ,m$ are disjoint cycles.
Here we assume that cycles of length 1 are also included and all integers
$1,2,\cdots,n$ appear in the representation (\ref{e71}). Further, as any two
disjoint cycles
 commute we can assume that if $l_i$ is the length of $\tau_i$, then,
 $1 \leq l_1\leq l_2\leq \ldots \leq l_m.$  We shall call the sequence of
integers $\{l_1,l_2, \ldots ,l_m \}$, the cycle pattern of $\sigma.$ From
\ref{e71}, it is clear that $$l_1+l_2+ \ldots +l_m = n,$$ i.e.,the cycle
pattern determines a partition of the integer $n$. Further, one can easily
see that given any partition of $n$, there exists a permutation with that
cycle pattern. If $e_i$ denotes the number of cycles of length $i$ in the
decomposition \ref{e71}, then we write that the permutation $\sigma$ is of
the type $1^{e_1}2^{e_2} \ldots n^{e_n}.$ Note that here $e_i\geq 0.$
\bex\label{ex72} For the identity permutation $Id$, $Id = (1)(2) \cdots (n).$
Thus in this case $ e_1 = n,$ and $e_i = 0$ for all $i > 1.$
\eex \bx\label{ex73}
 Let
 $$\sigma = \left(\begin{array}{ccccccc}
                 1 & 2 & 3 & 4 & 5 & 6 & 7\\
                 2 & 7 & 5 & 3 & 4 & 6 & 1
         \end{array}\right) \in S_7.$$
Then $\sigma = (127)(354)(6)$ is the decomposition of $\sigma$ into
disjoint cycles. Hence for $\sigma$, $e_1 = 1, e_2 = 0, e_3 = 2$, and
 $e_i = 0$ for all $i > 3$.
\ex
\bl\label{l77}
Any two permutations with same cycle pattern are conjugate and conversely.
\el
\pf Let $\sigma_1, \sigma_2\in S_n$ be two permutations with same cycle
pattern $\{m_1, m_2,\ldots, m_s\}$ and let \bea
  \sigma_1 & = & \omega_1 \omega_2 \ldots \omega_s\\
           & = & (x_1x_2\ldots x_{m_1})(y_1y_2\ldots y_{m_2})\ldots
                    (u_1u_2\ldots u_{m_s})\\
  \sigma_2 & = & \tau_1 \tau_2 \ldots \tau_s\\
           & = & (x_1'x_2'\ldots x_{m_1}')(y_1'y_2'\ldots y_{m_2}')\ldots
                    (u_1'u_2'\ldots u_{m_s}')
\eea be the decompositions of $\sigma_1$ and $\sigma_2$ into disjoint cycles.
 Consider the permutation $$ \theta =
\left(\begin{array}{ccccccccccccc}
 x_1 & x_2 & \ldots & x_{m_1} & y_1 & y_2 & \ldots & y_{m_2} & \ldots
                   & u_1 & u_2 & \ldots & u_{m_s} \\
 x_1' & x_2' & \ldots & x_{m_1}' & y_1' & y_2' & \ldots & y_{m_2}' & \ldots
                   & u_1' & u_2' & \ldots & u_{m_s}' \\
                  \end{array}\right)$$
in $S_n$.Then \bea
 \theta^{-1} \sigma_2 \theta & = & \theta^{-1}(\tau_1 \tau_2 \ldots \tau_s)\theta\\
     & = & \theta^{-1}(\tau_1)\theta \theta^{-1}(\tau_2)\theta \ldots
                     \theta^{-1}(\tau_s)\theta\\
     & = & (\omega_1 \omega_2 \ldots \omega_s)
\eea by the lemma \ref{l76}. Hence $\sigma_1$ and $\sigma_2$ are conjugate. The converse follows
easily from the lemma \ref{l76}. Hence the result is proved.
$\blacksquare$
\bt\label{t78} {\bf (Cauchy)} Suppose $\sigma \in S_n$\  has cycle
pattern $1^{e_1} 2^{e_2} \ldots n^{e_n}.$ Then the number of
permutations conjugate to $\sigma$ is $$ \frac{n!}{ e_1!e_2!
\ldots e_n! 1^{e_1}2^{e_2} \ldots n^{e_n}}.$$ \et \pf Let us
indicate the cycle pattern of $\sigma$ as follows \be\label{e72}
 \sigma = (.)(.)\cdots(.)(..)(..)(..) \ldots \\
\ee where (.) indicates a 1-cycle, (..) indicates a 2-cycle and so
on. All the $n$ symbols $\{1,2,\cdots,n\}$ are involved in the presentations
(\ref{e72}) of $\sigma$. We know that any two permutations are conjugate
if and only if they have the same cycle pattern (Lemma \ref{l77}). Hence,
any permutation of the $n$ symbols shall give conjugate to $\sigma$ and
conversely. Thus apparently $\sigma$ has $n!$ conjugates; the number of
permutations on $n$ symbols. Note that two distinct permutations of the
symbols in (\ref{e72}) may give rise same element in $S_n$. This can
happen in two ways.\\ (a) $j$-cycles are permuted among themselves\\ (b)
each $j$-cycle is written in different order.\\ By way of (a), as there
are $e_j$ distinct $j$-cycles, $e_j!$ repetitions occur i.e., $e_j!$
permutations in $S_n$ give same element on operating on $\sigma$.\\
Further, note that any $j$-cycle can be written in $j$ different ways as
$$ (a_1a_2 \ldots a_j) = (a_2 \ldots a_j a_1)
                        = \ldots = (a_ja_1 \ldots a_{j-1})$$
Hence there are $j^{e_j}$ different ways of writing one element by
changing writing pattern of $e_j$ $j$-cycles involved in  the
representation of $\sigma$. Hence while counting the number of conjugates
of $\sigma$ as $n!$ we actually counted each element $ e_1!e_2! \ldots
e_n! 1^{e_1}2^{e_2} \ldots n^{e_n}$ times. Therefore the number of
distinct conjugates of $\sigma$ is $$\frac{n!}{ e_1!e_2! \ldots e_n!
1^{e_1}2^{e_2} \ldots n^{e_n}}.  \blacksquare$$ \bco\label{c71} For
$\sigma = (1 2....n)$ in $S_n = G$, $C_G(\sigma) = gp\{\sigma \}$ i.e.,
the centralizer of $\sigma$ in $G$ is the cyclic group generated by
$\sigma$. \eco \pf As $\sigma = (1 2....n)$ is a $n$-cycle we have $e_n =
1$ and $e_i = 0$ for all $i < n.$ Hence by Cauchy's theorem, the number of
conjugates of $\sigma$ is equal to $\frac{n!}{n} = (n-1)!.$ However, the
number of conjugates of $\sigma$ is equal to $$ [G:C_G(\sigma)] =
\frac{\abs{G}}{ \abs{C_G(\sigma)}} =
                          \frac{n!}{\abs{C_G(\sigma)}}.$$
Therefore \ba{rl}
             & (n-1)! = \frac{n!}{\abs{C_G(\sigma)}}\\
\Rightarrow  & n = \abs{C_G(\sigma)} \ea We know that $gp\{\sigma \}
\subseteq C_G(\sigma)$, and $\circ(\sigma) = n.$ Hence $gp\{\sigma \} =
C_G(\sigma).    \blacksquare$ \br The result in the corollary holds for
any n-cycle.
\er
\bt\label{t75.5t} Let $n \geq 4$. Then every element of $S_n$ is a product
of two elements of order 2.\et \pf For any transposition $\tau$ in $S_n$,
$\tau^2 = Id$. Hence the result holds for the identity element in $S_n$.
Now, let $\sigma(\neq Id) \in S_n$. By the theorem \ref{t77}, $$\sigma =
w_1 w_2 \cdots w_s,$$ $s\geq 1$, where $w_i's$ are pair-wise disjoint
cycles with $l(w_i)>1$ for each $i$. We can assume that $l(w_i) \leq
l(w_{i+1})$ for all $1\leq i\leq s-1$, since any two disjoint cycles
commute. We shall first consider the case when $l(w_i)=2$ for all
$i=1,2,\cdots,s$. If $s = 1$ then $\sigma = w_1$. Let $\sigma = (ij)$ for
some $1\leq i< j\leq n$. As $n\geq 4$, we can choose $1\leq k < l\leq n$
such that $\{k,l\} \cap \{i,j\} = \emptyset$. Clearly $$\sigma =
(ij)(kl)\;(kl)$$ where $(ij)(kl)$ and $(kl)$ are elements of order 2 in
$S_n$. Next, let $s>1$. Then $$\sigma = (w_1 w_2 \cdots w_{s-1})w_s$$
where $(w_1 w_2 \cdots w_{s-1})$ and $w_s$ are elements of order 2 in
$S_n$. Thus if $l(w_i)=2$ for all $1\leq i\leq s$, then $\sigma$ is a
product of two elements of order 2 in $S_n$. Now, let $l(w_i) = 2$ for all
$0\leq i\leq r< s$. If we prove that for any cycle $w\in S_n$ with $l(w)
\geq 3$, $w = cd$ where $c$ and $d$ are two elements of order 2 in $S_n$
which move no element left fixed by $w$, then let $w_j = c_j\, d_j$, for
all $r+1 \leq j \leq s$, where $c_j ,d_j$ have required properties, then
$$
\begin{array}{ccc}
  \sigma & = & w_1 w_2 \cdots w_r (c_{r+1}d_{r+1})\cdots(c_{s}d_{s}) \\
   & = & (w_1 w_2 \cdots w_r c_{r+1}\cdots c_{s})(d_{r+1}\cdots d_{s})
\end{array}$$ Note that $w_{1},w_{2},\cdots , w_{r}, c_{r+1},\cdots c_{s}$
and $d_{r+1},\cdots , d_{s}$ are pairwise commuting elements, each of
order 2, in $S_n$. Thus $c = w_{1}w_{2}\cdots w_{r}c_{r+1}\cdots c_{s}$
and $d = d_{r+1}\cdots d_{s}$ are elements of order 2, and $\sigma = cd$.
Now, to complete the proof, we shall show that for any cycle $w \in S_n$
with $l(w)\geq 3$, $w = cd$ where $c,d$ are elements of order 2 such that
$c,d$ move no element left fixed by $w$ and $d$ does not move all the
symbols moved by $w$. We shall do this by induction on $l(w)=k\geq 3.$ We
can assume, if necessary by change of notation, that $w = (1 2\, 3 \cdots
k)$. If $k=3$, then $w = (1 2 3 ) = (1 2)(2 3).$ Hence our assertion holds
in this case. Let $k\geq 4$. Then by induction hypothesis
\be\label{egroup}
  w_1  = (12\cdots k\!-\!1)= cd
\ee  where $c,d$ are elements of order 2 with required properties. Let for
$1\leq l\leq k-1$, $l$ is fixed under $d$. Then from the equation
(\ref{egroup}), $$(lk)(12\cdots k\!-\!1) = (lk)c(lk)\,(lk)d$$
\be\label{ering}\Rightarrow\;\;\;\;\; e = (12\cdots l\!-\!1 k l l\!+\!1\cdots k\!-\!1)
= c_1 d_1
\ee  where $c_1 = (lk)c(lk)$ each and $d_1 = (lk)d$ are elements of order 2
and $l$ is left fixed by $c_1$. Now, there are two possibilities : \\ {\bf
Case I} The element $d_1$ keeps fixed one of the symbols $1,2,\cdots,k$.
\\ In this case, observe that the cycle $w = (12\cdots k)$ is a conjugate
of the cycle $e$ (use; $l(w)=l(e)=k$). Hence we can write $w = cd$ where
$c,d \in S_n$ have order 2 such that $c,d$ move no element left fixed by
$w$ and $d$ fixes at-least one symbol from $1,2,\cdots,k$. \\ {\bf Case
II} The element $d_1$ move all the symbols $1,2,\cdots,k$.  \\  From the
equation (\ref{ering}), we have $$e^{-1} = d_{1}^{-1}c_{1}^{-1}$$ where
$d_{1}^{-1}$ and $c_{1}^{-1}$ are elements of order 2 which move no element
left fixed by $e^{-1}$ and $c_{1}^{-1}$ keeps fixed at-least one symbol,
infact $k$, out of the symbols $1,2,\cdots,k$. Now, we are in the case I.
Hence the assertion holds. Thus the result is proved.       $\blacksquare$
\bd\label{d75} Let $\Z[X_1,X_2, \ldots ,X_n]$ be the polynomial ring over
integers in $n$-indeterminates $X_1, X_2, \ldots ,X_n.$ Then for any
polynomial $f(X_1,X_2, \ldots ,X_n)$ in $\Z[X_1,X_2, \ldots ,X_n]$ and
$\sigma \in S_n$, we define: $$\sigma(f(X_1,X_2, \ldots ,X_n )) =
                     f(X_{\sigma(1)},X_{\sigma(2)}, \ldots ,X_{\sigma(n)}).$$ \ed
  \brs (i) The action of $\sigma$ defined in the definition gives an
  action of the group $S_n$ over $\Z[X_1 ,\cdots X_n]$.  \\  (ii) For
  $\Delta=\prod_{1\leq i<j\leq n}(X_{i}-X_{j})$, $\sigma(\Delta) = \pm
  \Delta$ for any $\sigma \in S_n$. \\ (use; For every pair $(i,j)$, $1
  \leq i<j \leq n$, the pair $(\sigma(i),\sigma(j))$ is unique and either
  $\sigma(i) < \sigma(j)$ or $\sigma(j) < \sigma(i)$)\ers  \bd A permutation
  $\sigma\in S_n$ is called even (odd) if $\sigma(\Delta) =
  \zeta(\sigma)\Delta$ where $\zeta(\sigma) = 1\, (-1)$. The integer
  $\zeta(\sigma)$ is called the {\bf sign} of the permutation $\sigma$.\ed
  \bexe Let $f(X_{1},X_{2},\cdots,X_{n}) \in \Z[X_{1},\cdots X_{n}]$.
  Prove that $H = \{\sigma \in S_n \st \sigma(f(X_{1},X_{2},\cdots X_{n}))
  = f(X_{1},X_{2},\cdots X_{n})\}$ is a subgroup of $S_n$.\eexe \brs (i)
  The group $H$ in the exercise is called the {\bf group of symmetries} /
  {\bf symmetry group} of $f(X_{1},X_{2},\cdots X_{n})$. \\ (ii) Even
  permutations form the group of symmetries of $\Delta = \prod_{1\leq i<j \leq
  n}(X_{i}-X_{j})$.\ers \bd A polynomial $f(X_{1},X_{2},\cdots X_{n}) \in \Z[X_{1},X_{2},\cdots X_{n}]$
  is called symmetric if $\sigma(f(X_{1},X_{2},\cdots X_{n})) = f(X_{1},X_{2},\cdots
  X_{n})$ for all $\sigma \in S_n$.\ed  \br A polynomial $f(X_{1},X_{2},\cdots
  X_{n}) \in \Z [X_{1},X_{2},\cdots X_{n}]$ is a symmetric polynomial if
  and only if every homogeneous component of $f(X_{1},X_{2},\cdots X_{n})$
  is symmetric.\er \bexe Let $p(X_{1},X_{2},\cdots X_{n}) =
  X_{1}+X_{2}+\cdots +X_{n-1}+X_{1}^{2}$ be a polynomial in$\Z [X_{1},X_{2},\cdots
  X_{n}]$, $n\geq 3$. Find the group of symmetries of $p(X_{1},X_{2},\cdots
  X_{n})$ in $S_n$, and prove that this group is isomorphic to $S_{n-2}.$
  \eexe
 Let us, now, make some observations regarding the signs of permutations.
\bl\label{l78}
Let $\sigma, \tau \in S_n$. Then\\ (i) $\zeta(\sigma \tau) =
\zeta(\sigma) \zeta(\tau).$\\ (ii) $\zeta(\sigma^{-1} \tau \sigma)
= \zeta(\tau)$\\ (iii) If $\sigma$ is a transposition,
$\zeta(\sigma) = -1.$
\el
\pf (i)  We have \bea
        (\sigma \tau)(\Delta) & = & \sigma(\tau(\Delta))\\
                              & = & \sigma(\zeta(\tau)\Delta)\\
                              & = & \zeta(\tau)(\sigma(\Delta))\\
                              & = & \zeta(\tau) \zeta(\sigma) \Delta
\eea Therefore $\zeta(\sigma \tau)\Delta = \zeta(\sigma)\zeta(\tau)
\Delta.$\\ (ii) Note that for the identity permutation $id \in S_n$
$$id(\Delta) = \Delta.$$ Therefore $\zeta(id) = 1.$ Hence from (i) we get
$$ \zeta(\sigma) \zeta(\sigma^{-1}) = \zeta(\sigma \sigma^{-1})
                                       = \zeta(id) = 1.$$
Therefore $$\zeta(\sigma^{-1}\tau \sigma) =  \zeta(\sigma^{-1})
\zeta(\tau)
                      \zeta(\sigma) = \zeta(\tau) $$
(iii) We have proved that any two cycles of same length are conjugate
(Lemma \ref{l76}). Hence, any two  transpositions  are  conjugate.
Therefore, in view of (ii), it suffices to prove the result for $\sigma =
(1 2).$ We have
\be\label{e73}
\begin{array}{rcr}
 \Delta & = & (X_1-X_2)(X_1-X_3)(X_1-X_4) \ldots (X_1-X_n)\\
        &   &                (X_2-X_3)(X_2-X_4) \ldots (X_2-X_n)\\
        &   &                          (X_3-X_4) \ldots (X_3-X_n)\\
        &   &                                \cdots \cdots \cdots\\
        &  &                                       (X_{n-1}- X_n)
\end{array}
\ee Hence \be\label{e74}
\begin{array}{rcr}
\sigma (\Delta)  & = & (X_2-X_2)(X_2-X_3)(X_2-X_4) \ldots
(X_2-X_n)\\
                 &  &        (X_1-X_3)(X_1-X_4) \ldots (X_1-X_n)\\
                &   &                (X_3-X_4) \ldots (X_3-X_n)\\
               & &                      \cdots\cdots\cdots \\
          & &                              (X_{n-1}- X_n)
\end{array}
\ee where the expression from third row onwards remain unchanged.
It is clear from expressions \ref{e73} and \ref{e74} that $\sigma(\Delta)
= -\Delta.$ Hence $\zeta(\sigma) = -1.$
\bco\label{c72} A permutation $\sigma \in S_n$ is even (odd) if
and only if $\sigma$ is a product of even (odd) number of transpositions.
\eco \pf The proof is straightforward. \brs\label{rs9.32rs} (1) By the lemma \ref{l78} (ii)
 it is clear that the conjugate
 of an even
(odd) permutation is
 even (odd).
(2) By the lemma \ref{l78} (i), the map $\zeta : S_n\to \{\pm1 \}$ is a
homomorphism of groups. \ers
\bl\label{l79}
 Let $G$ be a group of permutations on $n$ symbols. i.e., G is a subgroup of
 $S_n$. Then even permutations in $G$ form a normal subgroup (say) $H.$
 Further, $[G:H]\leq 2.$
\el
\pf Clearly $H\neq \emptyset$, since $id\in H.$ Let $\sigma,\tau
\in H.$ Then \bea
     \zeta(\sigma \tau^{-1}) & = & \zeta(\sigma) \zeta(\tau^{-1})\\
                             & = & \zeta(\sigma) (\zeta(\tau))^{-1}\\
                             & = & 1
\eea Hence $\sigma \tau^{-1}\in H.$  Thus $H$ is a
subgroup of $G.$ Further, as conjugate of an even permutation is even
(Remark : \ref{rs9.32rs}(1)), $H$ is normal in $G.$ For any even
permutation $\sigma \in G$, $\sigma H = H$, as $\sigma \in H.$ However, if
$\tau_1, \tau_2 \in G$ are odd permutations, then
\bea
       \zeta(\tau_1^{-1} \tau_2) & = & \zeta(\tau_1^{-1}) \zeta(\tau_2)\\
                                 & = & \zeta(\tau_1)^{-1} \zeta(\tau_2)\\
                                 & = & 1.
\eea Therefore $\tau_1^{-1} \tau_2  \in H$, and we have  \ba{rl}
               &   \tau_1^{-1} \tau_2H = H\\
\Rightarrow    &   \tau_2H = \tau_1H\\  \mbox{Consequently    }   & [G:H]
\leq 2.     \blacksquare \ea
\br From the above lemma, either $G = H$ or $\abs{H} = \frac{1}{2} \abs{G}.$
\er \bd\label{d710} The subgroup $A_n$ of all even permutations of $S_n$
is called the Alternating group of degree $n.$ \ed \bo $\abs{A_n} =
\frac{1}{2} \abs{S_n} = \frac{n!}{2}$ for all $n > 1,$ and $A_n$ is
normal in $S_n.$ \eo
\bt\label{t711} The Alternating group $A_n$ is generated by three
cycles for all $n\geq 3.$ \et \pf For any three cycle $(abc)$, we have
$(abc) = (ac)(ab)$. Hence, $(abc)$ is an even permutation, and all three
cycles lie in $A_n$. Next, we know that any element in $A_n$ is a product
of even number of transpositions. Hence to prove our claim, it is
sufficient to show that the product of any two transpositions is a product
of three cycles. Let $(ab),(cd)$ be two transpositions. We shall prove
that the product $(ab)(cd)$ is a product of three cycles. Let us
consider:\\ Case I : $\{a,b\} \cap \{c,d\} = \emptyset$\\ In this case
\bea
 (ab)(cd) & = & (ab)(ac)(ac)(cd)\\
          & = & (acb)(ca)(cd)\\
          & = & (acb)(cda).
\eea Case (ii) :$\{a,b\} \cap \{c,d\} \neq \emptyset$\\ If
$\{a,b\} =  \{c,d\}$, then clearly $(ab)(cd) = id.$ Hence, let
$\{a,b\} \neq \{c,d\}$. Then there is exactly one common element.
Let $b = c.$ Then \bea
              (ab)(cd) & = & (ab)(bd)\\
                       & = & (ba)(bd)\\
                       & = & (bda)
\eea Thus it follows that the product of any two transpositions is
a product of three cycles. This completes the proof.     $\blacksquare$
\\ We, now, prove a result giving a more elaborate set of generators for
$A_n$.
\bt\label{t712} The alternating group $A_n (n > 2)$ is generated by the
three cycles $(123),(124),....,(12n).$ \et \pf By the theorem \ref{t711},
it suffices show that a three cycle $(a_1a_2a_3)$ is a product of the
given three cycles. To prove this, we shall consider :\\ Case I : $a_i =
1$ for some $i$.\\ There is no loss of generality in supposing that $a_1 =
1.$ Then, if $a_i \neq 2$ for $i = 2,3,$ we have
\bea
             (1a_2a_3) & = & (1a_3)(1a_2)\\
                       & = & (1a_3)(12)(12)(1a_2)\\
                       & = & (12a_3)(1a_22)\\
                       & = & (12a_3)(12a_2)^2
\eea However, if $a_3 = 2$, then it is clear from above that $(1a_2 a_3) =
(12a_2)^{2}$ and for $a_2 = 2$, there is nothing to prove.\\ Case II :
$a_i \neq 1$ for any $i$.\\ In this case, we have
\bea
                (a_1a_2a_3) & = & (a_1a_3)(a_1a_2)\\
                            & = & (a_1a_3)(a_11)(a_11)(a_1a_2)\\
                            & = & (a_11a_3)(a_1a_21)
\eea By case I and case II, it is immediate that $(a_1a_2a_3)$ is $a$
product of the given cycles. Hence the result follows.
\bl\label{l713}
If a normal subgroup $N$ of $A_n (n > 4)$  contains a three cycle,
 then $N = A_n$.
\el
\pf If necessary, by change of notation, we can assume $(123)\in N.$ Then
for any $k > 3,$ and $l\neq 1,2,3,k$, we have $$(3l)(3k)(123)(3k)(3l) =
(12k)$$ Therefore, since $N$ is normal in $A_n$, $(12k) \in A_n$ for all
$k\geq 3.$ Hence by the theorem \ref{t712}, $N = A_n$. \bd\label{d714} A group
$G\neq \id$ is called a simple group if it has no non-trivial normal
subgroups i.e., it has no normal subgroups other than identity subgroup
and the whole group. \ed \bt\label{t715} The alternating group $A_n (n >
4)$ is simple. \et \pf Let $N (\neq \id)$ be a normal subgroup of $A_n$.
The result will follow if we prove that $N$ contains a three cycle
(Lemma\ref{l713}). Choose an element $\tau (\neq \id)$ in $N$ which leaves
maximum number of symbols fixed. We shall prove that $\tau$ is a three
cycle. This will follow if we prove that $\tau$ displaces exactly three
symbols. Assume this is not true. Let us write $\tau$ as a  product of
disjoint cycles. Then either there is a cycle of length atleast three in
the decomposition or else $\tau$ is a product of transpositions. (Note
that $\tau$ is not a transposition as it is an even permutation.) \\
 In the first case, if necessary by change of notation, we can write
\be\label{e75}
 \tau = (123...)
\ee Let us note that $\tau$ shall move atleast five symbols in
this case , as $A_n$ contains no cycle of length four (Use : a
four cycle is an odd
 permutation). We can assume that $\tau$ moves $1,2,3,4,5.$ Now, for
 $\sigma = (345)\in A_n$, we have
$$ \tau_1 = \sigma \tau \sigma^{-1} = (124...) \in N.$$ Clearly $\tau^{-1}
\tau_1$ is a non-identity element in $N$ and moves at-best the same
elements which $\tau$ does. We,however, have $\tau^{-1} \tau_{1}(1) = 1.$
This contradicts the choice of $\tau$ that $\tau (\neq Id)$ is an element
of $N$ which moves minimum number of symbols. Hence $\tau$ does not
 have the form (\ref {e75}).\\
The other possibility is that \be\label{e76}
 \tau = (12)(34).....
\ee Again for $\sigma = (345)$, we have $$ \tau_1 = \sigma \tau
\sigma^{-1} = (12)(45) \ldots$$ is in $N.$ Hence $\tau^{-1} \tau_1 \in N$
is a non-identity element and moves atmost one more symbol than the
symbols moved by $\tau.$ We, however, have that $\tau^{-1} \tau_1$ keeps
$1,2$ fixed. This again contradicts the choice of $\tau.$ Consequently
$\tau$ moves exactly three symbols, and hence is a three cycle. Thus the
result follows.         $\blacksquare$  \brs (i) The group $A_4$ is not simple.In fact
 $$ N  = \{ id,(12)(34),(13)(24),(14)(23) \}$$ is a nontrivial normal subgroup of
$A_4$.\\ (ii) The subgroup $N$ in (i) is abelian. Hence $$ C = \{
id,(12)(34) \} = gp\{(12)(34)\}$$ is normal in $N.$  It is , however, easy
to check that the subgroup $C$ is not normal in $A_4$. This shows that
normal subgroup of a normal subgroup need not be normal.\\ (iii) $A_n$ is
simple for $n \neq 3,4.$ \ers \bt\label{t716} Let $n > 4.$ Then $A_n$ is
the only proper normal subgroup of $S_n$. \et \pf Let $H$ be a proper
normal subgroup of $S_n$. We claim $\circ(H) > 2.$ If not, then $H = \{
id,\sigma \st \circ(\sigma) = 2 \}$. As $\circ(\sigma) = 2,$ $\sigma$ is
either a transposition or is a product of disjoint transpositions. If
$\sigma$ is a transposition, let $\sigma = (ab).$  Choose $c$ other than
$a$ and $b$ (Use: $n > 4$). Then as $H$ is normal in $S_n$, $$
(ac)(ab)(ac) = (cb) \in H.$$  This contradicts the assumption that
$\circ(H) = 2.$ Hence $\sigma \neq (ab).$ In the other case, let $$ \sigma
= (a_1a_2)(a_3a_4) \sigma_1$$ where $\sigma_1$ is a product of disjoint
transpositions not involving $ a_1,a_2,a_3,a_4.$ Then for $\tau =
(a_1a_2a_3)$, we have $$ \tau \sigma \tau^{-1} = (a_2a_3)(a_1a_4) \sigma_1
\in H$$ Clearly $\tau\sigma\tau^{-1} \neq Id$, and also $\tau \sigma
\tau^{-1} \neq \sigma.$ This again contradicts that $\circ(H) = 2.$ Hence
the claim that $\circ(H)
> 2$ holds. Now, let $D = A_n \cap H.$ Then $D$ is normal in $A_n$ (Use:
$H$ is normal in $S_n$), and is the subgroup of all even permutations in
$H.$ Hence, by the lemma \ref{l79}, $\abs{D} = \frac{1}{2} \abs{H}
> 1.$ Consequently, by the theorem \ref{t715}, $D = A_n$. Thus $H
\supseteq A_n.$  As $[S_n:A_n] = 2,$ $H = A_n$, since $H$ is a proper
subgroup of $S_n$. This completes the proof.     $\blacksquare$
\bt\label{t717} Let $n \geq 2$. Then there exists no outer automorphism of $S_n$,
for $n\neq 6$ i.e., $Aut(S_n) = IAut(S_n)$ for $n\neq 6.$ \et \pf Clearly
$Aut(S_2) = IAut(S_2)$. Hence, let $n\geq 3.$ Now, let $\alpha$ be any
automorphism of $S_n$. Note that\\ (i) Any conjugacy class of $S_n$  is
mapped to a conjugacy class of $S_n$ under $\alpha$, i.e., for any $\sigma
\in S_n$, image of the conjugacy class of $\sigma$ under $\alpha$ is the
conjugacy class of $\alpha(\sigma).$\\ (ii) Any two permutations of $S_n$
are conjugate if and only if they have same cycle pattern.\\ (iii) For any
$\sigma \in S_n$, order of $\alpha(\sigma)$ is  equal to the order of
$\sigma.$\\ Put
\bea
     B_k & = & \mbox {The set of all products of $k$ disjoint transpositions
               in $S_n$.}\\
         & = & \left \{ \prod_{i=1}^k (a_ib_i) \st (a_ib_i) \mbox { \& }
                    (a_jb_j) \mbox { are disjoint if }i\neq j \right\}
\eea By (ii), each $B_k$ is a conjugacy class for $k\geq 1.$
Further, an element of $S_n$ has order $2$ if and only if it is a
 product of disjoint transpositions. Hence, $\alpha(B_1) = B_k$ for
some $k\geq 1.$  We shall show that $\alpha(B_1) = B_1.$ First of all,
note that $$
 \abs{B_k} = \frac{1}{k!} \prod_{i=0}^{k-1} \left(\begin{array}{c}
                n-2i\\ 2 \end{array} \right) \mbox { for all } k\geq 1.$$
If $\alpha(B_1) = B_k,$ then $\abs{B_1} = \abs{B_k}.$ Therefore $$
       \frac{1}{k!} \prod_{i=0}^{k-1} \left(\begin{array}{c}
       n-2i\\ 2 \end{array} \right)  =
           =  \left(\begin{array}{c}
                n\\ 2 \end{array} \right)$$
Hence, we get \be\label{e77}
 (n-2)(n-3) \ldots (n-2k+2)(n-2k+1) = k! 2^{k-1}
\ee By the equation \ref{e77}, $n\geq 2k.$ Therefore\\ $(n-2)(n-3)
\ldots (n-2k+2)(n-2k+1)$ \bea
 & \geq & (2k-2)(2k-3) \ldots 2.1.\\
              & = & ((2k-2)(2k-4) \ldots 2)((2k-3)(2k-5) \ldots 1))\\
              & = & 2^{k-1} (k-1)!\cdot ((2k-3)(2k-5) \ldots 1))
\eea Note that $2k-3 > k$ if and only if $k >3.$ Hence for $k\geq
4,$ the equation \ref{e77} does not hold for any $n\geq 2k.$  We shall now
consider the cases $k = 2, k = 3.$\\ Case I :  $k = 3$   \\  As $n\geq
2k$, we have $n\geq 7$ since $n \neq 6$. As $5.4.3.2.1.
> 3!.2$, the equation \ref{e77} does not hold true for any $n\geq 7.$\\ Case II : $k
= 2$ \\  Here $n\geq 4.$  It is easy to see that the equation \ref{e77}
does not hold for
 any $n\geq 4.$\\   \\
In view of the above we conclude that $\alpha(B_1) = B_1.$ We shall now
show that there exists an inner automorphism $\theta$ of $S_n$ such that
$$ \alpha(1,i) = \theta(1,i) \mbox { for all } i = 1,2,\ldots n.$$ As $\{
(1,i) \st i = 1,2,\ldots n. \}$ generate $S_n$, this will complete the
proof. Let $$ \alpha(1,2) = (i,j)$$ By (ii), $(1,2)$ and $(i,j)$ are
conjugate to each other. Thus there exists an inner automorphism $\theta_1
$ of $S_n$  such that $$  \alpha(1,2) = \theta_1(1,2).$$ Put $\alpha_1 =
\theta_1^{-1} \alpha$. Then $\alpha_1(1,2) = (1,2).$ Let $$ \alpha_1(1,3)
= (i,j)$$ Then \ba{rl}
               & \alpha_1((1,2)(1,3)) = (1,2)(i,j)\\
\Rightarrow    & \alpha_1(1,3,2) = (1,2)(i,j) \ea As $(1,3,2)$ is an
element of order $3$, order of the element $(1,2)(i,j)$ is also $3$. Hence
$\{1,2\} \cap \{i,j \} \neq \emptyset$, since otherwise order of
$(1,2)(i,j)$ shall be $2$. We can assume, without any loss of generality,
that $i = 1$. Then $$ \alpha_1(1,3) = (1,j)$$ where $j \geq 3.$  If $j >
3$, then $$ \alpha_1(1,3) = (3,j)(1,3)(3,j)$$ i.e., $\alpha_1(1,3)$ is the
conjugate of $(1,3)$ by $(3,j).$ Thus, if $\theta_2$ denotes the
conjugation by $(3,j)$, then for $\alpha_2 = \theta_2^{-1} \alpha_1$, we
have $$ \alpha_2(1,2) = (1,2) \mbox{ and } \alpha_2(1,3) = (1,3).$$ Thus
it is clear that, there exists an inner automorphism $\gamma$ of $S_n$
such that $$ \alpha(1,2) = \gamma(1,2) \mbox{ and } \alpha(1,3) =
\gamma(1,3).$$ Now, suppose for an inner automorphism $\theta$ of $S_n$,
$$ \theta(1,i) = \alpha(1,i)$$ for all $i = 1,2,\ldots ,r$ where $2 < r
<n$, and let $\beta = \theta^{-1} \alpha.$  Then $\beta(1i) = (1i)$ for all $1\leq i\leq
r$. If
 $$ \beta(1,r+1) = (i,j)$$ Then \bea
         \beta((1,2)(1,r+1)) & = & (1,2)(i,j)\\
         \beta((1,3)(1,r+1)) & = & (1,3)(i,j)
\eea Hence, as seen above $$  \{1,2\}\cap \{i,j \}\neq \emptyset,
\{1,3\}\cap \{i,j \}\neq \emptyset.$$ Note that $$\beta(23) =
\beta((13)(12)(13)) = (23)$$ Hence, as $\beta$ is one-one, $\beta(1 r+1) \neq (23)$ i.e.,
$(ij) \neq (23).$ Thus we can assume $i = 1$, and therefore $$
\beta(1,r+1) = (1,j).$$ Clearly $j\geq r+1$. If $j > r+1$, then $$
\beta(1,r+1) = (1,j) = (r+1,j)(1,r+1)(r+1,j)$$ i.e., $\beta(1,r+1)$ is a
conjugate of $(1,r+1)$ by $(r+1,j)$. Let us denote by $\theta_{r+1}$ the
conjugacy map of $S_n$ by $(r+1,j)$. Then \ba{rl}
            & \beta(1,r+1) = \theta_{r+1}(1,r+1)\\
\Rightarrow & \theta_{r+1}^{-1} \beta(1,r+1) = (1,r+1). \ea Hence,
 for all $i = 2,3,\cdots,r+1$, we have \ba{rl} & \theta_{r+1}^{-1} \beta(1,i) =
 (1,i)\\
   \Rightarrow & \beta(1,i) =
\theta_{r+1}(1,i)\\ \Rightarrow & \alpha(1i) = \theta\theta_{r+1}(1i) \ea
Therefore by induction there exists an inner
 automorphism $\theta$ of $S_n$  such that
$$ \alpha(1,i) = \theta(1,i)$$ for all $i = 2,3,\ldots ,n.$ Consequently,
$\alpha$ is an inner automorphism of $S_n.$ Hence the result follows.
$\blacksquare$  \clearpage
 {\bf EXERCISES}
\begin{enumerate}
  \item
  (i) For the element $(432156)$ in $S_7$, write $(432156)(25)$ and
  $(432156)(27)$ as product of disjoint cycles. \\ (ii) Prove that the
  product of a cycle and a transposition in $S_n$ is either a cycle or
  product of two disjoint cycles or identity.   \\  (iii) Prove that any
  three cycle is a commutator in $S_n.$
  \item
  (i) For the element $(125678)$ in $S_8$ write $\sigma^2$ and $\sigma^3$
  as product of disjoint cycles. \\ (ii) Let $\sigma(\neq Id) \in S_n$ be
  a cycle of length $t$. Prove that for any divisor $d \geq 1$ of $t$,
  $\sigma^d$ can be expressed as a product of $d$ disjoint cycles each
  of length $t/d$.
  \item
  (i) Write all elements of order $\leq 2$ in $S_5$ and $S_4$. \\ (ii) Let
  $\alpha(n)$ denotes the number of elements of order $\neq 2$ in $S_n$.
  Prove that $$\alpha(n+1) = \alpha(n)+n\alpha(n-1).$$
  \item
  Prove that any cycle of length $k$ has order $k$. Further, show that if
  for $\sigma \in S_n$, $\sigma = c_1c_2 \cdots c_m$, where $c_i's$ are
  pairwise disjoint cycles with $l(c_i) = l_i$ for all $1 \leq i \leq m$,
  then $\circ(\sigma) = l.c.m.$ of $l_1,l_2,\cdots l_m.$
  \item
  Write the cycle $(123456789)$ as the product of two elements of order 2
  in $S_9$.
  \item
  Let $A$ be a set of $1\leq r<n$ distinct natural numbers in
  $\{1,2,\cdots,n\}$. Prove that $H = \{\sigma \in S_n \st
  \sigma(A) \subset A\}$ is a subgroup of $S_n$, and show that $H$ is
  isomorphic to the direct product $S_{n-r}\times S_r$.
  \item
  Prove that the transposition $(12),(23),\cdots,(n-1\; n)$ together
  generate $S_n$.
  \item
  Let $\sigma \in S_n$ be an $n$-cycle and $\tau$ be a transposition.
  Prove that $S_n$ is generated by $\tau$ and $\sigma$.
  \item
  Show that for any odd integer $n$, (123)and $(123\cdots n)$ generate
  $A_n$, and if $n$ is an even integer than (123) and $(23\cdots n)$
  generate $A_n$.
  \item
  Prove that $A_n$ is the derived group of $S_n$ for all $n$.
  \item
  Find all the subgroups of order 6 in $S_4$.
  \item
  Let $\alpha,\beta \in S_n$, and $\alpha\beta = \beta\alpha.$ Prove that for $A
  = \{i \in A \st \alpha(i) = i\}$, $\beta(A) = A$.
  \item
  Prove that any symmetric polynomial of degree $t =$ 2 or 3 in
  $\Z[X_1,X_2,X_3]$ is a polynomial in $s_1 = X_1+X_2+X_3$, $s_2 =
  X_1X_2+X_1X_3+X_2X_3$ and $s_3 = X_1X_2X_3$ over $\Z$.
  \item
  Find the symmetry group of the polynomials $XYZ+XY+ZX+X^2$ and
  $(X+wY+w^2Z)^3$ in $\Z[X,Y,Z]$ where $w = (-1+\sqrt{3})/2.$
  \item
  Let $n \geq 2$. Show that there exists a subspace $W$ of the metric space
  $\R^{n-1}$ with $\abs{W} = n$ such that $S_n$ is isomorphic to $I(W)$,
  the group of isometrics of $W$.
  \item
  Prove that $Gl_2(\Z_2)$ is isomorphic to $S_3$.
  \item
  Let $n\geq 2$, and let $$g = l.c.m.\{\circ(\sigma) \st \sigma \in S_n\}$$
  $$t = l.c.m.\{\circ(\tau) \st \tau \in A_n\}.$$ Prove that $2t = g$ if
  $n = 2^k$ or $2^{k}+1$ and $t = g$ otherwise.
  \item
  Prove that $S_{999}$ contains an cyclic subgroup of order 1111, and this
  subgroup fixes atlast one $1\leq i \leq 999.$
  \item
  Let $n>1$, and let $e_{ij}\; i=1,2,\cdots,n$, be the transpose of the
  row vector $(0,\cdots,1^{i^{th}},0,\cdots,0).$ Prove that $\sigma \in
  S_n$ is even or odd if and only if
  $det(e_{\sigma(1)},e_{\sigma(2)},\cdots,e_{\sigma(n)})=1$ or -1
  respectively.
  \item
  Let $n>1$. Prove that for $\sigma,\tau \in S_n,$ $\sigma\tau$ and
  $\tau\sigma$ have same cycle pattern.
  \item

  Let $n\geq 3.$ Write the cycle (123) as product of two $n$-cycles in
  $S_n$.
  \item
  Let $n\geq 1.$ Prove that : \\  (i) The map, \[\begin{array}{rcl}
    \theta_{n,n+1}:S_n & \rightarrow & S_{n+1} \\
    \sigma & \mapsto & \widehat{\sigma} \
  \end{array}\] where $\widehat{\sigma}(i)=\sigma(i)$ foe all $1\leq i\leq
  n,$ and $\widehat{\sigma}(n+1)=n+1$ is a monomorphism of groups which
  maps $A_n$ to $A_{n+1}$. \\  (ii) Identity the group $S_n$ to the
  subgroup $\theta_{n,n+1}(S_n)$ of $S_{n+1}$. Prove that
  $G=\bigcup_{n\geq 1}S_n$ is a group and $A=\bigcup_{n\geq 1}A_n$ is a
  normal subgroup of $G$ which is simple.

   \end{enumerate}

\chapter{Sylow Subgroups}
 \footnote{\it contents group10.tex }
 We have proved that order of any
subgroup of a finite group divides the order of the group. It is natural
to ask if the converse is true. The Sylow theorems are in consequence to
this enquiry. We shall also note that the converse is not true in general.
\bd\label{d81} Let $G$ be a finite group of order $n$. Let for a prime
$p$, $n = p^r\cdot m,$ where $(p,m) = 1.$ Then any subgroup of order $p^r$
is called a Sylow $p$-subgroup of $G.$ \ed We shall prove that any finite
group contains a Sylow $p$-subgroup for every prime $p$ dividing the order
of the group. To prove this we shall prove a
 more general result.
\bt\label{t82} Let $G$ be a finite group of order $n$. Let for $r
> 0$, $p^r$ divides $n$. Then $G$ contains a subgroup of order
$p^r.$ Further, the number of subgroups of order $p^r$ is congruent to 1
mod $p$. \et \pf We shall prove this result in steps.\\ {\bf Step I} : Let
$\mathcal{S}$ be the set of all subsets of $G$ with $p^r$ elements. We
know that for any subset $A$ of $G$, $\abs{A} = \abs{xA}$ for any $x\in
G.$ Thus for any $x\in G,$ \bea
          \hat{x} : \mathcal{S} & \to & \mathcal{S}\\
                         A & \mapsto & xA
\eea is a transformation on $\mathcal{S}$. (Use : $x^{-1}(xA) = A$). For
any $x,y$ in $G,$ and $A \in \mathcal{S},$ $$ (\hat{x}\hat{y})(A) =
\hat{(xy)}(A)$$ Thus $G$ is a transformation group of {\it S} with respect
to this action. From the lemma \ref{l510}, for any $A \in \mathcal{S}$,
$$\abs{\orb(A)} = [G:\stab(A)]$$ Let $H = \stab(A).$ Then $HA = A.$ Hence
if $Ha_1,Ha_2,\ldots ,Ha_t$  are all distinct right cosets of $H$ with
respect to elements of $A$, we have $$ A = Ha_1\cup Ha_2\cup \ldots \cup
Ha_t.$$ and this is a disjoint union. This gives \ba{rl}
            & \abs{A} = t\abs{H}\\
\Rightarrow & \abs{H} \mbox { divides } \abs{A}\\ \ea Hence, $\abs{H} =
p^k$ for some $k \leq r$.   \\  {\bf Step II} : If for $A\in \mathcal{S}$,
$\stab(A) = H$ has $p^r$ elements, then $\orb(A)$ has exactly one subgroup
of order $p^r$. \\  We have $HA = A$, and $\abs{H} = \abs{A} = p^r$. Hence,
for any $a\in A$ \ba{rl}
                  &  Ha = A\\
\Rightarrow       & a^{-1}Ha = a^{-1}A \in \orb(A) \ea Note that
$a^{-1}Ha$ is a subgroup of $G$ with $\abs{a^{-1}Ha} = \abs{H} = p^r.$ Put
$K = a^{-1}Ha.$ Then $K\in \orb(A)$, and hence $\orb(A) = \orb(K).$ As any
element in the $\orb(K)$ is of the form $gK (g \in G),$ it
 is either equal to $K$ or is not a subgroup  of $G$ (Use: For any
 subgroup $K$ of $G$, and $x\in G, xK$ is a subgroup of $G$ if and only if $x\in
 K.$) Therefore the assertion holds.  \\
{\bf Step III}: If $\stab(A) = H$ is a subgroup of $G$ with $p^k
(k < r)$
 elements, then $\orb(A)$ has no subgroup.\\
Let us assume the contrary i.e., let $\orb(A)$ has a subgroup $K.$ Then $K
= xA$ is a subgroup for some $x\in G.$ We have \ba{rl}
                      & K = KK \\
 \Rightarrow           & xA =
KxA\\ \Rightarrow           & A = x^{-1}KxA\\ \Rightarrow & x^{-1}Kx
\subseteq \stab(A) = H\\ \Rightarrow           & \abs{x^{-1}Kx} = \abs{K}
\leq \abs{H}\\ \Rightarrow           & p^r \leq \abs{H} (\mbox{ since
}\abs{K} = p^r) \ea This contradicts our assumption on $H$. Hence
$\orb(A)$ has no subgroup.\\ {\bf Step IV :} There is atleast one $A\in
\mathcal{S}$, such that $\stab(A)$ has $p^r$ elements.\\ \hspace*{.2in}
First of all, note that the number of elements in $\mathcal{S}$ depends on
the
  order of the group and not on the group. In fact
$$ \abs{\mathcal{S}} = \left(\begin{array}{c}
          n\\
          p^r
     \end{array}\right) = N.$$
We have seen that for any $A\in \mathcal{S}$, $\abs{\stab(A)} = p^k (k
\leq r),$ and $\abs{\orb(A)} = [G:\stab(A)].$  Hence if $\abs{\stab(A)} =
p^r$, then $\orb(A)$ is minimal. Let $M_G$ denotes the number of distinct
minimal orbits in $\mathcal{S}$. As $\abs{\orb(A)} = [G:\stab(A)],$ if
$\abs{\stab(A)} = p^k, k
< r$, then $\abs{\orb(A)}$ is a multiple of $pl$, where $l =
\frac{n}{p^r},$ and if $\abs{\stab A} = p^r$, then $\abs{\orb(A)} = l$.
Therefore
\be\label{e81}
 N = plz + M_Gl
\ee If $G = C_n$, a cyclic group of order $n$, then there is
exactly one subgroup
 of order $p^r$ in $G$ (Theorem\ref{t28}). Hence, in this case, there is
 exactly one minimal orbit in $\mathcal{S}$ i.e., $M_{C_n} = 1$. Therefore we get
\be\label{e82}
 N =  plz_1 + l
\ee Now, from (\ref{e81}) and (\ref{e82}), we get \ba{rl}
            &    plz + M_G l =   plz_1 + l\\
\Rightarrow &  p z + M_G  =  p z_1 + 1 \\ \Rightarrow & M_G \equiv 1 (mod
p) \ea Thus $G$ has at least one subgroup of order $p^r$ and the number of
 such subgroups is congruent to $1 (\mbox{mod }p).$      $\blacksquare$
\bt\label{t83} {\bf (Sylow's first theorem)} Let G be a finite
group of order $n$, and $p$, a prime. Let $p^m, m\geq 1,$ be the highest
power of $p$ dividing $n$. Then there exists a subgroup of order $p^m$ in
$G.$
\et \pf This is a special case of theorem \ref{t82}.        $\blacksquare$
\bt\label{t84} Let $G$ be a finite abelian group and $p$, a prime. Let
 $H = \{ x\in G \st \circ(x) = p^{\alpha}, \alpha \geq 0 \}$. Then $H$ is
 the only Sylow $p$-subgroup of $G.$ \et \pf Clearly, $H \neq \emptyset$, as $e\in
H.$ Let $x,y \in H$, where $\circ(x) = p^{\alpha}, \circ(y) = p^{\beta}$.
Then \bea
   (xy^{-1})^{p^{\alpha + \beta}} & = & (x)^{p^{\alpha +\beta}}
                                            (y^{p^{\alpha +\beta}})^{-1}\\
                                    & = &
                                    (x^{p^{\alpha}})^{p^{\beta}}\cdot
                                    ((y^{-1})^{p^{\beta}})^{p^{\alpha}} \\
                                    & = & e
\eea Hence $xy^{-1} \in H$ for any $x,y\in H.$  Thus $H$ is a
subgroup of $G.$ Further, by Cauchy's theorem, $ \circ(H) = p^m.$ For any
Sylow $p$-subgroup $S$ of $G$, order of any element of $S$ is a power of
$p$ (Use : Order of an element divides the order of the  group ). Hence $S
\subseteq H$. Now, by definition of Sylow $p$-subgroup, $S = H.
\blacksquare$ \\  \hspace*{.2in} For any subgroup $H$ of a group $G$,
 and $x\in G$, the conjugate
of H with respect to $x$, i.e., $x^{-1}Hx$, is a subgroup of $G$. Further
$\abs{H} = \abs{x^{-1}Hx}.$ Hence if $H$ is a Sylow $p$-subgroup of $G$,
so is its conjugate $x^{-1}Hx.$ We shall, now, show that any two Sylow
$p$-subgroups of a finite group are conjugate. To prove this we introduce
:
\bd\label{d85} Let $A,B$ be two subgroups of a group $G$. For any element
$x\in G$, the set $$ AxB = \{ axb \st a\in A, b\in B \}$$ is called the
double coset of the pair $(A,B)$ with respect to $x.$ \ed
\bl\label{l86}
Let $A,B$ be two subgroups of a group $G$. Then any two double
cosets of the pair $(A,B)$ in $G$ are either identical or
disjoint.
\el
\pf Take $x,y\in G.$ If $$ (AxB) \cap (AyB) \neq \emptyset ,$$ then there
exist $a_1, a_2\in A$ and $b_1 ,b_2\in B$ such that \ba{rl}
                   & a_1xb_1 = a_2yb_2\\
  \Rightarrow      & Aa_1xb_1B = Aa_2yb_2B\\
  \Rightarrow      & AxB = AyB.
\ea Hence the result follows.      $\blacksquare$
\bl\label{l87}
Let $A,B$ be two subgroups of a group $G$. For $x,y\in G$, define $x\sim
y$  if and only if  $y = axb$ for some $a\in A, b\in B$. Then  $\sim$ is
an equivalence relation over $G$ and the equivalence class of an element
$x\in G$ is $AxB.$
\el
\pf The proof is routine and is left as an exercise.       $\blacksquare$
 \br Let for $x\in G$,
$[x]$ denote the equivalence class of $x$ with respect to $\sim$. If $\{
[x_i]\}_{i\in I}$ be the set of all distinct equivalence classes of $G$,
then $$ G = \cup_{i\in I} Ax_i B  \qquad (\mbox{ Use }: [x_i] = Ax_i B)$$
a disjoint union. This is called the double coset decomposition of $G$
with respect to $(A,B).$ \er \bt\label{t88} {\bf ( Sylow's second
theorem)} Let $G$ be a finite group. Then any two Sylow $p$-subgroups of
$G$ are conjugate. Thus the number of Sylow $p$-subgroups of $G$ is equal
to the number of conjugates of any Sylow $p$-subgroup of $G$.
\et \pf Let $H,K$ be two Sylow $p$-subgroups of $G$, and let
\be\label{e83} G = Hx_1K \cup Hx_2K \cup \ldots \cup Hx_tK
\ee be a double coset decomposition of $G$. Assume $\abs{H} =
\abs{K} = p^m.$ From the equation (\ref{e83}), we have $$\abs{G} =
\abs{Hx_1K}+ \abs{Hx_2K}+ \ldots +\abs{Hx_tK}.$$ Note that for any $x\in
G$,
\bea \abs{HxK} & = & \abs{(x^{-1}Hx)K}\\
          & = & \frac{\abs{(x^{-1}Hx)} \abs{K}}{\abs{(x^{-1}Hx)\cap K}}
                                       \mbox{  (Theorem \ref{t210})}\\
          & = &  \frac{\abs{H} \abs{K}}{\abs{(x^{-1}Hx)\cap K}}\\
          & = &  \frac{p^{2m}}{\abs{(x^{-1}Hx)\cap K}}
\eea Hence \be\label{e84}
 \abs{G} = p^{2m} \sum_{i=1}^t \frac{1}{d_i}
\ee where $d_i = \abs{(x^{-1}_{i}Hx_{i})\cap K} \leq \abs{K} = p^m.$ As
$(x^{-1}_{i}Hx_{i})\cap K$ is a subgroup of $K$, $d_i$ divides $p^m$. Thus
$d_i = p^{\alpha_i}, \alpha_i \geq 0,$ for all $i$. If $\alpha_i < m$ for
all $i$, then equation \ref{e84} gives that $p^{m+1}$ divides $\abs{G}.$
This contradicts that $p^m$ is the highest power of $p$ dividing the order
of $G$ (Use: Sylow $p$-subgroup of $G$ has order $p^m$). Hence $d_{i_0} =
p^m$ for some $i_0$. Thus $$ x_{i_0}^{-1}Hx_{i_0} = K
             (\mbox{Use} : \abs{x_{i_0}^{-1}Hx_{i_0}} = \abs{K} = p^m)$$
Consequently, $K$ is a conjugate of $H.$ The rest of the statement is
immediate.       $\blacksquare$
\bco\label{c81} A finite abelian group has unique Sylow $p$-subgroup  for
any prime $p$.
\eco \pf Clear from theorem \ref {t88}     $\blacksquare$
 \br The corollary also follows
from theorem \ref{t84}
\er \bco\label{c82} Let $P$ be a Sylow $p$-subgroup of a finite group $G$.
Then $P$ is unique
 if and only if $P$ is normal in $G.$ Thus $P$ is unique Sylow
 $p$-subgroup of $N_{G}(P).$
 \eco
\pf We know that the number of distinct conjugates of $P$ in $G$ is
$[G:N_G(P)].$ Hence the first part of the result follows from the theorem.
For the other part, note that $P\subset N_G(P) \subset G.$ By the
Lagrange's theorem (Theorem \ref{co22}), $\circ(P)$ divides $\circ(N_G(P))$
and $\circ(N_G(P))$
divides $\circ(G)$. Therefore it is clear that $P$ is a Sylow $p$-subgroup
of $N_G(P).$ Now, as $P\lhd N_G(P)$, the result is immediate from the
theorem. $\blacksquare$

\bt\label{t89} {\bf( Sylow's third theorem )} Let $G$ be a finite
group, and $H$ a Sylow $p$-subgroup of $G$. Let $r$ be the number of
distinct conjugates of $H$ in $G$. Then $r \equiv 1$ (mod $p$) and $r$
divides $\circ(G).$ \et \pf By theorem \ref{t88}, $r$ is equal to the
number of Sylow $p$-subgroups of $G$. It is already proved in theorem
\ref{t82} that $r \equiv 1$ (mod $p$). Now, as \bea
                \abs{G} & = & \abs{N_G(P)}\cdot \abs{[G:N_G(P)]}\\
                        & = & \abs{N_G(P)}\cdot r,
\eea  $r$ divides $\circ(G).$ Hence the result follows.
$\blacksquare$ \bt\label{t10.16t} Let $G,H$ be two finite groups, and let
$\varphi : G \rightarrow H$ be an epimorphism. Then for any Sylow
$p$-subgroup $S$ of $G$, $\varphi(S)$ is a Sylow $p$-subgroup $H$.
Conversely, if $A$ is a Sylow $p$-subgroup of $H$, then there exists a
Sylow $p$-subgroup $B$ of $G$ such that $\varphi(B)=A.$\et \pf Put
$\varphi(S)=T.$ By lemma \ref{l31}, $\circ(T)$ divides $\circ(S)$. Hence
$\circ(T)=p^{\beta}$ for some $\beta \geq 0.$ Let $K$ be the kernel of
$\varphi$. For the subgroup $KS=SK$ of $G$, let \be\label{e10.3e} G =
\bigcup_{i=1}^{t} x_i SK \ee be a left coset decomposition of $G$ with
respect to $SK$. From the equation (\ref{e10.3e}), we get \be\label{e10.4e}
H = \varphi(G) = \bigcup_{i=1}^{t} \varphi(x_i)T \ee Note that
$$\begin{array}{crcl}
   & \varphi(x_i)T & = & \varphi(x_j)T \\
  \Leftrightarrow & \varphi(x^{-1}_{j}x_i) \in T & = & \varphi(S) \\
  \Leftrightarrow & x_{j}^{-1}x_{i} \in \varphi^{-1}(\varphi(S)) & = & SK \\
  \Leftrightarrow & x_i SK & = & x_j SK
\end{array}$$ As $x_i SK \neq x_j SK$ for $i \neq j$, $\varphi(x_i)T \neq \varphi(x_j)T$
 for all $i\neq j$. Therefore (\ref{e10.4e}) is a coset decomposition of $H$
 with respect to $T$, and $[G:SK] = [H:T].$ We have $$[G:S] = [G:SK][SK:S]
 = t[SK:S]$$ As $\circ(S)$ is the highest power of $p$ dividing $\circ(G)$,
 $[G:S]$ is relatively prime to $p$. Therefore $(p,t) = 1$. Consequently
 $\circ(T) = p^{\beta}$ is the highest power of $p$ dividing $\circ(H).$
 Thus $T$ is Sylow $p$-subgroup of $H$. For the other part, let
 $C=\varphi^{-1}(A).$ Then $C$ is a subgroup of $G$. Let
 \be\label{ee10.5ee} G = \bigcup_{i=1}^{m} y_iC \ee be a left coset
 decomposition of $G$ with respect to $C$. From the equation
 (\ref{ee10.5ee}), we get
 \be\label{e10.5e}
 H = \varphi(G) = \bigcup_{i=1}^{m}\varphi(y_i)A.\ee Note that $$\begin{array}{cl}
    & \varphi(y_i)A = \varphi(y_j)A \\
   \Leftrightarrow & \varphi(y_{j}^{-1}y_i) \in A \\
   \Leftrightarrow & y_{j}^{-1}y_i \in \varphi^{-1}(A) = C \\
   \Leftrightarrow & y_j C = y_i C \
 \end{array}$$ As $y_i C \neq y_j C$ for all $i\neq j$, (\ref{e10.5e}) is a
 coset decomposition of $H$ with respect to $A$. Thus $[G:C] = [H:A]=m$.
 Since $A$ is a Sylow $p$-subgroup of $H$, and $[H:A]=m, (p.m)=1$. Now, as
 $m=[G:C]=\circ(G)/\circ(C)$, and $(p,m)=1$ highest power of $p$ dividing
 $\circ(G)$ is the same as highest power of $p$ dividing $\circ(C)$. Thus
 if $B$ is a Sylow $p$-subgroup of $C$, it is also a Sylow $p$-subgroup of
 $G$, and hence by the direct part $\varphi(B)$ is Sylow $p$-subgroup of
 $H$. As $\varphi(C)=A$, $\varphi(B) \subset A$, and hence $\varphi(B) =
 A.$ Thus the proof is complete.     $\blacksquare$
 \bl\label{l10.17l} Let $S$ be a Sylow $p$-subgroup of a finite group
 $G$, and let $H$ be a subgroup of $G$ containing $N_G(S)$. Then $N_G(H)=H.$
  \el \pf We have $H\subset N_G(H).$
 Therefore, if $H \neq N_G(H),$ choose $x \in N_G(H), x\not\in H.$ As
 $H\supset N_G(S), x \not\in H, xSx^{-1} \neq S.$ However $xHx^{-1}=H.$ Thus $S$ and
 $xSx^{-1}$ are two Sylow $p$-subgroups of $H.$ By theorem \ref{t88}, there exists
 an element $h\in H$ such that $$\begin{array}{crl}
    & h^{-1}Sh & =x^{-1}Sx \\
   \Rightarrow & xh^{-1}S(xh^{-1})^{-1} & =S \\
   \Rightarrow & xh^{-1} & \in N_G(S) \subset H \\
   \Rightarrow & x & \in H \
 \end{array}$$ A contradiction to our choice of $x$. Consequently $N_G(H)=H.$
 $\blacksquare$      \bt\label{t10.14t} Let $G$ be a finite group and $p$, a prime.
 Then any subgroup $K$ of $G$ with $\circ(K)=p^r, r\geq 1,$ is contained
 in a Sylow $p$-subgroup.\et  \pf Let $\mathcal{S}=\{H_1,H_2,\cdots,H_t\}$
  be the set of all Sylow $p$-subgroup of $G$. By theorem \ref{t82},
  $t\equiv1(mod\, p).$ As conjugate of a Sylow $p$-subgroup is a Sylow
  $p$-subgroup, $K$ operates by conjugation on $\mathcal{S}$. By lemma
  \ref{l510}, any orbit of $\mathcal{S}$ with respect to $K$ action has
  $p^{\alpha}, \alpha\geq 0,$ elements. Hence as $t\equiv 1(mod\,p)$,
   at-least one orbit shall have
  1-element. We can assume without any loss of generality that $\abs{\orb
  (H_1)}=1.$ Then by theorem \ref{t82}, $\stab (H_1) = G.$ Thus
  $x^{-1}H_1x = H_1$ for all $x \in K$. Hence $H_1K=KH_1$ is a subgroup of
  $G$ containing $K$ as well as $H_1$. By theorem \ref{t3.10},
  $$\circ(H_1K)= \frac{\circ(H_1)\circ(K)}{\circ(H_1\bigcap K)}.$$
  Therefore $\circ(H_1K) = p^r$ for some $r\geq 1.$ As $H_1$ is Sylow
  $p$-subgroup of $G$ and $H_1 \subset H_1K, H_1=H_1K.$ Thus $K\subset
  H_1.$         $\blacksquare$  \bco\label{co10.15co} Let $G$ be a finite
  group and $p$, a prime. Let for a subgroup $K$ of $G$, $p\mid\circ(K).$
  Then there exists a Sylow $p$-subgroup $P$ of $G$ such that $P\cap K$ is
  a Sylow $p$-subgroup of $K$.\eco \pf Let $B$ be a Sylow $p$-subgroup of
  $K$. Then $\circ(B)=p^r$ for some $r\geq 1.$ By the theorem $B\subset P$
  for a Sylow $p$-subgroup $P$ of $G$. Thus $P\cap K \supset B,$ and
  $\circ(P\cap K)$ is a power of $p$. Now, as $B$ is a Sylow $p$-subgroup
  of $K$, $B=P\cap K.$ Hence the result follows.      $\blacksquare$
  \bco\label{co10.16co} Let $G$ be a finite group and $p$, a prime. Let
  $H$ be a normal subgroup of $G$ such that $p$ divides $\circ(H)$. Then
  for any Sylow $p$-subgroup $P$ of $G$, $H\cap P$ is a Sylow $p$-subgroup
  of $H$.\eco \pf By corollary \ref{co10.15co} and theorem \ref{t88},
  there exists an element $x \in G$ such that $x^{-1}Px \cap H$ is a Sylow
  $p$-subgroup of $H$. Note that $$x^{-1}Px \cap H = x^{-1}(P\cap H)x
  \mbox{    since }H\lhd G.$$ Thus for $P_1 = x^{-1}Px,$  $$P\cap H =
  x(P_1 \cap H)x^{-1}$$ i.e., $P\cap H$ is a conjugate of $P_1\cap H,$ a
  Sylow $p$-subgroup $H$. Hence $P\cap H$ is a Sylow $p$-subgroup of $H$.
      $\blacksquare$   \\    \\

  {\bf APPLICATIONS}
\bt\label{t810} Let $A$ be a finite abelian group. Then $A$ is the direct
product of it's
 Sylow subgroups.
\et \pf Let $\abs{A}=l = p_1^{\alpha_1}\cdot p_2^{\alpha_2}\ldots
p_m^{\alpha_m},$ where $p_i; i = 1,2,\ldots ,m$ , are distinct primes and
$\alpha_i
> 0.$ Let $H_i$ be the Sylow $p_i$-subgroup of $A$ for $i = 1,2,\ldots
,m.$ Then consider the map \bea
      \Phi : H_1\times H_2\times \ldots \times H_m & \to & A\\
        \underline{a}=  (a_1,a_2, \ldots ,a_m) & \mapsto & a_1a_2 \ldots a_m.
\eea It is easy to check that $\Phi$ is a homomorphism. If
\ba{rl}
               & \Phi(\underline{a}) = e\\
\Rightarrow    & a_1a_2 \ldots a_m = e\\ \Rightarrow    & a_1
\ldots a_{i-1}a_{i+1}\ldots a_m = a_i^{-1}
                           \mbox{ for all } 1\leq i \leq m.
\ea Now, let $\circ(a_j) = (p_j)^{n_j}$ for $j = 1,2,\ldots ,m.$ Then, if
\\
$n = {p_1}^{n_1}\ldots {p_{i-1}}^{n_{i-1}}{p_{i+1}}^{n_{i+1}}
                  \ldots {p_m}^{n_m}$, we have
\be\label{e85}
 {(a_i^{-1})}^n = {(a_1 \ldots a_{i-1}a_{i+1}\ldots a_m)}^n = e
\ee Also \be\label{e86}
 \circ(a_i^{-1}) = \circ(a_i) = {p_i}^{n_i}
\ee Hence, as $(p_i,n) = 1$, by theorem \ref{t25}, we get
$a_i = e.$  Since $i$ is arbitrary, $\underline{a}=id.$ Therefore
 $\Phi$ is a monomorphism.\\
Finally, as $\abs{H_1\times\cdots\times H_m} = l = \abs{A}$,
 $\Phi$ is an isomorphism.
$\blacksquare$
\bt\label{t811} Let $G$ be a group of order $pq,$ where $p,q$ are primes,
$p>q,$ $p\not\equiv 1$(mod $q$). Then $G$ is cyclic. \et \pf Let $n_p$ be
the number of Sylow $p$-subgroups of $G$ and $n_q$, the number of Sylow
$q$-subgroups of $G.$ By Sylow's third  theorem, $n_p$ and $n_q$ divide
$\circ(G)$, $n_p \equiv 1$ (mod $p$) and $n_q\equiv 1$ (mod $q$). Thus
$(n_p,p) = 1 = (n_q,q)$ and $n_p$ and $n_q$ divide $pq.$ Now, as $q<p,$
and $n_p \equiv 1(mod\, p), n_p = 1$ and, as $n_q \equiv 1(mod\, q),$ and
$p \not\equiv 1(mod\, q), n_q = 1$. Let $H$ be the Sylow $p$-subgroup of
$G$ and $K$, the Sylow $q$-subgroup of $G$. By the
corollary \ref{c82} of Sylow's second theorem $H$ and $K$ are normal in G.
As $H, K$ are subgroups of prime orders these are cyclic. Moreover,
 using Lagrange's theorem and that $\circ(H) = p, \circ(K) = q$
and $(p,q) = 1$, we get $H\cap K =\{e\}$. For any $x\in H, y\in K,$
 we have $xyx^{-1}y^{-1} \in H\cap
K = \{e\}$, since $H$ and $K$ are normal in $G$. Therefore $ xy = yx$ for
all $x\in H, y\in K.$ As $H$ and $K$ are normal subgroups of $G$, $HK =
KH$ is a subgroup of $G$ with $$ \abs{HK} = \frac{\abs{H} \abs{K}}{\abs
{H\cap K}} = pq
                     \mbox { (Theorem \ref{t210})}$$
Hence $HK = G.$  Therefore by corollary \ref{co6.3co}, $G= H\times K$
i.e.,
$G$ is internal direct product of $H$ and $K$. As
$H$ and $K$ are cyclic groups of orders $p$ and $q$ respectively, and
$(p,q)=1,$ it follows from theorem \ref{t43} that $G$ is cyclic of order
$pq$.      $\blacksquare$
\br The conditions on the primes imposed in the theorem are
essential. e.g. take $G = S_3$. Then $\circ(G) = 6 = 2\times 3.$ Here for
$p = 3, q = 2, p > q, p \equiv 1$ (mod $q$), and the group $G$ is
non-abelian. Thus $p \not\equiv 1$ (mod $q$) is essential. \er We shall,
now, identify the Sylow subgroups of some familiar groups.  \\  1.  Let
$C_n = gp\{a\}$ be a cyclic group of order $n$, and let $p$ be a prime. If
$p^{\alpha}\, (\alpha \leq 1)$ divides $n$, and $m=n/p^{\alpha}$, then
$C_m = gp\{a^m\}$ has order $p^{\alpha}.$ Thus $C_m$ is the Sylow
$p$-subgroup of $C_n.$ \\ 2.  Let $n>1$ be an odd integer. Consider the
Dihedral group (Example \ref{th11}):  $$D_{n} = \{Id,
\tau,\sigma,\sigma^2,\cdots,\sigma^{n-1},\tau\sigma,\cdots,\tau\sigma^{n-1}
\st \tau^2 = \sigma^n = Id, \sigma\tau = \tau\sigma^{n-1}\}$$ As $n$ is
odd, 2 is the highest power of 2 dividing $\circ(D_n)=2n.$ Thus cyclic
subgroups of order 2 are precisely the Sylow 2-subgroups of $D_n$. One can
easily check that $gp\{\tau\sigma^i\}, 0\leq i\leq n-1,$ are all the Sylow
2-subgroups of $D_n.$  \\  3.  Let $\F =\Z_p$ a field with $p$-elements.
The group $Gl_n(\F_p)$ has order $\prod_{k=0}^{n-1}(p^n-p^k)$ (Use:
Exercise 12, Chapter 1). We have $$\circ(Gl_n(\F_p)) =
p^{n(n-1)/2}\prod_{k=0}^{n-1}(p^{n-k}-1)$$ Thus $p^{n(n-1)/2}$ is the
highest power of $p$ dividing $\circ(Gl_n(\F_p)).$ Now, check that
$$\begin{array}{cl}
  {\bf UT_n}(\F_p) = &  \{A=(a_{ij}): n\times n-\mbox{matrix over } \F_p \st
a_{ij}=0  \\
   & \mbox{ for all }i>j, a_{ii}=1,\mbox{ for all }1\leq
i\leq n\}
\end{array} $$  is a subgroup of $Sl_n(\F_p)$ of order $p^{n(n-1)/2}$
(Use : For all $a_{ij}$, $j>i$, we can choose any element of $\F_p$. Thus
for all $1\leq i<j\leq n,$ $a_{ij}$ has $p$-possibilities and the number of
elements in the set $\{1 \leq i<j\leq n\}$ is $n(n-1)/2.$). \\  \\ We
shall, now, consider a case complimentary to that in theorem \ref{t811} of
particular interest. \bt\label{t10.18} Let $G$ be a group of order $2p$,
where $p>2$ is a prime. Then $G$ is either the Dihedral group $D_p$ (upto
isomorphism) or is cyclic.\et \pf As in theorem \ref{t811}, $G$ has
exactly one Sylow $p$-subgroup $H$ which is cyclic of order $p$. By
corollary \ref{c82}, $H$ is normal in $G$. Let $H=gp\{y\}$. Further, by
theorem \ref{t82}, $G$ has an element $x$ of order 2. As $H\lhd G,$ $$xyx
= x^{-1}yx = y^i \mbox{    for some }1\leq i\leq p-1$$ If $i=p-1,$ then
$yx=xy^{p-1}.$ In this case it is easy to check that $G$ is isomorphic to
the Dihedral group $D_p$. Now, let $1\leq i<p-1.$ \\  Then  $$(xy)^2 =
y^{i+1} \neq e$$ and $$(xy)^p = (y^{i+1})^{\frac{p-1}{2}}y^{i}x \neq
e$$ since $H\cap gp\{x\}=e.$ Thus it follows that $\circ(xy)=2p$, and
$G=gp\{xy\}.$ Hence the result is proved.      $\blacksquare$
\bt\label{t10.19t} Let $G$ be a finite group with $\circ(G)=pqr,$ where
$p>q>r$ are primes. Then \\ (i)$G$ has only one Sylow $p$-subgroup.  \\
(ii) Let $p\not\equiv 1\, mod \, q.$ Then $G$ has only one Sylow
$q$-subgroup. If $H$ is the Sylow $p$-subgroup of $G$ and $K$ is the Sylow $q$-subgroup of $G$
 then $HK\lhd G$ and
is a cyclic group of order $pq$.\et \pf (i) Let $n_i, i=p,q,r,$ be the
number of Sylow $i$-subgroups in $G.$ By theorem \ref{t82}, $n_p\equiv
1(mod\, p).$ Further, by theorems \ref{t88} and \ref{t89}, $n_p$ divides
$\circ(G)=pqr.$ Hence $n_p=1$ or $qr$ (Use: $p>q>r>1$). If $n_p=1,$
then (i) holds. Assume $n_p=qr.$ In this case we shall show that either
$n_q=1$ or $n_r=1.$ If not, then as above we get that $n_q=p$ or $pr$, and
$n_r=p$ or $q$ or $pq.$ Thus $n_q\geq p,$ and $n_r\geq q.$ Note that any
Sylow subgroup of $G$ has prime order, hence is cyclic. Moreover, any two
distinct Sylow subgroups of $G$,being of prime orders, intersect in the
identity subgroup. Now, counting the elements in all the Sylow subgroups
of $G$, we get $$\begin{array}{crl}
   & (p-1)qr+(q-1)p+(r-1)q+1 & \leq pqr \\
  \Rightarrow & pqr+(q-1)(p-1) & \leq pqr
\end{array}$$ This, however, is not true. Therefore either $n_q=1$ or
$n_r=1.$ Assume $n_q=1,$  and let $K$ be the Sylow $q$-subgroup
of $G$. Let $H=gp\{a\}$ be a Sylow $p$-subgroup of $G$ and let
$K=gp\{x\}.$ \\ Then \begin{eqnarray}
  axa^{-1} & = & x^i \;\;\; \mbox{ for some }1\leq i< q \nonumber \\
 \Rightarrow  a^p x a^{-p} & = & x^{ip}   \nonumber \\
 i.e.,  x & = & x^{ip} \;\;\;  \mbox{ since } \circ(a)=p \nonumber \\
 \Rightarrow  i^p & \equiv & 1(mod\, q) \label{eq10.3}
   \end{eqnarray} We also have
   \be\label{eq10.11}
    i^{q-1}  \equiv  1(mod\, q)
    \ee Moreover, it is clear that $(p,q-1)=1.$ Hence from the equations (\ref{eq10.3}) and
    (\ref{eq10.11}), we conclude that $i\equiv 1(mod\,q).$ Therefore
    $axa^{-1}=x$ i.e., $ax=xa.$ Thus $x \in N_G(H).$ By theorem \ref{t88},
    we have $$\begin{array}{cl}
       & n_p = qr = \frac{\circ(G)}{\circ(N_G(H))} = \frac{pqr}{\circ(N_G(H))} \\
      \Rightarrow & \circ(N_G(H)) = p \\
\end{array}$$ As $H \subset N_G(H)$, and $\circ(H)=p, H=N_G(H).$ Thus
$x\not\in N_G(H)$, since $K\cap H = \{e\}.$ Hence $n_q\neq 1.$ Similarly we can see
that $ n_r \neq 1.$ Consequently $n_p= 1.$
 \\ (ii) As $H$ is unique Sylow $p$-subgroup of $G$, $H\lhd G$ (Corollary \ref{c82}).
Hence for a Sylow $q$ subgroup $K$ of $G$, $HK=KH$ is a subgroup of order $pq$ in $G$.
 Now, $[G:HK]=r$ is the smallest prime dividing the order of $G$.
 Thus by corollary \ref{co7.32co} $HK\lhd G$ and by theorem \ref{t811} $HK$ is cyclic.
  By theorem \ref{t28}, $K$ is unique subgroup of order $q$ in $HK$. Hence $K\lhd G.$
$\blacksquare$
\bt\label{t10.20t} Let $p,q$ be two, not necessarily distinct, primes.
Then any group $G$ of order $p^nq(n\geq 1)$ is not simple i.e., $G$
contains a proper normal subgroup.\et \pf If $p=q$, then
$\circ(G)=p^{n+1}.$ By theorem \ref{t82}, $G$ contains a subgroup of order
$p$, and hence by corollary \ref{co8.24co}, $G$ contains a proper normal
subgroup. Now, let $p\neq q$, and let $r$ be the number of Sylow
$p$-subgroups of $G$. By theorems \ref{t88} and \ref{t89}, $r\equiv
1(mod\,p)$ and $r$ divides $\circ(G)=p^nq.$ Therefore $r=1$ or $q$. If
$r=1$, then the Sylow $p$-subgroup (say) $P$ of $G$ is a normal subgroup
of order $p^n$ in $G$ (Corollary \ref{c82}). Hence the result holds in
this case. Next, let $r=q$. If any two distinct Sylow $p$-subgroups of $G$ intersect
in the identity subgroup, then the number of elements of order
$p^{\alpha}(\alpha\geq 1)$ in $G$ is $(p^n-1)q+1.$ Hence, $G$ has only one
Sylow $q$-subgroup which is a proper normal subgroup of $G$. Thus, to
complete the proof, we assume that not any two distinct Sylow $p$-subgroups of $G$
intersect in identity group. Let for two distinct Sylow $p$-subgroups $P_1,P_2$
of $G$, $\abs{P_1\cap P_2}$ is maximal. Put $A=P_1\cap P_2.$ Then $A$ is a
proper subgroup of $P_i, i=1,2.$ By corollary \ref{co8.24co},
$N_i=N_{{P_i}}(A)\supsetneqq A,$ for $i=1,2$. Let $H=gp\{N_1\cup N_2\}.$
Note that $H$ is not a $p$-group, since otherwise there exists a Sylow
$p$-subgroup $P'$ of $G$, $P' \supset H$ (Theorem \ref{t810}), and $P'\cap
P_1\supset N_1 \supsetneqq A.$ This contradicts our assumption on
maximality of $\abs{P_1\cap P_2}$. Thus $\circ(H)=p^bq$ where $b\leq n.$
As $N_i \subset N_{G}(A)$ for $i=1,2$, $H\subset N_{G}(A).$ Let $Q$ be a
Sylow $q$-subgroup of $H$, then $\abs{P_1Q}=\abs{P_1}\abs{Q}=p^nq.$ Thus
$P_1Q=G.$ Therefore any element of $G$ can be written as $ya$ where $y \in
P_1, a \in Q.$  We have $$\begin{array}{ccl}
  yaA(ya)^{-1} & = & yaAa^{-1}y^{-1} \\
   & = & yAy^{-1}\subset P_1,
\end{array}$$ since $a \in Q\subset H\subset N_{G}(A).$ Hence all
conjugates of $A$ in $G$ are contained in $P_1$ and are of the for
$yAy^{-1}, y\in P_1.$ Let $$K=gp\{yAy^{-1}\st y \in P_1\}.$$ Then clearly
$K\lhd G$ and $K\subset P_1.$ Thus $K$ is a proper normal subgroup of $G.$
Hence the proof is complete.           $\blacksquare$
\bt\label{t10.21} Let $G$ be a finite group and $p$, a prime dividing
$\circ(G)$. If a Sylow
$p$-subgroup $P$ of $G$ is contained in $Z(G)$, then $G$ isomorphic to the
direct product $P\times K$ for a subgroup $K$ of $G$.\et \pf Let
$\circ(G)=p^{\alpha}m$ where $(p,m)=1,$ and let $$G=Pa_1\cup Pa_2\cup
\cdots \cup Pa_m$$ be a right coset decomposition of $G$ with respect to
$P$. Then, for any \\ $x\in G, xa_i$ can be uniquely written as
$$xa_i=p_{x(i)}a_{i(x)}$$ where $p_{x(i)}\in P$ and $1\leq i(x)\leq m.$ We
have $$G=xG=\bigcup_{i=1}^mPxa_i=\bigcup_{i=1}^mPa_{i(x)}$$ Hence,
$$\left(\begin{array}{cccc}
  1 & 2 & \cdots & m \\
  1(x) & 2(x) & \cdots & m(x)
\end{array}\right)$$ is a permutation on the symbols $1,2,\cdots,m.$ \\
For any $x,y\in G,$ \begin{eqnarray} yxa_i & = &
yp_{x(i)}a_{i(x)}\nonumber \\
 & = & p_{x(i)}p_{y(i(x))}a_{i(x)(y)}
\label{line2} \
\end{eqnarray}  Define : $$\begin{array}{crcl}
  \varphi & :\, G & \rightarrow & P \\
   & x & \mapsto & \prod_{i=1}^m p_{x(i)}.
\end{array}$$ Then for any $x,y \in G,$ using \ref{line2}, we get $$\begin{array}{cclc}
  \varphi(yx) & = & \prod_{i=1}^m p_{x(i)} p_{y(i(x))} &  \\
   & = & \prod_{i=1}^mp_{y(i)} \prod_{i=1}^m p_{x(i)} & (\mbox{use : }P\subset Z(G)) \\
   & = & \varphi(y)\varphi(x) &
\end{array}$$ Therefore $\varphi$ is a homomorphism. Clearly, if $p\in P$,
then $p_{p(i)}=p$ for all $i=1,2,\cdots,m.$ Hence $\varphi(p)=p^m.$ As
$(m,p^{\alpha})=1,$ there exists an integer $n$ such that $mn \equiv
1(mod\,p^{\alpha}).$ As $P$ is an abelian group, $$\begin{array}{crcl}
  \theta_n & :\,P & \rightarrow & P \\
   & a & \mapsto & a^n
\end{array}$$ is an endomorphism of $P$ and for $\psi = \theta_n \circ \varphi
$, $\psi(x) = \varphi(x)^n$ for all $x\in G.$ If $p\in P,$ then
$\psi(p)=p^{mn}=p$ since $mn \equiv 1(mod\,p^{\alpha})$ and
$\circ(P)=p^{\alpha}.$ Therefore $\psi : G \rightarrow P$ is a
homomorphism such that $\psi(p)=p.$ \\ Let $K=\ker\psi .$ We shall prove
that $G$ is isomorphic to the direct product $P\times K.$ Define :
$$\begin{array}{rcl}
  \alpha : P\times K & \rightarrow & G \\
  (p,a) & \mapsto & pa
\end{array}$$ For $(p_1,a_1), (p_2,a_2)$ in $P\times K,$ we have $$\begin{array}{ccl}
  \alpha((p_1,a_1)(p_2,a_2)) & = & \alpha(p_1p_2, a_1a_2) \\
   & = & (p_1a_1)(p_2a_2) \mbox{  \;\;\;\; since }P\subset Z(G) \\
   & = & \alpha(p_1,a_1)\alpha(p_2,a_2)
\end{array}$$ Thus $\alpha$ is a homomorphism. If for $(p,a) \in P\times
K,$ $$\begin{array}{crlc}
   & \alpha(p,a) & = e &  \\
  \Rightarrow & pa & = e &  \\
  \Rightarrow & \psi(pa) & = \psi(p) = e &  \\
  \Rightarrow & p & = e & \mbox{ \;\;\; since }\psi(p) = p \\
  \Rightarrow & p & = a = e &
\end{array}$$ Therefore $\alpha$ is a monomorphism. Now, note that $G/K$
is isomorphic to $P$ (Corollary \ref{c43a}). Therefore
$$\begin{array}{ccl}
   & \circ(G/K) & = \circ(G)/\circ(K) = \circ(P) \\
  \Rightarrow & \circ(G) & = \circ(P)\circ(K) \\
   &  & =\circ(P\times K).
\end{array}$$ Hence $\alpha$ is an isomorphism.     $\blacksquare$
   \vspace{.7in}\\
 {\bf EXERCISES}
\begin{enumerate}
  \item
  Find all Sylow 3-subgroups of $Gl_2(\Z_3).$
  \item
  Prove that every subgroup of order 11 in a group of order 77 is normal.
  \item
  Let $G$ be a finite group with order $n$ and let $p$ be a prime such
  that $p$ divides $n$. If $1\leq n/p<p,$ then show that $G$ contains a
  normal subgroup of order $p$.
  \item
  Let $G$ be a non-abelian group of order 343. Prove that $\circ(G')=7.$
  \item
  Let $H$ be a normal subgroup of a finite group $G$ and $p$, a prime.
  Prove that if $\circ(H)=p^{\alpha}\,(\alpha\geq 1),$ then $H$ is
  contained in all Sylow $p$-subgroups of $G$.
  \item
  Find all Sylow subgroups of $D_{12}.$
  \item
  Write all Sylow subgroups of $S_5$.
  \item
  Let $\F$ be a finite field with $p^{\alpha}$ ($p:$ prime, $\alpha\geq
  1$) elements. Find a Sylow $p$-subgroup in $Sl_n(\F).$
  \item
  Let $K$ be a field and $n>1.$ Prove that
  $T_n(K)=N_{Gl_{n}(K)}(UT_n(K)).$
  \item
  Let $\F$ be finite field with $q=p^m$ ($p$: prime, $m\geq 1$) elements.
  Find the number of Sylow $p$-subgroup of $Gl_n(\F).$
  \item
  Let for a prime $p$, $G$ be a group of order $p^n.$ Prove that $G$
  contains a normal subgroup $H$ of order $p^k$ for each $0\leq k\leq n.$
  \item
  Let for a group $G$, $H$ be a finite normal subgroup of $G$. Prove that
  if $S$ is a Sylow subgroup of $H$, then $G=HN_G(S).$
  \item
  Let $H$ be a normal subgroup of a group $G$. Prove that for any Sylow
  $p$-subgroup $P$ of $G$, $HP/H$ is a Sylow $p$-subgroup of $G/H$.
  \item
  Let $p,q$ be distinct primes such that $p^2\not\equiv 1(mod\,q)$ and
  $q^2\not\equiv 1(mod\,p)$. Prove that any group of order $p^2q^2$ is
  abelian.
  \item
  Prove that the number of Sylow $p$-subgroups in $S_p$ is $(p-2)!.$
  \item
   Prove that the order of the Sylow $p$-subgroup of $S_{p^{m}}$ is
  $p^{(1+p+\cdots+p^{m-1})}.$
  \item
  Let $G$ be a group of order 21. Prove that \\ (i) There exists $x,y$ in
  $G$ such that $\circ(x)=7$, $\circ(y)=3$ and $G=gp\{x,y\}$. \\ (ii) We
  have $y^{-1}xy=x^r$ where $r=1$ or 2 or 4. \\ (iii) There are two
  classes of groups of order 21 upto isomorphism.
  \item
  Let $p\geq 3$ be a prime. Prove that any non-abelian group of order $2p$
  is Dihedral group.
  \item
  Let $p\geq 3$ be a prime. Prove that the Dihedral group $D_{2p}$ has $p$
  Sylow 2-subgroups.
  \item
  Let $G$ be a finite group, and $\circ(G)<60.$ Prove that if $G$ is
  simple then $G$ is cyclic of prime order.
  \item
  Let $G$ be a group such that $\Phi(G)$ is finitely generated. Prove that
  any Sylow $p$-subgroup $P$ of $\Phi(G)$ is normal in $G$.
  \end{enumerate}

\chapter{Finitely Generated And Divisible Abelian Groups}
 \footnote{\it contents group11.tex }
 In this chapter we shall study the structure of finitely generated
 abelian groups and divisible abelian groups. All groups shall be assumed additive
 unless the binary operation is specified.
\bd\label{d91}
An abelian group $F$ is called a finitely generated free abelian group if
there exist $f_i \in F,\, i=1,2,\cdots,n$, $n\geq 0$, such that\\ (i) $F =
gp \{f_1,f_2,\ldots ,f_n \}$\\ and \\  (ii) The set $\{f_1,f_2,\ldots ,f_n
\}$ is linearly independent i.e., for $n\geq 1$ and $a_i; i= 1,2,\ldots, n ,$ in
$\Z$, $a_1f_1+a_2f_2+\ldots +a_nf_n = 0,$ if and only if $a_i = 0$ for all $i =1,2,....,n .$
\ed
\brs
(i) If we take $n=0$ in the definition then $F$ is generated by the null
set $\emptyset$, and hence $F=\{0\}.$ \\ (ii) The set $\{f_1,f_2,\ldots
,f_n \}$ is called a basis of $F$. If $n\geq 1,$ then we shall denote the
fact that $F$ is a free abelian group with basis $\{f_1,f_2,\ldots ,f_n\}$
by writing $$ F = \langle f_1,f_2,\ldots ,f_n \rangle.$$ (ii) The abelian
group (0) is free with basis consisting of the empty set $ \emptyset$.
\ers
\bos
(1) For any non-zero element of a free abelian group $F$, the subgroup $G
= gp\{f\}$ is free.\\ Let $F = \langle f_1,f_2,\ldots ,f_n \rangle$. Then
$f = a_1f_1+a_2f_2+\ldots +a_nf_n,$ for some $a_i\in \Z$. Since $f\neq 0$,
all $a_i$'s can not be zero. Let for $t\in \Z$, \ba{rl}
               & tf = 0\\
\Rightarrow    & ta_1f_1+ta_2f_2+\ldots +ta_nf_n = 0\\ \Rightarrow
& ta_i = 0\;\;\; \forall i = 1,2,\ldots ,n.\\ \Rightarrow    & t =
0 \mbox{ since } a_i\neq 0 \mbox{ for some i} \ea Hence $G$ is a
free abelian group with basis $\{f\}$.\\ (2) Let $F = \langle
f_1,f_2,\ldots ,f_n \rangle$ be a free abelian group. Then any
$x\in F$ can be uniquely expressed as $$ x = a_1f_1+a_2f_2+\ldots
+a_nf_n, where \ a_i\in \Z\ for\ all\ i=1,2,\ldots \ n.$$ Let $$ x
= a_1f_1+a_2f_2+\ldots +a_nf_n = b_1f_1+b_2f_2+\ldots +b_nf_n$$
Then \ba{rl} & (a_1-b_1)f_1+(a_2-b_2)f_2+\ldots +(a_n-b_n)f_n =
0\\ \Rightarrow    & a_1-b_1 = a_2-b_2 = \ldots = a_n-b_n = 0\\
\Rightarrow     & a_i = b_i \;\;\; \forall i \ea (3) Let $F =
\langle f_1,f_2,\ldots ,f_n \rangle$, and $G = \langle
g_1,g_2,\ldots ,g_k \rangle$ be two finitely generated free
abelian groups. Then $F\oplus G$ is a free abelian group with
basis $$\{ (f_1,0),\ldots ,(f_n,0),(0,g_1), \ldots ,(0,g_k) \}.$$
Let $x\in F, y\in G$, then \bea
 x & = & a_1f_1+a_2f_2+\ldots +a_nf_n\\
 y & = & b_1g_1+b_2g_2+\ldots +b_kg_k
\eea
for some $ a_1,a_2,\ldots ,a_n, b_1,b_2,\ldots ,b_k$ in $\Z$.
Hence for $(x,y)\in F\oplus G$,
\bea
 (x,y) & = & \left(\sum_{i=1}^n a_if_i , \sum_{j=1}^k b_jg_j \right)\\
       & = & \sum_{i=1}^n (a_if_i,0) + \sum_{j=1}^k (0,b_jg_j)\\
       & = & \sum_{i=1}^n a_i(f_i,0) + \sum_{j=1}^k b_j(0,g_j)
\eea Further, \ba{rl}
            & \sum_{i=1}^n a_i(f_i,0) + \sum_{j=1}^k b_j(0,g_j) = 0\\
\Rightarrow & \left(\sum_{i=1}^n a_if_i , \sum_{j=1}^k b_jg_j \right) =
0\\ \Rightarrow & \sum_{i=1}^n a_if_i = 0,\mbox{ and } \sum_{j=1}^k b_jg_j
= 0\\ \Rightarrow & a_i = 0, b_j = 0 \forall 1\leq i\leq n, 1\leq j\leq k.
\ea Hence the assertion.
\eos
\bexe\label{ex91}
Let $\Z^n = (\Z^n,+)$, and $f_i = (0,0,...,1,..0 )$ for $i = 1,2,...n$,
where $1$ is in the $i^{th}$ place. Prove that
 $$ \Z^n  = \fab{f}{n}.$$
 i.e., $(\Z^n,+)$ is a free abelian group with basis $\{ f_1,f_2,\ldots ,f_n \}.$
\eexe
\bl\label{l91}
Let $F  = \fab{f}{n}$, $n\geq 1$, be a finitely generated  free abelian group.
Then $F \simeq \Z^n.$
\el
\pf Define
\bea
     \alpha : \Z^n & \to &  F\\
     (b_1,b_2,\ldots ,b_n) & \mapsto &  b_1f_1+b_2f_2+\ldots +b_nf_n.
\eea
For $\underline{b} = (b_1,b_2,\ldots ,b_n)$ and
$ \underline{c} = (c_1,c_2,\ldots ,c_n)$ in $\Z^n$, we have
\bea
\alpha(\underline{b}+ \underline{c}) & = &
                           \alpha(b_1+c_1,b_2+c_2,\ldots ,b_n+c_n)\\
              & = &  (b_1+c_1)f_1+(b_2+c_2)f_2+\ldots +(b_n+c_n)f_n \\
              & = & (b_1f_1+b_2f_2+\ldots +b_nf_n) +
                     (c_1f_1+c_2f_2+\ldots +(c_nf_n)\\
              & = & \alpha(\underline{b}) + \alpha(\underline{c}).
\eea
Therefore, $\alpha$ is a homomorphism.As $$ F = \Z f_1+ \Z f_2+ \ldots +\Z
f_n,$$ $\alpha$ is clearly onto. Finally, if $\alpha(\underline{b}) = 0,$
then \ba{rl}
             &  b_1f_1+b_2f_2,\ldots ,b_nf_n = 0\\
\Rightarrow  &  b_1 = b_2 = \ldots = b_n = 0, \ea since $f_1,f_2,\ldots
,f_n$ are linearly independent.Therefore $\underline{b} = 0$ which implies
$\alpha$ is one-one. Hence $\alpha$ is an isomorphism, and $F\simeq \Z^n.$
$\blacksquare$

\bexe\label{92}
  Prove that the converse of the lemma is true i.e.,if an abelian group
  $F$ is isomorphic to the free abelian group $\Z^n,$ then $F = \fab{f}{n}$
  for some $ f_i \in F,\,i=1,2,\cdots,n.$
\eexe
\bexe \label{93}
     Let $m\geq 1$, $n\geq 1$ be two integers. Prove that for the free
     abelian groups $\Z^n=(\Z^n, +)$ and $\Z^m=(\Z^m,+)$, $\Z^n\bigoplus
     \Z^m$ is isomorphic to $\Z^{m+n}=(\Z^{n+m},+).$ Then deduce that the
     direct sum of two free abelian groups of ranks $n$ and $m$ is a free
     abelian group of rank $n+m$.
     \eexe
\bl\label{l92}
Let $F = \fab{f}{n}$ be a free abelian group, and let  $A$ be any abelian group.
Then for any $a_i, 1\leq i \leq n$,in $A$ there exists a unique homomorphism
$$\theta : F \to A$$
such that $\theta(f_i) = a_i$.
\el
\pf Since $F$ is a free abelian group with basis $\{f_1,f_2,\ldots ,f_n \}$,
any $x\in F$ can be uniquely expressed as
$$ x =  \lambda_1f_1+ \lambda_2f_2+\ldots + \lambda_nf_n, (\lambda_i\in \Z).$$
Thus the map
\bea
              \theta : F & \to & A\\
       \sum_{i=1}^n \lambda_if_i & \mapsto & \sum_{i=1}^n \lambda_ia_i
\eea
is well defined. Now, let
$$ x = \sum_{i=1}^n \lambda_if_i , y = \sum_{i=1}^n \mu_if_i $$
be any two elements of $F$. Then
\bea
    \theta(x+y) & = & \theta (\sum_{i=1}^n (\lambda_i + \mu_i)f_i)\\
                & = & \sum_{i=1}^n (\lambda_i + \mu_i)a_i\\
                & = & \sum_{i=1}^n \lambda_ia_i + \sum_{i=1}^n \mu_ia_i\\
                & = & \theta(x) + \theta(y)
\eea
Hence, $\theta$ is a homomorphism from $F$ to $A$ with $\theta(f_i) =
a_i.$ The uniqueness of $\theta$ is clear as any element of $F$ is of the
form $$\sum_{i=1}^n \lambda_if_i, \lambda_i \in \Z.\, \blacksquare$$
\bl\label{l93}
For any two integers $1\leq m \neq n < \infty$, $ \Z^m \not\simeq \Z^n.$
\el
\pf Define a relation  $\sim$ on $\Z^n (\Z^m)$ as follows :\\
For $a = \fab{a}{n}$ and $b=\fab{b}{n}$ in $\Z^n$,
$$ a\sim b  \Leftrightarrow  a-b\in 2\Z^n \mbox {i.e., }a_i-b_i \in 2\Z
                         \mbox { for all } i  = 1,2,...n.$$
We shall show that $\sim$ is an equivalence relation.\\ (i) Reflexivity :
$a\sim a,$ as $a-a =  0 \in 2\Z^n.$ \\ (ii) Symmetry : Let $a\sim b.$ Then
$$ a-b \in 2\Z^n $$  $$ \Rightarrow b -a \in 2\Z^n$$ Hence $b\sim a.$\\ (iii)
Transitivity : Let for $a,b,c \in \Z^n, a\sim b, b\sim c.$ Then
\begin{eqnarray}
             & &   a-b\in 2\Z^n \mbox{ and } b-c\in 2\Z^n\\
\Rightarrow  & & a-c = (a-b) - (c-b) \in 2\Z^n \nonumber \\
\Rightarrow  & & a\sim c.\nonumber
\end{eqnarray}
Hence $\sim$ is an equivalence  relation\\
Now, suppose $\alpha : \Z^n \to \Z^m$ is an isomorphism.Then, for $a,b\in
\Z^n$,
$$ \alpha(a-b) = \alpha(a)-\alpha(b).$$
If $a\sim b,$ then $a-b = 2c$ for some $c\in \Z^n.$  Hence
\ba{rl}
             &  \alpha(a-b)  = \alpha(2c)\\
\Rightarrow  &  \alpha(a)- \alpha(b) = 2 \alpha(c)\\
\Rightarrow  &  \alpha(a)\sim \alpha(b).
\ea
Conversely, if $\alpha(a) \sim \alpha(b)$, then
\ba{rl}
            & \alpha(a)- \alpha(b) = 2d \mbox { for some }d\in \Z^m\\
\Rightarrow & \alpha(a-b) = 2 \alpha(c) \mbox{ for some }c\in \Z^n, \mbox{
as }
                        \alpha \mbox{ is onto.}\\
\Rightarrow & \alpha(a-b) = \alpha(2c)\\
\Rightarrow & (a-b) = 2c \mbox { as }\alpha \mbox{ is one-one.}\\
\Rightarrow & a\sim b.
\ea
Hence $\alpha$ maps bijectively  $\Z^n/\sim$ to $\Z^m/\sim.$

We  shall  now show that $\abs{\Z^n/\sim} = 2^n$. Note that, if $a =
(a_1,a_2,\ldots ,a_n)$\\ $ \in \Z^n,$  then $$  a\sim
(\delta(a_1),\delta(a_2),\ldots ,\delta(a_n)),$$ where $\delta(a_i) = 0$
if $a_i$ is even and, $\delta(a_i) = 1$ if $a_i$ is odd. Further, clearly
no two elements of the form $(c_1,c_2,\cdots,c_n)$ in $\Z^n$, where
$c_i=0$ or 1, can be equivalent under the relation $\sim$
 unless they are equal. Hence,
$$\abs{\Z^n/\sim} = 2^n.$$
As $\alpha$ maps bijectively $\Z^n/\sim$ to $\Z^m/\sim$, we get $2^n = 2^m.$
Therefore, $m = n$. Hence, if $m \neq n$, $\Z^n$ is not isomorphic to
$\Z^m.$ $\blacksquare$
\bl\label{l94}
Any two bases of a finitely generated free abelian group have  same number
of elements.
\el
\pf Let $F$ be a finitely generated free abelian group, and let
\bea
            F & = & \fab{f}{n}\\
              & = & \fab{g}{m}
\eea
Then, by lemma \ref{l91}, $F \simeq \Z^n$ and $F \simeq \Z^m.$  Thus
 $\Z^n \simeq \Z^m.$  Hence, by the lemma \ref{l93}, $m = n.$
\bexe\label{ec95}
Show that in a finitely generated abelian group $A$, every infinite subset
is linearly dependent.
\eexe
\bd\label{d97}
The number of elements in a basis of a finitely generated free abelian group
is called the rank of the free abelian group.
\ed
\br
A free abelian group $F$ has rank $0$ if and only if $F = {0}.$
\er
\bl\label{l99}
Let $A$ be a finitely  generated abelian  group. Let $$ A^t = \{ a\in A
\st a \mbox{ is a torsion element i.e., }\circ (a)< \infty \}.$$
Then $A^t$ is a subgroup of $A.$
\el
\pf Clearly $0\in A^t.$  Let $a,b\in A^t$. Then there exist $m \geq 1$ and $n\geq 1$
 such that $na = 0$ and $mb = 0.$ Hence,
\bea
          mn(a-b) & = & mna - mnb\\
                  & = & m(na)-n(mb) = 0.
\eea
$\Rightarrow$ $a-b \in A^t.$ \\  Hence, $A^t$ is a subgroup.$\blacksquare$
\br\label{d910}
An abelian group A is a torsion group if and only if  $A^t = A$.
\er
\bl\label{l911}
For an abelian group $A$ if $A \neq A^t$, then the factor group $A/A^t$ is
a torsion free group.
\el
\pf Let $a\in A-A^t,$ and let $n(a+A^t) = A^t$ for some $n\geq 1.$ Then
\ba{rl}
            &  na+A^t = A^t\\
\Rightarrow & na\in A^t \\ \Rightarrow & mna = 0 \mbox{ for some }m\geq
1\\ \Rightarrow & a\in A^t. \ea This contradicts the assumption that $a
\not\in A^t.$ Therefore,in case $A \neq A^t,$ $A/A^t$ is a torsion free
group.         $\blacksquare$
\bl\label{l912}
Let $F = \fab{f}{n},\,n\geq 0,$ be a finitely generated free abelian
group. Then any subgroup $H$ of $F$ is free of rank $\leq n.$
\el
\pf If $F=(0),$ there is noting to prove. Hence, let $F\neq (0).$ \\ If $H
= (0)$, then the proof is clear since in this case $H$ is free of rank
$0$. Hence, let $H\neq (0)$. We shall, now, prove the result by
induction on $n.$\\
 $n = 1$ : In this case $F = \langle f_1 \rangle = \Z f_1.$
Put $$ I = \{ m\in \Z \st mf_1\in H \}.$$ As $H \neq (0)$, $I \neq
(0).$ Further, if $m$ and $n$ $\in I$, then $(m-n)f_1 = mf_1-nf_1
\in H \mbox{ as }H \mbox { is a subgroup.}$ Hence $m-n \in I.$
Therefore $(I,+)$ is a subgroup of $(\Z,+)$. As $(\Z,+)$ is a
cyclic group, $(I,+)$ is a cyclic group. Let $I = \Z k.$ Then
$k\neq 0$ and $H = \Z (kf_1)$. Thus, clearly $H = \langle kf_1
\rangle.$ Hence, the result holds for $n = 1.$\\ $n > 1$
: In this case  $F = \fab{f}{n}$ where $n > 1.$ For any $a_1,a_2,\ldots
,a_n$ in an abelian group $A$, the map  \[\begin{array}{rl}
  \theta : F & \to A  \\
  \sum_{i=1}^n m_if_i & \mapsto \sum_{i=1}^n m_ia_i
\end{array}\]

is a homomorphism (Lemma \ref{l92}). For $F_1 = gp\{f_1 \},$ consider the
homomorphism,
\bea
       p_1 : F & \to & F_1 \\
  \sum_{i=1}^n m_if_i & \mapsto & m_1f_1.
\eea
Then $p_1$ is clearly onto and $\ker(p_1) = K = \langle f_2,\ldots ,f_n
\rangle.$ Thus, $K$ is free of rank $n-1$. By induction $H \cap K$ is
free of rank $\leq n-1.$ If $p_1(H) = (0)$, then $H = H\cap K$ and hence
is free of rank $\leq n-1.$ If $p_1(H) \neq (0)$, then by the case $n =
1$, $p_1(H)$ is free of rank $1.$ Let $p_1(H) = \langle m_1f_1 \rangle.$
Then there exists an element $$  h = m_1f_1+ m_2f_2+\ldots +m_nf_n \in
H.$$ Claim : $H = (H \cap K) \oplus \langle h \rangle$.\\ Let $x\in H$,
then \ba{rl}
              &  p_1(x) \in p_1(H) = \langle m_1f_1 \rangle\\
\Rightarrow   & p_1(x) = p_1(th) \mbox{ for some } t\in \Z\\
\Rightarrow & p_1(x-th) = 0\\ \Rightarrow   & x-th \in K.\\
\Rightarrow   & x-th \in H\cap K \mbox { since } x, th \in H.\\
\Rightarrow   & x = (x-th)+th \in (H\cap K) + \langle h \rangle.
\ea Next, as for any $t (\neq 0)\in \Z, p_1(th) = tm_1f_1 \neq 0.$
Therefore, $ K \cap \langle h \rangle = 0.$ Hence, $ H = (H\cap K)
\oplus \langle h \rangle$. Now, as $H\cap K$ is free of rank $\leq
n-1$  and $\langle h \rangle$ is free of rank $1$, $H$ is free of
rank $\leq n (Exercise \ref {93}).$ \bl\label{l913} Let $A$ be a
finitely generated abelian group. Then any subgroup of $A$ is
finitely generated. \el \pf Let $B$ be a subgroup of $A =
gp\{a_1,a_2,\ldots ,a_n \}$. Consider the homomorphism .
\[\begin{array}{rl}
  \theta : \Z^n & \to A  \\
  \sum_{i=1}^n \lambda_ie_i & \mapsto \sum_{i=1}^n \lambda_ia_i
\end{array}\]
 which is
clearly onto. Put $H = \theta^{-1}(B)$. Then $H$ is free of rank
$\leq n,$ i.e., $H$ is finitely generated. Now, as $\theta$ is
onto, $\theta(H) = B$. Hence $B$ is generated by images of
generators of $H$ under $\theta$ and hence is finitely generated.
$\blacksquare$ \bl\label{l915} Let $A$ be a finitely generated
torsion free abelian group. Then $A$ is free. \el \pf Let $A = gp
\{a_1,a_2,\ldots ,a_n,a_{n+1},a_{n+2},\ldots ,a_k \}.$ If
necessary, by change of order, we can assume that
$\{a_1,a_2,\ldots ,a_n \}$ is maximal linearly independent subset
of $\{a_1,\ldots ,a_n,a_{n+1},\ldots ,a_k \}$. Then
$\{a_1,a_2,\ldots ,a_n,a_i \}$ is linearly dependent for each
$i\geq n+1$. Hence there exist $t_i(\neq 0) \in \Z$ such that
$t_ia_i \in F = \fab{a}{n}$ for all $i\geq n+1.$ Put $t =
t_{n+1}\cdot t_{n+2}\cdots t_k \in \Z$, then, $ta_i \in F$ for all
$i.$ Hence $tA \subseteq F$. Now, consider the homomorphism $\tau:
A \to F$ defined by $\tau(a) = ta.$ As $A$ is torsion free, $\tau$
is a monomorphism. Thus $A \simeq tA \subseteq F.$  As $F$ is
free, $tA$ is free of rank $\leq n$ (Lemma \ref{l912}). Hence, $A$
is free of rank $\leq n.$ \bexe\label{ec916} Show that the torsion
free abelian group $(\Q,+)$ is not finitely generated. \eexe
\bd\label{d917} Let $F$ be a free abelian group of rank $n$. Then
for any $f(\neq 0)\in F = \fab{f}{n}$, if $f =
b_1f_1+b_2f_2+\ldots +b_nf_n$, we define, $$\mbox{ Content of }f =
\mbox {cont}(f) = \mbox {g.c.d.}(b_1,b_2,\ldots ,b_n).$$ \ed \br
The definition of content of an element uses a basis of $F$. We
shall see that in fact it is independent of the basis chosen. \er
\bl\label{l918} Let $F \neq (0)$ be a finitely generated free
abelian group. Suppose $$ F = \fab{f}{n} = \fab{g}{n}. $$ Let for
$x (\neq 0)\in F$, $$ x = a_1f_1+a_2f_2+\ldots +a_nf_n =
b_1g_1+b_2g_2+\ldots +b_ng_n ,( a_i, b_i \in \Z ). $$ Then $$
g.c.d.(a_1,a_2, \ldots ,a_n) = g.c.d.(b_1,b_2, \ldots ,b_n)$$
i.e., content of an element in $F$ is independent of the basis of
$F$. \el \pf Let $g.c.d.(a_1,a_2,\ldots ,a_n) = d$ and
$g.c.d.(b_1,b_2,\ldots ,b_n) = e.$ Since $ F = \fab{f}{n} =
\fab{g}{n}$, we can write $$ f_i = \sum_{j=1}^n
\lambda_{ij}g_j,\mbox{  where  }\lambda_{ij}\in \Z.$$ for all $i =
1,2,\ldots ,n.$ Hence \bea
        x & = & a_1f_1+a_2f_2+\ldots +a_nf_n\\
& = & a_1(\sum_{j=1}^n \lambda_{1j}g_j)+\ldots
                              +a_n(\sum_{j=1}^n \lambda_{nj}g_j)\\
& = & (\sum_{k=1}^n a_k \lambda_{k1})g_1+ \ldots +
                                (\sum_{k=1}^n a_k \lambda_{kn})g_n.
\eea As any element of $F$ can be expressed uniquely as linear
combination of basis elements with respect to any given basis, we
conclude, \ba{rl} & b_i = \sum_{k=1}^n a_k \lambda_{ki}\;\;\;\;
\forall i = 1,2,\ldots ,n.\\ \Rightarrow & d \mid b_i \;\;\;\;
\forall i = 1,2,\ldots ,n.\\ \Rightarrow & d\mid e. \ea Similarly,
we can show that $e\mid d$. Hence $ e = d.$ $\blacksquare$
\bl\label{l919} Let $F = \fab{f}{n}$ and $g (\neq 0)\in F$ be such
that $cont(g) = 1$. Then, $g$ is part of a basis of $F.$ \el \pf
Let $g = a_1f_1+a_2f_2+\ldots +a_nf_n,$ where $ a_i\in \Z$ for all
$i=1,2\cdots,n.$ As $cont(g) = 1$,  $g.c.d.(a_1,a_2,\ldots ,a_n) =
1.$ Thus there exist $ b_1,b_2,\ldots ,b_n \in \Z$ such that
$a_1b_1+a_2b_2+\ldots +a_nb_n = 1.$ Define : \bea
         \theta : F & \!\!\!\!\!\!\! \to  \Z  \\
 \sum t_if_i & \mapsto \sum t_ib_i
\eea
Then $\theta$ is a group homomorphism. Further, \ba{rl}
            & \theta(g) = \theta(\sum a_if_i) = \sum a_ib_i = 1\\
\Rightarrow & \theta \mbox { is onto.} \ea Let $H = \ker(\theta).$ Then
$H$ is free (Lemma \ref{l912}). We shall prove that $F = H \oplus \langle
g\rangle$. Clearly, $gp\{g\}=\langle g\rangle$ and $H + \langle g\rangle
\subseteq F.$ Let for $x\in F,$ \ba{rl}
             & \theta(x) = m \\
\Rightarrow  & \theta(x) = m \theta(g) \mbox { since }\theta(g) = 1.\\
\Rightarrow  & \theta(x-mg) = 0\\ \Rightarrow  & x-mg \in \ker(\theta) =
H\\ \Rightarrow  & x\in H + \langle g\rangle \\
    \Rightarrow & H+\langle g\rangle=F. \ea  Next, if $y\in H\cap \langle g\rangle,$
     then \ba{rl}
               & y = tg \mbox{ and }y = h,\mbox{ where }t\in \Z \mbox { and }
                                    h\in H.\\
\Rightarrow & \theta(y) = \theta(tg) \mbox{ and }\theta(y) = \theta(h) =
0\\ \Rightarrow & \theta(y) = t = 0, \mbox{ since }\theta(g) = 1.\\
\Rightarrow & y = 0\\ \Rightarrow & H\cap \langle g\rangle = 0\\ \ea
Hence, $F = H \oplus \langle g\rangle.$ Now, $g$ together with a basis of
$H$ gives a basis of $F.$ Consequently $g$ is part of a basis of $F$.
 $\blacksquare$
\bt\label{t920}
Let $F = \fab{f}{n}$ be a free abelian group of rank $n$ and $H(\neq 0)
\subseteq F$, a subgroup. Then there exists a basis $\{g_1,g_2,\ldots ,g_n
\}$ of $F$ and non-zero integers $a_1,a_2,\ldots ,a_k$ $(0 < k \leq n)$,
$a_i \mid a_{i+1}$, such that $H$ is a
 free abelian group with basis $\{a_1g_1,a_2g_2,\ldots ,a_kg_k \}.$
\et
\pf Let $h (\neq 0)\in H$ be an element of least content in $H$ and let
 $a_1 = cont(h).$ Then we can write $h = a_1h_1,$ where $h_1\in H,$ and
 $cont(h_1) = 1.$ By lemma \ref{l919}, $h_1$ is part of a basis of $F$. Therefore, we can
 choose $ h_2',h_3',\dots ,h_n'$ in $F$ such that
$F = \langle h_1,h_2',h_3',\dots ,h_n'\rangle$. Note that $F = \langle
h_1\rangle \oplus \langle h_2',h_3',\dots ,h_n'\rangle$. We shall prove
that $H
=
<a_1h_1> \oplus H\cap \langle h_2',h_3',\dots ,h_n'\rangle.$ Let $x
(\neq 0)\in H.$ Then $$ x = \lambda_1h_1 +\lambda_2h_2'+ \ldots
+\lambda_nh_n',\mbox{ for some } \lambda_i \in \Z.$$ Let $\lambda_1 =
a_1q_1 + r_1, 0\leq r_1\lvertneqq a_1.$ Then $$ x-q_1h = x-q_1a_1h_1 =
r_1h_1+ \lambda_2h_2'+\ldots +\lambda_nh_n'.$$ If $r_1\neq 0,$ then,
$cont(x-q_1h) \leq r_1 < a_1.$ This contradicts the choice of $h.$
Therefore $r_1 = 0,$ and we get \ba{rl}
         & x-q_1h = x-q_1a_1h_1\\
         & = \lambda_2h_2'+\ldots +\lambda_nh_n'\in
                         H\cap \langle h_2',h_3',\dots ,h_n'\rangle.\\
\Rightarrow & x\in \langle a_1h_1 \rangle \oplus
                                  H\cap \langle h_2',h_3',\dots
                                  ,h_n'\rangle \\
\Rightarrow  & H=\langle a_1h_1 \rangle \oplus
                                  H\cap \langle h_2',h_3',\dots
                                  ,h_n'\rangle. \ea  If $H\cap
\langle h_2',h_3',\dots ,h_n'\rangle = (0),$ then the proof is finished.
However, if $H\cap \langle h_2',h_3',\dots ,h_n'\rangle \neq (0),$ then
 by induction, there exists a basis $\{h_2,h_3,\dots ,h_n \}$ of
 $ \langle h_2',h_3',\dots ,h_n'\rangle$, and non-zero integers
$a_2,a_3,\ldots ,a_k (2\leq k \leq n),$  $a_i \mid a_{i+1},$ for $i\geq 2$,
such that
$ H\cap \langle h_2',h_3',\dots ,h_n'\rangle
                = \langle a_2h_2,\ldots ,a_kh_k\rangle.$
In this case, $F = \fab{h}{n}$ and $H = \langle a_1h_1,a_2h_2,\ldots
,a_kh_k\rangle.$ The proof will be  complete if we show that $a_1 \mid a_2.$ If
$a_1$ does not
 divide $a_2$, then
$ cont(a_1h_1+a_2h_2) = g.c.d.(a_1,a_2) < a_1$, and $(a_1h_1+a_2h_2)\in H.$
This contradicts the  choice of $h.$ Hence, $a_1$ divides $a_2$ and the proof
 follows.      $\blacksquare$
\bl\label{l921}
Let $A$ be an abelian group and $F = \fab{f}{n}$ be a free abelian group.
If $\theta : A \to F$ is an onto homomorphism, then, $$ A \simeq F\oplus
\ker \theta$$
\el
\pf Since $\theta$ is onto, there exist $a_i\in A$ such that
$ \theta(a_i) = f_i$ for $ i = 1,2,\ldots ,n.$ Define
\bea
           \psi : F & \to & A\\
         \sum t_if_i & \mapsto &  \sum t_ia_i.
\eea
Then $\psi$ is a homomorphism, and \ba{rl}
     & (\theta \circ \psi)(\sum t_if_i) = \theta(\sum t_ia_i)
                = \sum t_i\theta(a_i) = \sum t_if_i.\\
\Rightarrow &  \theta \circ \psi = Identity.
\ea
Note that,
\ba{rl}
                    & \psi(\sum t_if_i) = 0\\
\Rightarrow & (\theta \circ \psi)(\sum t_if_i) = 0\\ \Rightarrow & \sum
t_if_i = 0\\ \Rightarrow &  \psi \mbox { is one-one.} \ea Hence, $\psi$ is
an isomorphism onto  $\psi(F).$ We shall now prove that
$A\simeq\psi(F)\bigoplus \ker\theta$. Clearly, $ \psi(F) + \ker \theta \subseteq
A.$ Let $x\in A, x\neq 0.$ Then, \ba{rl}
      & \theta(x-(\psi \circ \theta)(x))\\
      & = \theta(x) - (\theta \circ \psi)(\theta(x))\\
      & = \theta(x)- \theta(x)\\
      & = 0.\\
\Rightarrow  & x - (\psi \circ \theta)(x) \in \ker \theta\\
\Rightarrow  & x\in \psi(F) + \ker \theta\\
\Rightarrow  & A = \psi(F) + \ker \theta.\\
\ea
Now, let $y\in \psi(F) \cap \ker \theta.$ Then
\ba{rl}
            & y = \psi(f) = a \mbox { where }f\in F \mbox { and }
                               a\in \ker\theta.\\
\Rightarrow & \theta(y) = \theta (\psi(f)) = \theta(a)\\
\Rightarrow & \theta(y) = f = 0\\
\Rightarrow & y = \psi(f) = 0\\
\Rightarrow & (0) = \psi(F) \cap \ker \theta.\\
\ea
Hence,
\ba{rl}
            & A = \psi(F) \oplus \ker \theta\\
\Rightarrow & A \simeq F \oplus \ker \theta. \blacksquare \ea
\br
In the above proof freeness of $F$ is used only to construct $\psi.$
However, if $F$ is not necessarily free and there exists an onto
homomorphism $\theta
: A \to F$ and a homomorphism $\psi : F \to A$ such that
$(\theta \circ \psi) = Identity,$ then $A = \psi(F) \oplus \ker \theta.$
\er  \vspace{.25in}
{\bf STRUCTURE OF FINITELY GENERATED ABELIAN  GROUPS  } \\ Let $A$ be a finitely
generated abelian group and let $A^t$ be the torsion subgroup of $A.$ If
$A\neq A^t,$ then $A/A^t$ is torsion free, finitely  generated,
 abelian group and hence is free by lemma \ref{l915}. Consider the natural
 homomorphism
$$ \eta : A \to A/A^t$$ Since $\ker\eta = A^t,$ by the lemma \ref{l921},
$A\simeq  A^t \oplus A/A^t.$ As $A$ is finitely generated, $A^t$ is a
finitely generated abelian group (Lemma \ref{l913}). Hence, any finitely
generated  abelian group is the direct sum of a finitely generated free
group and a finitely generated  torsion group. Clearly, a finitely
generated torsion abelian group is finite. Hence as a finitely generated
free abelian group is isomorphic to $\Z^n$ for some $n\geq 0,$ to know the
structure of a finitely generated abelian group, we just have to know the
structure of a finite abelian group. We first, prove:
\bt\label{t922}
Let $A$ be a finitely generated abelian group, and let
$$ A\simeq F\oplus B, \mbox { and } A\simeq F_1\oplus C$$
where $ F,F_1 $ are free abelian groups of finite ranks and $ B,C$ are finite
abelian groups. Then $F\simeq F_1$ and $B\simeq C.$
\et  \pf Let
$$\theta_1 : A \simeq F\oplus B \mbox{ and }\theta_2 : A \simeq F_1 \oplus C$$
be isomorphisms of groups. Then
$$ \theta_1(A^t) = (F\oplus B)^t = B$$
and
$$ \theta_2(A^t) = (F_1\oplus C)^t = C$$
Therefore $B\simeq C$. Further, $\theta_1$ and $\theta_2$ shall induce the
isomorphisms
\bea
          \bar{\theta_1} : A/A^t & \simeq & (F\oplus B)/B = F\\
          a+A^t & \mapsto & \theta_1(a)+B
\eea
and
\bea
          \bar{\theta_2} : A/A^t & \simeq & (F_1 \oplus C)/C = F_1\\
          a+A^t & \mapsto & \theta_2(a)+C
\eea
Hence $F$ is isomorphic to $F_1$.     $\blacksquare$
\bl\label{l923}
Let $B (\neq (0))$ be a finite abelian group. Then,
$$ B \simeq \Z/(n_1) \oplus \Z/(n_2) \oplus \ldots \oplus \Z/(n_k)$$
 where $1\leq n_1$  and $n_i$ divides $n_{i+1}$  for all $1\leq i \leq k-1.$
\el \pf Let $k$ denotes the number of elements in a minimal set of
generators of $B$. As $B (\neq (0))$, $k\geq 1.$ Thus there exists
an onto homomorphism $$ \psi : \Z^k \to B$$ Put $F = \Z^k$, and
$F_1 = \ker \psi$. Clearly $F_1 \neq (0)$, since $F/F_1 $ is
isomorphic to $B$, which is a torsion group. By theorem
\ref{l915} there exists a basis $\{g_1,g_2,\ldots ,g_k \}$ of $F$
and integers $1 < n_1 \leq n_2\leq \ldots \leq n_l $, $1\leq l\leq
k$, such that $n_i \mid n_{i+1}$ for all $i\geq 1,$ and $$ F_1 =
\langle n_1g_1, n_2g_2,\ldots ,n_lg_l\rangle.$$ Then \bea F/F_1 &
\simeq & \frac{\Z g_1 \oplus \Z g_2 \oplus \ldots \oplus \Z g_k}
               {\Z n_1g_1 \oplus \Z n_2g_2 \oplus \ldots \oplus \Z n_lg_l}\\
  & \simeq & \frac{\Z g_1}{\Z n_1g_1} \oplus \ldots
              \oplus \frac{\Z g_l}{\Z n_lg_l} \oplus \Z g_{l+1}\oplus
                         \ldots \oplus \Z g_k
\eea
As $F/F_1 \simeq B$ is a torsion group, $l = k$. Hence
\bea
F/F_1 & \simeq & \frac{\Z g_1}{\Z n_1g_1} \oplus \ldots
                         \oplus \frac{\Z g_l}{\Z n_kg_l}\\
      & \simeq & \frac{\Z}{(n_1)} \oplus \ldots
              \oplus \frac{\Z}{(n_k)}
\eea
where $1 < n_1 \leq n_2\leq \ldots \leq n_l$, $n_i \mid n_{i+1}$
 for all $ i\geq 1.$     $\blacksquare$\\
 \\  {\large\bf An alternate proof of lemma \ref{l923},} \\
 \bl\label{l923l.1} Let $A(\neq 0)$ be a finite abelian group, and let
 $x\in A$ be an element of maximal order in $A$. Then the cyclic group
 $gp\{x\}$ is a direct summand of $A$. \el  \pf Let $\circ(x)=d$. We shall
 first show that for any $y\in A$, $\circ(y)$ divides $d$. If not, then
 for a prime $p$, and $\alpha >0$, $p^{\alpha}$ divides $\circ(y)$, but,
 does not divide $d$. Let $\circ(y)=mp^{\alpha}$. Then
 $\circ(my)=p^{\alpha}$, and $p^{\alpha}\nmid d$. Let $d=p^{\beta}n$
 where $\beta\geq 0$, and $(p,n)=1.$ As $p^{\alpha}\nmid d$, $\alpha
 >\beta$, and clearly $\circ(p^{\beta}x)=n$. Now, for $x_1=p^{\beta}x$,
 $y_1=my$, $(\circ(x_1),\circ(y_1))=(n,p^{\alpha})=1.$ Therefore
 $\circ(x_1+y_1)=p^{\alpha}n>d=p^{\beta}n$ (Theorem 2.35 ). This contradicts
 the choice of $x$. Hence $\circ(y)$ divides $\circ(x)$ for all $y\in A$.
 \\ Choose a maximal subgroup $B$ of $A$ such that $B\cap gp\{x\}=(0)$.
 Put $C=B+gp\{x\}$. To prove the result it suffices to show that $A=C$. If
 not, then the factor group $A/C$ is non-zero. Choose an element of prime
 order (say) $q$ in $A/C$. Thus there exists $z\in A-C$, such that \[qz\in
 C = B+\{x\}\] \be\label{epr}
   \Rightarrow\;\;qz = b+ mx \mbox{ for some }m\in\Z\mbox{ and }b\in B\ee
   Let $\circ(z)=k$. We claim that $q|k$. If not, then $(q,k)=1.$ Hence,
   as $kz=0$ and $qz\in C$, $z\in C$. A contradiction to our assumption.
   Thus $q|k$, and hence $q|d$ since $k|d$. From the equation
   (\ref{epr}),
   we get \[\begin{array}{cl}
      & 0=dz=d/q\cdot b+d/q\cdot mx \\
     \Rightarrow & d/q\cdot b= d/q\cdot mx=0\;\;\mbox{ since }B\cap gp\{x\}=(0) \\
     \Rightarrow & d\;\mbox{ divides }d/q\cdot m \\
     \Rightarrow & q|m. \\
   \end{array}\] Thus the equation (\ref{epr}), can be written as
   \[\begin{array}{cll}
      & qz=b+qm_1x & \mbox{ where }m=qm_1 \\
     \Rightarrow & q(z-m_1x)=b &  \\
   \end{array}\] Put $z_1=z-m_1x.$ Then $z_1\not\in B+gp\{x\}$ as
   $z\not\in B+gp\{x\}$, moreover, $qz_1=b\in B.$ As $z_1\not\in
   C=B+gp\{x\}$, and $B$ is maximal disjoint from $gp\{x\}$,
   \[\begin{array}{clc}
      & (B+gp\{z_1\})\bigcap gp\{x\}\neq 0 &  \\
     \Rightarrow & b_1+sz_1=tx\neq 0 & \mbox{ for some }s,t\in\Z \mbox{ and } b_{1}\in B  \\
   \end{array}\] Let us note that $q\nmid s$, since otherwise
   \[\begin{array}{clr}
      & b_1+sz_1=tx\in B\bigcap gp\{x\}=0 & \mbox{  since } qz_1\in B \\
     \Rightarrow & tx=0 &  \\
   \end{array}\] A contradiction to our assumption that $tx\neq 0$. Now,
   as $sz_1=tx-b_1\in B+gp\{x\}$, and $qz_1\in B$, $z_1\in B+gp\{x\}$
   since $(q,s)=1.$ This contradicts the fact that $z_1\not\in B+gp\{x\}$.
   Hence $A=C=B+gp\{x\}$, and the result follows.       $\blacksquare$
   \bco Let $A(\neq 0)$ be a finite abelian group. Then $A$ is isomorphic
   to a direct sum $C_{n_1}\bigoplus C_{n_2}\bigoplus\cdots \bigoplus C_{n_k}$
   ($C_{n_i}$ : cyclic group of order $n_i$) where $0<n_1\leq
   n_2\leq\cdots\leq n_k$,and $n_i\mid n_{i+1}$ for all $1\leq i\leq k-1$.
   \eco \pf The proof is straight forward deduction from the theorem. \br
   The above corollary is a restatement of lemma \ref{l923}.\er
Before proceeding further, we shall prove:
\bl\label{l924}
Let $m,n$ be two positive integers. Then
$$n\Z_m = \Z_d$$
where $d = \frac{m}{(m,n)}$ and $\Z_m$ denotes $(\Z/(m))$.
\el
\pf The group $\Z_m$ is cyclic generated by the element $1+(m)$. Hence,
$n\Z_m$ is cyclic generated by $n(1+(m))$. By theorem \ref{t27}, order of the element
$n(1+(m))$ is $\frac{m}{(n.m)} = d$ (say). Hence the result follows.
$\blacksquare$
\bco\label{co925}
 In the lemma, $m \mid n$ if and only if $n(\Z_m) = (0)$.
\eco
\bl\label{l926}
Let $m,m_1$ and $n$ be three positive integers where $m$ divides $m_1$. Then
$\frac{m}{(n,m)}$ divides $\frac{m_1}{(n,m_1)}$.
\el
\pf Let $m_1 = ma$. Then
\bea
               (n.m_1) & = & (n,ma)\\
                       & = & (n,m)\left(\frac{n}{(n,m)},\frac{m}{(n,m)}a\right)\\
                       & = & (n,m)\left(\frac{n}{(n,m)},a\right), \mbox { since }
                                  \left(\frac{n}{(n,m)},\frac{m}{(n,m)}\right) = 1.
\eea
Hence
\bea
              \frac{m_1}{(n,m_1)} & = & \frac{ma}{(n,m)(\frac{n}{(n,m)},a)}\\
              & = & \frac{m}{(n,m)} \frac{a}{(\frac{n}{(n,m)},a)}\\
  \Rightarrow & & \frac{m}{(n,m)} \mid \frac{m_1}{(n,m_1)}
\eea
Hence the result.      $\blacksquare$
\bt\label{t927}
 Let  $n_i, 1\leq i\leq k$, and $m_j, 1\leq j\leq l$, be positive integers such that
$ n_1> 1, m_1> 1$, $n_i\mid n_{i+1}, m_j\mid m_{j+1}$ for all $1\leq i\leq
k-1$ and $1\leq j\leq l-1.$ Then if for an abelian group $A$,
\be\label{e95}
 A \simeq \Z/(n_1) \oplus \Z/(n_2) \oplus \ldots \oplus \Z/(n_k)
\ee
and also
\be\label{e95e}
    A \simeq \Z/(m_1) \oplus \Z/(m_2) \oplus \ldots \oplus \Z/(m_l),
\ee
then $l = k$ and $m_i = n_i$ for all $i.$
\et
\pf Choose a prime $p$ such that $p \mid n_1.$ . Note that whenever $p$
divides an integer $m$,
\bea
                    p \Z_m & = & \{p(k+m\Z) \st k\in \Z \}\\
                           & = & \{(pk+m\Z) \st k\in \Z \}\\
                           & = & p\Z/m\Z
\eea
Further, if $p\notdivide m$, then there exist integers $a,b$ such that
$pa+mb = 1$. Hence
\bea
             1+m\Z & = & pa+mb+m\Z\\
                   & = & p(a+m\Z) \in p\Z_m\\
\Rightarrow  p\Z_m & = & \Z_m
\eea
Thus, from equation (\ref{e95}) \ba{rl}
          &  A/pA \simeq \Z/(p) \oplus \Z/(p) \oplus \ldots
                                \oplus \Z/(p) \mbox{ $k$-times}\\
\Rightarrow & \abs{A/pA} = p^k \ea Next, from equation (\ref{e95e}),
clearly $$ \abs{A/pA} \leq p^l$$ Therefore \ba{rl}\label{e96}
           & p^k\leq p^l\\
\Rightarrow & k\leq l. \ea Similarly, we can show that $l\leq k$. Hence $l
= k.$ For any fixed $1\leq t\leq k$,
using the lemma \ref{l924}, from equations (\ref{e95}) and (\ref{e95e}),
we get
\be\label{e97}
m_tA \simeq \Z/(d_1) \oplus \Z/(d_2) \oplus \ldots \oplus \Z/(d_k)
\ee
where $d_j = \frac{n_j}{(m_t,n_j)}$,\\
and
\be\label{e98}
 m_tA \simeq \Z/(e_1) \oplus \Z/(e_2) \oplus \ldots \oplus \Z/(e_k)
\ee
where $e_j = \frac{m_j}{(m_t,m_j)}.$ Now as $m_j$ divides $m_t$ for $j\leq
t$, we get $e_j = 1$ for all $j\leq t$. Therefore, from equation
(\ref{e98})
\be\label{e99}
 m_tA \simeq \Z/(e_{t+1}) \oplus \Z/(e_{t+2}) \oplus \ldots \oplus \Z/(e_k)
\ee
By the lemma \ref{l926}, $e_i\mid e_{i+1}$ and $d_i\mid d_{i+1}$ for all
$i$. Hence from equations (\ref{e99}) and (\ref{e97}), we get two
decompositions of $m_tA$ which satisfy the hypothesis of the theorem.
Therefore, from the first part of the theorem, we get that the two
decompositions (\ref{e99}) and (\ref{e97}) have same number of terms. Hence
\ba{rl}
                  & d_j = 1 \mbox{ for all } j\leq t.\\
\Rightarrow       & n_t \mid m_t \ea Similarly we can show that
$m_t$ divides $n_t$. Therefore $n_t = m_t \mbox{ for all } t.$
Hence the proof of the theorem is complete.        $\blacksquare$
\bd An additive (multiplicative) abelian group $G$ which is a
finite direct sum (product) of cyclic groups of order $p$ for a
prime $p$ is called an elementary abelian $p$-group.\ed \bl Let
$G$ be a finite group with order $n>1.$ If $Aut\,G$ acts
transitively over the non-identity elements of $G$, i.e., given
any two non-identity elements $x,y$ in $G$ there exists $\alpha\in
Aut\,G$ such that $\alpha(x)=y,$ then $G$ is elementary abelian
$p$-group for a prime $p$.\el  \pf First of all, note that if
$\alpha\in Aut\,G$ and $x\in G$ then $\circ(x)=\circ(\alpha(x)).$
Hence if $x$ is any non-identity element in $G$, then as $Aut\,G$
acts transitively over non-identity element of $G$, for any
$1<i<\circ(x),$ there exists $\beta_i\in Aut\; G$ such that
$\beta_i(x)=x^i$. Therefore $\circ(x)=p,$ a prime. Further, for
any two non-identity elements $x,y$ in $G$, there exists $\beta\in
Aut\,G$ such that $\beta(x)=y$. Hence $\circ(g)=p$ for all $g\in
G.$ Now, to complete the proof we have to show that $G$ is
abelian. By Cauchy's theorem (Theorem \ref{CAUCHY}),
$\circ(G)=p^m$ for some $m\geq 1.$ Further, as $G$ is finite
$p$-group $Z(G)\neq e$ (Theorem \ref{t614}). Take any $x(\neq
e)\in Z(G).$ Then for any $y\in G$ and $\alpha\in Aut\,G$, \\
\[\begin{array}{rl}
   & \alpha(xy)=\alpha(yx) \\
  \Rightarrow & \alpha(x)\alpha(y)=\alpha(y)\alpha(x) \\
  \Rightarrow & \alpha(x)g=g\alpha(x)\,\,\,\mbox{  for all }g\in
  G,\\

  \Rightarrow & \alpha(x)\in Z(G)\,\,\,\mbox{  for all }\alpha\in Aut\,G \\
  \Rightarrow & Z(G)=G.\
  \end{array}\]
   since $Aut\,G$ acts transitively over non-identity elements
in $G$. \\ Hence $G$ is elementary abelian $p$-group.
$\blacksquare$ \hspace{.5in}  \\     \\
      {\bf\Large Divisible groups} \\   \\  In this section, we shall
      analyse divisible (Definition \ref{d26.06d}) abelian groups. We know that
      : \\  (i) Any divisible subgroup of an abelian group
      is its direct factor (Theorem \ref{t6.40}). \\  (ii) Any homomorphic
      image of a divisible group is divisible (Exercise \ref{p12}). \\
      By (ii), any direct factor of a divisible group is divisible. We
      shall see below that the abelian groups $(\Q,+)$ and
      $\C_{p^{\infty}}$ are the only indecomposable divisible abelian
      groups (upto isomorphism).    \bt\label{th11.266th} Let $A$ be an
      additive abelian divisible group. Then  \\  (i) If $A$ contains an
      element of infinite order, $A$ contains a torsion free divisible
      subgroup.  \\   (ii) If $A$ has an element $(\neq 0)$ of finite
      order, then $A$ contains a subgroup isomorphic to $\C_{p^{\infty}}$
      for some prime $p$. \et \pf (i) Let $x_1(\neq 0)\in A$ be an element
      of infinite order. As $A$ is divisible we can choose $x_n\in A$,
      for all $n\geq 2$, such that $nx_n=x_{n-1}$. Let $H=gp\{x_n\st n\geq
      1\}$. We shall prove that $H$ is a divisible torsion free
      subgroup. First of all, we shall show that $ H$ is torsion free.
      Any element $h(\neq 0)$ of $H$ can be expressed as : \[h=m_1x_1+\cdots +
      m_tx_t\] where $ t\geq 1 $ and $0\leq m_i<i$ for all $i\geq 2.$ If $ t=1$ then clearly
      order of $h $ is infinite. Assume $t\geq 2.$  We have
      \[\begin{array}{cl}
        t!h & =t!m_1x_1+\frac{t!}{2!}m_2x_1+\cdots +m_tx_1 \\
         & =(t!m_1+t(t-1)\cdots 3m_2+\cdots +tm_{t-1}+m_t)x_1 \
      \end{array}\] We claim: \[m_1t!+m_2t(t-1)\cdots 3++\cdots +tm_{t-1}
      +tm_{t-1}+m_t\neq 0 \] If not, then \[\begin{array}{cl}
         & t!m_1+m_2t(t-1)\cdots 3+\cdots +tm_{t-1}=-m_t \\
        \Rightarrow & t \mid m_t \\
        \Rightarrow & m_t=0\;\;\;\mbox{ since by assumption } m_i<i \mbox{ for } i\geq 2 \\
        \Rightarrow & m_1t!+m_2t(t-1)\cdots 3+tm_{t-1}=0 \\
        \Rightarrow & (t-1)\mid m_{t-1} \\
        \Rightarrow & m_{t-1}=0 \;\;\;\mbox{ as above, if }  t-1 \geq 2\
      \end{array}\] Continuing the argument we get $m_1=m_2=\cdots m_t =0.$
      This gives $h=0$, which is not true. Hence the claim follows.
      Consequently any $h\neq 0$ in $H$ has infinite order i.e., $H$ is
      torsion free. We shall, now, prove that $H$ is divisible. Note that
      for any $i\geq 1$, \[(i+1)(i+2)\cdots (i+m)x_{i+m}=x_i\] As $m$
      divides $(i+1)(i+2)\cdots (i+m)$, it follows that $my=x_i$ has a
      solution in $H$ for all $m\geq 1$ and $i\geq 1$. Therefore it is
      clear that $H$ is divisible. \\  (ii) Let $x(\neq 0)\in A$ be an
      element of finite order $m$ (say). If $p$ is a prime dividing $m$,
      then for $x_1=m_1 x$ where $m_1p=m$ , $\circ(x_1)=p$. As $A$ is divisible, we can
      choose $x_1,x_2,\cdots x_n,\cdots$ in $A$ such that $px_{n+1}=x_n$
      for $n\geq 1$. Put $H=gp\{x_1,x_2,\cdots,x_n,\cdots\}.$ Then it is
      easy to check that the map: \[\begin{array}{ccl}
        H & \rightarrow & \C_{p^{\infty}} \\
        x_m & \mapsto & e^{2\pi/p^m} \
      \end{array}\] is an isomorphism. The fact that $H$ is divisible
      follows by observing that $qH=H$ for all primes $q$    $\blacksquare$
      \bt\label{27.06th} Let $A(\neq 0)$ be an additive abelian torsion
      free divisible group. Then $A$ contains a subgroup isomorphic to
      $(\Q,+)$, the additive group of rational numbers. \et \pf Let
      $x_1(\neq 0)\in A$. For all $n\geq 1$,choose $x_{n+1}\in A$ such
      that $(n+1)x_{n+1}=x_n$. Then, as in theorem \ref{th11.266th}, the
      subgroup of $H=gp\{x_m\st m\geq 1\}$  is a torsion free divisible
      subgroup of $A$. Define: \[\begin{array}{rcl}
        \phi :(\Q,+) & \rightarrow & H \\
        m/n & \mapsto & y \
      \end{array}\] where $ny=mx_1$. As $H$ is divisible, existence of
      such a $y$ is guaranteed. If $m/n=0$, then $m=0$, hence $ny=0.$
      Therefore $y=0$. Further, let \[\begin{array}{clc}
         & \frac{m}{n}=\frac{k}{l}(\neq 0) & \mbox{ where }m,n,k,l\in \Z \\
        \Rightarrow & ml=kn &  \
      \end{array} \] Let \[\begin{array}{cll}
         & ny=mx_1, & lz=kx_1 \\
        \Rightarrow & nly=mlx_1, & nlz=knx_1 \\
        \Rightarrow & nly=nlz & \mbox{ since }ml=kn \\
        \Rightarrow & y=z & \mbox{ since }A\mbox{ is torsion free } \
      \end{array}\] Hence $\phi$ is well defined. Now, let for $m/n$,
      $m_1/n_1\in \Q$, \[\begin{array}{cl}
         & \phi(m/n)=y, \phi(m_1/n_1)=z \\
        \Rightarrow & mx_1=ny, m_1x_1=n_1z \\
        \Rightarrow & (mn_1+nm_1)x_1= nn_1(y+z) \\
        \Rightarrow & \phi(\frac{mn_1+nm_1}{nn_1})=\phi(\frac{m}{n})+\phi(\frac{m_1}{n_1}) \
      \end{array}\]  Therefore $\phi$ is a homomorphism. Clearly
      $\phi(m/n)=0$ if and only if $m=0$. Hence $\phi$ is a monomorphism
          $\blacksquare$    \brs (i) The groups $(\Q,+)$ and
          $\C_{p^{\infty}}$ are not finitely generated, hence a divisible
          abelian group is not finitely generated. \\ (ii) The divisible
          groups $(\Q,+)$ and $\C_{p^{\infty}}$ are the only
          indecomposable, divisible abelian groups (upto isomorphism). \\
          (iii) Using above result we can prove that any divisible
          group$(\neq id)$ is a direct sum of groups each isomorphic
          either to the additive group of rationals or to the group
          $C_{p^{\infty}}$ for some prime $p$ (different summands may
          correspond to different primes). The proof involves transfinite
          induction and is beyond our scope. \ers

      \vspace{.6in}

  {\bf EXERCISES}
\begin{enumerate}
  \item
  Prove that $\Q_p = \{\frac{m}{p^k} \st m\in \Z, p\mbox{ a prime, }  k\geq
1.\}$ is an additive subgroup of $(\Q,+)$. Show  that $\Q_p$ is not
isomorphic to $\Q_q$ if $p$ and $q$ are distinct primes.
\item
Prove  that $\Q/\Z$ is a torsion group and is not  finite.
\item
Prove that $\Q/\Z$ is the union of subgroups of prime power orders. \item
Show  that $\abs{Aut(\Z)} = 2$ and $\abs{Aut(\Z_n)} = \phi(n).$
\item
Let $(A,+)$ be a torsion free abelian group. Prove that if $\{B_i\}_{i\in
I}$ is a family of divisible subgroups of $A$, then $\bigcap_{i\in I} B_i$
is a divisible subgroup of $A$.
\item
Show that every cyclic group is isomorphic to a subgroup of a divisible
group.
\item
Prove that any finitely generated abelian group is isomorphic to a
subgroup of a divisible abelian group. \\ The result in the above exercise
is true in general. To arrive at that, we define: \\ An abelian group
$(F,+)$ is called free if there exists a subset $S=\{f_i\}_{i\in I}$ in
$F$ such that: \\ (i) $F=gp\{S\}$\\ (ii) $S$ is linearly independent i.e.,
given any finite subset $\{f_{i_1}, \cdots,f_{i_k}\}$ of $S$ and $\lambda_{i}\in Z $,
 $1\leq i\leq k$
 if $\lambda_{1}f_{i_1}+\lambda_{2}f_{i_2}+\cdots+\lambda_kf_{i_k}=0,$ then
$\lambda_{1}=\lambda_{2}=\cdots =\lambda_{k} =0$
\item
Show that any abelian group is homomorphic image of a free abelian group.
\item
Prove that any abelian group is isomorphic to a subgroup of a divisible
group.
\item
 Prove that any non-identity element of a free abelian group has infinite
 order, but, the converse is not true.
 \item
 Let $A,B,C$ be three finitely generated abelian groups. Use structure
 theorem of finitely generated abelian groups to prove: \\ (i) $A\bigoplus
 B \cong A\bigoplus C$ implies $B\cong C$ \\ (ii) $A\bigoplus A \cong
 B\bigoplus B$ implies $A\cong B$ \\ (iii) If $A\bigoplus B$ is a free
 abelian group, then $A$ and $B$ are free.
 \item
 Let $A$ be an abelian group, and $B$, a subgroup of $A$. Prove that if
 $A$ is generated by $n$ elements then $B$ can be generated by $m$
 elements where $m\leq n.$
 \item
 Find upto isomorphism all abelian groups of order 48,56, and 120.
 \item
 Let $A(\neq 0)$ be a finite abelian group. Prove that if $A$ is not
 cyclic then there exists a prime $p$ such that $p^2$ divides $\circ(A).$
 \item
 Prove that upto isomorphism the number of cyclic subgroups of the group
 $H=\Z_6\bigoplus \Z_{12}$ is 5. Further, find all distinct cyclic subgroups
 of $H$.
 \item
 Let $m\geq 1, n\geq 1$ be two integers. Prove that $\Z_m\bigoplus \Z_n$
 is isomorphic to $\Z_d\bigoplus \Z_l$ where $d=g.c.d.(m,n)$ and
 $l=l.c.m.(m,n).$
 \item
 Let $G$ be a finite abelian group. Prove that for any subgroup $H$ of
 $G$, there exists a subgroup $K$ of $G$ such that $K$ is isomorphic to
 $G/H.$
 \item
 Let $p$ be a prime. Prove that for the abelian group $G=\Z_p\bigoplus
 \cdots \bigoplus \Z_p$ ($n$-copies), the group of automorphisms $G$ i.e.,
 $Aut G$ is isomorphic to $Gl_n(\Z_p).$
 \item
 Let $p$ be a prime. Prove that for the additive abelian group
 $$A=\Z_p\bigoplus \Z_p\bigoplus \cdots \bigoplus \Z_p \mbox{ ($n$-times)
 }$$ the number of distinct subgroups of order $p^m\,(m\leq n)$ is:
 \[\frac{(p^n-1)(p^n-p)\cdots (p^n-p^{m-1})}{(p^m-1)(p^m-p)\cdots
 (p^m-p^{m-1})}.\]
 \item
 Let $A$ be a finitely generated infinite abelian group. Prove that if for
 every non-identity subgroup $H$ of $A$, $A/H$ is finite cyclic, then $A$
 is infinite cyclic.
 \item
 Determine the number of non-isomorphic abelian groups of order 56.
 \item
 Prove that there are exactly eleven abelian groups of order 729 upto
 isomorphism.
 \item
 Let $G$ be a finite abelian group. Prove that $Aut\,G$ is abelian if and
 only if $G$ is cyclic.
 \item
 Prove that any finite abelian group $A$ is Frattini subgroup of an
 abelian group $B$.
\end{enumerate}

\chapter{Chains Of Subgroups}
 \footnote{\it contents group12.tex }
 In this chapter, we shall study some special chains of subgroups of a
 group.
\bd\label{d10.1} Let $G$ be a
group. Then\\
 (a) A chain of subgroups
 $$G=G_0\rhd G_1\rhd G_2 \rhd \cdots\rhd G_{n-1}\rhd G_{n}=\{e\}$$
 (repetitions allowed) is called a subnormal
chain of the group $G$,and the factor groups $G_i/G_{i+1}$, $i=0,1,\ldots
,n-1$, are called the factors of the subnormal chain. Further, if $G_i \neq
G_{i+1} $ for all $0 \leq  i \leq n-1$, then n is defined as the length of
the subnormal chain.\\ (b) A subnormal chain $$G=H_0\rhd H_1\rhd H_2\rhd
\cdots\rhd H_{m-1}\rhd H_{m}=\{e\}$$ of $G$ with length m, i.e., $H_i \neq
H_{i+1}$, for $0 \leq i
 \leq m-1$, is called a composition series of $G$ if all the factors
 of the chain are simple groups i.e., $H_i/H_{i+1}$ is a simple group
 for all $i \geq 0$.
 \ed
 \bex\label{e10.2}
1. For the symmetric group $S_n,$ $S_n \rhd A_n \rhd {e}$ is a subnormal
chain of length 2 for all $n\geq 3$.\\ 2. For an abelian group $G$, any
chain of subgroups of $G$ starting from $G$
 and ending at the identity subgroup is subnormal chain of $G$. Thus
 $$ \Z \supset n\Z \supset m\Z \rhd \{0\}$$
 is a subnormal chain of the group $(\Z,+)$ if and only if $n$ divides
 $m$.\\
3. For the group $Gl_n(\R),$ $Gl_n(\R) \rhd Sl_n(\R) \rhd {id}$ is a
subnormal chain.
\eex
\bl\label{l10.3}
Any finite group $G (\neq \id)$ has a composition series.
\el
\pf Let $G$ be a finite group with $o(G) = n > 1$. Choose a maximal normal
subgroup $H_1$ of $G$ i.e., $H_1 \lhd G$, $H_1 \neq G$ and whenever $K
\lhd G$ and $H_1\subset K$, then either $K = G$ or $K = H_1$. Clearly
$o(H_1) < o(G)$. Next choose a maximal normal subgroup of $H_2$ of $H_1$,
and continue. As $o(G) < \infty$, this process will terminate into
identity subgroup. If $H_n = {e}$, then $$G=G_0\rhd G_1\rhd G_2 \rhd
\cdots\rhd G_{n-1}\rhd G_{n}=\{e\}$$ is a composition series of $G$ since
$H_i/H{i+1}$ is a simple group for all $0\leq i \leq n-1$ (Lemma
\ref{l301a}).     $\blacksquare$ \\ {\bf Note:} Try to find an infinite
group which has a composition series.
\bexe\label{e10.4}
Show that an infinite abelian group has no composition series.
\eexe
\bd\label{d10.5}
Let $G$ be a group and let
\be\label{e10.6}
G=G_0\rhd G_1\rhd G_2 \rhd \cdots\rhd G_{n-1}\rhd G_{n}=\{e\}
\ee
and
\be\label{e10.7}
G=H_0\rhd H_1\rhd H_2 \rhd \cdots\rhd H_{m-1}\rhd G_{m}=\{e\}
\ee
be two subnormal chains of $G$ such that $G_i \neq G_j$, $H_i\neq H_j$
for all $i\neq j$. The chains (\ref{e10.6}) and (\ref{e10.7}) are said to be equivalent if
$n=m$ and there exists a permutation $\sigma$ of $\{0,1,\cdots n-1\}$ such
that
$$ G_i/G_{i+1} \simeq H_{\sigma(i)}\ H_{\sigma(i)+1}$$
for all $i\geq 0$.
\ed
\bx\label{e10.8}
Let $G = C_6 = gp\{a\}$ be a cyclic group of order 6. Then for the subnormal
chains $$G = G_0 \rhd G_1=gp\{a^2\} \rhd G_2 = \{e\}$$ and $$G = H_0 \rhd
H_1 = gp\{a^3\} \rhd H_2 = \{e\}$$ $G/G_1 \simeq H_1/H_2$ and $G_1/G_2\simeq
H_0/H_1$. ( Use: any two cyclic groups of same order are isomorphic).
Hence the two chains are equivalent.
\ex
\bd\label{d10.9}
Given two subnormal chains
\be\label{e10.10}
G=G_0\rhd G_1\rhd G_2 \rhd \cdots\rhd G_{n-1}\rhd G_{n}=\{e\}
\ee
and
\be\label{e10.11}
G=H_0\rhd H_1\rhd H_2 \rhd \cdots\rhd H_{m-1}\rhd G_{m}=\{e\}
\ee
of a group $G$, the chain (\ref{e10.11}) is called refinement of the chain
 (\ref{e10.10}) if
each $G_i$ is equal to some $H_j$ i.e., every member of the chain
 (\ref{e10.10})is a member of the chain (\ref{e10.11}).
 \ed
 \brs\label{r10.13}
 (i) We shall say that (\ref{e10.11}) is a proper refinement of
 (\ref{e10.10}) if (\ref{e10.11}) is a refinement of
 (\ref{e10.10}) and for some $0 \leq j_0 \leq m$, $H_{j_{0}} \neq
 G_{i}$ for any $i$.\\
 (ii) By the definition of a composition series, a composition
 series can not have a proper refinement.
 \ers
 \bex \label{e10.12}
 1. Let $G=S_n$, $n \geq 3$. Then
 $$G=S_n \rhd A_n \rhd \{e\}$$
 is a refinement of $S_n \rhd \{e\}.$\\
 2. For the abelian group $\Z=(\Z,+)$, the subnormal chain
 $$\Z \supset 2\Z \supset 4\Z \supset 8\Z \supset \{0\}$$ of $\Z$
 is a refinement of the subnormal chain
 $$\Z \supset 4\Z \supset \{0\}$$
 \eex
 \bt\label{t10.121}
 Let $A,B,H$ be subgroups of a group $G$ where $H \lhd G$ and $B
 \lhd A$. Then $BH$ is a normal subgroup of $AH$ and
 $$AH/BH \simeq A/B(H\cap A).$$
 In particular if $B=\{e\}$, then
 $$ AH/H \simeq A/(H\cap A).$$
 \et
 \pf Since $H \lhd G$, $AH = HA$ and $BH = HB$. Hence, by
  theorem \ref{t29}, $AH$ and $BH$ are subgroups of $G$. Clearly
 $BH \subset AH$. Further, for any $a \in A$, $h \in H$
 \ba{rcl}
 (ah)BH(ah)^{-1} & = & ahBHh^{-1}a^{-1}\\
                 & = & ahBHa^{-1}\\
                 & = & ahHBa^{-1}\\
                 & = & aHBa^{-1}\\
                 & = & aBa^{-1}H \quad (H\lhd G)\\
                 & = & BH \quad {\mbox{since} B \lhd A}
 \ea
 Hence $BH \lhd AH.$ Now, as $A$ is a subgroup of $AH$, the map
 \bea
 \phi : A & \to & AH/BH\\
 x & \mapsto & xBH
 \eea
 is a homomorphism. For any $a \in A$, $h \in H$, we have
 \ba{rcl}
 ahBH & = & ahHB\\
      & = & aHB\\
      & = & aBH
 \ea
 Hence $\phi$ is onto. Further, for any $x\in A$,
 \ba{rl}
             & \phi (x) = BH\\
 \Rightarrow & xBH = BH\\
 \Rightarrow & x = bh_1 \mbox{for some } b\in B, h_1 \in H\\
 \Rightarrow & h_1 = b^{-1}x \in A\cap H\\
 \Rightarrow & x \in B(A\cap H)
 \ea
 Conversely if $x \in B(A\cap H)$, then we can write $x=bh'$ for
 some $b\in B$ and $h' \in A\cap H$, moreover, we have
 \ba{rclll}
 xBH & = & bh'BH & = & bh'HB\\
     &   &       & = & bHB\\
     &  &       & = & bBH\\
     &  &       & = & BH
 \ea
 Therefore $\Ker \theta = B(A\cap H)$, and by  the I
 isomorphism theorem (Theorem \ref{t33}),
 $$AH/BH \simeq A/B(H\cap A)$$
 In case $B=\{e\}$, it is immediate that
$$AH/H \simeq A/(H\cap A)       \blacksquare$$
 \bl\label{l10.14}
   Let
$$G=G_0\rhd G_1\rhd G_2 \rhd \cdots\rhd G_{n-1}\rhd G_{n}=\{e\}$$
be a subnormal chain of a group $G$. Then\\
 (a) For any subgroup $H$ of $G$,
$$H=H_0\rhd H_1\rhd H_2 \rhd \cdots\rhd H_{n-1}\rhd H_{n}=\{e\}$$ is a
subnormal chain of $H$ where $H_i = H\cap G_i$.\\ (b) If for a group $A$,
$\phi
: G \to A$ is an onto homomorphism, then for $A_i = \phi (G_i)$, $i \geq
0$, $$A=A_0\rhd A_1\rhd A_2 \rhd \cdots\rhd A_{n-1}\rhd A_{n}=\{e\}$$ is a
subnormal chain of $A$.
\el
\pf (a) As $G_{i+1} \lhd G_i$ for all $i\geq 0$, $$ H_{i+1} = H \cap
G_{i+1} \lhd H_i = H\cap G_i.$$ Hence, clearly $$H=H_0\rhd H_1\rhd H_2 \rhd
\cdots\rhd H_{n-1}\rhd H_{n}=\{e\}$$ is a subnormal chain of $H$.\\ (b) Since
$\phi$ is onto homomorphism $\phi (G_{i+1}) \lhd \phi (G_i)$ for all
$i\geq 0$ (Theorem \ref{t301a}(i)). Hence the proof is clear.
$\blacksquare$
\bco\label{c10.15}
Let $H$ be a normal subgroup of a group $G$. If
$$G=A_0\rhd A_1\rhd A_2 \rhd \cdots\rhd A_{m-1}\rhd A_{m}=\{e\}$$
is a subnormal chain of $G$, then
$$G/H\rhd A_1H/H\rhd A_2H/H \rhd \cdots\rhd A_{m-1}H/H\rhd H$$
is a subnormal chain of $G/H$.
\eco
\pf
Let $$\eta : G \rightarrow G/H$$
be the natural homomorphism. Then
\ba{rcl}
 \eta (A_i) & = & \{xH \mid x\in A_i\}\\
            & = & \{ xhH \mid x\in A ,h\in H\}\\
            & = & A_iH/H
\ea for all $i\geq 0$. Hence, by the lemma, the result follows.
$\blacksquare$
\bl\label{l10.16}
If a group $G$ has a composition series of length $n$, then\\
(a) Every normal subgroup $H$ of $G$ has a composition series of
length $\leq n$.\\
(b) Any homomorphic image of $G$ has a composition series of
length $\leq n$.
\el
\pf Let
\be\label{e10.17}
G=G_0\rhd G_1\rhd G_2 \rhd \cdots\rhd G_{n-1}\rhd G_{n}=\{e\}
\ee
    be a composition series of $G$. We, now, prove: \\
(a) By the lemma \ref{l10.14},
\be\label{e10.18}
H=H_0\rhd H_1\rhd H_2 \rhd \cdots\rhd H_{n-1}\rhd H_{n}=\{e\}
\ee
is a subnormal chain of $H$ where $H_i=H\cap G_i$. Since $H\lhd G$, $H_i =
H\cap G_i \lhd G_i$ for all $i\geq 0$. Consider the natural homomorphism:
\ba{rcl}
    H_i/H_{i+1} & \rightarrow & G_i/G_{i+1}\\
    xH_{i+1}   & \mapsto & xG_{i+1}
\ea This is a monomorphism and as $H_i \lhd G_i$, the image of
$H_i/H_{i+1}$ is a normal subgroup of $G_i/G_{i+1}$. Since $G_i/G_{i+1}$
is a simple group, $im (H_i/H_{i+1}) = G_i/G_{i+1}$ or ${id}$. i.e.,$im(
H_i/H_{i+1}) = G_i/G_{i+1}$ or $ H_i = H_{i+1}$. Hence after deleting
repetitions, if any, in (\ref{e10.18}), we get a composition series of $H$
with length $\leq n$.\\ (b) By lemma \ref{l10.14} if $\phi : G \to A$ is
an onto homomorphism, then
\be\label{e10.19}
A=A_0\rhd A_1\rhd A_2 \rhd \cdots\rhd A_{n-1}\rhd
A_{n}=\{e\}
\ee
is a subnormal chain of $A$ where $A_i = \phi(G_i)$. Consider the natural
map: \ba{rcl} G_i/G_{i+1} &\to  A_i/A_{i+1}\\ xG_{i+1} & \mapsto & \phi
(x)A_{i+1} \ea This is an onto homomorphism. As $G_i/G_{i+1}$ is a simple
group, by the lemma \ref{l302a}, either $A_i = A_{i+1}$ or $G_i/G_{i+1}
\simeq A_i/A_{i+1}$. Hence, after deleting repetitions if any in
(\ref{e10.19}), we get a composition series of $A$ with length $\leq n$.
$\blacksquare$
\bl\label{l10.20}
If for a normal subgroup $H$ of a group $G$, $H$ and $G/H$ admit
composition series of length $k$ and $l$ respectively, then $G$
admits a composition series of length $k+l$.
\el
\pf Let
\be\label{e10.21}
H=H_0\rhd H_1\rhd \cdots\rhd H_{k-1}\rhd H_{k}=\{e\}
\ee
and
\be\label{e10.22}
G/H=A_0/H \rhd A_1/H \rhd A_2/H \rhd \cdots\rhd A_{l-1}/H \rhd H
\ee
be composition series of $H$ and $G/H$ respectively. Then, using corollary
\ref{c301a}, $$G=A_0\rhd A_1\rhd A_2 \rhd \cdots\rhd A_{l-1}\rhd H\rhd
H_1\rhd$$
\be\label{e10.23}
\cdots\rhd H_{k-1}\rhd H_{k}=\{e\}
\ee
is a subnormal chain of $G$.We have $$(A_i/H)/(A_{i+1}/H)\simeq A_i/A_{i+1}$$
for all $i\geq 0$. Thus $A_i/A_{i+1}$ is a simple group for all $i\geq 0$.
Further, as $H_j/H_{j+1}$ is a simple group for all $j\geq 0$, (\ref{e10.23})
is a composition series of $G$ with length $l+k$. Hence the result.
$\blacksquare$
\bt\label{t10.24} ({\bf Schreier})
Any two subnormal chains of a group $G$ have equivalent refinements.
\et
\pf
Let
\be\label{e10.25}
G=G_0\rhd G_1\rhd G_2 \rhd \cdots\rhd G_{n-1}\rhd G_{n}=\{e\}
\ee
and
\be\label{e10.26}
G=H_0\rhd H_1\rhd H_2 \rhd \cdots\rhd H_{m-1}\rhd H_{m}=\{e\}
\ee
be two subnormal chains of a group $G$. Put $$G_{ij} = (G_i \cap
H_j)G_{i+1}$$ and $$H_{ji} = (H_j \cap G_i)H_{j+1}$$ for all $i=0,1,
\ldots , n$ and $j=0,1, \ldots , m$ assuming that $G_{n+1} = H_{m+1} =
\{e\}$. We have \ba{rcl} G_{i_0} & = & (G_i \cap H_0)G_{i+1}\\ & = &(G_i
\cap G)G_{i+1}= G_i \ea and \ba{rcl} G_{i_m} & = & (G_i \cap H_m)G_{i+1}\\
& = &(G_i \cap \{e\})G_{i+1}= G_{i+1} \ea Further taking $G= G_i$,
$H=G_{i+1}$, $A=G_i \cap H_j$ and $B=G_i \cap H_{j+1}$ in the theorem
\ref{t10.121}, we get
\be\label{e10.27}
G_{i,j+1} = (G_i \cap H_{j+1}) G_{i+1} \lhd (G_i \cap H_{j}) G_{i+1}
=G_{ij}
\ee
\be\label{e10.28}
\frac{G_{ij}}{G_{i,j+1}} \simeq \frac{G_i \cap H_j}{(G_i\cap
H_{j+1})(G_{i+1}\cap H_j)}
\ee
Similarly, using the subnormal chain \ref{e10.26}, we get
\be\label{e10.29}
H_{ji+1}\lhd H_{ji}
\ee
and
\be\label{e10.30}
\frac{H_{ji}}{H_{j,i+1}} \simeq \frac{G_i \cap H_j}{(G_i\cap
H_{j+1})(G_{i+1}\cap H_j)}
\ee
Thus using (\ref{e10.27}), (\ref{e10.28}), (\ref{e10.29}) and
(\ref{e10.30}) we conclude that
\be\label{e10.31}
G = G_{00} \rhd G_{01} \rhd \cdots \rhd G_{0m}=G_{10}\rhd G_{11}\rhd
\cdots \rhd G_{nm}=\{e\}
\ee
and
\be\label{e10.32}
G = H_{00} \rhd H_{01} \rhd \cdots H_{0n}=H_{10}\rhd H_{11}\rhd
\cdots \rhd H_{mn}=\{e\}
\ee
are subnormal chains for $G$ which are refinements of the subnormal chains
(\ref{e10.25}) and (\ref{e10.26})  respectively and are equivalent to each
other. Hence the result follows.    $\blacksquare$
\bco\label{c10.33} ({\bf Jordan- Holder Theorem})
Any two composition series of a group are equivalent.
\eco
\pf By definition, a composition series has no proper refinement. Further,
by the theorem, any two normal chains have equivalent refinements. Hence,
 any two composition series are equivalent.      $\blacksquare$
\br
If group $G$ has a composition series then any two composition
series of $G$ have same length.
\er
We give below an independent proof of the corollary \ref{c10.33}
in case $G$ is finite.
\bt\label{t10.34}
Let $G$ be a finite group. Then any two composition series of $G$
are equivalent.
\et
\pf Let
\be\label{e10.35}
G=G_0\rhd G_1\rhd G_2 \rhd \cdots\rhd G_{n-1}\rhd G_{n}=\{e\}
\ee
and
\be\label{e10.36}
G=H_0\rhd H_1\rhd H_2 \rhd \cdots\rhd H_{m-1}\rhd H_{m}=\{e\}
\ee
be two composition series  of $G$.\\
We shall prove the result by induction on $n$. If $n=1$, then $G$
is simple and the result is clear in this case. Now, let $n>1$. We
shall consider two cases.
Case I. $H_1 = G_1$ \\
In this case
\be\label{e10.37}
 G_1\rhd G_2 \rhd \cdots\rhd G_{n-1}\rhd G_{n}=\{e\}
\ee and \be\label{e10.38} G_1\rhd H_2 \rhd \cdots\rhd H_{n-1}\rhd
H_{n}=\{e\} \ee are two composition series of $G_1$, and length of
(\ref{e10.37}) is $n-1$. Hence, by induction, the subnormal chain (\ref{e10.37}) of
the group $G_1$ is equivalent to the subnormal chain (\ref{e10.38}) and
consequently the subnormal chain (\ref{e10.35}) of the group $G$ is equivalent
to the subnormal chain (\ref{e10.36}).\\ Case II\;\; $H_1 \neq G_1$\\ We have
$H_1\lhd G$, and $G_1 \lhd G$, hence $G_1 H_1\lhd G$. As (\ref{e10.35}) and
(\ref{e10.36}) are composition series of $G$, $G_1$ and $H_1$ are maximal
normal subgroups of $G$. Clearly $H_1 \subset {G_1}{H_1}$, $G_1 \subset
{G_1}{H_1}$. Hence if ${G_1}{H_1} \neq G$, ${G_1}{H_1} = G_1 = H_1$, which
is not true. Therefore ${G_1}{H_1} = G$. Put $D = H_1 \cap G_1$. By lemma
\ref{e10.10}, $D$ has a composition series (say) \be\label{e10.39}
D=D_0\rhd D_1\rhd D_2 \rhd \cdots\rhd D_{t-1}\rhd D_{T}=\{e\} \ee unless
$D = \{e\}.$ Now, consider the normal chains \be\label{e10.40} G=G_0\rhd
G_1\rhd D \rhd D_1\cdots\rhd D_{t-1}\rhd D_{t}=\{e\} \ee and
\be\label{e10.41} G=H_0\rhd H_1\rhd D \rhd D_1 \cdots\rhd
D_{t-1}\rhd D_{t}=\{e\} \ee of $G$. We have $$\frac{G}{G_1} = \frac
{{G_1}{H_1}}{G_1} \simeq \frac {H_1}{H_1\cap G_1} = \frac{H_1}{D}$$ and
$$\frac{G}{H_1} = \frac {{G_1}{H_1}}{H_1} \simeq \frac {G_1}{H_1\cap G_1} =
\frac{G_1}{D}$$ Hence $\frac{G_1}{D}$ and $\frac{H_1}{D}$ are simple
groups. Therefore, as (\ref{e10.39}) is a composition series,
(\ref{e10.40}) and (\ref{e10.41}) are composition series of $G$ and are
equivalent to each other. However, by the case I, the composition series
(\ref{e10.40}) and (\ref{e10.41}) are equivalent to the composition series
(\ref{e10.35}) and (\ref{e10.36}) respectively. Hence the result follows.
$\blacksquare$
\bd\label{d10.42} If a group $G$ has a composition series of length $n$,
then $n$ is defined as length of $G$ and we write $l(G)= n$. If $G =
\{e\}$, then length of $G$ is defined to be zero. Further, if $G\neq \{e\}$ has
no composition series, then we define length of $G$ to be $\infty$ i.e.,
$l(G)= \infty$. \ed
\bd\label{d10.43} Let $G$ be a group. A subgroup $H$ of $G$ is
called sub-normal if $H$ is a member of a subnormal chain of $G$. \ed
\bl\label{l10.44} Let $G$ be a group which admits a composition
series. Then any sub-normal subgroup $H \neq \{e\}$ of $G$ admits a
composition series and $l(H) \leq l(G).$ \el \pf As $H$ is sub-normal
subgroup of $G$, there exists a subnormal chain of $G$ with $H$ as it's
member. By Schreier's theorem the subnormal chain admits a refinement
equivalent to any given composition series of $G$. Hence the result is
immediate.   $\blacksquare$ \bl\label{l10.45} Let $G$ be a group which
admits a composition series and let $l(G)= n$. Then for any normal
subgroup $H$ of $G$,
 $l(G/H) + l(H) = l(G).$
 \el
 \pf If $H = \{e\}$, the result is clear. Hence let $H\neq \{e\}$.
 Let
\be\label{e10.46}
G=G_0\rhd G_1\rhd G_2 \rhd \cdots\rhd G_{n-1}\rhd G_{n}=\{e\}
\ee
be a composition series of $G$. Since $H\lhd G$, $G\rhd H \rhd \{e\}$ is a
subnormal chain of $G$. By Schreier's theorem the subnormal chain  $G\rhd H \rhd
\{e\}$ admits a refinement equivalent to the composition series
(\ref{e10.46}) of $G$. Let this be the subnormal chain :
\be\label{e10.47}
G=A_0\rhd A_1\rhd \cdots\rhd A_t=H\rhd A_{t+1}\rhd \cdots \rhd A_n=\{e\}
\ee
Clearly
\be\label{e10.48}
G/H=A_0/H\rhd A_1/H \rhd \cdots\rhd A_{t-1}/H\rhd H
\ee
is a subnormal chain of $G/H$ where
$$\frac{A_i/H}{A_{i+1}/H} \simeq \frac{A_i}{A_{i+1}}$$
for all $0\leq i \leq t-2$. Thus (\ref{e10.48}) is a composition
series of $G/H$. Further, as (\ref{e10.47}) is a composition
series of $G$.
\be\label{e10.49}
H\rhd A_{t+1}\rhd  \cdots\rhd A_n=\{e\}
\ee
is a composition series of $H$. Hence $$l(G) = n = (n-t)+t = l(H) + l(G/H)
\blacksquare$$
\bco\label{c10.50}
Let $G$ be a group with $l(G) = n< \infty$. Then any homomorphic
image of $G$ has length $m\leq n$.
\eco
\pf If $n=0$, then $G=\{e\}$, and the result is clear. Hence, let $n\geq
1$,
and $$\phi : G \to A$$ be an onto homomorphism. Let $K=\Ker \phi$. Then by
the corollary \ref{c43a}, $G/K$ is isomorphic to $A$. Hence the result
follows from the lemma.       $\blacksquare$ \\ {\bf NOTE:} We have not
proved the natural result that for any subgroup $H$ of a group $G$, $l(H)
\leq l(G)$. Analyse this.  \\ Apart from the concepts of subnormal chain and composition
series of a group studied above, the following definitions are also of
interest.
 \bd\label{d12.24dd} Let $G$ be a group. Then \\ (a) A subnormal
chain $$G=G_0\rhd G_1\rhd \cdots\rhd G_{n-1}\rhd G_n=\{e\}$$ of $G$ is
called a normal chain if $G_i\lhd G$ for all $i\geq 0.$ Further, if
$G_i\neq G_{i+1}$ for any $i$, then $n$ is called the length of the normal
chain. The factor groups $G_i/G_{i+1}$ are called factors of the normal
chain. \\ (b) A normal chain $$G=G_0\supsetneqq G_1\supsetneqq \cdots
\supsetneqq G_{n-1}\supsetneqq G_n=\{e\}$$ of $G$ is called a chief series
or principal series of length $n$ if $G_{i+1}$ is maximal in $G_i$ for all
$i\geq 0.$\ed  \bd\label{d12.25d} Two normal chains of a group are called
equivalent if those are equivalent subnormal chains.\ed \brs (i) A chief
series of a group $G$ admits no proper refinement into a normal chain.
\\ (ii) The factors of
a chief series of a group are not necessarily simple groups. However, if
the group is abelian, then all the factors of a chief series are simple,
hence are cyclic of prime order.\\  Following the proof of corollary
\ref{c10.33}, it is easy to prove :\ers \bexe Any two principal series of
a group $G$ are equivalent.\eexe  \bt\label{t12.27t} Let
\[G=G_0\supsetneqq G_1\supsetneqq G_2\supsetneqq\cdots \supsetneqq
G_{n-1}\supsetneqq G_n=\{e\}\] be a chief series of a finite group $G$
with length $n\geq 1.$ Then $G_{n-1}$ is isomorphic to a subgroup of a
finite direct product of simple groups all of which are isomorphic to each
other. \et  \pf If $G_{n-1}$ is simple, there is nothing to prove. Let
$G_{n-1}$ be non-simple and let $K$ be a maximal normal subgroup of
$G_{n-1}$. Then $G_{n-1}/K$ is a simple group. Let $K=K_1, K_2,\cdots
,K_t$ be all distinct conjugates of $K$ in $G$. Clearly, $G_{n-1}$ is a
minimal normal subgroup of $G$. Hence $K$ is not normal in $G$. Therefore
$t>1.$ Note that $H=\cap_{i=1}^{t}K_i$ is a normal subgroup of $G$
contained in $K \subsetneqq G_{n-1}.$ Therefore $H=\{e\}.$ If
$K_i=x_iKx_{i}^{-1}$ for some $x_i\in G$, then the inner automorphism:
\[\begin{array}{rcl}
  \tau_{x_i}:G & \rightarrow & G \\
  g & \mapsto & x_igx_{i}^{-1}
\end{array}\] induces the isomorphism. \[\begin{array}{rcl}
  \overline{\tau_{x_i}}:G_{n-1}/K & \rightarrow & G_{n-1}/K_i \\
  gK & \mapsto & \tau_{x_i}(g)K_i
\end{array}\] Therefore each $G_{n-1}/K_i$ is a simple group isomorphic to
$G/K$. As $H=\{e\},$ it is easy to see that the map: \[\begin{array}{rcl}
  G_{n-1} & \rightarrow & \prod_{i=1}^{t}G_{n-1}/K_i \\
  x & \mapsto & (xK_1,xK_2,\cdots,xK_t)
\end{array}\] is a monomorphism. Hence the result follows.
$\blacksquare$ \bco\label{12.28cc} In the theorem if $G$ has a composition
series with abelian factors, then $G_{n-1}$ is an abelian group of order
$p^{\alpha}$ ($p$: prime, $\alpha\geq 1$) where $\circ(x)=p$ for all
$x(\neq e)$ in $G_{n-1}.$\eco \pf By theorem \ref{t10.24}, given a
composition series of $G$, the subnormal chain \be\label{111.}G=G_0\supsetneqq
G_1\supsetneqq \cdots\supsetneqq G_{n-1}\supsetneqq K\supsetneqq \{e\}\ee
has a refinement equivalent to the composition series of $G$. As $K$ is a
maximal normal subgroup of $G_{n-1}$, $G_{n-1}/K$ is a factor group of any
refinement of the subnormal chain (\ref{111.}). Further, as $G$ has a
composition series with abelian factor groups, by theorem \ref{t10.24}, any
composition series of $G$ has abelian factor groups. Therefore $G_{n-1}/K$
is abelian. Now, as $G_{n-1}/K$ is abelian simple group, it is cyclic of
order $p$ for a prime $p$. The result, now, follows from the proof of the
theorem.         $\blacksquare$   \vspace{.25in}
\\ {\bf EXERCISES}
\begin{enumerate}
  \item
  Find the length of the symmetric group $S_n$ for $n\geq 2.$
  \item
  Write composition series for the groups $C_{124}, D_{8},$ and $S_4$.
  \item
  Give two non-isomorphic groups with equal length.
  \item
  Give two non-isomorphic, non-abelian, groups with equal length.
  \item
  Find length of the group $Gl_2(\Z_3).$
  \item
  Let $G$ be a finite $p$-group with $\circ(G)=p^n, n\geq 1.$ Prove that
  $G$ has a composition series of length $n$.
  \item
  Let $C_n$ be a cyclic group of order $n$. Prove that if
  $n=p_{1}^{\alpha_1}p_{2}^{\alpha_2}\cdots p_{k}^{\alpha_k}$ where
  $\alpha_i\geq 1$ and $p_1,p_2,\cdots,p_k$ are distinct primes, then
  $l(C_n)=\sum_{i=1}^n \alpha_i.$
  \item
  Prove that if a finite group $G$ has a subnormal chain with all its factor
  groups abelian, then $G$ has a subnormal chain with all factor groups
  cyclic.
  \item
  Show that a group is finite if and only if it has a subnormal chain with
  all factor groups finite.
  \item
  Prove that the group $G=\bigcup_{n\geq 1}S_n$ (Exercise 23, Chapter 9)
  has a composition series, but has a subgroup $H$ which has no
  composition series nor a principal series. \\
   If $H$ is a subnormal subgroup of a group $G$, then the smallest $t$
  for which there exists a chain \[G=H_0 \rhdn H_1\rhdn \cdots \rhdn
  H_{t-1}\rhdn H_t=H\] of subgroups of $G$ is denoted by $m(G,H)$.
  \item
  Let $H$ be a subnormal subgroup of a group $G$. For $i\geq 0$, define $H_0=G$, and
  $H_i$, the normal closure of $H$ in $H_{i-1}$ for all $i\geq 1.$ Prove
  that \[G=H_0 \rhdn H_1\rhdn \cdots \rhdn
  H_{m-1}\rhdn H_m=H\] where $m=m(G,H).$ Further, show that if $\alpha\in
  Aut\,G$ such that $\alpha(H)=H$, then $\alpha(H_i)=H_i$ for all $i\geq
  0.$
  \item
  Prove that for any subnormal subgroup $H$ of a group $G$ and $\alpha\in
  Aut\,G$, $\alpha(H)$ is subnormal in $G$ and $m(G,H)=m(G,\alpha(H)).$
  \item
  Let $A,B$ be two subnormal subgroups of a group $G$ such that $A\subset
  N_G(B).$ Prove that $AB$ is a subnormal subgroup of $G$. \\ (Hint : use
  induction on $m(G,B)$)
  \item
  Let $H$ be a subgroup of a group $G$. Let $\mathcal{S}=\{A\subset H\st
  A:\mbox{ subnormal subgroup of }G\}$ has a maximal subgroup $M$ with
  respect to containment. Prove that $M$ is normal in $H$.
  \item
  Let the group $G$ admits a composition series. Prove that : \\ (i) Any
  non-empty set of subnormal subgroups of $G$ contains a maximal element.
  \\ (ii) Show that if $A,B$ are subnormal subgroups of $G$, then
  $gp\{A,B\}$ is a subnormal subgroup of $G$.
  \item
  Let $A,B$ be two subgroups of a group $G$ such that $A$ is subnormal in
  $G$. Prove that : \\ (i) If $m(G,A)=t$, then
  $m(\overline{N_G}(A),A)=t-1.$  \\  (ii) If $A,B$ are finite, then
  $gp\{A,B\}$ is a finite subgroup of $G$. \\ (Hint : use induction on
  $m(G,H)$)
  \item
  Let for a subgroup $A$ of a group $G$, $[G:H]<\infty$. Prove that if $A$
  commutes with all its conjugates then $A$ is subnormal in $G$.
  \item\label{1a}
  Any factor group of a chief series of a finite group $G$ is isomorphic
  to a subgroup of a finite direct product of simple groups all of which
  are isomorphic to each other.
  \item
  In the exercise \ref{1a} if we assume that $G$ has a composition series
  with abelian factor groups, then any factor group of a chief series is
  an abelian $p$-group with all its non-identity elements with order $p$
  (the prime $p$ need not be same for all factor groups).
  \item
  Let $M$ be a subnormal subgroup of a finite group $G$. Prove that if $M$
  is maximal with respect to the properties that it is subnormal in $G$
  and contains a unique Sylow-$p$ subgroup then $M\lhd G$.
  \end{enumerate}

\chapter{Solvable And Nilpotent Groups}
 \footnote{\it contents group13.tex }
 In this chapter, we introduce the concepts solvable and nilpotent
 groups using the idea of subnormal chain.
We define:
\bd\label{d10.51}
A group $G$ is called Solvable if there exists a subnormal chain $$G=G_0\rhd
G_1\rhd G_2 \rhd \cdots\rhd G_{n-1}\rhd G_{n}=\{e\}$$ of $G$ with all its
factors $G_0/G_1 ,G_1/G_2 , \ldots G_{n-2}/G_{n-1},$ and $G_{n-1}$ abelian
groups.
\ed
\bex\label{ex10.52}
1. Any abelian group is solvable.\\ 2. The symmetric groups $S_2$,$S_3$
and $S_4$ are solvable. The group $S_2$ is cyclic of order $2$, moreover,
$$S_3 \rhd A_3 \rhd \{e\}$$ is a normal chain of $S_3$ where $\frac
{S_3}{A_3}$ is cyclic of order $2$ and $A_3$ is cyclic of order $3$.
Further, $$S_4 \rhd A_4 \rhd H \rhd \{e\}$$ where $H =
\{(12)(34),(13)(24),(14)(23), e \}$. Here  $\frac {S_4}{A_4}$ is cyclic
of order $2$, $\frac {A_4}{H}$ is cyclic of order $3$ and $H$ is abelian.
\eex
\bexe\label{exe 10.53}
Show that $S_n$ is not solvable for $n >4$.
\eexe
\bd\label{d10.531} Let $G$ be a group. Then\\ (i) $G' = [G,G] =
gp\{[a,b]|a,b \in G\}$ is called the $1^{st}$ derived group of $G$. In
general, for $n\geq 2$, we define $G^{(n)} = (G^{(n-1)})'$, the $1^{st}$
derived group of $G^{n-1}$, as the $n^{th}$ derived group of $G$. The
$0^{th}$ derived group $G^{0}$ of $G$ is defined to be $G$.\\ (ii) For any
two subgroups $A$, $B$ of a group $G$, we write $$[A,B] = gp\{[x,y]|x\in A,
y \in B\}.$$ \ed \br For any two subgroups $A,B$ of a group $G$,
$[A,B]=[B,A].$\er
\bd\label{d10.532} A subgroup $H$ of a group $G$ is called a
characteristic subgroup of $G$ if $\phi (H) \subset H$ for all $\phi \in
Aut G$.
\ed
 \brs\label{r10.533}
 (i) For a group $G$, the identity subgroup and $G$ are
 characteristic subgroups of $G$.\\
 (ii) Let $H$ be a characteristic subgroup of a group $G$. Let for
 $g \in G$, $\tau_g$ denotes the inner automorphism determined by
 $g$. Then $\tau_g(H) \subset H$ i.e., $gHg^{-1} \subset H$. Hence
 $H\lhd G$.\\
 (iii) If $H$ is a characteristic subgroup of $G$ and $\phi \in
 Aut G$. Then, as $\phi^{-1} \in Aut G,$ $\phi^{-1}(H) \subset H$.
  Hence $H \subset \phi (H)$, and
 consequently $\phi (H) = H$ for all $\phi \in Aut G$.
 \ers
 \bl\label{l10.534}
 Let $K \subset H$ be subgroups of a group $G$ such that $H$ is a
 characteristic subgroup of $G$ and $K$ is a  characteristic
 subgroup of $H$. Then $K$ is a  characteristic subgroup of $G$.
 \el
\pf Let $\phi \in AutG$. As $H$ is a characteristic subgroup of $G$,
$\phi(H)=H.$ Hence $\phi \in Aut H.$
 Now,
as $K$ is a characteristic subgroup of $H$, $\phi(K) \subset K.$ Thus $K$
is a characteristic subgroup of $G$.           $\blacksquare$
 \bl\label{l10.535}
 Let $H$ be a characteristic subgroup of a group $G$ , and $K/H$
 be a characteristic subgroup of $G/H$, then $K$ is a
 characteristic subgroup of $G$.
 \el
 \pf
 Since $H$ is a characteristic subgroup of $G$, $H\lhd G$. Hence
 $G/H$ makes sense. Now, let $\alpha \in Aut G$. Then $\alpha (H) =
 H$. Therefore $\alpha$ induces the automorphism:
 \bea
 \bar{\alpha} : G/H & \to & G/H\\
 xH & \mapsto & \alpha (x) H
 \eea
 on the factor group $G/H.$ Hence $\bar{\alpha} (K/H) = K/H$, and
 consequently $\alpha (K) \subset K$.
 Thus $K$ is a characteristic subgroup of $G$.         $\blacksquare$
 \bl\label{l10.536}
 Center of a group $G$ is a characteristic subgroup of $G$.
 \el
 \pf
 Let $\alpha \in Aut G$, and $x \in Z(G)$. Then, for any $a \in G$
 \ba{rl}
             & \alpha (xa) = \alpha (ax)\\
 \Rightarrow & \alpha (x) \alpha (a) = \alpha (a) \alpha (x)\\
 \Rightarrow & \alpha (x) \in Z(G) \mbox{ since } \alpha (G) = G\\
 \Rightarrow & \alpha (Z(G)) \subset Z(G)
 \ea
 Thus $Z(G)$ is a characteristic subgroup of $G$.        $\blacksquare$
 \bl \label{l10.537}
 If $A$, $B$ are characteristic subgroups of a group $G$ then so
 is $[A,B]$.
 \el
 \pf
 Let $\alpha \in Aut\;G$. Then for any $x \in A$, $y \in B$,
 \ba{rl}
         & \alpha [x,y] = [\alpha (x), \alpha (y)] \in [A,B]\\
\Rightarrow & \alpha [A,B] \subset [A,B]
 \ea
Hence $[A,B]$ is a characteristic subgroup of $G$.       $\blacksquare$
\bl\label{l10.538} Let $K\subset H$ be two subgroups of a group
$G$, where $H\lhd G$ and $K$ is a characteristic subgroup of $H$. Then $K
\lhd G$.
\el
 \pf As $H \lhd G$, for any $x \in G$, for the
inner automorphism $\tau_x$ of $G$ determined by $x,$ $\tau_x (H) =H$.
Hence $\tau_x \in Aut H$ for all $x \in G$. As $K$ is a characteristic
subgroup of $H$, $\tau_x (K) \subset K$ i.e., $xKx^{-1} \subset K$. Hence
$K \lhd G$.            $\blacksquare$
 \bexe\label{ex10.539}
 Prove that for a subgroup $H$ of a group $G$, $H \lhd G$ if and
 only if $[G,H] \subset H$.
 \eexe
 \bl\label{l10.540}
Let $A,B,H$ be subgroups of a group $G$ such that $B\subset A$ and $[G,A]
\subset B$. Then $HB$ is a normal subgroup of $HA$
\el
 \pf
 By assumption,
 \ba{rl}
 & [G,A] \subset B \subset A\\
 \Rightarrow & A \lhd G
 \ea
 Further,
 \ba{rl}
             &  [G,B] \subset [G,A] \subset B\\
 \Rightarrow & B \lhd G
 \ea
 Thus $HB \subset HA$ are subgroups of $G$ (Lemma \ref{w}). For any
 $a \in A$, $h \in H$,
 $$ ahHBh^{-1}a^{-1}=BaHa^{-1}$$
 and
 \bea
 Baha^{-1} & = & Baha^{-1}h^{-1}h\\
           & = & B[a,h]h\\
           & = & Bh \mbox{ since }[a,h]\in [G,A] \subset B
 \eea
 Hence
 \ba{rl}
 & ahHBh^{-1}a^{-1} \subset BH=HB\\
 \Rightarrow & HB \lhd HA  \blacksquare \ea
 \bl\label{l10.54}
 All derived groups of a group $G$ are
characteristic subgroups of $G$.
\el
\pf By definition of characteristic subgroup, $G^0=G$ is a characteristic
subgroup of $G$. For $n\geq 1,$ we shall prove the result by induction on
$n$. Let $n=1$.
 If $\phi \in Aut\,G$ and $a,b\in G$, then \ba{rcl} \phi [a,b] & = &
\phi (aba^{-1}b^{-1})\\
           & = & \phi (a)\phi (b) {\phi (a)}^{-1} {\phi
           (b)}^{-1}\\
           & = & [\phi (a),\phi (b)]\\
\Rightarrow & \phi (G^{\prime}) \subset & G^{\prime} \ea Hence $G'$ is a
characteristic subgroup of $G$. Thus the result holds for $n=1$.
 Now, let $n\geq 1$ and $G^n$ is a
characteristic subgroup of $G$. By the case $n=1$, $G^{n+1} =
{(G^n)}^{\prime}$ is a characteristic subgroup of $G^n$. Thus,
using lemma \ref{l10.534}, $G^n$ is a characteristic subgroup of
$G$ for all $n\geq 1$. $\blacksquare$ \br In view of the above
lemma, for any group $G$, $$ G\rhd G^{\prime} \rhd G^2 \rhd \cdots
\rhd G^n \rhd G^{n+1}\rhd \cdots$$ is possibly an infinite chain
of normal subgroups of $G$. This is called the derived chain of
$G$. \er \bl\label{l10.55} Let $H$ be a normal subgroup of a group
$G$. Then $\frac{G}{H}$ is an abelian group if and only if
$G^{\prime}\subset H$. \el \pf If $\frac{G}{H}$ is abelian, then
for any $a,b \in G$, \ba{rcl} aH.bH & = & bH.aH\\ \Rightarrow abH
& = & baH\\ \Rightarrow Haba^{-1}b^{-1} & = & H\\ \Rightarrow
[a,b] & = & aba^{-1}b^{-1} \in H\\ \Rightarrow gp\{[a,b] \mid
a,b\in G\} & = & G^{\prime}\subset H \ea Conversely if $G^{\prime}
\subset H$, then for any $a,b \in G$, \ba{rcl} [a,b] & = &
aba^{-1}b^{-1} \in H\\ \Rightarrow Haba^{-1}b^{-1} & = & H\\
\Rightarrow Hab & = & Hba\\ \Rightarrow (Ha)(Hb) & = & (Hb)(Ha)
\ea Hence $\frac{G}{H}$ is an abelian group. $\blacksquare$
\bl\label{l10.56} A group $G$ is solvable if and only if $G^t$,
the $t^{th}$ derived
 group of $G$, is identity for some $t\geq 1$.
 \el
 \pf Let $G$ be a solvable group. Then there exists a subnormal chain
 $$G=G_0\rhd G_1\rhd G_2 \cdots\rhd G_{t-1}\rhd G_{t}=\{e\}$$
 such that $\frac{G_i}{G_{i+1}}$ is abelian for all $i\geq 0.$
 We shall show that $G_i\supset G^i$, the $i^{th}$ derived group,
 for all $i\geq 0.$ Clearly $G^0 = G = G_0$. Let for $i\geq 0$, we
 have proved that $G^i \subset G_i$. As $\frac{G_i}{G_{i+1}}$ is
 abelian, by the lemma \ref{l10.55}, ${G_i}^{\prime} \subset G_{i+1}.$
Hence $G^{i+1} = (G^i)^{\prime} \subset {G_i}^{\prime} \subset G_{i+1}.$
Therefore by induction $G^i \subset G_i$ for all $i\geq 0.$ Thus, in
particular, $G^t \subset G_t = \{e\}$. Hence $G^t = \{e\}$. Conversely, if
$G^t = \{e\}$ for some $t\geq 1$, then $$G=G^0\rhd G^1\rhd G^2 \cdots\rhd
G^{t-1}\rhd G^{t}=\{e\}$$ is a normal chain of $G$ and
$\frac{G_i}{G_{i+1}} = \frac{G_i}{(G_i)^{\prime}}$ is abelian group for
all $i\geq 0$ (Lemma \ref{l10.55}). Hence the result follows. $\blacksquare$
 \bl\label{l10.57} Let $G$ be a solvable group. Then any
subgroup as well as homomorphic image of $G$ is solvable.
 \el
 \pf
As $G$ is solvable, by lemma \ref{l10.56}, $G^t = \{e\}$ for some $t\geq
1.$ Let $H$ be a subgroup of $G$. As $H\subset G$, it is clear that$$
H^{\prime} = gp\{[x,y] \mid x,y\in H\} \subset G^{\prime}$$. Let for
$i\geq 1$, we have $H^i\subset G^i$. Then $$ H^{i+1} = {(H^i)}^1 \subset
{(G^i)}^1 = G^{i+1}.$$ Hence by induction, $H^i\subset G^i$ for all $i\geq
0.$ Now, as $G^t = \{e\}$, $H^t = \{e\}$.Thus $H$ is solvable. Next, let
$\phi : G\to A$ be an onto homomorphism. As for any $x,y\in G$, $\phi
[x,y] = [ \phi(x), \phi (y) ]$, it is clear that $\phi (G^{\prime}) =
A^{\prime}.$ Now, let for $i\geq 1$, $\phi (G^i) = A^i$. Then, as above
$$\phi (G^{i+1}) = \phi ({(G^i)}^{\prime})  = {(A^i)}^{\prime} =
A^{i+1}.$$ Hence by induction, $\phi (G^i) = A^i$ for all $i\geq 1$.
Therefore $ A^t = \phi (G^t) = \{e\}.$ Hence $A$ is solvable.
$\blacksquare$
\bexe\label{a1} Prove that $Gl_n (\R)$, and $Sl_n (\R)$ are not solvable for any
$n\geq 3.$
 \eexe
\bl\label{l10.58} A finite direct product of solvable groups is
solvable.
\el
\pf Let $G_1, G_2,\cdots G_n$ be $n$ solvable groups and let $$ G = G_1
\times G_2 \times \cdots \times G_n$$ be the direct product of $G_i's$. By
definition of derived groups it is clear that $$ G^{\prime} =
{G_1}^{\prime} \times {G_2}^{\prime} \times \cdots \times {G_n}^{\prime}$$
Using induction, one can easily see that $$ G^i = {G_1}^i \times {G_2}^i
\times \cdots \times {G_n}^i$$ for all $i\geq 0$ i.e., the $i^{th}$
derived group of $G$ is the product of $i^{th}$ derived groups of $G_l,
1\leq l\leq n.$ As each $G_l$ is solvable, by lemma \ref {l10.56}, there
exists
 $t\geq 1$ such that ${G_l}^t = \{e\}$ for all $l\geq 1.$
 Hence
$$ G^t = {G_1}^t \times {G_2}^t \times \cdots \times {G_n}^t = \{id\}$$
Therefore $G$ is solvable.      $\blacksquare$ \br From the proof of the
lemma it is clear that if $G$ is solvable then so is each $G_i$. This also
follows from the fact that each $G_i$ is a homomorphic image of $G$ (Lemma
\ref{l10.57}).
\er \bl\label{l10.59} Let $H$be a normal subgroup of a group $G$. If $H$
and $\frac{G}{H}$ are solvable then $G$ is solvable.
\el
 \pf
  Let $\eta : G\to \frac{G}{H}$ be the natural homomorphism. As seen in lemma
\ref{l10.57}, $\eta (G^i) = (\frac{G}{H})^i$ for all $i\geq 0.$ As
$\frac{G}{H}$ is solvable, there exists $l\geq 1$ such that
$(\frac{G}{H})^l = id$. Hence \ba{rcl} \eta (G^l) & = & id\\
\Rightarrow G^l & \subset & H\\ \Rightarrow {(G^l)}^i & \subset & H^i
\;\;\; \forall i\geq 0 \ea As $H$ is solvable, there exists $s\geq 1$ such
that $H^s = id$. Hence \ba{rcl} G^{l+s} = {(G^l)}^s & \subset H^s =
id\\ \Rightarrow G^{l+s} & =  id \ea Therefore $G$ is solvable.
$\blacksquare$
\bco\label{c10.60} Any finite $p$-group is solvable. \eco \pf Let $G$ be a
finite $p$-group with $o(G) = p^s$. We shall prove the result by induction
on $s.$ If $s = 0$, $G = id$. Hence $G$ is clearly solvable. Now,
assume $s\geq 1.$ By theorem \ref{t614}, $o(Z(G)) = p^m, m\geq 1.$ Now,
$Z(G) \lhd G$ and is solvable since it is abelian. If $Z(G) = G$, then the
result clearly holds. If $Z(G) \subsetneqq G$, then $o(\frac{G}{Z(G)}) =
p^{s-m} < p^s$. Hence $\frac{G}{Z(G)}$ is solvable by induction. Now, $G$
is solvable by lemma \ref{l10.59}.         $\blacksquare$
\bl\label{13.25l} A finite group $G$ is solvable if and only if $G$ has a
composition series with all its factor groups cyclic of prime orders. \el
\pf Let $G$ be solvable, then $G$ has a subnormal chain with all the
factor groups abelian.
Therefore any subnormal chain which is a refinement of this chain has
abelian factor groups. Since $G$ is finite, any subnormal chain can be
refined to a composition series. Therefore $G$ has a composition series
with abelian factor groups. As factor groups of a composition series are
simple, and a simple abelian group is cyclic of prime order, the direct
part of the lemma follows. The converse of the statement is clear by
definition of solvable group.      $\blacksquare$
 \bt\label{t13.25tt} Let $G$
be a finite group and $p,q,r$ be primes. Then \\ (i) If $\circ(G)=pq$ or
$p^2q,$ $G$ is solvable. \\ (ii) If $p>q>r$ and $\circ(G)=pqr,$ $G$ is
solvable. \et \pf (i) We can assume that $p,q$ are distinct primes, since
otherwise the result follows by corollary \ref{c10.60}. By theorem \ref{t10.20t},
 $G$ contains a proper normal subgroup $H$(say). If $\circ(G)=pq,$ then $\circ(H)=p$ or $q$.
 Therefore by corollary \ref{c10.60}, $H$ and $G/H$ are solvable, and
 hence $G$ is solvable (Lemma \ref{l10.59}). Next, let $\circ(G)=p^2q.$
 Then $\circ(H)=p$ or $q$ or $pq$ or $p^2$. Hence again $H$ and $G/H$ are
  solvable by the corollary \ref{c10.60} in case $\circ(H)=p^2$ or $q$ and
  by the corollary \ref{c10.60} and the case $\circ(G)=pq$ in case
  $\circ(H)=pq$ or $p$. Thus $G$ is solvable by lemma \ref{l10.59}. \\ (ii)
  By theorem \ref{t10.19t}, $G$ has only one Sylow $p$-subgroup $H$(say).
  Then $H$ is normal in $G$ by corollary \ref{c82} and is a solvable group
  (Corollary \ref{c10.60}). Further by (i), $G/H$ is solvable. Hence $G$
  is solvable (Lemma \ref{l10.59}).     $\blacksquare$
  \bt\label{13.27tr} Let $G$ be a finite solvable group with $\circ(G)=mn$
  where $m\geq 1, n\geq 1$ are integers such that $(m,n)=1.$ Then \\ (i)
  $G$ contains a subgroup of order $m$. \\ (ii) Any two subgroups of $G$ with  order
  $m$ are conjugate. \\ (iii) If $H$ is a subgroup of $G$ with $\circ(H)$
  dividing $m$, then $H$ is contained in a subgroup of order $m$ in
  $G$.\et \pf If $m=1$ or $n=1,$ the result is trivial. Hence assume $m>1,
  n>1.$ Let $m=p^a$ where $p$ is a prime and $a\geq 1.$ Then the result
  follows from Sylow's first and second theorems (Theorems \ref{t83} and
  \ref{t88}) and theorem \ref{t10.14t}. Thus we can assume that $m>1,
  n>1$, and $G$ is not a $p$-group. We shall use induction on $\circ(G)$
  to prove the result in this case. Let us consider the two cases:
   \\ {\bf Case I:} The group $G$ contains a proper
  normal subgroup $K$ with $\circ(K)=m_1n_1$ such that $m=m_1m_2,$
  $n=n_1n_2$ and $n_2>1.$ \\  As $(m,n)=1,$ $(m_i,n_j)=1$ for all $1\leq
  i,j\leq 2.$ Further, $G/K$ is a solvable group (Lemma \ref{l10.57}), and
  $\circ(G/K)=m_2n_2<\circ(G)$ since $\circ(K)>1.$ We shall first prove
  (i). By induction, there exists a subgroup $B/K$ of $G/K$ of order $m_2$
  where $B$ is a subgroup of $G$ containing $K$. We have $\circ(B)=mn_1<mn,$
  since $n_2>1.$ As $B$ is solvable (Lemma \ref{l10.57}), by induction,
  $B$ contains a subgroup of order $m$. Thus (i) holds. To prove (ii), let
  $M,N$ be two subgroups of $G$ of order $m$. As $K\lhd G$, $MK$ is a
  subgrop of $G$ containing $M$ and $K$. Therefore
  $\circ(M)=m$ and $\circ(K)=m_1n_1$ divide $\circ(MK).$ Consequently
  $mn_1$ divides $\circ(MK).$ Next, note that $\circ(MK)$ divides
  $\circ(M)\circ(K)=mm_1n_1$ as well as $\circ(G)=mn.$ Hence $\circ(MK)$
  divides $mn_1$. Consequently $\circ(MK)=mn_1.$ Similarly, we get
  $\circ(NK)=mn_1.$ Thus $\circ(MK/K)=m_2=\circ(NK/K).$ The factor group
  $G/K$ is solvable of order $m_2n_2<n$ where $(m_2,n_2)=1,$ and $MK/K$
  and $NK/K$ are subgroups of $G/K$ each of order $m_2.$ Hence, by
  induction, $MK/K$ and $NK/K$ are conjugate in $G/K$. Thus there exists
  an element $g\in G$ such that $gNg^{-1}K=MK.$ Now, $gNg^{-1}$ and $M$
  are two subgroups of $MK$ of order $m$, and $MK$ is solvable group of
  order $mn_1<\circ(G).$ Therefore, by induction, $gNg^{-1}$ is conjugate
  of $M$ in $MK$. Consequently $M$ and $N$ are conjugate in $G$. This
  proves (ii). To prove (iii), let $M'$ be a subgroup of $G$ with
  $\circ(M')=m',$ a divisor of $m$. Clearly $\circ(M'K/K)$ divides
  $\circ(M'),$ and hence $m=m_1m_2.$ Also $\circ(M'K/K)$ divides
  $\circ(G/K)=m_2n_2.$ Therefore $\circ(M'K/K)$ divides $m_2.$ Hence, by
  induction, $M'K/K$ is contained in a subgroup $B/K$ of $G/K$ with order
  $m_2$. We have $\circ(B)=mn_1<\circ(G)$ where $B$ is a subgroup of $G$
  containing $M'K.$ Now, $M'$ is a subgroup of $B$ with $\circ(M')=m'$
  dividing $m.$ Hence, as $B$ is solvable group with $\circ(B)<\circ(G)$,
  by induction, $M'$ is contained in a subgroup $D$ of $B$ of order $m$.
  Thus (iii) is proved. \\ We, now, consider the situation when $G$ does
  not contain any proper normal subgroup as described in case I. As $G$ is
  solvable, $G$ has a composition series with abelian factor groups (Lemma
  \ref{13.25l}). Hence, by corollary \ref{12.28cc}, the smallest normal
  subgroup $K$ in a chief series of $G$ is an abelian group of order $p^a$
  where $p$ is a prime and $a\geq 1.$ If $n\neq p^a,$ then $K$ will be a
  proper normal subgroup of $G$ as in case I. Hence assume $n=p^a.$
  Clearly, we can assume that any minimal normal subgroup of $G$ has order
  $p^a$ since otherwise we fall in case I. As $K\lhd G$, by Sylow's second
  theorem, $K$ is unique minimal normal subgroup of $G$ (use: $K$ is
  Sylow $p$-subgroup of $G$). We, now, consider the following: \\ {\bf Case
  II}: $G$ has unique minimal normal subgroup $K$ with $\circ(K)=n=p^a$
  where $p$ is a prime and $a\geq 1.$ \\ Let $L$ be a minimal normal
  subgroup of $G$ containing $K$. Then $L/K$ is a minimal normal subgroup
  of the solvable group $G/K.$ Hence, as seen above, by corollary
  \ref{12.28cc}, $L/K$ is an abelian group of order $q^b$, $b\geq 1,$
  where $q$ is a prime dividing $\circ(G/K)=m.$ Clearly $q\neq p,$ and
  $\circ(L)=p^aq^b.$ Let $Q$ be a Sylow $q$-subgroup of $L$. Then
  $L=QK.$ Put $E=N_G(Q)$, and $C=E\cap K.$ If $C=\{e\},$ then all
  conjugates of $Q$ with respect to elements of $K$ are distinct. Thus the
  number of conjugates of $Q$ in $G$ i.e., $[G:E]\geq p^a=\circ(K).$
  Further, as $L\lhd G$, the number of conjugates of $Q$ in $G$ equals to
  the number of conjugates of $Q$ in $L$ by Sylow's second theorem. Thus
  $[G:E]\leq [L:Q]=p^a.$ Consequently $[G:E]=p^a,$ and $E$ is a subgroup
  of order $m$ in $G$. Therefore if $C=id,$ then statement (i) is proved.
  Now, let $C\neq id$. We have $(\circ(K),\circ(Q))=1,$ hence $K\cap
  Q=id$. Thus as $C\subset K, C\cap Q=id.$ We have $K\lhd G$, hence $C\lhd
  E.$ Further, by the definition of normaliser, $Q\lhd E.$ Hence as $C\cap
  Q=id$, elements of $C$ and $Q$ commute (use: $[x,y]\in C\cap Q=\{e\},$
  for all $x\in C, y\in Q$). As $K$ is abelian, and $L=QK$, it is clear
  that $C\subset Z(L).$ Note that $Z(L)$ is a characteristic subgroup of
  the normal subgroup $L$ of $G$. Therefore $Z(L)\lhd G$ (Lemma
  \ref{l10.538}). As $K\lhd G, Z(L)\cap K\lhd G.$ We have $Z(L)\subset
  N_G(Q)=E,$ and $C\subset Z(L).$ Hence $E\cap K=C\subset Z(L)\cap
  K\subset K.$ As $K$ is a minimal normal subgroup of $G$ and $Z(L)\cap K$
  is a non-trivial normal subgroup of $G$, $K=Z(L)\cap K.$ Thus
  $Z(L)\supset K$ i.e., elements of $K$ and $Q$ commute. Hence $L=K\times
  Q.$ This shows that $Q\lhd L.$ Thus $Q$ is unique Sylow $q$-subgroup of
  $L$. Now, as $L\lhd G, Q\lhd G.$ Note that any minimal normal subgroup $G$
  contained in $Q$ has order $q^{\alpha}$ for some $\alpha\geq 1.$ This
  contradicts the fact that $K$ is unique minimal normal subgroup of $G$.
  Thus the case $C\neq id$ does not arise. Hence (i) is proved in case II.
  In fact $E=N_G(Q)$ is a subgroup of order $m$ in case II. Now, we shall
  prove (ii). Let $N$ be any subgroup of $G$ with order $m$. Then \[\begin{array}{rrl}
     & NL & \supset NK  \\
    \Rightarrow & \circ(NL) & \geq \circ(NK)=mp^a \\
    \Rightarrow & NL & =G \\
    \Rightarrow & G/L & \cong N/N\cap L \\
    \Rightarrow & \circ(G)/\circ(L) & =\circ(N)/\circ(N\cap L) \\
    \mbox{ i.e., } & \frac{mp^a}{p^aq^b} & =\frac{m}{\circ(N\cap L)} \\
    \Rightarrow & \circ(N\cap L) & =q^b \\
  \end{array}\] Let $Q_2$ be a Sylow $q$-subgroup of $N$ and let $Q_1$ be
  a Sylow $q$-subgroup of $G$ containing $Q$. Since $\circ(G)=mp^a, q\neq
  p,$ and $(m,p)=1$, clearly $Q_2$ is a Sylow $q$-subgroup of $G$ as well.
  Hence by Sylow's second theorem, $Q_1$ and $Q_2$ are conjugates in $G$.
  Therefore $Q_2$ contains a conjugate of $Q$. As $L\lhd G,$ $N\cap L\lhd
  N.$ If for $g\in G, g^{-1}Qg\subset Q_2$, then $g^{-1}Qg \subset N.$
  Also $g^{-1}Qg\subset L.$ Therefore \[g^{-1}Qg\subset N\cap L\]
  \[\Rightarrow g^{-1}Qg=N\cap L \mbox{ since }\circ(Q)=\circ(N\cap
  L)=q^b\] As $N\cap L\lhd N,$ we have \[\begin{array}{crl}
     & N \subset  N_G(g^{-1}Qg) & =g^{-1}N_G(Q)g \\
     &  & =g^{-1}Eg \\
    \Rightarrow & N = & g^{-1}Eg\;\;\mbox{  since }\circ(N)=m=\circ(E) \\
  \end{array}\] This proves (ii) in case II. Finally, we shall prove
  (iii). Let $M'$ be a subgroup of $G$ with $\circ(M')=m'$ such that $m'$
  divides $m$.
   If $m=m',$
  there is nothing to prove. Hence let $m'<m.$ Clearly $EM'K=G$. Therefore
   \[\begin{array}{crl}
     & \circ(G)= & \frac{\circ(E)\circ(M'K)}{\circ(E\cap M'K)} \\
    \Rightarrow & mp^a= & \frac{mm'p^a}{\circ(E\cap M'K)}\;\;\mbox{  since }\circ(M'K)=m'p^a \\
    \Rightarrow & m'= & \circ(E\cap M'K) \\
  \end{array}\] Note that $M'$ and $E\cap M'K$ are subgroups of $M'K$ with
  order $m'$, where $M'K$ is a solvable group of order
  $m'p^a<mp^a=\circ(G).$ Hence, by induction, $M'$ and $E\cap M'K$ are
  conjugate in $M'K$. If for $g\in G,$ \[\begin{array}{crl}
     & g^{-1}(E\cap M'K)g & =M' \\
    \Rightarrow & g^{-1}Eg & \supset M' \\
  \end{array}\] Thus $g^{-1}Eg$ is a subgroup of order $m$ in $G$
  containing $M'$. Thus the proof is complete.           $\blacksquare$
  \br The first property of the theorem characterizes finite
  solvable groups. The proof of this fact is beyond our scope of study.\er

 \vspace{.25in}   {\large\bf Nilpotent
groups
:}  \\   \\
 We first define {\bf upper central chain} of a group $G$. Put
 $Z_0(G) = \{e\}$, $Z_1(G) = Z(G)$, the centre of $G$, and for
 $i\geq 1$,
 $$Z(\frac{G}{Z_i(G)}) = \frac{Z_{i+1}(G)}{Z_i(G)}.$$
 By lemma \ref{l10.536}, $Z(G)$ is a characteristic subgroup of
 $G$. Hence, using induction and lemma \ref{l10.535}, $Z_i(G)$ is a
 characteristic subgroup of $G$ for all $i\geq 1.$ The chain
 $$\{e\}=Z_0(G) \subset Z_1(G) \subset Z_2(G) \subset \cdots $$
 of characteristic subgroups of $G$ is called the upper central
 chain of $G$. Note that by definition, $[G,Z_{i+1}(G)] \subset
 Z_i(G)$ for all $i\geq 0.$
 \bd\label{d10.61}
 A group $G$ is said to have upper central series/ chain of length
 $m$ if
$$\{e\}=Z_0(G)\subsetneqq Z_1(G)\subsetneqq \cdots
Z_{m-i}(G)\subsetneqq Z_m(G) = G. $$ \ed We,now define the {\bf
lower central chain} of a group $G$. Put $\Gamma_0 = G,
\Gamma_1(G) = [G,\Gamma_0(G)]$, and for all $i\geq 1$, $
\Gamma_{i+1}(G) = [G,\Gamma_i(G)]$. Note that $\Gamma_1(G) = [G,G]
= G^{\prime}$ is a characteristic subgroup of $G$ (Lemma
\ref{l10.54}). Further, if $\Gamma_i(G)$ is a characteristic
subgroup of $G$, then for any $\alpha \in AutG$, \ba{rcl} \alpha
(\Gamma_{i+1}(G)) & = & \alpha [G,\Gamma_i(G) ]\\
               & = & [ \alpha (G), \alpha (\Gamma_i(G)) ]\\
               & = & [ G, \Gamma_i(G) ] = \Gamma_{i+1}(G)
\ea Therefore $\Gamma_{i+1}(G)$ is a characteristic subgroup of $G$.
Hence,
by induction, $\Gamma_i(G)$ is a characteristic subgroup of $G$ for all
$i\geq 0.$ The chain $$ G = \Gamma_0(G)\supset \Gamma_1(G)\supset
\Gamma_2(G) \cdots$$ is called the lower central chain/series of $G$.
\bd\label{d10.62} A group $G$ is said to have lower central series/ chain
of length $m$ if $$ G = \Gamma_0(G)\supsetneqq \Gamma_1(G)\supsetneqq
\cdots \Gamma_{m-1}(G) \supsetneqq \Gamma_m(G) =\{e\}.$$ \ed Generalising the
concepts of upper and lower central chains of a group $G$, we define:
\bd\label{d13.28dd} Let $G$ be a group. A chain $$G=H_0\supset H_1\supset
\cdots \supset H_l=\{e\},$$ of normal subgroups of $G$ is called a central
chain / series of $G$ if $[G,H_i]\subset H_{i+1}$ for all
$i=0,1,\cdots,l-1.$ \ed \brs (i) By definition, it is clear that
$$\frac{H_i}{H_{i+1}}\subset Z(\frac{G}{H_{i+1}})$$ for all $i$. \\ (ii)
Any chain of normal subgroups of $G$ which refine a central chain is a central
chain of $G$.\ers
\bt\label{t10.63} A group $G$ has upper central chain of length $m$ if and
only if $G$ has lower central chain of length $m$. \et \pf Suppose $G$ has
upper central chain of length $m$. Then $$\{e\}=Z_0(G)\subsetneqq
Z_1(G)\subsetneqq \cdots Z_{m-i}(G)\subsetneqq Z_m(G) = G. $$ We shall
prove by induction that $\Gamma_i(G) \subseteq Z_{m-i}(G)$ for all $i\geq
0.$ If $i=0$, then $G = \Gamma_0(G) \subseteq Z_m(G) = G.$ Now, let $t\geq
1$, and $\Gamma_j(G) \subseteq Z_{m-j}(G)$ for all $0\leq j\leq t-1$. Then
$$
\begin{array}{cl}
   & \Gamma_{j+1}(G) = [G,\Gamma_j(G) ] \subseteq [G,Z_{m-j}(G) ] \subseteq
Z_{m-j-1}(G) \\ \Rightarrow & \Gamma_t(G)\subset Z_{m-t}(G)
  \end{array}.$$ Hence, by induction, $\Gamma_i(G) \subseteq Z_{m-i}(G)$ for
all $i\geq 0$. In particular, $\Gamma_m(G) \subseteq Z_0(G) = \{e\}.$
Therefore $G$ has lower central chain of length $\leq m.$\\ Next, let $G$
has lower central chain of length $m$. i.e., $$ G = \Gamma_0(G)\supsetneqq
\Gamma_1(G)\supsetneqq \cdots \Gamma_{m-1}(G) \supsetneqq \Gamma_m(G)
=\{e\}.$$ We shall prove by induction that $\Gamma_{m-i}(G) \subseteq
Z_i(G)$ for all $i\geq 0.$ If $i=0$, then $\{e\} = \Gamma_m(G) \subseteq
Z_0(G) = \{e\}.$ Now, let $t\geq 1$, and let $\Gamma_{m-j}(G) \subseteq
Z_j(G)$ for all $j < t$. Then, by our assumption $$\begin{array}{ccc}
   & \Gamma_{m-t+1}(G)=[G,\Gamma_{m-t}(G)]\subset Z_{t-1}(G) &  \\
  \Rightarrow & [g,x]\in Z_{t-1}(G) & \mbox{ for all $g \in G$ and }x \in \Gamma_{m-t}(G) \\
  \Rightarrow & Z_{t-1}(G)[g,x]=Z_{t-1}(G) &  \\
  \Rightarrow & Z_{t-1}(G)gZ_{t-1}(G)x=Z_{t-1}(G)xZ_{t-1}(G)g &  \\
  \Rightarrow & Z_{t-1}(G)x\in Z(G/Z_{t-1}(G))=Z_t(G)/Z_{t-1}(G) &  \\
  \Rightarrow & x\in Z_t(G) &  \\
  \Rightarrow & \Gamma_{m-t}(G) \subset Z_t(G) &
\end{array}$$ Hence, by induction,
$\Gamma_{m-i}(G) \subseteq Z_i(G)$ for all $i\geq 0.$ In particular, $G =
\Gamma_0(G) \subseteq Z_m(G).$ Hence $G$ has upper central chain of length
$\leq m.$ Thus the result is proved.          $\blacksquare$
\bd\label{d10.64}
A group $G$ is called Nilpotent if $G$ has an upper (lower)
central chain of finite length. We say that $G$ is Nilpotent of class $m$
if its upper (lower) central chain has length $m$.
\ed
\brs
(i) A nilpotent group is solvable.\\
 If $G$ is a nilpotent group of class $m$, then $G$ has upper
 central chain of length $m.$ i.e.,
 \be\label{e10.65}
\{e\}=Z_0(G)\subsetneqq Z_1(G)\subsetneqq \cdots Z_{m-i}(G)\subsetneqq
Z_m(G) = G
\ee As $Z_i(G)$ is a characteristic subgroup of $G$ for all $i\geq 0$,
(\ref{e10.65}) is a normal chain of $G$. By definition
$$\frac{Z_{i+1}}{Z_i(G)} = Z(\frac {G}{Z_i(G)})\;\;\;\;\; \forall i\geq
0.$$ Hence all factor groups of the normal chain (\ref{e10.65}) are
abelian. This proves $G$ is solvable.\\ (ii)If $G = \{e\}$, then $G$ is
nilpotent of class 0. Further, if $G$ is an abelian group of order $> 1,$
then $Z_1(G) = G.$ Therefore $G$ is nilpotent of class 1.\\ (iii) A
nilpotent group need not be abelian. e.g., if $G$ is the group of
quaternions (Example \ref{ex1.10}), then $$ Z_1(G) = \left\{ \left[
\begin{array}{cc} 1 & 0 \\  0 & 1
                        \end{array}\right],
                  \left[\begin{array}{cc}  -1 & 0 \\  0 & -1
                         \end{array}\right]
             \right\}$$
 and $o(\frac{G}{z_1(G)}) = 4$, which is abelian (Corollary \ref{c615}). Hence
  $$ \{e\} = Z_0(G) \subsetneqq Z_1(G) \subsetneqq Z_2(G) = G$$
  is the upper central chain of $G$. Thus $G$ is nilpotent of
  class 2.\\
  (iv) A solvable group need not be nilpotent.\\
  Let $G = S_3$, the symmetric group on three symbols. We have
 observed that $S_3$ is solvable, and also noted that $Z(S_3) = \{id\}.$
Hence $S_3$ does not have an upper central chain and thus is not
nilpotent. \ers \bexe Let $G$ be a group and $x_i\in G$ for
$i=1,2,\cdots,n.$ A commutator of weight $n$ in
$x_1,x_2,\cdots,x_n$ is defined to be $x_1$ if $n=1,$ $[x_1,x_2]$
if $n=2$ and $[x_1,\cdots,x_n]=[x_1,[x_2,\cdots,x_{n-1},x_n]]$ for
$n\geq 3.$ Prove that $G$ is nilpotent of class $m$ if and only if
$[a_1,a_2,\cdots,a_{m+1}]=e$ for $a_i's$ in $G$.\eexe  \bexe Prove
that for a group $G$, $\Gamma_i(\Gamma_j(G))\subset
\Gamma_{i+j}(G)$ for all $i\geq 0, j\geq 0.$ Then deduce that if
$G$ is nilpotent of class $m$, then for any $i\leq m,$
$\Gamma_i(G)$ is nilpotent of class $\leq m-i.$\eexe
\bl\label{l10.66} Let $G$ be a nilpotent group. Then\\ (i) Every
subgroup of $G$ is nilpotent.\\ (ii) Any homomorphic image of $G$
is nilpotent. \el \pf (i) Let $H$ be a subgroup of $G$. We have $H
= \Gamma_0(H) \subset \Gamma_0(G) = G$. Further, if $\Gamma_i(H)
\subset \Gamma_i(G)$, then $$\Gamma_{i+1}(H) = [H,\Gamma_i(H)]
\subset [G,\Gamma_i(G)] = \Gamma_{i+1}(G).$$ Hence by induction,
$\Gamma_i(H) \subset \Gamma_i(G)$ for all $i\geq 0.$ If $G$ is
nilpotent of class $m$, then $\Gamma_m(H) \subset \Gamma_m(G) =
\{e\}$. Therefore $\Gamma_m(H) = \{e\}.$ Thus $H$ is nilpotent of
class $\leq m.$\\ (ii) Let for a group $K$, $f:G\to K$ be an onto
homomorphism. Then $$ f(\Gamma_o(G)) = f(G) = K = \Gamma_0(K)$$
Further, let for some $i\geq 0$, $f(\Gamma_i(G)) = \Gamma_i(K).$
Then \ba{rcl} f(\Gamma_{i+1}(G)) & = & f [G,\Gamma_i(G)]\\
                  & = & [f(G),f(\Gamma_i(G))]\\
                  & = & [K,\Gamma_i(K)]\\
                  & = & \Gamma_{i+1}(K)
\ea Hence, by induction $f(\Gamma_i(G)) = \Gamma_i(K)$ for all $i\geq 0.$
If $G$ is nilpotent of class $m$, then $\Gamma_m(G) = \{e\}$. Therefore
$\Gamma_m(K)=f(\Gamma_m(G))=\{e\}.$ Thus $K$ is nilpotent of class $\leq
m.$        $\blacksquare$
\bco\label{c10.67}
Any factor group of a nilpotent group is nilpotent.
 \eco
\pf As a factor group of a group is its homomorphic image, the result
follows from (ii).
 \br
If $G$ is a nilpotent group and $H\lhd G$, then from the lemma $H$ and
$\frac{G}{H}$ are nilpotent. The converse is not true e.g., if $G = S_3$,
then $G$ is not nilpotent, but, for $H = A_3 \lhd S_3 = G$, $H$ is cyclic
of order $3$ and $\frac{G}{H}$ is cyclic of order $2$, hence both $G$ and
$\frac{G}{H}$ are nilpotent. \er
\bl\label{l10.68} Any finite $p$-group is nilpotent.
\el
 \pf Let
$G$ be a finite $p$-group. Then $o(G) = p^m, m\geq 0.$ If $m = 0$, $G =
\{e\}$ and hence is nilpotent of class $0$. Let $m\geq 1$, and let
\be\label{e10.69}
\{e\}=Z_0(G)\subset Z_1(G)\subset \cdots
Z_i(G)\subset \cdots
\ee
be the upper central chain of $G$. If $Z_i(G)\neq G$ for some $i$, then as
$Z_i(G)$ is a subgroup of $G$, $o(Z_i(G)) = p^{\alpha_i}$ for some $0\leq
\alpha_i < m.$ Note that $$o(\frac{G}{Z_i(G)}) = p^{m-\alpha_i} > 1.$$
Hence by theorem \ref{t614}, $$\frac{Z_{i+1}(G)}{Z_i(G)} =
Z(\frac{G}{Z_i(G)})\neq \{id\}$$ $$\Rightarrow Z_{i+1}(G)\supsetneqq
Z_i(G).$$ Thus in the upper central chain (\ref{e10.69}) of $G$,
$Z_i(G)\neq G$ implies $ Z_{i}(G)\subsetneqq Z_{i+1}(G).$ As $G$ is finite
we must have $Z_k(G) = G$ for some $k$. Hence $G$ is nilpotent.
$\blacksquare$
\bl\label{l10.70}
A finite direct product of nilpotent groups is nilpotent.
\el
\pf Let $G_i$, $i =1,2,\cdots ,t$ be nilpotent groups of class $m_i$. Then
for $m = max \{m_i\mid 1\leq i\leq t\}$, $\Gamma_m(G_i) = \{e\}$ for all
$i\geq 1.$ Put $$ G = G_1\times G_2\times \cdots \times G_t$$ Then
\ba{rcl} \Gamma_0(G) = G & = & G_1\times G_2\times \cdots \times G_t\\ & =
&\Gamma_0(G_1)\times \Gamma_0(G_2)\times \cdots \times \Gamma_0(G_t) \ea
If $$ \Gamma_i(G) = \Gamma_i(G_1)\times \Gamma_i(G_2)\times \cdots \times
\Gamma_i(G_t)$$ then \ba{rcl} [G,\Gamma_i(G)] & = &
[G_1,\Gamma_i(G_1)]\times [G_2,\Gamma_i(G_2)]\times
                       \cdots \times [G_t,\Gamma_i(G_t)]\\
\Rightarrow \Gamma_{i+1}(G) & = & \Gamma_{i+1}(G_1)\times \Gamma_{i+1}(G_2)
                      \times \cdots \times \Gamma_{i+1}(G_t)
\ea Hence, by induction $$ \Gamma_t(G) = \Gamma_t(G_1)\times
\Gamma_t(G_2)\times \cdots \times \Gamma_t(G_t)$$ for all $t\geq
0.$ Therefore $$ \Gamma_m(G) = \Gamma_m(G_1)\times
\Gamma_m(G_2)\times \cdots \times \Gamma_m(G_t) = \{e\}.$$
Consequently, $G$ is nilpotent of class $m$.        $\blacksquare$
\bl\label{l10.71} Let $G$ be a nilpotent group and $H\subsetneqq
G$, a subgroup of $G$. Then $H\subsetneqq N_G(H).$ \el \pf Let $G$
be nilpotent of class $m.$ As $H\subsetneqq G$, $G\neq \{e\}.$
Hence $m > 0.$ Choose $r$ maximal such that $Z_r(G) \subset H,
Z_{r+1}(G) \notsubset H.$ Clearly $r < m$, as $Z_m(G) = G\;
\notsubset H.$ Now, take an element $x\in Z_{r+1}(G), x\notin H.$
Then, as \ba{rcl} [G,Z_{r+1}(G)] & \subset & Z_r(G) \subset H\\
\Rightarrow [h,x] & \in & H \;\;\;\;\forall h\in H  \\
i.e.,hxh^{-1}x^{-1} & \in & H\\ \Rightarrow xh^{-1}x^{-1} & \in &
H\;\;\;\; \forall h\in H\\ \Rightarrow x & \in & N_G(H) \ea
Hence$H \subsetneqq N_G(H).$         $\blacksquare$
\bco\label{c10.72} Every subgroup $H$ of a finite nilpotent group
$G$ is subnormal. \eco \pf Clearly $G$ and the identity subgroup
are subnormal in $G$. Hence, let $H$ be a proper subgroup of $G$.
By the lemma $H \subsetneqq N_G(H) = H_1$(say). If $H_1\neq G$,
then as above, $H_1 \subsetneqq N_G(H_1) = H_2$. If we define
$N_G(H_i)=H_{i+1}$ for all $i\geq 1$, then as $G$ is finite this
chain shall terminate in $G$. If $H_r = G$, then $$\{e\} \lhd
H\lhd H_1\lhd \cdots \lhd H_{r-1}\lhd H_r = G $$ is a subnormal
chain of $G$. Hence $H$ is subnormal.$\blacksquare$ \\ We shall,
now, prove that the corollary is true for any nilpotent group $G$
i.e., $G$ need not be finite. \bl\label{l10.73} Any subgroup $H$
of a nilpotent group $G$ is subnormal. \el \pf Let $G$ be
nilpotent of class $m$. We shall prove the result by induction on
$m.$ If $m = 0$, $G=\{e\}$, and if $m =1$, $G$ is abelian. Hence
the result is trivial for $0\leq m\leq 1.$ Let $m > 1.$ Clearly,
for $Z_1 = Z(G)$, $Z_1H$ is a subgroup of $G$ such that $H\lhd
Z_1H$, moreover, $\frac{G}{Z_1(G)}$ is a nilpotent group of class
$m-1$. Hence by induction $\frac{Z_1H}{Z_1}$ is a subnormal
subgroup of $\frac{G}{Z_1}$. Thus there exists a normal chain of
$G$ of which $Z_1H$ is a member(Use : corollary \ref{c301a}). Now,
as $H\lhd Z_1H$, the result follows.$\blacksquare$\\ The next
lemma gives an alternative way of looking at the above statement.
\bl\label{l10.74} Let $G$ be a nilpotent group of class $m$ and
$H$, a subgroup of $G$. Put $N_0 = H$, $N_1 = N_G(H)$ and $N_{i+1}
= N_G(N_i)$ for all $i\geq 1.$ Then $N_m = G$. \el \pf The proof
is clear for $m=0$. Let $m\geq 1$. Clearly $Z_1(G) \subset N_G(H)
= N_1.$ Suppose we have proved $Z_k(G)\subset N_k.$ Then \ba{rcl}
[G,Z_{k+1}(G)] &\subset & Z_k(G)\subset N_k\\ \Rightarrow
[x^{-1},y^{-1}] &\in & N_k \;\;\;\;\forall x\in N_k, y\in
                                                    Z_{k+1}(G)\\

\Rightarrow y^{-1}xy &\in & N_k \;\;\;\;\forall x\in N_k, y\in
                                                    Z_{k+1}(G)\\
\Rightarrow y &\in & N_G(N_k) = N_{k+1} \\ \Rightarrow Z_{k+1}(G)
& \subset &  N_{k+1} \ea  Hence, by induction, $G = Z_m(G) \subset
N_m.$ Thus $N_m = G.       \blacksquare$ \br From the above
statement it is clear that $H$ is subnormal in $G$. \er
\bt\label{t13.41tt} Let $G$ be a finite group of order $p^n$ where
$n\geq 1$ and $p$ is a prime. Then $G$ is cyclic of order $p^n$ if
and only if $G$ has unique composition series of length $n$. \et
\pf By theorem \ref{t28}, any cyclic group of order $m$ has
exactly one subgroup of order $d$ for each divisor $d\geq 1$ of
$m.$ Now, if $G$ is cyclic of order $p^n$, let $H_0=G;H_1$, the
cyclic group of order $p^{n-1}$ in $H_0$, and in general,
$H_{i+1}$ the cyclic group of order $p^{n-(i+1)}$ in $H_{i}$ for
all $i\geq 0.$ Then $H_i/H_{i+1}$ is a cyclic group of order $p$,
hence is simple. Thus, it is clear that $$G=H_0\rhd H_1\rhd \cdots
\rhd H_{n-1}\rhd H_n=\{e\}$$ is composition series of $G$ of
length $n$. To prove that the composition series is unique note
that for any two subgroups $H\supset K$ in $G$, $H/K$ is simple if
and only if $\circ(H/K)=p.$ Hence if $$G=A_0\rhd A_1\rhd A_2\rhd
\cdots \rhd A_{s-1}\rhd A_s=\{e\}$$ is a composition series of
$G$, then $\circ(A_i)=p^{n-i}.$ Thus $s=n$, and the uniqueness of
the composition series follows by theorem \ref{t26}. We shall,
now, prove the converse by induction on $n$. If $n=1$, then $G$ is
simple. Hence, by theorem \ref{t614}, $\circ(G)=p$.
 Therefore $G$ is cyclic of order $p$. Now, let $n>1$
and let $$G=K_n\rhd K_{n-1}\rhd\cdots \rhd K_1\rhd K_0=\{e\}$$ be
the unique composition series of length $n$ for $G$. Clearly
$K_{n-1}$ has unique composition series of length $n-1$ and is a
group of order $p^{\alpha}.$ Therefore by induction hypothesis
$K_{n-1}$ is cyclic of order $p^{n-1}.$ As $G/K_{n-1}$ is simple
$p$-group, $\circ(G/K_{n-1})=p.$ Therefore $\circ(G)=p^n.$
Further, let us note that $K_{i+1}/K_i$ is simple $p$-group for
each $0\leq i\leq n-1.$ Hence $\circ(K_i)=p^i$ for all $i\geq 0$.
Now, let $x\in G-K_{n-1}.$ Then $\circ(x)=p^{\beta}$ for some
$\beta\leq n.$ Assume $\beta <n,$ and put $C=gp\{x\}$. By lemma
\ref{l10.68} and corollary \ref{c10.72}, $C$ is subnormal subgroup
of $G$. Therefore using theorem \ref{t10.24} there exists a
composition series for $G$ with $C$ as its member. Now, by
uniqueness of the composition series $C=K_{\beta}.$ This
contradicts our assumption that $x\in G-K_{n-1}.$ Hence $\beta
=n,$ and $G$ is cyclic group of order $p^n.$     $\blacksquare$
\bt\label{t10.75} Any non-trivial normal subgroup $H$ of a
nilpotent group $G$ intersects the center $Z_1(G)$ of $G$
non-trivially i.e., $H\cap Z_1(G)\neq \{e\}.$ \et \pf We shall
prove the result by induction on the nilpotency class $m$ of $G$.
Clearly $m\geq 1$. If $m = 1$, then $G$ is abelian and the proof
is clear. Let $m>1$. If $H\subset Z_1(G)$, the result is
immediate. Now, let $H \notsubset Z_1(G) = Z_1.$ Then
$Z_1H\supsetneqq Z_1$ and $Z_1H\lhd G.$ Hence $\frac{Z_1H}{Z_1}$
is a nontrivial normal subgroup of $\frac{G}{Z_1}$. Note that
$\frac{G}{Z_1}$ is nilpotent of class $m-1$. Hence by induction
$Z_1H\cap Z_2\supsetneqq Z_1.$ Choose $x\in Z_1$, $h\in H$ such
that $xh\in Z_2$, and $xh\notin Z_1$. Then $h\notin Z_1$ and $h\in
Z_2.$ As $h\notin Z_1$ there exists an element $g\in G$ such that
$[g,h]\neq e$. As $[G,Z_2]\subset Z_1$ we get $[g,h]\in Z_1\cap
H$. (Use: $H\lhd G$). Hence $Z_1\cap H\neq \{e\}.$$\blacksquare$
\bl\label{l10.76} Let $G$ be a nilpotent group of class $m\geq 2.$
Then for any $a\in G$, the subgroup $H = gp \{a,[G,G]\}$ is
nilpotent of class $\leq m-1.$ \el \pf First of all, we shall show
that $$Z_i(G)\cap H \subset Z_i(H)$$ for all $i\geq 0.$ As $Z_0(G)
= Z_0(H) = \{e\}$, the claim is clear for $i = 0.$ Now, let $t\geq
0$, and $Z_t(G)\cap H \subset Z_t(H).$ Then \ba{rcl}
         [H,Z_{t+1}(G)\cap H] &\subset & [H,Z_{t+1}(G)]\cap H\\
\Rightarrow [H,Z_{t+1}(G)\cap H] &\subset & Z_t(G)\cap H \subset Z_t(H)\\
\Rightarrow Z_{t+1}(G)\cap H &\subset & Z_{t+1}(H) \ea Hence by induction
the claim follows. Now, for $i = m-1$, \ba{rcl} Z_{m-1}(G)\cap H &\subset
& Z_{m-1}(H)\\ \Rightarrow [G,G] &\subset & Z_{m-1}(H) \ea Thus
$\frac{H}{Z_{m-1}(H)}$ is cyclic. Note that $$ Z (\frac{H}{Z_{m-2}(H)}) =
\frac{Z_{m-1}(H)}{Z_{m-2}(H)}$$ and $$
\frac{\frac{H}{Z_{m-2}(H)}}{\frac{Z_{m-1}(H)}{Z_{m-2}(H)}}
                \cong \frac{H}{Z_{m-1}(H)}$$ is cyclic. Hence
$\frac{H}{Z_{m-2}(H)}$ is abelian. This implies $Z_{m-1}(H) = H.$
Therefore the result follows. $\blacksquare$ \bt\label{t10.77} Let $G$ be
a nilpotent group. Then $$ G^t = \{x\in G\mid o(x)< \infty \}$$ is a
subgroup of $G.$ This is called the torsion subgroup of $G.$
\et \pf We shall prove the result by induction on the nilpotency
class $m$ of $G.$ If $m=0$, then $G=\{e\}.$ Hence the result is trivial.
Now, let $m\geq 1.$ If $m=1,$ then $G$ is abelian. Therefore the result
follows by lemma \ref{l99}. Now, let $m>1$ and let the result holds for
every nilpotent group of class $<m$. Choose $a,b\in G^t.$ By lemma
\ref{l10.76} for $H = gp \{a,[G,G]\}$ and $K = gp \{b,[G,G]\}$, $H^t$ and
$K^t$ are subgroups of $H$ and $K$ respectively. Now, note that for any
endomorphism of a group, image of an element of finite order is of finite.
Hence $H^t$ and $K^t$ are characteristic subgroups of $H$ and $K$
respectively. Further, as $[G,G]\subset H$, $H\lhd G$ (Lemma \ref{l101l})
  For the same reasons $K\lhd G$. Therefore
$H^t \lhd G$ and $K^t\lhd G$ (Lemma \ref{l10.538}). Now, for any $x\in
H^t$, if $o(x) = l$, then \ba{rcl}
 {(xK^t)}^l & = & K^t\\
 \Rightarrow {(xy)}^l & \in & K^t \;\;\;\;\forall y\in K^t \\
\Rightarrow o(xy) & < & \infty \;\;\;\;\forall x\in H^t, y\in K^t \\
\Rightarrow o(ab^{-1}) & < & \infty \ea Therefore $G^t$ is a subgroup of
$G$. $\blacksquare$
\bl\label{10.78} If for a subgroup
$H$ of a nilpotent group $G$, $H[G,G] = G$, then $H = G$.
\el
\pf Let $G$ be nilpotent of class $m$ and let $r$ be maximal such that
$H\Gamma_r (G) = G$. By assumption $H\Gamma_1 (G) = H[G,G] = G$.
 Hence $r \geq 1$. We claim $r=m$. If
not, then by lemma \ref{l10.540}, $$H\Gamma_{r+1} (G) \lhd H\Gamma_r
(G) =G.$$
 As $\Gamma_r (G) / \Gamma_{r+1} (G)$ is abelian, we get
 $$G/(H\Gamma_{r+1} (G)) = H\Gamma_r (G)/H\Gamma_{r+1} (G)$$
 is
abelian. Hence $[G,G] \subset H\Gamma_{r+1} (G)$. Thus $G = H[G,G] =
H\Gamma_{r+1} (G)$. This contradicts the maximality of $r$ and
consequently $r=m$. Thus $H = H\Gamma_m (G) = G$. $\blacksquare$
\bco\label{co10.79co} If $G$ is a nilpotent group, then $[G,G]\subset
\Phi(G)$, the Frattini subgroup of $G$. \eco \pf Let $x\in [G,G]$, and let
for a subset $S$ of $G$, $$\begin{array}{cl}
   & G=gp\{x,S\} \\
  \Rightarrow & G=[G,G]H
\end{array}$$ where $H=gp\{S\}$. Hence by the lemma, $G=gp\{S\}.$ Thus $x$
is a non-generator of $G$, and consequently $[G,G]\subset
\Phi(G).$ $\blacksquare$ \bl\label{l10.79} Let $G$ be a torsion
free nilpotent group i.e., $G^t = \{e\}$. Then if $a^n = b^n$ for
$a,b \in G$ and $n>0$, $a=b$. \el \pf Let $G$ be nilpotent of
class $m$. We shall prove  the result by induction on $m$. If
$m=0$, then $G$ is identity. Hence there is nothing to prove. In
case $m=1$, $G$ is abelian. Therefore $(ab^{-1})^n = a^n
(b^n)^{-1} = e$. Hence as $G$ is torsion free $a=b$. Now, let $m >
1$. By the lemma \ref{l10.76}, $H = gp\{a, [G,G]\}$ is nilpotent
of class $\leq m-1$. Clearly $bab^{-1} \in H$ since $H\lhd G$.
Moreover \ba{rcl} (bab^{-1})^n & = & ba^nb^{-1}\\
             & = & bb^n b^{-1} = b^n = a^n
\ea Hence, by induction, $bab^{-1} = a$. This gives $ab = ba$, and
consequently ${(ab^{-1})}^n = a^n {(b^n)}^{-1} = e$. Therefore $ab^{-1} =
e$ i.e., $a = b.$$\blacksquare$
\bt\label{t13.44tt} If for a nilpotent group $G$, $\circ(G/G')<\infty$,
then $G$ is finite.\et \pf Let $G$ be nilpotent group of class $n$. We
shall prove the result by induction on $n$. If $n\leq 1$, then
 $G'=id,$
and hence the result is trivial. Now, let $n\geq 2.$ Consider the natural
epimorphism: \[\eta : G \rightarrow G/\Gamma_{n-1}(G)=H.\]  Then
$\eta$ induces the epimorphism: $$\begin{array}{crcl}
  \overline{\eta}: & G/G' & \rightarrow & H/H' \\
   & xG' & \mapsto & \eta (x)H'
\end{array}$$ As $\circ(G/G')<\infty, \circ(H/H')<\infty.$ Note that $H$
is nilpotent of class $n-1$. Hence, by induction, $\circ(H)<\infty.$ As
$\Gamma_{n-1}(G)\subset Z(G), \circ(H)<\infty$ implies $[G:Z(G)]<\infty.$
Therefore by theorem \ref{t4.14tt}, $G'$ is a finite group. Now, as
$\circ(G)=\circ(G/G')\cdot \circ(G')$, $G$ is a finite group.
$\blacksquare$    \bl\label{l13.46ll} Let $G$ be a nilpotent group. Then
any finitely generated subgroup of $G^t$ is finite.\el \pf Let $S$ be a
finite subset in $G^t$ and $H=gp\{S\}.$ Then $H$ is a nilpotent group
(Lemma \ref{l10.66}). By corollary \ref{b101z}, the factor group $H/H'$ is
an abelian group. Clearly $H/H'$ is generated by the finite set $\{xH'\st
s\in S\}$ where $xH'$ of finite order for all $x\in S.$ Thus $H/H'$ is a
finite group. The result, now, follows by theorem \ref{t13.44tt}.
$\blacksquare$  \\    \\ {\bf STRUCTURE THEOREM FOR FINITE NILPOTENT
GROUPS}
\bt\label{t10.80}
Let $G$ be a finite group and $$o(G) =
{p_1}^{\alpha_1}{p_2}^{\alpha_2}\cdots {p_k}^{\alpha_k}$$ where $p_i$'s
are primes, $\alpha_i > 0$ and $k\geq 1.$ Then $G$ is nilpotent if and only
if $G$ is a direct product of its Sylow $p$-subgroups.
\et
\pf Let $P_i, i = 1,2,\cdots ,k$ denotes Sylow $p_i$-subgroup of $G$. Then
$o(P_i) = {p_i}^{\alpha_i}.$ Let $$ G = P_1\times P_2\times \cdots \times
P_k$$ As a finite $p$-group is nilpotent(Lemma \ref{l10.68}) and any
finite direct product of nilpotent groups is nilpotent (Lemma
\ref{l10.70}), it follows that $G$ is nilpotent. Conversely, let $G$ be
nilpotent and $k>1$, since otherwise there is nothing to prove. Let $P =
P_i$ for an $i$ and $H = N_G(P)$. Then, using lemma \ref{l10.17l} and lemma
\ref{l10.71}, we get $H=G$. Thus $P_i \lhd G$ for each $i = 1,2,\cdots
,k$. Now, note that $A = P_1P_2\cdots P_k \subset G$ is a subgroup of $G$.
As $P_i \subset A$, $o(P_i) = {p_i}^{\alpha_i}$, $1\leq i \leq k$,
divides, $o(A)$. Consequently \ba{rcl} o(P_1P_2\cdots P_k) & \geq &
       {p_1}^{\alpha_1}{p_2}^{\alpha_2}\cdots {p_k}^{\alpha_k} =
       o(G)\\
\Rightarrow G & = & P_1P_2\cdots P_k \ea We have
$({p_i}^{\alpha_i},{p_j}^{\alpha_j}) = 1$ for all $i\neq j$. Hence for
$1\leq i\neq j\leq k$, $P_i \cap P_j = \{e\}.$ Now, for any $1\leq i\neq
j\leq k$, as $P_i$ and $P_j$ are
normal in $G$ and $P_i \cap P_j = \{e\}$, $ab = ba$ for all $a\in P_i$,
$b\in P_j$(Lemma \ref{w}(b)). Note that $$ x  \in  P_i\cap (P_1P_2\cdots
P_{i-1}P_{i+1}\cdots P_k)$$
\be\label{*}
 \Rightarrow x  =  a_i = a_1a_2\cdots a_{i-1}a_{i+1}\cdots a_k
\ee
where $a_t \in P_t$ for $1\leq t\leq k$. As $a_t \in P_t$, $o(a_t) =
{p_t}^{\beta_t}$ for all $1\leq t\leq k.$ If $n_i = {p_1}^{\beta_1} \cdots
{p_{i-1}}^{\beta_{i-1}} {p_{i+1}}^{\beta_{i+1}} \cdots
{p_k}^{\beta_k}$, then by the equation \ref{*} $$ x^{n_i} = {(a_1a_2\cdots
a_{i-1}a_{i+1}\cdots a_k)}^{n_i} = e.$$ Further, $o(x) = o(a_i) =
{p_i}^{\beta_i}.$ Therefore ${p_i}^{\beta_i}$ divides $n_i$. Since $(n_i ,
p_i) = 1$, we conclude $o(x) = 1$. Thus $P_i\cap (P_1P_2\cdots
P_{i-1}P_{i+1}\cdots P_k) = \{e\}$ for all $1\leq i\leq k.$ Consequently
$$  G = P_1\times P_2\times \cdots \times P_k. \;\;\;\; \blacksquare$$
\bt\label{t13.45tt} For a finite group $G$, the following statements are
equivalent: \\ (i) The group $G$ is nilpotent. \\ (ii) Every subgroup of
$G$ is subnormal. \\ (iii) The group $G$ is direct product of its Sylow
$p$-subgroups. \\ (iv) $[G,G]\subset \Phi(G).$ \et \pf (i) $\Rightarrow$
(ii): (Corollary \ref{c10.72}).
\\ (ii) $\Rightarrow$ (iii): By (ii), for any proper subgroup $H$ of $G$,
$H\subsetneqq N_{G}(H).$ Hence using lemma \ref{l10.17l}, if $P$
is a Sylow $p$-subgroup of $G$, then $N_{G}(P)=G.$ Thus any Sylow
$p$-subgroup of $G$ is normal in $G$. Now, as in theorem
\ref{t10.80}, $G$ is direct product of its Sylow $p$-subgroups. \\
(iii) $\Rightarrow$ (i): (Theorem \ref{t10.80}). \\ To complete
the proof, we shall now prove that (i) is equivalent to (iv). \\
(i) $\Rightarrow$ (iv): \: This follows by corollary
\ref{co10.79co}. \\ (iv) $\Rightarrow$ (i):\: By theorem
\ref{t10.80}, it is sufficient to show that every sylow
$p$-subgroup of $G$ is normal in $G$. Let $P$ be a sylow
$p$-subgroup of $G$, and let $H=N_G(P).$ If $H\neq G$, take a
maximal subgroup $M$ of $G$, $M\supset H.$ As $N_G(P)\subset M,$
$N_G(M)=M$ by lemma \ref{l10.17l}. Now,as, $[G,G]\subset
\Phi(G)\subset M$ (Theorem \ref{t2.09t}), $[G,M]\subset M.$ Hence
$M\lhd G$. Thus $N_G(M)=G.$ This implies $H=G$ i.e., $P\lhd G.$
$\blacksquare$    \\  \\  {\large\bf Finitely Generated Nilpotent
Groups:} \\ \bl\label{l13.48ll} Let $G$ be a nilpotent group, and
$x,y\in G$ be elements of finite order. If $\circ(x)=m,
\circ(y)=n,$ and $(m,n)=1$ then $xy=yx.$\el \pf Let $K=gp\{x,y\}$.
By lemma \ref{l13.46ll}, $K$ is a finite group. Next, by lemma
\ref{l10.66}, $K$ is nilpotent. By the structure theorem of finite
nilpotent groups, $K$ is direct product of its Sylow subgroups.
Let $\circ(K)=p_1^{\alpha_1}\cdots p_t^{\alpha_t}$ where $p_i's$
are primes and $\alpha_i>0$ for all $1\leq i\leq t.$ Let $S_i$
denotes the Sylow $p_i$-subgroup of $K$ for $i=1,2,\cdots,t.$ As
$K=S_1S_2\cdots S_t,$ we can write $$\begin{array}{cl}
  x & =a_1a_2\cdots a_t \\
  y & =b_1b_2\cdots b_t
\end{array}$$ where $a_i,b_i\in S_i$ for all $i=1,2,\cdots t$. Note that
for all $1\leq i\neq j \leq t;$
$(\circ(b_i),\circ(b_j))=1=(\circ(a_i),\circ(a_j))$ and $ab=ba$
whenever $a\in S_i,$ and $b\in S_j.$ Therefore
$\circ(x)=\circ(a_1)\circ(a_2)\cdots \circ(a_t),$ and
$\circ(y)=\circ(b_1)\circ(b_2)\cdots \circ(b_t).$ As $(m,n)=1$, it
is immediate that for no $1\leq i\leq t,$ $a_i\neq e,$ $b_i\neq
e.$ Therefore, as $ab=ba$ for all $a\in S_i,$ $b\in S_j$, $i\neq
j,$ we get $xy=yx.$           $\blacksquare$ \bt\label{t13.47tt} A
nilpotent group $G(\neq id)$ is finitely generated if and only if
every subgroup of $G$ is finitely generated.\et  \pf It is enough
to prove the direct part as the converse is trivial. Let $G$ be
nilpotent of class $n$. We shall prove the result by induction on
$n$. If $n=1$, then $G$ is abelian. Hence the result follows by
lemma \ref{l913}. Now, let $n>1.$ As $G$ is finitely generated,
$G/\Gamma_{n-1}(G)$ is a finitely generated nilpotent group of
class $n-1$. Hence by induction $\Gamma_{n-2}(G)/\Gamma_{n-1}(G)$
is finitely generated group. Let for $b_1,b_2,\cdots,b_k$ in
$\Gamma_{n-2}(G),$ $b_i\Gamma_{n-1}(G),\; i=1,2,\cdots,k,$
generate $\Gamma_{n-2}(G)/\Gamma_{n-1}(G).$ Then any element of
$\Gamma_{n-2}(G)$ can be written as
$b_1^{\alpha_1}b_2^{\alpha_2}\cdots b_k^{\alpha_k}z$ where
$\alpha_i\in \Z$ and $z\in \Gamma_{n-1}(G).$ Further, let
$G=gp\{x_1,x_2,\cdots,x_t\}.$ Then any element of $G$ is a finite
product of $x_i^{\epsilon}, \epsilon=\pm 1$ (repetitions allowed).
Now, as $[G,\Gamma_{n-2}(G)]=\Gamma_{n-1}(G)\subset Z(G),$ and for
any $a,b,c$ in $G$, $[ab,c]=a[b,c]a^{-1}[a,c]$,
$[a,bc]=[a,b]b[a,c]b^{-1},$ $[a^{-1},b]=a^{-1}[b,a]a$ and
$[a,b^{-1}]=b^{-1}[b,a]b,$ it follows that
$\Gamma_{n-1}(G)=[G,\Gamma_{n-2}(G)]$ is generated by the finite
set of elements $[x_i,b_j], 1\leq i\leq t, 1\leq j\leq k.$ As
$\Gamma_{n-1}(G)\subset Z)G),$ it is abelian. Now, let $H(\neq
id)$ by any subgroup of $G$. Then $H/\Gamma_{n-1}(G)\cap H$ is a
subgroup of the finitely generated nilpotent group
$G/\Gamma_{n-1}(G)$ of class $n-1.$ Hence by induction $H/H\cap
\Gamma_{n-1}(G)$ is finitely generated. Further as
$\Gamma_{n-1}(G)$ is finitely generated abelian group, $H\cap
\Gamma_{n-1}(G)$ is finitely generated (Lemma \ref{l913}).
Therefore $H$ is finitely generated.        $\blacksquare$
\bt\label{t13.53} Let $G$ be a finitely generated nilpotent group.
If $Z(G)$ is finite then $G$ is finite.\et \pf Assume $G$ is
non-abelian since otherwise $G=Z(G)$ is finite. Let $G$ be
nilpotent of class $m.$ We shall first prove, by induction on $m$,
that $G$ contains a torsion free, normal subgroup of finite
index.As $G$ is non-abelian, $m\geq 2.$ If $m=2,$ then
$$G\supsetneqq \Gamma_1(G)\supsetneqq \Gamma_2(G)=id$$ is the
lower central chain of $G$. As $[G,\Gamma_1(G)]=\Gamma_2(G)=id,$
$\Gamma_1(G)\subset Z(G)$. Thus $\Gamma_1(G)$ is a finite abelian
group.
Now,note that $G/\Gamma_1(G)$ is a non-identity, finitely
generated abelian group. Therefore by structure theorem of
finitely generated abelian groups, there exists a chain
$$H_{r+1}=G\supset H_r\supset\cdots\supset H_1\supset
\Gamma_1(G)=H_0$$ of subgroups of $G, r\geq 0,$ such that
$H_{i+1}/H_i$ is an infinite cyclic group for all
$i=0,1,\cdots,r-1$ and $G/H_r$ is a finite abelian group. As
$H_i\supset \Gamma_1(G)$, $H_i\lhd G$ for all $i=0,1,\cdots,r$
 (Lemma \ref{l101l}). If $r=0,$ then $G/\Gamma_1(G)$ is
finite. Hence $G$ is finite as $\Gamma_1(G)$ is finite.
 Thus the claim
follows in this case. Now, let $r\geq 1$. By theorem
\ref{t8.22tt}, there exists a torsion free, normal subgroup $M$ of
$H_1$ with $[H_1:M]<\infty$. Hence for some $\alpha\geq 1,$
$H^{\alpha}_{1} =gp\{x^{\alpha} \st x\in H_1\} \subset M$ is a
torsion free normal subgroup of $G$ (use:$H_1\rhd G$). Now, note
that $H_1/H_{1}^{\alpha}$ is a nilpotent group (Lemma
\ref{l10.66}) and is torsion. Hence, as $H_1/H_{1}^{\alpha}$ is
finitely generated (Theorem \ref{t13.47tt}), it is finite (Lemma
\ref{l13.46ll}). A repeated use of theorem \ref{t8.22tt} gives a
torsion free normal subgroup $K$ of $G$ such that $K\subset H_r$
and $[H_r:K]<\infty.$ As $[G:H_r]<\infty,$ we get
$[G:K]=[G:H_r][H_r:K]<\infty$. Hence the result follows in case
$m=2.$ If $m>2,$ then since $\Gamma_1(G)$ is nilpotent of class
$m-1$, by induction, there exists a torsion free normal subgroup
$N$ of $\Gamma_1(G)$ with $[\Gamma_1(G):N]<\infty.$ Now, as in
case $m=2$, we get a torsion free, normal subgroup $K$ of $G$ with
$[G:K]<\infty$. If $G$ is not finite then $K\neq id.$ By theorem
\ref{t10.75}, $K\cap Z(G)\neq id.$ Hence, as $K$ is torsion free,
$\circ(Z(G))=\infty.$ This contradicts our assumption on $Z(G).$
Consequently $\circ(G)<\infty.$             $\blacksquare$

\vspace{1in}
 {\bf EXERCISES}
\begin{enumerate}
  \item Let $A$ and $B$ be solvable subgroups of a group $G$ and $A\lhd
  G.$ Prove that $AB$ is a solvable group.
  \item
  Let $p,q$ be, not necessarily distinct, primes and let $G$ be a group of
  order $p^nq, n\geq 1.$ Prove that $G$ is solvable.
  \item
  Prove that a group of order 105 is solvable.
  \item
  Let $p,q,r$ be primes. Prove that any group of order $pqr$ is solvable.
  \item
  Let $H,K$ be two solvable, normal subgroup of a group $G$. Prove that $HK$
  is a solvable, normal subgroup of $G$.
  \item
  Let $G$ be a finite group. Prove that $G$ contains a solvable, normal
  subgroup $H$ such that for any solvable, normal subgroup $K$ of $G$,
  $K\subset H.$
  \item
  Let $G$ be a finite group. Prove that $G$ is solvable if and only if for
  any subgroup $H(\neq e)$ of $G$, $H^{\prime}\neq H.$
  \item
  Let $G$ be a solvable group, and $H$ its maximal subgroup such that
  $[G:H]<\infty .$ Prove that $[G:H]=p^{\alpha}$ where $p$ is a prime and
  $\alpha \geq 1.$
  \item
  Let $G(\neq id)$ be a solvable group. Prove that $G$ contains an
  abelian, characteristic subgroup $H(\neq id).$
  \item
  Let $N$ be a minimal normal subgroup of a finite solvable group $G$.
  Prove that $N$ is an elementary abelian $p$-group for a prime $p$.
  \item
  Let $G$ be a finite nilpotent group and $\circ(G)=n$. Prove that for any
  $m\geq 1, m\mid n,$ there exists a normal subgroup $H$ of $G$ with
  $\circ(H)=m.$
  \item
  Show that derived group of the symmetric group $S_4$ is not nilpotent.
  \item
  Let for a group $G, G/Z(G)$ is nilpotent. Prove that $G$ is nilpotent.
  \item
  Let $G$ be a finite nilpotent group. Prove that for any maximal subgroup
  $H$ of $G,$ $[G:H]$ is prime.
  \item
  Let for a group $G$ and a subgroup $H$ of $Z(G), G/H$ be nilpotent.
  Prove that $G$ is nilpotent.
  \item
  Let $G$ be a nilpotent group and $H$, a simple normal subgroup of $G$.
  Prove that $H\subset Z(G)$.
  \item
  Let $G$ be a torsion free nilpotent group and let $m\geq 1, n\geq 1.$
  Prove that if for
  $x,y\in G$, $x^my^n=y^nx^m,$ then $xy=yx.$
  \item
  Prove that if $G$ is a torsion free nilpotent group, then $G/Z(G)$ is
  torsion free nilpotent. Further, show that $G/Z_i(G)$ is torsion free
  nilpotent for all $i\geq 2.$
  \item
  Let $n\geq 2$ be a fixed integer. Prove that the group
  \[G=\left\{\left(\begin{array}{ccc}
    1 & n^k & a \\
    0 & 1 & b \\
    0 & 0 & 1 \
  \end{array}\right)\st a,b,k\in \Z\right\}\] is nilpotent of class 2 and
  $G/\Gamma_{\!1}(G)$ is not torsion free.
  \item
  Let $G$ be a finitely generated nilpotent group. Prove that there exists
  a central chain of $G$ with all its factor groups cyclic.
  \item
  Prove that the torsion subgroup of a finitely generated nilpotent group
  is finite.
  \item
  Let $G$ be a group.
  Prove that $\Gamma_n(G)=gp\{[a_1,a_2,\cdots,a_{n+1}]\st a_i\in G \mbox{
  for all }1\leq i \leq n+1\}$ for all $n\geq 0.$
  \item
Let $G$ be a nilpotent group. Prove that if $A$ is a maximal, normal
abelian subgroup of $G$, then $C_G(A)=A.$

  \item
  Let $G=D_{2n},$ $(n\geq 1),$ be the Dihedral group of order $2^{n+1}$
  (Example \ref{th11}). Prove that $G$ is nilpotent of class $n+1.$
  \item
  Prove that a finite group $G$ is nilpotent if and only if every maximal
  subgroup of $G$ is normal.
  \item
    Let $ G$ be a nilpotent group , and let $A$, $B$ be two normal
    subgroups of $G$ where $B\subsetneqq A$. Prove that there exists a
    normal subgroup $K$ of $G$, $B\subsetneqq K\subset A$ such that $K/B$ is cyclic.
  \item
  Let $N$ be a normal subgroup of a nilpotent group $G$ such that
  $\circ(N)=p^n$ ($p$ : prime, $n\geq 1$). Prove that $N\subset Z_n(G)$.
  \item
  Prove that every finite group $G$ contains a nilpotent subgroup $A$ such
  that $\overline{N_G}(A)=G$.
  \item
  Let $G$ be a finite group such that for any subgroup $H$ of $G$, $HG'=G$
  implies $H=G.$ Prove that $G$ is nilpotent.
  \item
  Prove that for any finite group $G$, $\Phi(G)$ is nilpotent.
  \item
  Let $G$ be a finite $p$-group. Prove that $G/\Phi(G)$ is an elementary
  abelian $p$-group.
\end{enumerate}

\chapter{Free Groups}
 \footnote{\it contents group14.tex }
If for a subset $X$ of a group $G$, $gp\{X\} = G$, then any $g(\neq e)\in
G$ can be expressed as \be\label{wet} g = x_1^{n_1}x_2^{n_2}\cdots
x_k^{n_k} \ee where $x_i \in X$, $n_i\neq 0$ and $x_i\neq x_{i+1}$ for any
$1\leq i < k.$ In general, such representation is not unique e.g., if for
an element $x\in X$, $o(x) = m <\infty$, then $$g=x_1^{n_1}x_2^{n_2}\cdots
x_k^{n_k} =
                      x_1^{n_1}x_2^{n_2}\cdots
                      x_k^{n_k}x^m$$
Hence, in this case $g(\neq e) \in G$ shall admit more than one
representation of the type (\ref{wet}). Therefore, if we wish to have a
generator set $X$ for a group $G$ with the property that for any $g(\neq
e) \in G$, the representation of $g$ of the type(\ref{wet}) is unique then
all the elements in $X$ have to have order $\infty$. This, however, is not
sufficient. To see this, note that all the non-zero elements of the group
$G = (\Q , +)$  have order $\infty$ and $G$ is not cyclic. If $r,s$ be any
two non-zero elements in a generator set of $(\Q , +)$, then there exist
non-zero integers $m$ and $n$ such that $mr = ns$. Thus $(\Q , +)$ can not
have such a set of generators. We, therefore, define:
\bd\label{d11.1} A subset $X$ of a group $G\neq id$ is called a free
set of generators for $G$ if every $g(\neq e) \in G$ can be
uniquely expressed as $$g = x_1^{n_1}x_2^{n_2}\cdots x_k^{n_k}$$
where $x_i \in X$, $n_i\neq 0$, $k>0$  and $x_i\neq x_{i+1}$
 for any $1\leq i<k.$
\ed
\bd\label{d11.2}
A group $G\neq id$ is called a {\bf free group} if $G$ has a free set of
generators $X$, and the set $X$ is called its basis.
\ed
\brs
(i) An infinite cyclic group $C = gp\{a\}$ is free and $\{a\}$ is a
free set of generators for $C$.\\
(ii) By convention the identity group is assumed free with
$\emptyset$(null set) its free set of generators.
\ers
\bexe\label{exe11.3}
Prove that if for a free group $G$ with a basis $X$, $\abs{X} >
1$, then $G$ is non-abelian.
\eexe
\bexe\label{exe11.4}
Let $X$ ($ \neq \emptyset$ ) be a set of generators for a group $G$. Prove
that $X$ is a free set of generators for $G$ if and only if the following
condition holds:\\ If for $z_i \in X$, $1\leq i\leq l$, $z_i \neq
z_{i+1}$, $$ {z_1}^{m_1}{z_2}^{m_2}\cdots {z_l}^{m_l} = e \;\;\;\;( m_i
\in \Z),$$ then $m_i = 0 $ for all $1\leq i\leq l$.
\eexe
We, now, give a concrete example of a free group with a free set of
generators having two elements.
\bt\label{t11.5}
Let $n\geq 2$ be an integer. The subgroup
$$H =  gp \left\{\left[\begin{array}{cc}
  1 & n \\
  0 & 1
\end{array}\right],\left[\begin{array}{cc}
  1 & 0 \\
  n & 1
\end{array}\right] \right\}$$
of $Sl_2(\Z)$ is a free group and
$$\left\{\left[\begin{array}{cc}
  1 & n \\
  0 & 1
\end{array}\right],\left[\begin{array}{cc}
  1 & 0 \\
  n & 1
\end{array}\right] \right\}$$
is a free set of generators of $H$.
\et
\pf Put
$x = \left[\begin{array}{cc}
  1 & n \\
  0 & 1
\end{array}\right]$ and $y = \left[\begin{array}{cc}
  1 & 0 \\
  n & 1
\end{array}\right]$. To prove the result we have to show that
there is no non-trivial relation among $x$ and $y$ i.e., any
alternating product $w$ of non-zero powers of $x$ and $y$ is
non-identity. For any $k\neq 0$, we have
$$x^k = \left[\begin{array}{cc}
  1 & nk \\
  0 & 1
\end{array}\right] \mbox{ and }y^k = \left[\begin{array}{cc}
  1 & 0 \\
 nk & 1
\end{array}\right]$$ Thus $x^k \neq id$ and $y^k\neq id$. To
complete the proof, we, now consider:\\
Case I : $$ w = x^{k_1}y^{k_2} \cdots z^{k_s}$$
where $s\geq 2$, $k_i \neq 0$ for all $i\geq 1$ and $z = x$ or
$y$. Let $r_i$ denotes the first row of the matrix
$x^{k_1}y^{k_2} \cdots z^{k_i}.$ Then $r_1 = (1, nk_1)$,
$r_2 = (1+nk_2(nk_1), nk_1)$ and if $r_{2t-1} = (a_{2t-1},a_{2t})$,
then
\ba{rcl}
r_{2t} = r_{2t-1}y^{k_{2t}} & = & (a_{2t-1},a_{2t})\left[\begin{array}{cc}
                                                      1 & 0 \\
                                                   nk_{2t} & 1
                                                     \end{array}\right]\\
 & = & (a_{2t-1}+nk_{2t}a_{2t}, a_{2t})\\
 & = & (a_{2t+1},a_{2t}) \mbox{(say)}
 \ea
 and
 \ba{rcl}
 r_{2t+1} = r_{2t}x^{k_{2t+1}} & = & (a_{2t+1},a_{2t})\left[\begin{array}{cc}
                                            1 & nk_{2t+1} \\
                                                   0 & 1
                                                   \end{array}\right]\\
  & = & (a_{2t+1}, a_{2t}+nk_{2t+1}a_{2t+1})\\
  & = & (a_{2t+1}, a_{2t+2}) \mbox{(say)}
  \ea
Thus we have
$$ a_{i+2} = a_i + nk_{i+1}a_{i+1}$$
for all $i = 1,2,\cdots ,s-1$. We have seen that $a_1 = 1$, $a_2 =
nk_1$ and $a_3 = a_1 + nk_2a_2$. As $n\geq 2$, and $k_i\neq 0$,
$\abs{a_2} = n\abs{k_1} > 1= \abs{a_1}$, and
\ba{rcl}
\abs{a_3} = \abs{a_1 + nk_2a_2} & \geq & n\abs{k_2a_2} -
                                          \abs{a_1}\\
    & \geq & 2\abs{a_2} - \abs{a_1}\\
    & \geq & \abs{a_2} + 1 > \abs{a_2}
\ea
Now, if $i\geq 2$, and $\abs{a_{i+1}} > \abs{a_i}$, then
\ba{rcl}
\abs{a_{i+2}} & = & \abs{a_i + nk_{i+1}a_{i+1}}\\
             & \geq & n\abs{a_{i+1}} - \abs{a_i}\\
             & \geq & \abs{a_{i+1}} + 1 > \abs{a_{i+1}}
\ea Hence $w\neq id.$ \\  Case II : $$ w = y^{k_1}x^{k_2} \cdots z^{k_s}$$
where $s\geq 2$, $k_i \neq 0$ for all $i\geq 1$ and $z = x$ or $y$.\\ In
this case $$ w_1 = y^{-k_1}wy^{k_1} = x^{k_2} \cdots z^{k_s}y^{k_1}$$ Thus
$w_1$ is a product of the form in case I, and hence $w_1\neq id.$
Consequently $w\neq id.$ This completes the proof.     $\blacksquare$

Let $X$ be a non-empty set. We shall construct a group with $X$ as
a free set of generators. Let us define a word in the alphabet $X$
a finite product $$ w = x_1^{m_1}x_2^{m_2}\cdots x_k^{m_k}$$
where $x_i\in X$ and $m_i\in \Z$ for all $1\leq i \leq k.$ The
integer $\sum_{i=1}^k \abs{m_i}$ is called the length of the word $w$
and is denoted by $l(w)$. The word $w$ is called reduced if $m_i\neq 0$
, and $x_i \neq x_{i+1}$ for all $i\geq 1.$ To make matters
precise, if $X = \{x,y,z\}$, then $$ w = x^{-3}x^2y^3z^2$$ is not
a reduced word as first two alphabets of the word are equal, but
$$ w = x^2y^3zx^4$$ is a reduced word. It is clear that each word
can be brought to a reduced word by collecting powers of alphabets
when adjacent alphabets are equal and by deleting zeroth powers
and continuing this process till we obtain a reduced word, e.g.,
if $$ w = x^3x^2z^0yz^4z^{-2}x^0z^{-2}y^{-2}$$
then
\ba{rcl}
w & = & x^5yz^2z^{-2}y^{-2}\\
 & = & x^5yz^0y^{-2}\\
 & = & x^5yy^{-2}\\
 & = & x^5y^{-1}
 \ea
 is now a reduced word or is in reduced form. The steps undertaken
 to bring a word in reduced form can be implemented in more than
 one way, e.g., we can also do
 \ba{rcl}
 w & = &(x^3x^2z^0yz^4)(z^{-2}x^0z^{-2}y^{-2})\\
   & = & (x^5yz^4)(z^{-2}z^{-2}y^{-2})\\
   & = &(x^5yz^4)(z^{-4}y^{-2})\\
   & = & x^5yz^4z^{-4}y^{-2}\\
   & = & x^5yz^0y^{-2}\\
   & = & x^5yy^{-2} = x^5y^{-1}
\ea In the above example the end result is same. Note that reducing the
word $x^0y^0$ gives a word with no symbols. We call a word with no symbols
as an empty word. As seen above each word in the alphabet $X$ can be
brought to reduced form. This, however, can be done in more than one ways.
We shall prove below that any word reduces to a unique reduced word by
above process of simplification. We shall multiply two words $w_1$ and
$w_2$ just by writing one after the other. Clearly even if $w_1$ and $w_2$
are reduced the product $w_1w_2$ may not be reduced because the last
alphabet of $w_1$ may be equal to the first alphabet of $w_2$. However,
$w_1w_2$ can be brought to reduced form by simplification explained above.
We, now, prove:
\bt\label{t11.6}
Each word in the alphabet $X$ can be simplified to exactly one
reduced word i.e., the process of simplification is immaterial.
\et
\pf For each $x\in X$, define a function $\phi_x$ on the set of
reduced words $\overline{W}$ in the alphabet $X$ as below:
\ba{rcl}
\phi_x : \overline{W} &\to & \overline{W}\\
           w & \mapsto & \overline{xw}
\ea Here $\overline{xw}$ is a reduced word obtained by simplification of
the word $xw$. Clearly $\overline{xw}$ is uniquely determined. Thus
$\phi_x$ is well defined. It is easy to check that $\phi_x$ is one-one and
onto from $\overline{W}$ to $\overline{W}$. Thus $\phi_x$ defines a
permutation on $\overline{W}$. Now, for any word $$ u =
x_1^{m_1}x_2^{m_2}\cdots x_k^{m_k}$$ in the alphabet $X$, define $$ \phi_u
= \phi_{x_1}^{m_1}\phi_{x_2}^{m_2}\cdots
                       \phi_{x_k}^{m_k}$$
the composite of the permutations $\phi_{x_i}$ and of their inverses. Thus
for each word $u$, $\phi_u$ is a permutation of $\overline{W}$. Clearly if
a word $u$ reduces in some way to a word $v$ then $\phi_u = \phi_v$. Hence
 $$ \begin{array}{cl}
    & \{\phi_u \mid u: \mbox{a word in the alphabet }X \} \\
   = & \{\phi_v \mid v:
\mbox{a reduced word in the alphabet } X \}
 \end{array} $$ is a subgroup of the group of
permutations (transformations)
 of the set $\overline{W}$. Now, let a word $w$ reduces by
two different ways to two reduced words $w_1$ and $w_2$, then $$ \phi_w =
\phi_{w_1} = \phi_{w_2}$$ Now, as $\phi_{w_i}$ sends the empty word to
$w_i$, $i =1,2$, we get $w_1 = w_2$. Hence the result follows.
$\blacksquare$ \\
\\ {\bf NOTATION}
: Let $X (\neq \emptyset)$ be a subset of a group $G$. Given two words
$w_1,w_2$ in the alphabet $X$, we choose to write $w_1\equiv w_2$ if $w_1$
and $w_2$ are same as written. Clearly, even if $w_1\not\equiv w_2$, we
may have $w_1 = w_2$ in $G$.
\bt\label{t11.7}
Let $X$ be a non-empty set. Then there exists a free group with
$X$ as a free set of generators.
\et
\pf By the theorem \ref{t11.6} each word in the alphabet $X$ can
be simplified to a unique reduced word. Let for a word $w$ in the
alphabet $X$, $\overline{w}$ denotes the unique reduced word
obtained by simplifying $w$. Then for any two words $w_1$ and
$w_2$, we have $\overline{\overline{w_1}}=\overline{w_1}$ and
$\overline{\overline{w_1}.\overline{w_2}}=\overline{w_1.\overline{w_2}}
 =\overline{w_1w_2}$. Now, consider the set
$$ F(X) = \{w\mid w: \mbox{ a reduced word in the alphabet} X\}$$ of all
reduced words in the alphabet $X$ and define a binary operation $\star$ on
$F(X)$ as follows : $$ w_1 \star w_2 = \overline{w_1.w_2}$$ For any
$w_1,w_2$ and $w_3$ in $F(X)$, we have \ba{rcl} (w_1\star w_2)\star w_3 &
= & \overline{w_1.w_2}\star w_3\\
                        & = & \overline{\overline{w_1.w_2}.w_3}\\
                        & = & \overline{(w_1w_2)w_3}\\
                       & = & \overline{w_1(w_2w_3)}\\
                       & = & \overline{w_1\overline{(w_2.w_3)}}\\
                       & = & w_1\star (w_2\star w_3)
\ea Thus $\star$ is an associative binary operation. Note that for any
reduced word $ u = x_1^{m_1}x_2^{m_2}\cdots x_k^{m_k}$, $ v =
x_k^{-m_k}\cdots x_2^{-m_2} x_1^{-m_1}$ is a reduced word and $u\star v =
v\star u$ is the empty word which is the identity element of $F(X)$ with
respect to the binary operation $\star$. Hence $F(X)$ is a group. The
words with single alphabet in $X$ clearly form a set of generators of
$F(X)$.
If $$x_1^{m_1}\star x_2^{m_2}\star \cdots \star x_k^{m_k} = id,$$ where
$m_i \in \Z$ and $x_i\neq x_{i+1}$ for any $i\geq 1$, then as
$x_1^{m_1}x_2^{m_2}\cdots x_k^{m_k}$ is a reduced word, it is the empty
word. Hence $m_i = 0$ for all $1\leq i\leq k$. Thus $X$ is a free set of
generators of the group $F(X)$.      $\blacksquare$
\br
The group $F(X)$ is free group generated by the set $X$. In fact
$X$ is a basis of $F(X)$. If $X = \{x_1,x_2,\cdots ,x_n \}$, then
we write $F_n(X)$ for $F(X)$ and if $X = \{x_1,x_2,\cdots ,x_n,\cdots  \}$
, then $F(X)$ is denoted by $F_{\infty}(X)$. In case it is not
necessary to specify $X$, we simply write $F_n$ and $F_{\infty}$.
\er
\bexe\label{exe11.8}
Let $H$ be the subgroup of $F(x,y)$ generated by the set
$$ S = \{xyx^{-1},x^2yx^{-2},\cdots ,x^nyx^{-n} \}, n\geq 1.$$
Prove that $S$ is a free set of generators for $H$.
\eexe
\br
Any free group $F_n, n\geq 1$ can be embedded in $F_{\infty}$. As for free
group $F(x,y)$, $H = gp\{ x^{-n}yx^n \mid n\geq 1 \}$ is a free subgroup
of $F(x,y)$ (Exercise \ref{exe11.8}), $F_{\infty}$ can be embedded in $F_2$. Thus, as
from the theorem \ref{t11.5}, $Sl_2(\Z)$ contains $F_2$, $F_{\infty}$ and
hence $F_n$ for all $n\geq 1$ can be embedded in $Sl_2(\Z)$.
\er
\bt\label{t11.9}
Any free group $G(\neq id)$ is isomorphic to the free group $F(X)$
for some $X$.
\et
\pf Let $A = \{a_i \mid i\in I \}$ be a free set of generators for
$G$ and let $X = \{x_i \mid i\in I \}$ be a set of alphabets with
$\abs{A} = \abs{X}$. Consider the map
\ba{rcl}
\phi: F(X) & \to & G\\
u = {x_{i_1}}^{m_1}{x_{i_2}}^{m_2}\cdots {x_{i_k}}^{m_k}
    & \mapsto & {a_{i_1}}^{m_1}{a_{i_2}}^{m_2}\cdots {a_{i_k}}^{m_k}
\ea from the group $F(X)$ to $G$. By definition of $F(X),$ and free set of
generators of a group, it is clear that $\phi$ is an isomorphism.
$\blacksquare$
\bco\label{co11.10}
If $G$ is a free group of rank $n$, then $G$ is isomorphic to
$F_n$.
\eco
\br
In the above theorem, if $G$ is not necessarily a free group and we
take $A$ to be a set of generators of $G$ in the proof, then
$\phi$ is an epimorphism from $F(X)$ to $G$. Thus any group is
homomorphic image of a free group. The kernel $N$ of the
epimorphism $\phi : F(X)\to G$ is called the set of relations of
the group $G$ in the alphabet $X$. If for a subset $R$ of $N$, $N$
is the smallest normal subgroup of $F(X)$ containing $R$, then $R$
is called a set of defining relations of $G$ in the alphabet $X$.
Since $G\cong \frac{F(X)}{N}$, $G$ can be completely determined
with the help of $X$ and the subset $R$ of $F(X)$. We call the
pair $<X\mid R>$, a presentation of $G$ in terms of generators and
relations. This fact is expressed by writing $G = <X\mid R>$ (the
equality is upto isomorphism). If $X = \{x_1,x_2,\cdots ,x_n \}$ is
a finite set and $R = \{w_1,w_2,\cdots ,w_m \}$ is also a finite
set of reduced words in the alphabet $X$, then
$$ G = <X\mid R> = <x_1,x_2,\cdots ,x_n \mid w_1,w_2,\cdots ,w_m
>$$
is called a finitely presented group.
\er
We shall give below presentations of some familiar groups, to
elaborate the above discussion, in terms of generators and
relations.
\bex\label{ex11.11}
$S_3 = <x,y \mid x^2,y^3, {(xy)}^2 >$
\eex
\pf For the set of generators $\{(12), (123)\}$ of $S_3$, consider
the epimorphism :
$$ \phi : F(x,y) \to S_3$$ for which $\phi(x) = (12)$ and
$\phi(y) = (123)$. Now, let $M = \ker \phi$ and let $N$ be the
normal closure of the set $R = \{x^2,y^3, {(xy)}^2\}$ in $F(x,y)$.
Clearly $x^2,y^3, {(xy)}^2 \in M$. Hence $N\subset M$. Consider
the natural epimorphism :
\ba{rcl}
\frac{F(x,y)}{N} & \to & \frac{F(x,y)}{M}\\
wN &\mapsto & wM
\ea
As $S_3$ is isomorphic to $\frac{F(x,y)}{M}$,
$o(\frac{F(x,y)}{N}) \geq o(\frac{F(x,y)}{M}) = 6.$ Clearly
$xN,yN$ generate $\frac{F(x,y)}{N}$ and
\be\label{**}
 {(xN)}^2 = N, {(yN)}^3 = N, \mbox{ and } xyN = y^2xN
\ee
since $x^2$, $y^3$ and ${(xy)}^2$ belong to $N$. Using the
relations \ref{**}, we can easily see that
$$\{ N,xN,yN,y^2N,xyN,xy^2N\}$$
is the group $\frac{F(x,y)}{N}$ i.e., it  has six elements. Hence
$\frac{F(x,y)}{N}$ is isomorphic to $\frac{F(x,y)}{M}$ which is
isomorphic to $S_3$. Thus $M=N$ and consequently
$S_3 = <x,y \mid x^2,y^3, {(xy)}^2 >.$
\bex\label{ex11.12}
$\Z_n = <x \mid x^n>$
\eex
\bex\label{ex11.13}
$\Z_6 = <x \mid x^6> = <x,y \mid x^3,y^2,xyx^{-1}y^{-1} >$
\eex
\bexe\label{exe11.14}
Prove that the presentation
$$ <x_1,x_2,\cdots ,x_n \mid x_ix_j{x_i}^{-1}{x_j}^{-1}, 1\leq
i<j\leq n >$$
gives a free abelian group of rank $n$.
\eexe
The example \ref{ex11.13} above gives two different presentations
of $\Z_6$. Thus, in general, a group may have more than one
presentation.
\bt\label{t11.15}
For the free group $F_n$, the quotient group
$\frac{F_n}{{F_n}^{\prime}}$ is isomorphic to the free abelian
group $\Z^n$.
\et
\pf Let $\{x_1,x_2,\cdots ,x_n\}$ be a free set of generators for
$F_n$. Then the map
\ba{rcl}
\phi : F_n & \to & \Z^n\\
u = {x_{i_1}}^{m_1}{x_{i_2}}^{m_2}\cdots {x_{i_k}}^{m_k}
    & \mapsto & m_1e_{i_1}+m_2e_{i_2}+\cdots +m_ke_{i_k}
\ea where $e_i = (0,0,\cdots ,1^{ith},0\cdots ,0)$, $i= 1,2,\cdots ,n,$ is
the canonical basis of $\Z^n$, is an epimorphism from $F_n$ to $\Z^n.$ As
$\Z^n$ is an abelian group ${F_n}^{\prime}\subset \ker\phi$. Hence $\phi$
induces the epimorphism: $$\overline{\phi} : \frac{F_n}{{F_n}^{\prime}}
\to \Z^n$$ such that for any $u = {x_{i_1}}^{m_1}{x_{i_2}}^{m_2}\cdots
{x_{i_k}}^{m_k}$ in $F_n$, $\overline{\phi}(u{F_n}^{\prime}) =
m_1e_{i_1}+m_2e_{i_2}+\cdots +m_ke_{i_k}.$ As $\frac{F_n}{{F_n}^{\prime}}$
is an abelian group , any non-identity element of
$\frac{F_n}{{F_n}^{\prime}}$ can be written as
${x_{i_1}}^{m_1}{x_{i_2}}^{m_2}\cdots {x_{i_k}}^{m_k}{F_n}^{\prime}$ where
$i_1<i_2<\cdots i_k$, $k>0$, and $m_t(\neq 0)\in \Z$ for all $1\leq t\leq
k.$ Thus it is clear that no non-identity element of
$\frac{F_n}{{F_n}^{\prime}}$ is mapped to identity under
$\overline{\phi}$. Hence $\overline{\phi}$ is an isomorphism.
$\blacksquare$
\br
If $F(X)$ is a free abelian group on the set of alphabet
$X(\neq \emptyset)$ where $X$ is not necessarily finite, then from
the proof of the theorem it is clear that
$$\frac{F(X)}{{F(X)}^{\prime}}\cong \bigoplus_{x\in X} \Z_x$$
where $\Z_x = \Z$ for all $x\in X.$
\er
\bexe\label{exe11.16}
Show that $F_n$ contains a normal subgroup of index $a^m$ for all $a\geq
1$, and $0\leq m\leq n.$\\ $$\left(\mbox{ Hint : Use the composite epimorphism: } F_n
\to \frac{F_n}{{F_n}^{\prime}}\cong {\Z}^n \to \left({\frac{\Z}{a}}\right)^m\right)$$
\eexe
\bt\label{t11.17}
If $F_m$ is isomorphic to $F_n$, then $m = n.$
\et
\pf Let $$\psi : F_m\to F_n $$ be an isomorphism. Then $\psi
({F_m}^{\prime}) = {F_n}^{\prime}$, and hence $\psi$ induces the
isomorphism \ba{rcl} {\displaystyle\overline{\psi}:
\frac{F_m}{{F_m}^{\prime}}} & \to &
{\displaystyle\frac{F_n}{{F_n}^{\prime}}}\\ u{F_m}^{\prime} &\mapsto &
\psi(u){F_n}^{\prime} \ea where $u\in F_m$. Thus by the theorem
\ref{t11.15}, we get that $\Z^m$ is isomorphic to $\Z^n$, and hence by
lemma \ref{l93}, $m = n.$
\bco\label{co11.18}
If a free group $G(\neq id)$ has a basis $X$ with $\abs{X} = n <
\infty$, then any other basis of $G$ has $n$ elements.
\eco
\pf Let $Y$ be a basis of $G$. Then, any $x\in X$ can be expressed as $$ x
= y_1^{n_1}y_2^{n_2}\cdots y_k^{n_k}$$ where $y_i \in Y$ and $n_i\neq 0$
for some $1\leq i\leq k$. Thus there exists a finite subset
$\{y_1,y_2,\cdots ,y_l \}$ of $Y$ such that $G = gp\{y_1,y_2,\cdots ,y_l
\}$. As $Y$ is a basis of $G$ it is clear that $\abs{Y} < \infty$. Now,
the result follows by the corollary \ref{co11.10} and the theorem.
\br
If $F(X) = G$ is a free group on the set of alphabet $X\neq \emptyset$,
and $Y$ is any other free set of generators of $G$, then
$\abs{X}=\abs{Y}$. \\ We have
$$\frac{G}{G^{\prime}} = \frac{F(X)}{F^{\prime}(X)} \cong \bigoplus_{x\in X}
\Z_x$$ and also $$ \frac{G}{G^{\prime}} = \frac{F(Y)}{F^{\prime}(Y)} \cong
\bigoplus_{y\in Y} \Z_y .$$ Thus the free abelian group $A =
\bigoplus_{x\in X} \Z_x$ is isomorphic to the free abelian group $B =
\bigoplus_{y\in Y} \Z_y $. Consequently $$ \frac{A}{2A} \cong
\frac{B}{2B}$$ As $\frac{A}{2A}$ and $\frac{B}{2B}$ are vector spaces over
$\Z_2$ of dimensions $\abs{X}$ and $\abs{Y}$ respectively, $\abs{X} =
\abs{Y}$.\er We, now, define:

\bd\label{d11.19}
If $X$ is a free set of generators for a group $G$, then the cardinal
number $\abs{X}$ is called the rank of the free group $G$.
\ed
\bt\label{t11.20}
A set of generators $X$ of a group $G$ is a free set of generators
of $G$ if and only if given any group $H$ and a function $\alpha$
from $X$ to $H$ there exists a unique homomorphism $\phi$ from $G$
to $H$ such that $$\phi(x) = \alpha(x)$$ for all $x\in X$.
\et
\pf We shall first prove the direct part. Clearly, the map \ba{rcl}
\hat{\alpha}: F(X) & \to \ & H\\ w = x_1^{m_1}x_2^{m_2}\cdots x_k^{m_k} &
\mapsto & {\alpha (x_1)}^{m_1}{\alpha (x_2)}^{m_2}\cdots {\alpha
(x_k)}^{m_k} \ea is a homomorphism from the free group $F(X)$ to $H$.
Further, as $X$ is a free set of generators of $G$, the natural map
\ba{rcl}
\eta : F(X) & \to & G\\
w = x_1^{m_1}x_2^{m_2}\cdots x_k^{m_k} & \mapsto &
x_1^{m_1}x_2^{m_2}\cdots x_k^{m_k} \ea which sends every reduced word in
$F(X)$ to the corresponding product of elements in $G$ is an isomorphism.
Hence $$ \phi = \hat{\alpha}{\eta}^{-1} : G \to H $$ is a homomorphism
such that $$\phi (x) = \hat{\alpha}{\eta}^{-1} (x) = \alpha (x) $$ for all
$x\in X$. For uniqueness of $\phi$, let $\psi$ be another homomorphism
from $G$ to $H$ such that $\psi (x) = \alpha (x)$ for all $x \in X$. As
$X$ is a set of generators of $G$ any element $g\in G$ can be expressed as
$g = x_1^{m_1}x_2^{m_2}\cdots x_t^{m_t},$ where $x_i \in X$ for $1\leq
i\leq t.$ Hence \ba{rcl} \phi (g) & = & {\phi (x_1)}^{m_1}{\phi
(x_2)}^{m_2}\cdots {\phi (x_t)}^{m_t}\\ & = & {\alpha (x_1)}^{m_1}{\alpha
(x_2)}^{m_2}\cdots {\alpha (x_t)}^{m_t}\\ & = & {\psi (x_1)}^{m_1}{\psi
(x_2)}^{m_2}\cdots {\psi (x_t)}^{m_t}\\ & = & \psi (g) \ea Thus $\phi$ is
uniquely determined. To prove the converse, take $H = F(X)$ and $\alpha
(x) = x$ for all $x\in X$. Then $\alpha$ extends to a unique homomorphism
$$ \phi : G \to F(X).$$ If $X$ is not a free set of generators for $G$,
then for some $t>1$, there exist $x_i \in X$, $x_i\neq x_{i+1}$ and
$m_i\neq 0$ for $1\leq i\leq t$ such that $x_1^{m_1}x_2^{m_2}\cdots
x_t^{m_t} = id$. Therefore \ba{rcl} \phi (x_1^{m_1}x_2^{m_2}\cdots
x_t^{m_t}) & = & id.\\ \Rightarrow x_1^{m_1}x_2^{m_2}\cdots x_t^{m_t} & =
& id \ea in $F(X)$. This, however, is not possible as $X$ is a free set of
generators for $F(X)$. Hence $X$ is a free set of generators for
$G$.$\blacksquare$\\ We shall, now, prove that if $G$ is a free group of
rank $n$, $n < \infty$, then any set of $n$ generators of $G$ is a free
set of generators of $G$. To prove this we need the following result.
\bl\label{l11.21}
Every free group is residually finite.
\el
\pf Let $G (\neq id)$ be a free group with $X$, a free set of generators
for $G$. For each $a\in X$, let $\alpha_a$ denote the function from $X$ to
$C_a = gp\{a\}$ such that $\alpha_a (a) = a$, and $\alpha_a (x) = e $ for
all $x\neq a$. Then, by the theorem \ref{t11.20}, $\alpha_a$ extends to
the unique homomorphism $$ \phi_a : G \to C_a.$$ Now, for each $n\geq 1$,
consider the composite homomorphism $$ \phi_{a,n} : G
\stackrel{\phi_{a}}\rightarrow C_a \stackrel\eta\rightarrow \frac
{gp\{a\}}{gp\{a^n\}}$$ where $\eta$ is the natural map. Let $K_{a,n} =
\ker \phi_{a,n}.$ Then $[G:K_{a,n}] = n $. If $w\in G$ lies in all
subgroups of $G$ with finite index, then $\phi_{a,n} (w) = id$ for all $a$
and $n\geq 1$. This clearly implies that $w = id.$ Hence $G$ is a
residually finite group. $\blacksquare$
\bt\label{t11.22}
If $G$ is a free group of rank $n < \infty,$ and $$S =\{a_1,a_2,\cdots
,a_n \}$$
 is a set of generators of $G$, then $S$
is a free set of generators of $G$.
\et
\pf Let $X = \{x_1,x_2,\cdots ,x_n \}$ be a free set of generators of $G$.
By the theorem \ref{t11.20}, there exists a unique homomorphism $\phi : G
\to G$ such that $\phi (x_i) = a_i$ for all $1\leq i \leq n$. As $S$ is
set of generators of $G$. $\phi$ is onto. Let ${\cal A} = \{ H\subset G :
subgroup \mid [G:H] < \infty\}$. Then, by the corollary \ref{l7.1l},
\ba{rcl} \hat{\phi} : {\cal A} & \to & {\cal A}\\ H & \mapsto &
{\phi}^{-1}(H) \ea is one-one from ${\cal A}$ onto ${\cal A}$. Hence, if
$K = \ker\phi$, then $K\subset H$ for all $H\in {\cal A}$. By lemma
\ref{l11.21}, $$ \bigcap_{H\in {\cal A}}H = \{e\}.$$ Hence $K = \{e\}$.
Thus $\phi$ is an automorphism of $G$. Consequently $S$ is a free set of
generators of $G$.$\blacksquare$
\bexe\label{exe11.23}
Show that a free group of rank $k$ can not be generated by less
than $k$ elements.
\eexe
\bexe\label{exe11.24}
Show that $F_m\times F_n$ is not a free group for all $m\geq 1, n\geq 1.$
\eexe

 Let $G_i, 1\leq i\leq n,$ be groups, and let
 $$X = \coprod_{1\leq i\leq n} G_i$$ be the disjoint union of
 $G_i$'s. We can form words $w = x_1x_2\cdots x_k$ from elements of
 $X$. A word $w = x_1x_2\cdots x_k$ is called reduced if $x_j$ and
 $x_{j+1}$ do not lie in the same group and no $x_j$ is equal to
 the identity element of any group $G_i$, $ 1\leq i\leq n.$ All
 the reduced words formed by elements of $X$ together with the
 empty word form a group if we define the product of two words
 $w_1,w_2$ as $w_1\star w_2$, the reduced word obtained after
 normal simplification of $w_1w_2$, the word obtained after
 writing $w_2$ adjacent to $w_1$. This group is denoted by
 $$G_1\star G_2\star \cdots \star G_n$$ and is called {\bf free product}
 of $G_i$'s.
 \bexe\label{exe11.25}
 For the free group $F(x_1,x_2,\cdots ,x_n) = F_n$,
 $$ F_n \cong C_1\star C_2\star \cdots \star C_n$$
 where $C_i = gp\{x_i \}$, $1\leq i\leq n$, is the cyclic group
 generated by $\{x_i \}$.
 \eexe

 Let $F(X)$ be a free group on a non-empty set $X$ of alphabets.
 We shall prove that every subgroup of $F(X)$ is free. To prove
 this, we need some new notions and results. Let us define:
 \bd\label{d11.26}
 A (left) Schreier system in a free group $F(X)$ is a set $S$ of
 reduced words in $F(X)$ such that if the reduced word
 $$w = x_1^{m_1}x_2^{m_2}\cdots x_k^{m_k}$$
 is in $S$, then $x_i^{m_i}x_{i+1}^{m_{i+1}}\cdots x_k^{m_k}\in S$
 for all $i\geq 2$.
 \ed
 \bd\label{d11.27}
Let $H$ be a subgroup of a free group $F(X)$. A left Schreier
transversal of $H$ in $F(X)$ is Schreier system which is also
a left transversal of $H$ in $F(X)$.
\ed
\bl\label{l11.28}
Let $H$ be a subgroup of a free group $F(X)$. Then there exists a (left)
Schreier transversal of $H$ in $F(X)$.
\el
\pf For any left coset $fH$ of $H$ in $F(X)$, put $$l(fH) = \min \{l(w)
\mid w\in fH \mbox{ is a reduced word } \}.$$ Clearly, $l(fH)$ is the
length of the shortest reduced word in $fH$. We shall call this the length
of the coset $fH$. We shall construct a (left) Schreier transversal of $H$
in $F(X)$ by choosing a coset representative $\Phi(f)$ from each coset
$fH$ of $H$ in $F(X)$ with $l(\Phi(f))=l(fH),$ using induction on the
lengths of the cosets. For the coset $H$, choose the empty word as its
representative. Note that $H$ is the only left coset of $H$ in $F(X)$ with
length $0$. For every left coset of $H$ with length $1$ choose a word of
length $1$ as its representative. Now, let for every left coset of $H$
with length $s<t$, $t\geq 2$, we have chosen a representative with length
$s$. Let $gH$ is any left coset of $H$ in $F(X)$ with $l(gH) = t$. Assume
$l(g) = t$, and $g = x_1^{e_1}x_2^{e_2}\cdots x_t^{e_t}$, where $x_i\in X
$, and $e_i = \pm 1$ for all $1\leq i\leq t$. We define $$\Phi(g) =
x_1^{e_1} \Phi({x_2^{e_2}}x_3^{e_3}\cdots x_t^{e_t})$$ Then $\Phi(g)$ is a
representative of the coset $gH$ and has length $t$. This completes the
induction step by choosing a coset representative for every left coset of
length $t$ as above. The fact that such a (left) transversal of $H$ in
$F(X)$ is a Schreier transversal is clear by construction. Hence the proof
follows. $\blacksquare$
\bt\label{t11.29}{\bf (Nielson- Schreier)}
Any subgroup $H$ of a free group $F(X)$ is free.
\et
\pf We can assume that $H\neq \{e\}$, since otherwise the result is
trivial.Let $T$ be a (left) Schreier transversal of $H$ in $F(X)$ and let
$\Phi$ be the corresponding coset representative function. By the remark
\ref{d23.02f}, $H$ is generated by the set $S$ of non-identity elements in
$\{ \Phi{(xt)}^{-1}xt \mid x\in X, t\in T\}.$ We
shall prove that $H$ is freely generated by the elements of $S$. This will
be done in steps.\\ {\bf Step I : } If $\Phi{(xt)}^{-1}xt \neq \{e\}$, then
it is a (non-empty) reduced word in the alphabet $X$ as written.\\
 As $\Phi{(xt)}^{-1}$ and $t$ are reduced words, any cancellation
 in the word $\Phi{(xt)}^{-1}xt$ has to begin with $x$. Therefore
 for the cancellation to occur either $t \equiv x^{-1}t_1$, in which
 case, $\Phi{(xt)}^{-1}xt = t_1^{-1}t_1 = e$, or $\Phi{(xt)}^{-1} \equiv
 t_2^{-1}x^{-1}$, which implies
 \ba{rcl}
 \Phi(xt)H & = & xt_2H\\
 \Rightarrow xtH & = & xt_2H\\
 \Rightarrow tH & = & t_2H\\
 \Rightarrow t = t_2\\
 \Rightarrow \Phi{(xt)}^{-1}xt & = & t_2^{-1}x^{-1}xt_2 = e
 \ea
 a contradiction to our assumption. Thus $\Phi{(xt)}^{-1}xt$, as written, is a
 non-empty reduced word in the alphabet $X$.\\
 Before moving to step II, we shall make an observation. First of
 all, note that $\Phi (x^{-1} \Phi (xt)) = t.$ Hence
 \ba{rcl}
 {(\Phi{(xt)}^{-1}xt)}^{-1} & = & t^{-1}x^{-1}\Phi(xt)\\
                            & = &
                            \Phi(x^{-1}\hat{t})^{-1}x^{-1}\hat{t}
 \ea
 where $\hat{t} = \Phi(xt)$. Thus any element of $S\cup S^{-1}$
 has the form $\phi(x^{e} t)^{-1}x^e t$ where $x\in X, t\in T$
  and $e = \pm 1$. We, now, state:\\
 {\bf Step II : } Let $u,v\in S\cup S^{-1}$, and $uv\neq e$. If
$u = \Phi(x^{\epsilon} t_1)^{-1}x^{\epsilon} t_1$ and $v = \Phi(y^{\delta}
t_2)^{-1}y^{\delta} t_2$ where $\epsilon, \delta = \pm 1$, $t_1,t_2 \in T$
and $x,y\in X$, then $x^{\epsilon}$ and $y^{\delta}$ do not cancel out in
the reduced form of $uv$ as a word in the alphabet $X$.\\ By the step I,
$u,v$ are non-empty reduced words in the alphabet $X$. Therefore any
cancellation in the word $uv$ in the alphabet $X$ has to begin at the
interface of $u$ and $v$. We claim that any cancellation in $uv$ in the
alphabet $X$ shall stop before reaching $x^{\epsilon}$ in $u$ and
$y^{\delta}$ in $v$. Let $x^{\epsilon}$ cancels earlier than $y^{\delta}$
then we must have $$ \Phi(y^{\delta} t_2)^{-1} \equiv
t_1^{-1}x^{-\epsilon}w^{-1}$$ Taking $s = \Phi(y^{\delta} t_2)$, we get
\ba{crl}  & wx^{\epsilon}t_1  = & s\in T\\ \Rightarrow  & x^{\epsilon}t_1  \in
& T\;\;\; \mbox{ since $T$ is a Schreier transversal}\\ \Rightarrow &
\Phi(x^{\epsilon}t_1)  = & x^{\epsilon}t_1\\ \Rightarrow & u  = & e \ea
This contradicts our assumption that $u\neq e$. Thus $x^{\epsilon}$ does
not cancel earlier than $y^{\delta}$. Next, let $y^{\delta}$ cancels
earlier than $x^{\epsilon}$, then we have \ba{rcl} t_1 & \equiv &
w_1y^{-\delta}\Phi(y^{\delta} t_2)\\
    & \equiv & w_1y^{-\delta} s^{\prime}
\ea
where $s^{\prime} = \Phi(y^{\delta} t_2)$. As $w_1y^{-\delta}
s^{\prime} = t_1\in T$, and $T$ is a Schreier transversal,
$y^{-\delta} s^{\prime}\in T$. Hence
\ba{crl}
 & \Phi(y^{-\delta} s^{\prime})  = & y^{-\delta} s^{\prime}\\
\Rightarrow & \Phi(y^{-\delta}\Phi(y^{\delta} t_2))  = &
             y^{-\delta}\Phi(y^{\delta} t_2)\\
\Rightarrow & t_2  = & y^{-\delta}s^{\prime}\;\;\;  \\  & i.e.,\;\; t_2  = &
y^{-\delta}\Phi(y^{\delta} t_2)\\ \Rightarrow & v  = & \Phi(y^{\delta}
t_2)^{-1}y^{\delta} t_2 = e \ea This is a contradiction to our assumption
that $v\neq e$. Thus $y^{\delta}$ does not cancel before $x^{\epsilon}$.
\\ Finally, if $y^{\delta}$ and $x^{\epsilon}$ cancel simultaneously then
we have \ba{crl}
 & x^{\epsilon}t_1 \equiv & y^{-\delta}\Phi(y^{\delta} t_2)\\
 \Rightarrow & \Phi(x^{\epsilon}t_1)  = &
                        \Phi(y^{-\delta}\Phi(y^{\delta} t_2))\\
                        & = & t_2\\
\Rightarrow & uv  = & e \ea which is not true by our assumption. Hence the
statement in the Step II is proved.\\ {\bf Step III : } The set $S$ is a free
set of generators of $H$.\\ We shall first prove that $S\cap S^{-1} =
\emptyset$. If not, then there exists $x,y\in X$, and $t_1,t_2\in T$ such
that \ba{rcl} \Phi(xt_1)^{-1}xt_1 & = & (\Phi(yt_2)^{-1}yt_2)^{-1}\\
\Rightarrow \Phi(xt_1)^{-1}xt_1\Phi(yt_2)^{-1}yt_2 & = & e\\ \ea and
$$\Phi(yt_2)^{-1}yt_2\Phi(xt_1)^{-1}xt_1 = e$$ Then by the Step II, we get
$$\Phi(xt_1) = t_2, \Phi(yt_2) = t_1 \mbox{ and } xt_1 =
y^{-1}\Phi(yt_2)$$ These relations imply $x = y^{-1}$, which is absurd.
Hence $S\cap S^{-1} = \emptyset$. Now, let $$ \alpha =
u_1^{m_1}u_2^{m_2}\cdots u_k^{m_k}$$ be a reduced word in the elements of
$S$ i.e., $k\geq 1$, $u_i\in S$, $m_i\neq 0$ and $u_i\neq u_{i+1}$ for all
$i = 1,2,\cdots ,k-1.$ To prove the result it suffices to show that
$\alpha \neq e$. This will follow if we prove that $\alpha$ as a reduced
word in the alphabet $X$ is non-empty. This is now immediate from the
Steps I and II since by the Step I each $u_i$ is a reduced word in the
alphabet $X$ and by the Step II any cancellation clearly begins at the
interface of adjacent $u_i$'s and stops before reaching the basic symbols
$x\in X$ (use:$S\cap S^{-1} = \emptyset$). Thus the theorem is proved.
$\blacksquare$
\bco\label{co11.30}
Let $H$ be a subgroup of a free group $F(X)$, and let $T$ be a
(left) Schreier transversal of $H$ in $F(X)$. Then the set of
non-identity elements of $T^{-1}XT\cap H$ is a free set of
generators of $H$.
\eco
\pf We shall use notations in the proof of the theorem. Clearly $$ \{
\Phi(xt)^{-1}xt \mid t\in T, x\in X\} \subset T^{-1}XT\cap H$$ Further, if
for any $s,t\in T$ and $x\in X$, \ba{rcl} s^{-1}xt & = & h\in H\\
\Rightarrow xt & = & sh\\ \Rightarrow \Phi(xt) & = & \Phi(sh) = s \ea
Hence $$ \{ \Phi(xt)^{-1}xt \mid t\in T, x\in X\} = T^{-1}XT\cap H$$
Therefore the proof is immediate from the theorem.
\bl\label{l11.31}
The derived group of the free group $F_2 = F(x,y)$ is free of
countable infinite rank.
\el
\pf As $\frac{F_2}{F_2^{\prime}}$ is free abelian group of rank 2
generated by $xF^{\prime}_{2}$ and $yF_{2}^{\prime}$ (Theorem
\ref{t11.15}), any element of this group has the form $x^my^nF_2^{\prime}$
$(m,n\in \Z)$. Further, by the theorem \ref{t11.15}, the map \ba{rcl} \phi
: \frac{F_2}{F_2^{\prime}} &\to & \Z^2\\
     x^my^nF_2^{\prime}&\mapsto & me_1+ne_2
\ea where $e_1 = (1,0)$ and $e_2 = (0,1)$, the canonical basis of the free
abelian group $\Z^2$, is an isomorphism. Hence it is immediate that the
set $T = \{x^my^n \mid m,n\in \Z \}$ is a Schreier transversal of
$F_2^{\prime}$ in $F_2$. Thus by the proof of the Nielsen-Schreier theorem
the set
 $ \{(x^my^{n+1})^{-1}yx^my^n \mid m,n\in \Z \}$ is a free set of
 generators of $F_2^{\prime}$. Hence the result
 follows.$\blacksquare$
 \bt\label{t11.32}
 Let $H$ be a subgroup of the free group $$F_n = F(x_1,x_2,\cdots
 ,x_n)$$ with $[F_n : H] = j < \infty$. Then $H$ is a free group of
 rank $j(n-1)+1.$
 \et
 \pf Put $X = \{x_1,x_2,\cdots ,x_n\}$. Let $T$ be a left Schreier
 transversal of $H$ in $F_n$ and let $\Phi$ be the corresponding
 coset representative function. For the collection $M$ of all formal
 expressions i.e., words
$\Phi(xt)^{-1}xt$ where $t\in T$ and $ x\in X$, consider the map \ba{rcl}
\alpha : T_1 = T-\{e\} & \to & M\\ t & \mapsto & \alpha(t)= \left\{
\begin{array}{c}
  \Phi(xt)^{-1}xt \mbox{ if } t\equiv x^{-1}t_1,\mbox{ for } x\in X \\
  \Phi(xt_1)^{-1}xt_1 \mbox{ if } t\equiv xt_1,\mbox{ for } x\in X
   \end{array}\right.
\ea We shall show that $\alpha$ is one-one. Let for $x,y\in X$, and $t,s \in
T_1$,  $t\equiv xt_1$, and $s\equiv y^{-1}s_1$. Then \ba{crl} & \alpha(t)  \equiv
& \alpha(s)\\ \Rightarrow & \Phi(xt_1)^{-1}xt_1  \equiv & \Phi(ys)^{-1}ys\\
\Rightarrow & t^{-1}xt_1  \equiv & s_{1}^{-1}ys \\ \Rightarrow & s \equiv
& xt_1 \\ \Rightarrow & s \equiv & t \\ \mbox{i.e.,} & s = &
t.\ea   In all other cases, it is clear that $\alpha(t)\neq \alpha(s)$
for $t\neq s$ in $T_1$. Therefore $\alpha$ is one-one. Now, let for $t\in
T_1, y\in X$, \ba{rcl} \Phi(yt)^{-1}yt & = & e \;\;\;\mbox{  in }
F_n\\ \Rightarrow \Phi(yt) & = & yt \ea If $\Phi(yt)\not\equiv yw$, then
$t\equiv y^{-1}t_1$ and $\alpha(t) = \Phi(yt)^{-1}yt.$ Further, if
$\Phi(yt) \equiv yw$, then $t = w$, $\Phi(yt)\in T_1$, and \ba{rcl}
\alpha(\Phi(yt)) & = & \Phi(yw)^{-1}yw\\
                 & = & \Phi(yt)^{-1}yt
\ea Finally, for $e\in T$ and $x\in X$, if \ba{rcl} \Phi(xe)^{-1}xe & = &
e\\ \Rightarrow \Phi(x) & = & x \in T_1 \ea Hence $$ \alpha(\Phi(x)) =
\Phi(xe)^{-1}xe$$ Therefore $\alpha(T_1)$ is the set of all expressions in
$M$ which represent identity element in $F_n$. Now, note that $M$ has $nj$
expressions and $\alpha(T_1)$ consists of exactly those expressions in $M$
which represent identity. By the theorem \ref{t11.29}, $H$ is freely
generated by the expressions in $M$ which are non-identity. Hence, $H$ is
free of rank $nj-(j-1) = j(n-1)+1$.       $\blacksquare$  \\
\vspace{.5in} \\
  {\bf A Graph theoretic proof of Nielsen-Schreir theorem} \\  \\ We shall
  first introduce some graph theoretic preliminaries. \bd A graph $\Gamma$
  consists of two sets, a non-empty set $X$ of elements called vertices of
  $\Gamma$ and a subset $A(\neq \emptyset)$ of $X \times X$ such that : \\
   (i) $(x,x)\not\in A$ for any $x\in X$. \\  (ii) $\alpha =(x_1,x_2)\in
   A$ implies $\overline{\alpha} =(x_2,x_1)\in A.$  \\ The elements of $A$
   are called edges of $\Gamma$. \ed {\bf Notation : } (i) We shall denote
   a graph $\Gamma$ by $\Gamma(X,A)$ where $X$ is the set of vertices of
   $\Gamma$ and $A$ is the set of edges. \\  (ii) If $\alpha =(x_1,x_2)\in
   A$ is an edge of $\Gamma$, then $i(\alpha)=x_1$ is called initial
   vertex of $\alpha$ and $t(\alpha)=x_2$ is called terminal vertex of
   $\alpha$.   \bd Let $\Gamma(X,A)$ be a graph. For any two vertices
   $x_1,x_2$ of $\Gamma$, a path from $x_1$ to $x_2$ is a finite sequence,
   $P : \alpha_1,\alpha_2,\cdots,\alpha_m$, of edges of $\Gamma$ such that
   $i(\alpha_1)=x_1, t(\alpha_m)=x_2$ and $t(\alpha_s)=i(\alpha_{s+1})$
   for all $1\leq s\leq m-1.$ The integer $m$ is called the length of the
   path $P$. \ed \bd Let $\Gamma(X,A)$ be a graph. A path of the form
   $\alpha\overline{\alpha}$ in $\Gamma$ is called a round trip.\ed
    \bd A graph $\Gamma(X,A)$ is called connected if given
   any two vertices $x_1, x_2$ of $\Gamma$ there is a path $P$ in $\Gamma$
   from $x_1$ to $x_2$. \ed \bd A graph $\Gamma(X,A)$ is called a tree if
   $\Gamma$ is a connected graph such that any path in $\Gamma$ which
   joins a vertex to itself contains a round trip. \ed \bd A path $P$ in a
   graph $\Gamma$ is called reduced if $P$ contains no round trip. \ed
   \bl Let $\Gamma(X,A)$ be a tree. Then any two distinct vertices in
   $\Gamma$ are connected by a unique reduced path. \el  \pf Let $x_1\neq
   x_2$ be two vertices of $\Gamma$, and let $P_1$, $P_2$ be two reduced
   paths connecting $x_1$ to $x_2$. Further, assume that $P_1$ is the
   reduced path in $\Gamma$ connecting $x_1$ to $x_2$ with smallest length
   any other reduced path connecting two vertices of $\Gamma$ with length
   less than that of $P_1$ is unique. Now, consider the path
   $P_1P_{2}^{-1}$ connecting $x_1$ to $x_1$. Since $\Gamma$ is a tree
   $P_1P_{2}^{-1}$
    has to contain a round trip. Thus either last edge
   occurring in $P_1$ equals last edge of $P_2$ or first edge of $P_1$
   equals the first edge of $P_2$. In the first case removing last edge
   $\alpha$ of $P_1$ and $P_2$ give two reduced paths from $x_1$ to
   $i(\alpha_1)=d$. This contradicts our assumption that $P_1$ is
   non-unique reduced path of smallest length connecting two distinct
   vertices of $\Gamma$. In the other case deleting first edge $\beta$
   from $P_1$ as well as $P_2$ give two distinct reduced paths from
   $t(\beta)$ to $x_2$. This again gives a contradiction as above. Hence
   the result follows. \br If $\Gamma(X,A)$ is a tree, then any path
   connecting two distinct vertices $x_1$, $x_2$ of $\Gamma$ is obtained
   from the unique reduced path connecting $x_1$ to $x_2$ by inserting a
   finite number of round trips. \er  \bd If $\Gamma(X,A)$ is a tree, then
   the unique reduced path connecting two distinct vertices $x_1$, $x_2$ of
   $\Gamma$ is called a Geodesic and is denoted by
   $\overrightarrow{x_1x_2}$. \ed \bl Let $\Lambda$ be a subtree of the
   tree $\Gamma(X,A)$. If $u,v$ are two distinct vertices of $\Lambda$,
   then $\overrightarrow{uv}$, the geodesic in $\Gamma$ connecting $u$ to
   $v$, is contained in $\Lambda$.\el \pf If $\overrightarrow{uv}$ is not
   contained in $\Lambda$ then $\overrightarrow{uv_{\Lambda}}$, the
   geodesic connecting $u$ to $v$ in $\Lambda$,
   is a path from $u$ to $v$ in $\Gamma$ distinct
   from $\overrightarrow{uv}$. Hence $\overrightarrow{uv_{\Lambda}}$ is
   obtained from $\overrightarrow{uv}$ by inserting finite number of round
   trip in $\overrightarrow{uv}$. This, however, is not true as
   $\overrightarrow{uv_{\Lambda}}$ is a geodesic. Hence
   $\overrightarrow{uv}$ is contained in $\Lambda$.  \bd Let
   $\Gamma(X,A)$ be a graph and $G$ be a group. We say that $G$ acts on
   $\Gamma$ if $G$ acts on $X$ as well as $A$ such that for any $g\in G$
   and $\alpha\in A$, $g(\overline{\alpha})=\overline{g(\alpha)}$,
   $g(i(\alpha))=i(g(\alpha))$, moreover, $g(\alpha)\neq
   \overline{\alpha}$. \ed \br For any $g\in G$, $\alpha\in A$,
   \[g(t(\alpha))=g(i(\overline{\alpha}))=t(g(\alpha)).\] Thus if
   $\alpha=(x_1,x_2)\in A$, then $g(\alpha)=(g(x_1),g(x_2)).$ \er \bd Let
   $G$ be a group acting on a graph $\Gamma(X,A)$. We say that $G$ acts
   freely on $\Gamma$ if the stabilizer of each vertex of $\Gamma$ is the
   identity subgroup of $G$. \ed Assume that a group $G$ is acting on a
   tree $\Gamma(X,A)$. Let $\mathcal{C}$ be the set of all sub-trees of $\Gamma$
   which contain atmost one vertex and one edge from each orbit of $X$ and
   $A$ under the given $G$-action. For any $\Lambda_1, \Lambda_2$ in $\mathcal{C}$,
   define: \[\Lambda_1\leq\Lambda_2 \Leftrightarrow
   \Lambda_1\subset\Lambda_2\] Then $(\mathcal{C},\leq)$ is a partially ordered set.
   One can easily check that any chain in $(\mathcal{C},\leq)$ has an upper bound in
   $\mathcal{C}$. Hence, by Zorn's lemma, $\mathcal{C}$ contains a maximal element called a
   {\bf reference tree} for the group action. We, now, prove. \bl Let $G$ be a
   group acting on a tree $\Gamma(X,A)$, and let $\Lambda$ be a reference
   tree in $\Gamma$ with respect to the group action. Then $\Lambda$
   contains precisely one vertex from each $G$-orbit. \el \pf Assume that
   the result is not true. Then there exists a vertex $z\in\Gamma$ such
   that $Gz\cap\Lambda \neq \emptyset$. Choose a vertex $x\in\Lambda$, and
   cosider the geodesic $\overrightarrow{zx}$. Let $\alpha$ be the first
   edge in $\overrightarrow{zx}$, and let $t(\alpha)=y$. If
   $Gy\cap\Lambda\neq\emptyset$, then there exists an element $g\in G$
   such that $g(y)\in\Lambda$. Now, adding $g(z)$, $g(\alpha)$ and
   $g(\overline{\alpha})$ to $\Lambda$ gives a tree in $\Gamma$
   larger than $\Lambda$ which contains at most one vertex and one edge
   from each orbit. This contradicts our assumption that $\Lambda$ is a
   reference tree. In case $Gy\cap\Lambda = \emptyset$, take $y$ for $z$
   and continue as above. We may end up in the situation when
   $\overrightarrow{zx} = (z,x)$, an edge $\beta$ in $\Gamma$. Here again,
   adding $z$, $\beta$ and $\overline{\beta}$ to $\Lambda$, we get a
   tree in $\Gamma$ larger than $\Lambda$ which contains atmost one vertex
   and one edge from each orbit. This gives a contradiction to the
   assumption that $\Lambda$ is a reference tree. Hence the result
   follows. \bt\label{end} If a group $G$ acts freely on a tree $\Gamma(X,A)$, then
   $G$ is a free group.\et \pf We can assume that $G\neq \{e\}$, since
   otherwise there is nothing to prove. Let $\Lambda$ be a reference tree
   in $\Gamma$ with respect to the group action. Put
   \[A^*=\{\alpha\in\Gamma : \mbox{ edge }\st \alpha \not\in \Lambda
   ,\mbox{ but } i(\alpha)\in\Lambda \}\] For each $\alpha\in A^*$,
   $t(\alpha)\not\in \Lambda$. Since $\Lambda$ contains exactly one edge
   from each $G$-orbit, there exists $g_{\alpha}(\neq e)\in G$ such that
   $g_{\alpha}^{-1}(t(\alpha))\in \Lambda$. Note that $g_{\alpha}$ is
   unique for each $\alpha\in A^*$ since if $h_{\alpha}^{-1}(t(\alpha))\in
   \Lambda$ for some $h_{\alpha}\in G$, then as $\Lambda$ is a reference
   tree  \[\begin{array}{cl}
      & g_{\alpha}^{-1}(t(\alpha))= h_{\alpha}^{-1}(t(\alpha)) \\
     \Rightarrow & h_{\alpha}g_{\alpha}^{-1}(t(\alpha))=t(\alpha) \\
     \Rightarrow & h_{\alpha}g_{\alpha}^{-1}=e \mbox{ i.e., } h_{\alpha}=g_{\alpha}. \\
   \end{array}\] Note that for any $\alpha \in A^*$,
   $g_{\alpha}^{-1}(t(\alpha)) \in \Lambda$, and
   $g_{\alpha}^{-1}(i(\alpha)) \not\in \Lambda$ since
   $i(\alpha)\in\Lambda$ and $g_{\alpha}\neq e$. Thus $\beta=
   g_{\alpha}^{-1}(\overline{\alpha})\in A^*$. It is easy to see that
   $g_{\beta}=g_{\alpha}^{-1}$. From each such pair $g_{\alpha},
   g_{\alpha}^{-1}$ select one element and let $S$ be the set of all such
   elements in $G$. We shall prove that $S$ is a free set of generators of
   $G$.  \\  {\bf Step I.} $S$ is a set of generators of $G$. \\ Take any
   $g(\neq e)\in G$ and choose a vertex $v\in \Lambda$. Consider the
   geodesic $\overrightarrow{vg(v)}$. As $G$-acts freely on $\Gamma$ and
   $\Lambda$ is a reference tree with respect to $G$-action
   $\overrightarrow{vg(v)}$ is not contained in $\Lambda$. Let $\alpha_1$
   be first edge in the geodesic $\overrightarrow{vg(v)}$ outside
   $\Lambda$, then clearly $\alpha_1\in A^*$. Apply $g_{\alpha_1}^{-1}$ to
   the geodesic $\overrightarrow{t(\alpha_1)g(v)}$. Then its image is the
   geodesic starting from $g_{\alpha_1}^{-1}(t(\alpha_1))$ and ending at
   $g_{\alpha_1}^{-1}g(v)$. If this lies in
   $\Lambda$, then $g_{\alpha_1}^{-1}g(v)\in \Lambda$ and hence
   $g_{\alpha_1}^{-1}g(v) = v$. Thus $g=g_{\alpha_1}$. However, if the
   geodesic $g_{\alpha_1}^{-1}(t(\alpha_1))$ to $g_{\alpha_1}^{-1}g(v)$ is
   outside $\Lambda$ then let $\alpha_2$ be the first edge outside
   $\Lambda$. Now, apply $g_{\alpha_2}^{-1}$ to the geodesic from
   $t(\alpha_2)$ to $g_{\alpha_1}^{-1}g(v)$. The image path starts from
   $g_{\alpha_2}^{-1}(t(\alpha_2))\in \Lambda$ and ends at
   $g_{\alpha_2}^{-1}g_{\alpha_1}^{-1}g(v)$. If
   $g_{\alpha_2}^{-1}g_{\alpha_1}^{-1}g(v)\in \Lambda$, then as above, $g
   = g_{\alpha_1}g_{\alpha_2}$. If not, then proceed as above. Note that
   at each stage of our operation we get a path with smaller length.
   Continue till the resultant geodesic lies in $\Lambda$. In that case if
   the end point in $g_{\alpha_k}^{-1}\cdots
   g_{\alpha_2}^{-1}g_{\alpha_1}^{-1}g(v)$, then \begin{eqnarray}
     & g_{\alpha_k}^{-1}\cdots g_{\alpha_2}^{-1}g_{\alpha_1}^{-1}g = e \nonumber \\
   \Rightarrow & g = g_{\alpha_1}g_{\alpha_2}\cdots g_{\alpha_k}\label{a+1}
   \end{eqnarray} Let us observe that if at some stage last edge is only
   one outside $\Lambda$. Then it lies in $A^*$. At that stage if the end
   point of the path is $g_{\alpha_{k-1}}^{-1}\cdots g_{\alpha_1}^{-1}(v)$
   and the last edge is $\alpha_k\in A^*$, then as $t(\alpha_k)=
    g_{\alpha_{k-1}}^{-1}\cdots g_{\alpha_1}^{-1}(v)$,
     $g_{k}^{-1}g_{\alpha_{k-1}}^{-1}\cdots g_{\alpha_1}^{-1}(v)=v$. Hence
     $g=g_{\alpha_1}\cdots g_{\alpha_k}$. Thus we conclude that
     $G=gp\{S\}$. Further, note that the expression for $g$ obtained above
     as a product of elements in $S$ is a reduced word $w(g)$ in symbols
     of $S$ i.e., an element $g_{\alpha}$ does not appear adjacent to
     $g_{\alpha}^{-1}$ at any place in the expression. Now, let we start
     with any other path $P$ from $v$ to $g(v)$ instead of the geodesic
     $\overrightarrow{vg(v)}$. Then $P$ is obtained from
     $\overrightarrow{vg(v)}$ by adding a finite number of round trips. We
     shall show that our construction even when we start from $P$ shall
     give same reduced word $w(g)$. For simplicity assume $P$ is obtained
     from $\overrightarrow{vg(v)}$ by adding a round trip
     $\gamma\overline{\gamma}$. If $\gamma$ is not the first edge outside
     $\Lambda$, then by our operation the round trip
     $\gamma\overline{\gamma}$ always result in a round trip. Clearly for
     any $\delta \in A$, $\delta \in \Lambda$ if and only if
     $\overline{\delta}\in \Lambda$. Thus during our process if $\gamma$
     changes to $\tau\not\in\Lambda$, and is the first edge outside
     $\Lambda$,
      then the resulting expression for
     $g$ shall be of the form $g_{\alpha_1}\cdots g_{\alpha_t}g_{\tau}
     g_{\tau}^{-1}g_{\alpha_{t+1}}\cdots g_k$. Thus in reduced form it is
     unchanged. Consequently our assumption holds.  \\  {\bf Step II} : $S$ is a
     free set of generator for $G$.  \\ To prove the assertion it is
     sufficient to show that any expression of the type $g=g_{\alpha_1}
     g_{\alpha_2}\cdots g_{\alpha_m}$, where
     $\alpha_1,\alpha_2,\cdots,\alpha_m$ are in $A^*$ is obtained from a
     suitable path from $v$ to $g(v)$. We shall prove this by induction on
     $m$. If $m=1$, consider the path \[P = \overrightarrow{vi(\alpha_1)}
     (i(\alpha_1),t(\alpha_1))\overrightarrow{t(\alpha_1)g(v)}\] In $P$,
     $(i(\alpha_1),t(\alpha_1))$ is the first edge outside $\Lambda$, and
     clearly $g_{\alpha_1}^{-1}(\overrightarrow{t(\alpha_1)g(v)})$ lies in
     $\Lambda$ since $g_{\alpha_1}^{-1}(t(\alpha_1))\in\Lambda$ and
     $g_{\alpha_1}^{-1}g(v)=v$. Thus the assertion holds for $m=1$. Now,
     let $g_1 = g_{\alpha_2}\cdots g_{\alpha_m}$ and let $Q$ be a path from
     $v$ to $g_1(v)$ which realises the decomposition for $g_1$. Consider
     the path \[P=(\overrightarrow{vi(\alpha_1)})(i(\alpha_1),t(\alpha_1))
     (\overrightarrow{t(\alpha_1)g_{\alpha_1}(v)})g_{\alpha_1}(Q)\] from
     $v$ to $g(v)$. As $v, i(\alpha_1)\in\Lambda$,
     $\overrightarrow{(vi(\alpha_1))}$ lies $\Lambda$. Thus $\alpha_1$ is
     the first edge in $P$ outside $\Lambda$. As per our procedure as a
     first step we consider the path \[P_1 =
     (\overrightarrow{g_{\alpha_1}^{-1}(t(\alpha_1)v)})Q\] Here as
     $g_{\alpha_1}^{-1}(t(\alpha_1))$, $v\in\Lambda$,
     $(\overrightarrow{g_{\alpha_1}^{-1}(t(\alpha_1)v)})$ lies in
     $\Lambda$. Thus $P_1$ realises $g_1$, and consequently $P$ realises
     $g=g_{\alpha_1}g_1$. Hence the result.      $\blacksquare$ \\  \\  We
     shall, now, use the above result to prove the Nielsen-Schreier
     theorem. \\ \hspace*{.3in} Let $G$ be a non-identity group and
     $S$,
     a set of generators of $G$ such that $z^2\neq e$ for any $z\in S$.
     Using this data, we shall define a graph $\Gamma$ as below:
     \[\begin{array}{cl}
       X & =G :\mbox{  set of vertices of }\Gamma \\
       A & =\{(g,gz)\st g\in G, \mbox{ and }z\mbox{ or }z^{-1}\in S\}, \\
        & \;\;\;\mbox{  the set of edges of }\Gamma. \\
   \end{array}\] Clearly $(g,g)\not\in A$ for any $g\in G$ and
   $(gz,g)=(gz,(gz)z^{-1})$ is in $A$. Thus $\Gamma = \Gamma(G,A)$ is a
   graph. \\ If $g_1, g_2\in G$ are any two distinct elements, and \[
   g_{1}^{-1}g_2=z_1z_2\cdots z_k\] where $z_i$ or $z_{i}^{-1}$ is in $S$
   for all $1\leq i\leq k$, then \[P=(g_1,g_1z_1)(g_1z_1,g_1z_1z_2) \cdots
   (g_1z_1\cdots z_{k-1},g_2)\] is a path from $g_1$ to $g_2$ in $\Gamma$.
   Thus $\Gamma$ is a connected graph. We, now, prove: \bl If $G(\neq e)$
   is a free group and $S$ is a free set of generators, then $\Gamma(G,A)$
   is a tree. \el \pf To prove the result it suffices to show that for
   any $g\in G$, any path in $\Gamma$ which joins $g$ to itself contains a
   round trip. Let \[ P : (g,gz_1)(gz_1,gz_1z_2)\cdots (gz_1\cdots
   z_{n-1},g)\] be a path in $\Gamma$ which joins $g$ to itself. Then there
   exists $z_n\in G$ such that \[\begin{array}{cl}
      & g=gz_1z_2\cdots z_{n-1}z_n\;\mbox{  where }z_n\mbox{ or } z_n^{-1}\in S \\
     \Rightarrow & z_1z_2\cdots z_{n-1}z_n = e \\
     \Rightarrow & z_i = z_{i+1}^{-1}\;\mbox{  for some }1\leq i\leq n-1 \\
   \end{array}\] since $S$ is a free set of generators for $G$. Thus for $\alpha
   = (gz_1\cdots z_{i-1},gz_1\cdots z_i)$, the round trip
   $\alpha\overline{\alpha}$ appears in $\Gamma$. Hence $\Gamma$ is
   a tree.            $\blacksquare$ \\   \\  We, now, define an action of
   the group $G$ on $\Gamma(G,A)$ as below: \\ Let $g_1\in G$, $x\in G$,
   and $(g,gz)\in A$. Then define: \[g_1(x)=g_1x\] and \[g_1(g,gz) =
   (g_1g,g_1gz).\] Clearly, this gives a free action of $G$ on $\Gamma$.
   We can, now, prove: \bt (Nielsen-Schreier): Let $G(\neq e)$ be a free
   group. Then any subgroup $H$ of $G$ is tree.\et \pf If $H=\{e\}$, then
   there is nothing to prove. Assume $H\neq \{e\}$, and let $S$ be a free
   set of generator of $G$. Clearly for any $z\in S$, $z^2\neq e$. Thus
   with notation and assumptions as above, $\Gamma(G,A)$ is a tree and $G$
   acts freely on $\Gamma$. Restriction of the action of $G$ on $\Gamma$
   to $H$, gives an action of $H$ on $\Gamma$. As $G$ acts freely action
   of $H$ is free on $\Gamma$. Hence the result follows from theorem
   \ref{end}.    $\blacksquare$

\begin{theindex}
\item{abelian group,\;4}
\item{action of a group on a set},\;{118}
\item{action of $G$ on a graph},\;{272}
\item{affine transformation},\;{10}
\item{automorphism},\;{64}
\item{Alternating group},\;{158}
\item{basis},\;{253}
\item{Burnside lemma},\;{127}
\item{canonical homomorphism},\;{68}
\item{center of group},\;{56}
\item{central chain},\;{235}
\item{centralizer},\;{134}
\item{centralizer},\;{135}
\item{centre},\;{134}
\item{characteristic subgroup},\;{225}
\item{chief series},\;{220}
\item{class equation},\;{140}
\item{collineation},\;{10}
\item{commutator subgroup},\;{56}
\item{commutator},\;{55}
\item{completely reducible},\;{102}
\item{conjugacy classes},\;{134}
\item{conjugate},\;{133}
\item{conjugate},\;{136}
\item{connected graph},\;{271}
\item{cycle of length $k$},\;{148}
\item{cyclic group},\;{29}
\item{Caley's Theorem},\;{121}
\item{Content of $f$},\;{193}
\item{derived group},\;{225}
\item{derived group},\;{56}
\item{diagonal group},\;{10}
\item{direct factor},\;{96}
\item{direct factors},\;{92}
\item{disjoint cycles},\;{148}
\item{divisible group},\;{100}
\item{double coset},\;{171}
\item{Dihedral Group},\;{9}
\item{edges},\;{270}
\item{eft Schreier transversal},\;{265}
\item{elementary abelian $p$-group},\;{203}
\item{elementary group},\;{28}
\item{elementary matrices},\;{28}
\item{endomorphism},\;{64}
\item{epimorphism},\;{64}
\item{equivalent chains},\;{210}
\item{equivalent normal chains},\;{220}
\item{even permutation},\;{155}
\item{exponent of group},\;{27}
\item{external direct product},\;{92}
\item{Euler's function},\;{35}
\item{factor group},\;{49}
\item{factors of the subnormal chain},\;{209}
\item{finitely generated free abelian group},\;{185}
\item{finitely generated group},\;{24}
\item{free group},\;{253}
\item{free product},\;{264}
\item{free set of generators},\;{252}
\item{FIRST ISOMORPHISM THEOREM},\;{73}
\item{Frattini subgroup},\;{26}
\item{Free action of $G$  on $\Gamma$},\;{272}
\item{general linear group},\;{7}
\item{graph},\;{270}
\item{group of transformations},\;{8}
\item{group,\;4}
\item{Geodesic},\;{271}
\item{homomorphism},\;{64}
\item{Holomorph},\;{108}
\item{indecomposable},\;{96}
\item{index of coset},\;{41}
\item{inner automorphism},\;{79}
\item{internal direct product},\;{104}
\item{internal direct product},\;{93}
\item{invariant group},\;{125}
\item{isomorphism},\;{64}
\item{Jordan- Holder Theorem},\;{216}
\item{kernel},\;{69}
\item{Klein's $4$-group},\;{6}
\item{left coset},\;{40}
\item{left transversal},\;{43}
\item{length of  the subnormal chain},\;{209}
\item{length of $G$},\;{218}
\item{linearly independent},\;{185}
\item{lower central chain},\;{235}
\item{Lagrange's Theorem},\;{42}
\item{maximal normal subgroup},\;{48}
\item{minimal normal subgroup},\;{48}
\item{minimal set of generators},\;{24}
\item{monomorphism},\;{64}
\item{natural homomorphism},\;{68}
\item{normal chain},\;{220}
\item{normal closure},\;{50}
\item{normal subgroup},\;{48}
\item{normalizer},\;{135}
\item{Nielson- Schreier Theorem},\;{266}
\item{Nilpotent group},\;{236}
\item{Non-Generator},\;{25}
\item{odd permutation},\;{155}
\item{orbit},\;{124}
\item{order of element},\;{27}
\item{order of group},\;{10}
\item{orthogonal group},\;{15}
\item{p-group},\;{144}
\item{path},\;{270}
\item{periodic},\;{114}
\item{permutation group},\;{8}
\item{principal series},\;{220}
\item{proper normal subgroup},\;{48}
\item{quotient group},\;{49}
\item{Quaternions},\;{6}
\item{rank of the free abelian group},\;{190}
\item{rank of the free group},\;{262}
\item{reduced},\;{271}
\item{reference tree},\;{272}
\item{refinement of the chain},\;{211}
\item{residually finite group},\;{124}
\item{right coset},\;{40}
\item{round trip},\;{270}
\item{semi-direct product},\;{107}
\item{semi-simple group},\;{102}
\item{simple group},\;{159}
\item{simple group},\;{49}
\item{special linear group},\;{7}
\item{special orthogonal group},\;{15}
\item{special unitary group},\;{16}
\item{split extension},\;{105}
\item{sub-normal},\;{218}
\item{subgroup},\;{10}
\item{subnormal chain},\;{209}
\item{symmetric group},\;{8}
\item{symmetric polynomial},\;{155}
\item{Schreier system},\;{265}
\item{Schreier Theorem},\;{215}
\item{Solvable group},\;{224}
\item{Stabilizer},\;{125}
\item{Sylow's first theorem},\;{170}
\item{Sylow's second theorem},\;{172}
\item{Sylow's third theorem},\;{173}
\item{Sylow $p$-subgroup},\;{168}
\item{torsion / periodic group},\;{27}
\item{torsion element},\;{27}
\item{torsion free element},\;{27}
\item{torsion free group},\;{27}
\item{transformation},\;{118}
\item{transitive},\;{127}
\item{tree},\;{271}
\item{trivial homomorphism},\;{68}
\item{trivial subgroups},\;{12}
\item{Transposition},\;{148}
\item{unitary group},\;{16}
\item{upper central chain},\;{234}
\item{vertices},\;{270}
\end{theindex}

\end{document}